\documentclass{agtart_a}
\pdfoutput=1

\usepackage[dvips]{graphicx}      
\usepackage{color}   

\title{Small genus knots in lens spaces have small bridge number}

\author{Kenneth L\,Baker} 
\givenname{Kenneth L\,}
\surname{Baker}
\address{School of Mathematics\\Georgia Institute of Technology\\\newline
Atlanta, GA 30332-0160, USA}
\email{kb@math.gatech.edu}
\urladdr{}

\volumenumber{6}
\issuenumber{}
\publicationyear{2006}
\papernumber{56}
\startpage{1519}
\endpage{1621}

\doi{}
\MR{}
\Zbl{}


\received{12 June 2005}
\revised{}
\accepted{16 August 2006}
\published{11 October 2006}
\publishedonline{11 October 2006}
\proposed{}
\seconded{}
\corresponding{}
\editor{}
\version{}

\arxivreference{math.GT/0612427}

\AtBeginDocument{\let\bar\wbar\let\tilde\wtilde\let\hat\what}
\def\hatT{\smash{\mskip2mu\what{\mskip-2mu T\mskip-1.5mu}\mskip.5mu}}
\def\hatA{\smash{\mskip3.5mu\what{\mskip-3.5mu A\mskip1mu}\mskip-1mu}}
\def\SetFigFont#1#2#3#4#5{\small}%
\def\tsty{\textstyle}

\makeatletter
\def\cnewtheorem#1[#2]#3{\newtheorem{#1}{#3}[section]
\expandafter\let\csname c@#1\endcsname\c@thm}
\makeatother

\newtheorem{thm}{Theorem}[section]
\cnewtheorem{lemma}[thm]{Lemma}
\cnewtheorem{Lemma}[thm]{Lemma}
\makeautorefname{Lemma}{Lemma}
\cnewtheorem{prop}[thm]{Proposition}
\cnewtheorem{conj}[thm]{Conjecture}
\cnewtheorem{question}[thm]{Question}
\cnewtheorem{problem}[thm]{Problem}

\theoremstyle{remark}
\newtheorem{remark}{Remark}[section]

\DeclareMathOperator{\Int}{Int}

\newcommand{\A}{\ensuremath{\mathcal{A}}}
\newcommand{\T}{\ensuremath{\mathcal{T}}}

\newcommand{\bdry}{\ensuremath{\partial}}
\newcommand{\cut}{\backslash}

\newcommand\kreis[1]{\smash{\ensuremath{\mathbin{\settowidth{\dimen7}{\mbox{$\bigcirc$}}\makebox[0pt][l]{$\bigcirc$}\makebox[\dimen7]{{\footnotesize #1}}}}}}

\newcommand\bigkreis[1]{\ensuremath{\mathbin{\settowidth{\dimen7}{\mbox{\large$\bigcirc$}}\makebox[0pt][l]{\large$\bigcirc$}\makebox[\dimen7]{{\footnotesize #1}}}}}

\begin{document}

\begin{asciiabstract}
In a lens space X of order r a knot K representing an element of the 
fundamental group \pi_1 X \cong \Z/r\Z of order s \leq r contains a connected 
orientable surface S properly embedded in its exterior X-N(K)
such that the boundary of S intersects the meridian of K minimally s times. 
Assume S has just one boundary component.  Let g be the minimal genus of 
such surfaces for K, and assume s \geq 4g-1.  Then with respect to the genus 
one Heegaard splitting of X, K has bridge number at most 1.
\end{asciiabstract}

\begin{htmlabstract}
In a lens space X of order r a knot K representing an element of the
fundamental group &pi;<sub>1</sub> X &asymp; <b>Z</b>/r<b>Z</b> of
order s &le; r contains a connected orientable surface S properly
embedded in its exterior X-N(K) such that &part; S intersects the
meridian of K minimally s times.  Assume S has just one boundary
component.  Let g be the minimal genus of such surfaces for K, and
assume s &ge; 4g-1.  Then with respect to the genus one Heegaard
splitting of X, K has bridge number at most 1.
\end{htmlabstract}

\begin{abstract}

In a lens space $X$ of order $r$ a knot $K$ representing an element of
the fundamental group $\pi_1 X \cong \mathbb Z/r\mathbb Z$ of order $s
\leq r$ contains a connected orientable surface $S$ properly embedded
in its exterior $X-N(K)$ such that $\partial S$ intersects the meridian
of $K$ minimally $s$ times.  Assume $S$ has just one boundary
component.  Let $g$ be the minimal genus of such surfaces for $K$, and
assume $s \geq 4g-1$.  Then with respect to the genus one Heegaard
splitting of $X$, $K$ has bridge number at most $1$.

\end{abstract}

\maketitle

\section{Statement of results}

Any knot $K$ in a lens space $X = L(r, q)$, $r>0$, is rationally nullhomologous, ie\ $[K] = 0 \in H_1(X;\Q) \cong 0$.  We say $r$ is the {\em order\/} of the lens space $X$, and we say the smallest positive integer $s$ such that $s[K] = 0 \in H_1(X;\Z) \cong \Z/r\Z$ is the {\em order\/} of the knot $K$.  Note $s \leq r$.  The exterior $X-N(K)$ of $K$ thus contains a connected properly embedded orientable surface $S$ such that when $S$ is oriented $\bdry S$ is coherently oriented on $\bdry \bar{N}(K)$ and intersects the meridian $\mu \subseteq \bdry \bar{N}(K)$ of $K$ minimally $s$ times, ie\ $|\mu \cdot \bdry S| = s$.  Such a surface $S$ is an analogue of a Seifert surface for a knot in $S^3$.  We refer to the genus of a knot $K$ in $X$ as the minimal genus of these ``rational'' Seifert surfaces for $K$.  For this article we will restrict our attention to knots with rational Seifert surfaces that have just one boundary component.

In this paper we prove the following theorem.
\begin{thm}\label{main}
Let $K$ be a genus $g$ knot of order $s$ in a lens space $X$ whose Seifert surfaces have one boundary component.  If $s \geq 4g-1$ then, with respect to the Heegaard torus of $X$, $K$ has bridge number at most $1$.
\end{thm}

\fullref{main} may be curiously rephrased as saying small genus knots in lens spaces have small bridge number.  

In \cite{berge:skwsyls} Berge shows that double-primitive knots (ie\ simple closed curves that lie on a genus $2$ Heegaard surface in $S^3$ and represent a generator of $\pi_1$ for each handlebody) admit lens space surgeries.  We refer to these knots as {\em Berge knots\/}.  Berge further shows that the corresponding knot in the resulting lens space is $1$--bridge.

\fullref{main} may be used to show the following theorem which one may care to compare with \fullref{conj:berge}.
\begin{thm}\label{thm:smallgenusbergeknots}
Let $K'$ be a genus $g$ knot in $S^3$.  If $K'$ admits a lens space surgery of order $r \geq 4g-1$ then $K'$ is a Berge knot.
\end{thm}

\begin{conj}[Berge~\cite{berge:skwsyls}] \label{conj:berge}
If $K$ is a knot in a lens space $X$ with an $S^3$ surgery, then with respect to the  genus one Heegaard splitting of $X$, $K$ has bridge number at most $1$.  In particular, the corresponding knot $K' \subseteq S^3$ is a Berge knot.
\end{conj}

\subsection{A quick overview of knots with lens space surgeries}
For coprime integers $p$ and $q$, a $p/q$ Dehn surgery on a knot in $S^3$ is the process of removing a solid torus neighborhood of the knot and attaching a new solid torus so that on the torus boundary of the knot exterior the new meridian is homologous to $p$ times the old meridian and $q$ times the old longitude.  The rational number $p/q$ (including $1/0$) is called the {\em slope\/} of the surgery.  We say a knot in $S^3$ {\em admits a lens space surgery\/} if some $p/q$ Dehn surgery with $q \neq 0$ yields a lens space.

As a consequence of Thurston's geometrization \cite{thurston:gt3m,thurston:survey}, a knot in $S^3$ is either a torus knot, a satellite knot, or a hyperbolic knot.  For any given torus knot, Moser describes all $p/q$ Dehn surgeries that produce lens spaces \cite{moser:esaatk}.  Satellite knots and the slopes along which they admit lens space surgeries are classified by Bleiler and Litherland \cite{bleilerlitherland:lsads}.  Hyperbolic knots admitting lens space surgeries have yet to be classified though they are conjectured to be Berge knots.  Because the nonhyperbolic knots that admit lens space surgeries are also Berge knots, \fullref{conj:berge} contains this conjecture, albeit in another guise.  Furthermore, by the Cyclic Surgery Theorem \cite{cgls:dsok}, if $p/q$ Dehn surgery on any nontorus knot produces a lens space, then $q=\pm1$, ie\ the surgery slope is {\em integral\/}.

\subsection{Berge's list of double-primitive knots}
Berge also gives a conjecturally complete list of double-primitive knots \cite{berge:skwsyls}.  Based on calculations from this list, Goda and Teragaito propose the following conjecture.
\begin{conj}[Goda--Teragaito~\cite{gt:dsokwylsagok}]\label{conj:gt}
If a hyperbolic knot in $S^3$ of genus $g$ admits a lens space surgery of order $r$, then $2g+8 \leq r \leq 4g-1$.
\end{conj}
If Berge's list is complete, then \fullref{thm:smallgenusbergeknots} not only would confirm this upper bound but also would reprove the following theorem of Rasmussen.
\begin{thm}[Rasmussen~\cite{rasmussen:lssaacogat}]\label{cor:rasmussenbound}
If a nontrivial knot in $S^3$ of genus $g$ admits a lens space surgery of order $r$, then $r \leq 4g+3$.
\end{thm}

\subsection{Towards a conjecture of Bleiler and Litherland}
Oszv\'ath and Szab\'o have shown that the Alexander polynomial $\Delta_K(T)$ of a knot $K$ with a lens space surgery determines $\smash{\widehat{HFK}}$ \cite{os:kfh+lss}.  This then implies that the degree of $\Delta_K(T)$ equals twice the genus of $K$ \cite{os:hdagb}.  At the end of \cite{os:agfh} for each lens space of order at most $26$ they list polynomials among which must be the Alexander polynomial of any knot that yields the lens space by positive integer surgery.  (Indeed they present an algorithm for listing such polynomials associated to any given lens space.)  From the results of \cite{os:kfh+lss} these lists may be refined by removing any polynomial whose nonzero coefficients do not alternate between $+1$ and $-1$.  
For each polynomial $A(T)$ listed for a lens space of order $r < 18$ with the exception of the polynomial listed for $L(14,11)$ we have $4(\frac{1}{2} \deg A(T)) - 1 \leq r$.  Thus by \fullref{thm:smallgenusbergeknots} if $K'$ is a nontrivial knot in $S^3$ admitting a lens space surgery of order $r < 18$ other than $L(14,11)$, then $K'$ is a Berge knot.  Computations discussed by Berge \cite{berge:skwsyls} confirm that $1$--bridge knots in lens spaces of order less than $1000$ with surgery yielding $S^3$ correspond to knots in Berge's list.  Thus we conclude the following theorem.
\begin{thm}
If a hyperbolic knot in $S^3$ admits a lens space surgery of order $r$, then $r=14$ or $r \geq 18$.
\end{thm}
This nearly obtains the following conjecture.
\begin{conj}[Bleiler--Litherland \cite{bleilerlitherland:lsads}] \label{conj:lensspaceorder}
If a hyperbolic knot in $S^3$ admits a lens space surgery of order $r$, then $r \geq 18$.
\end{conj}
The $(-2,3,7)$--pretzel knot which is genus $5$ and has a lens space surgery of order $18$ (and a second of order $19$) is hyperbolic and realizes this conjectured bound.  Note that $L(14, 11)$ is obtained by integral surgery on $T_{(3,5)}$, the $(3,5)$ torus knot, which has genus $4$.  To resolve \fullref{conj:lensspaceorder}, one must address whether $\pm 14$ surgery on a genus $4$ hyperbolic knot $K$ in $S^3$ with $\Delta_K(T) = \Delta_{T_{(3, 5)}}(T)$ can yield $L(14, 11)$.

Note that if for every polynomial $A(T)$ that Ozsv\'ath and Szab\'o's algorithm lists for $L(r,q)$ we have $4(\frac{1}{2} \deg A(T)) - 1 \leq r$, then the above method permits us to determine every knot in $S^3$ that has a surgery producing $L(r,q)$.  Since such knots must be Berge knots, one appeals to the last paragraph of \cite{berge:skwsyls} to determine the knots if $r < 1000$.  

If $r \geq 1000$ and there exists a nontrivial $A(T)$ in Ozsv\'ath and Szab\'o's list for $L(r,q)$, then since Berge's list is not yet known to be complete one must also examine all $1$--bridge knots in $L(r,q)$ for which an integer surgery may yield $S^3$.  This may be done in the manner in which the list at the end of \cite{berge:skwsyls} is produced.  Berge shows that there are finitely many such knots to consider \cite[Theorem~3]{berge:skwsyls}.  Since deciding whether a genus $2$ Heegaard diagram represents $S^3$ is algorithmic \cite{hot:aafrS3}, this is a finite process. 

Further note that Ozsv\'ath and Szab\'o do not (necessarily) list polynomials for $L(r,q)$ that may be Alexander polynomials of knots for which nonintegral surgery yields $L(r,q)$.  As mentioned previously, such knots are known to be only torus knots \cite[Cyclic Surgery Theorem]{cgls:dsok}, and Moser's classification \cite{moser:esaatk} gives a means of checking for these.

For example, for the lens space $L(19,11)$ Ozsv\'ath and Szab\'o list only the polynomial $T^{-5}-T^{-4}+T^{-2}-T^{-1}+1-T+T^2-T^4+T^5$.  Since $4 \cdot 5 - 1 =19$, we examine Berge's list to see that this corresponds only to the $(-2,3,7)$--pretzel knot.  In fact, this is the only knot in $S^3$ for which positive surgery yields $L(19,11)$, and moreover it is hyperbolic.  Indeed, $-L(19,11) = L(19,8)$ may be only obtained by $19/2$ surgery on $T_{(-2,5)}$ and by $19/3$ surgery on $T_{(-2,3)}$, a trefoil.  Thus $L(19,11)$ may also be obtained by $-19/2$ surgery on $T_{(2, 5)}$ and by $-19/3$ surgery on $T_{(2, 3)}$.

We further examine Ozsv\'ath and Szab\'o's list for lens spaces of order $19$.  Because it has genus $4$ we may conclude that $L(19,5)$ is obtained only by $19$ surgery on $T_{(-2,9)}$.  Although $L(19,16)$ may be obtained by $19$ surgery on $T_{(-4,5)}$, its genus is $6$ which is too large for our methods to determine whether there are other knots for which surgery may yield this lens space.  Aside from reflections of the above mentioned torus knots with the corresponding sign change for surgery slope, no other torus knots have lens space surgeries or order $19$.

\subsection[Proof of Theorem 1.2]{Proof of \fullref{thm:smallgenusbergeknots}}
\begin{proof}
If $K'$ admits a lens space surgery $X$ of order $r$, then the surgery slope crosses the $0$--slope minimally $r$ times.  Thus if $K$ is the corresponding knot in the lens space $X$, $K$ has order $r$.  

If $r \geq 4g-1$ then by \fullref{main}, $K$ is at most a $1$--bridge knot in $X$.  If $K$ is a  $0$--bridge knot, then $K'$ is a torus knot and hence a Berge knot.  If $K$ is not $0$--bridge then Theorem~2 of \cite{berge:skwsyls} shows that $K'$ is a double-primitive knot in $S^3$ and hence a Berge knot.  
\end{proof}

\subsection{Further questions}
The proof of \fullref{main} is dependent upon $K$ having a Seifert surface with just one boundary component.
\begin{question}
Do similar results hold if the Seifert surfaces for $K$ have more than one boundary component?
\end{question}

One wonders whether there are other relationships among knot order, genus, and bridge number (or width) for knots in lens spaces.  In a given lens space for any specified order and bridge number one may construct knots of arbitrarily large genus.  Indeed the collection of torus knots of a given order in a lens space ought to contain knots of arbitrarily large genus.  Thus we ask questions akin to \fullref{main}.  In the remainder of this section we restrict ourselves to considering knots in lens spaces whose Seifert surfaces have just one boundary component.

\begin{question}
If $K$ is a genus $g$ knot of order $s$ in a lens space such that $s=4g-2$, is the bridge number of $K$ bounded?
\end{question}

If we instead ask that $s=4g-3$ then nullhomologous genus $1$ knots satisfy the hypothesis.  Assumably Whitehead doubles of knots in lens spaces can be concocted to have arbitrarily large bridge number.  Further alterations to the question are then required such as restricting attention to knots whose order equals that of the lens space.

In another direction, \fullref{main} shows that any genus $1$ knot of order at least $3$ is $1$--bridge, any genus $2$ knot of order at least $7$ is $1$--bridge, etc.  Clearly genus $1$ knots of order $1$ contained in a ball in a lens space correspond to genus $1$ knots in $S^3$.  Indeed there are many other nullhomologous genus $1$ knots in a lens space.

\begin{question}
Is there a characterization of genus $1$ knots of order $2$ in lens spaces?
\end{question}

\begin{problem}
Classify genus $1$ knots of order at least $3$ in lens spaces.
\end{problem}

\begin{question}
If $K$ is a knot of order $s$ in a lens space, then for each integer $b \geq 2$ what is the minimal possible genus of $K$ with bridge number $b$?
\end{question}

\subsubsection{Acknowledgements}
The author would like to thank John Berge and especially John Luecke for many useful conversations.  The author is also indebted to the reviewer for the many constructive comments and suggestions.  This work was partially supported by a VIGRE postdoc under NSF grant number DMS-0089927 to the University of Georgia at Athens.

\section{Preliminaries}\label{sec:prelims}

The notion of thin position for knots in $S^3$ was first introduced by Gabai and employed in his proof of Property~R \cite{gabai:fatto3mIII}.  The theory of thin position for knots has since been greatly developed; see Scharlemann's survey \cite{scharlemann:tpittock}.  As one may observe, the idea of thin position naturally generalizes to knots in three manifolds with a given Heegaard splitting.  We discuss this below in \fullref{sec:thinpositioninlensspaces} for genus 1 Heegaard splittings.

The graphs of intersection arising from two intersecting properly embedded surfaces in the complement of a knot have a history of usefulness in questions about Dehn surgery.  The power of this idea was first demonstrated by Litherland~\cite{litherland:sokistII} and most notably used to great success in Gordon and Luecke's portion of the proof of the Cyclic Surgery Theorem \cite{cgls:dsok} and their solution to the knot complement problem \cite{gl:knotcomplement}.  If one of the two intersecting surfaces becomes a Heegaard surface in some Dehn filling of the knot complement, the theory of thin position may be used to establish some nice properties of the resulting graphs of intersection.  Gabai does this for the standard Heegaard splitting of $S^3$ \cite{gabai:fatto3mIII}, which Rieck promotes to Heegaard splittings in general \cite{rieck:hsomitdfs}.  See Gordon's survey \cite{gordon:cmids}.  We construct such graphs of intersection from a genus 1 Heegaard splitting of our lens space and a Seifert surface for our knot in \fullref{sec:graphsofintersection}.

\subsection{A few words about notation}  
By $A-B$ we denote the usual set difference, the set of points in $A$ that are not in $B$.  If $A$ and $B$ are two intersecting manifolds then by $A \cut B$ we denote $A$ ``cut along'' $B$, the closure of $A-B$ in the path metric on $A$.

A regular open neighborhood of $A$ is denoted $N(A)$.  Its closure is denoted $\bar{N}(A)$.

\subsection{Thin position in lens spaces}\label{sec:thinpositioninlensspaces}

Let $h \co X \to \R \cup \{\pm \infty\}$ be the height function induced by the genus one Heegaard splitting of the lens space $X = L(r,q)$ so that $\hatT_z = h^{-1} (z)$ is a torus for $z \in \R$ and $\hatT_{\pm \infty} = h^{-1}(\pm \infty)$ are circles.  The tori $\hatT_z$, $z \in \R$, are {\em level tori\/}, and $\hatT_{\pm \infty}$ are the {\em circles at infinity\/}.  Each level torus separates $X$ into two solid torus components; the component $X^+$ {\em above\/} containing $\hatT_{+\infty}$ and the component $X^-$ {\em below\/} containing $\hatT_{-\infty}$.

Let $K$ be a knot in $X$.   By an isotopy of $K$ we may assume $K \subseteq X - \{\hatT_{\pm \infty} \}$ and that $h|K$ is a Morse function so that $h|K$ has finitely many critical points, all of which are nondegenerate and have distinct critical values. 

Given such a Morse presentation of $K$, let $\hatT_{z_1}, \dots, \hatT_{z_n}$ be level tori so that exactly one is between each consecutive pair of critical levels.  Define the {\em width\/} of the Morse presentation to be $\sum_{i=1}^n |\hatT_{z_i} \cap K|$.  The {\em width of $K$\/} is the minimum of all widths of Morse presentations of $K$.  A Morse presentation of $K$ that realizes the width of $K$ is said to be a {\em thin presentation of $K$\/}, and $K$ itself is said to be in a {\em thin position\/}. 

If the first critical value of $h|K$ above the level torus $\hatT_{z_i}$ is a maximum and the first critical value below is a minimum, then $\hatT_{z_i}$ is a {\em thick level (torus)\/}.  Similarly, if the first critical value of $h|K$ above the level torus $\hatT_{z_i}$ is a minimum and the first critical value below is a maximum, then $\hatT_{z_i}$ is a {\em thin level (torus)\/}.  Every Morse presentation of $K$ must have a thick level torus, but not all Morse presentations of $K$ have a thin level torus.  

Over all Morse presentations of $K$ that have no thin level tori, one that minimizes the width of $K$ is said to be a {\em bridge presentation of $K$\/}, and $K$ itself is said to be in a {\em bridge position\/}.  If $K$ is in bridge position and $\hatT_0$ is the thick level torus, then the {\em bridge number of $K$\/} is $\frac{1}{2}|K \cap \hatT_0|$.  We say a knot with bridge number $n$ is an {\em $n$--bridge knot\/}.

In the event that $K$ may further be isotoped to lie as an embedded curve in a level torus, then either $K$ is a trivial knot (and hence bounds a disk in $X$) or $K$ is a {\em lens space torus knot\/}.  For both of these cases we will say that $K$ has width $0$ and bridge number $0$.

Let $K$ be in Morse position, and let $\hatT$ be a noncritical level torus.  Suppose that $\hatT$ contains an arc $\alpha$ with interior disjoint from $K$ that together with an arc $\beta$ of $K - \hatT$ bounds an embedded disk $\Delta$ with interior disjoint from $K$.  If $\beta$ lies above (resp.\ below) $\hatT$ then we say that $\Delta$ is a {\em high (resp.\ low) disk\/} for $\hatT$ or for the arc $K \cap \Delta$.  

\begin{Lemma}\label{highdisklowdisk}
Assume $K$ is not $0$--bridge.  If $\hatT$ is a noncritical level torus in a thin presentation of $K$, then $\hatT$ cannot have a high disk and a low disk whose boundaries have empty intersection in the complement of $K$.
\end{Lemma}

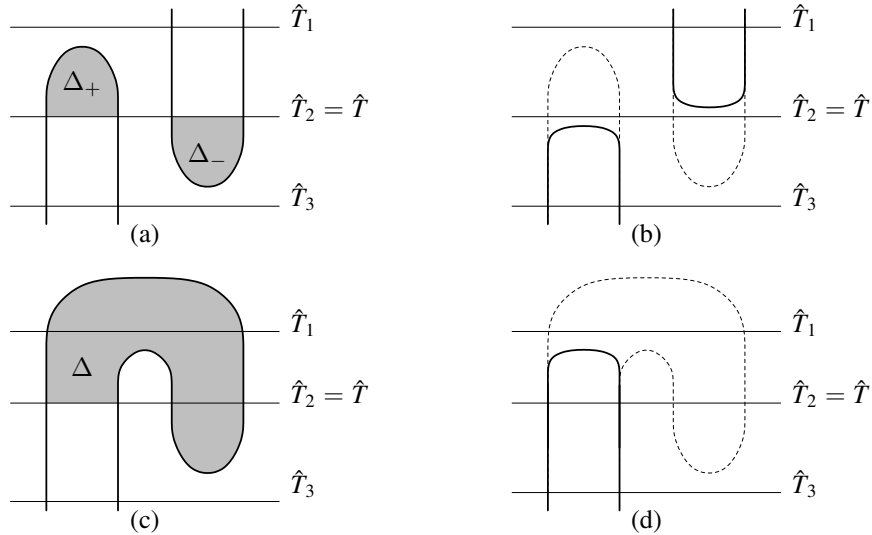
\begin{figure}[ht!]
\centering
\begin{picture}(0,0)%
\includegraphics{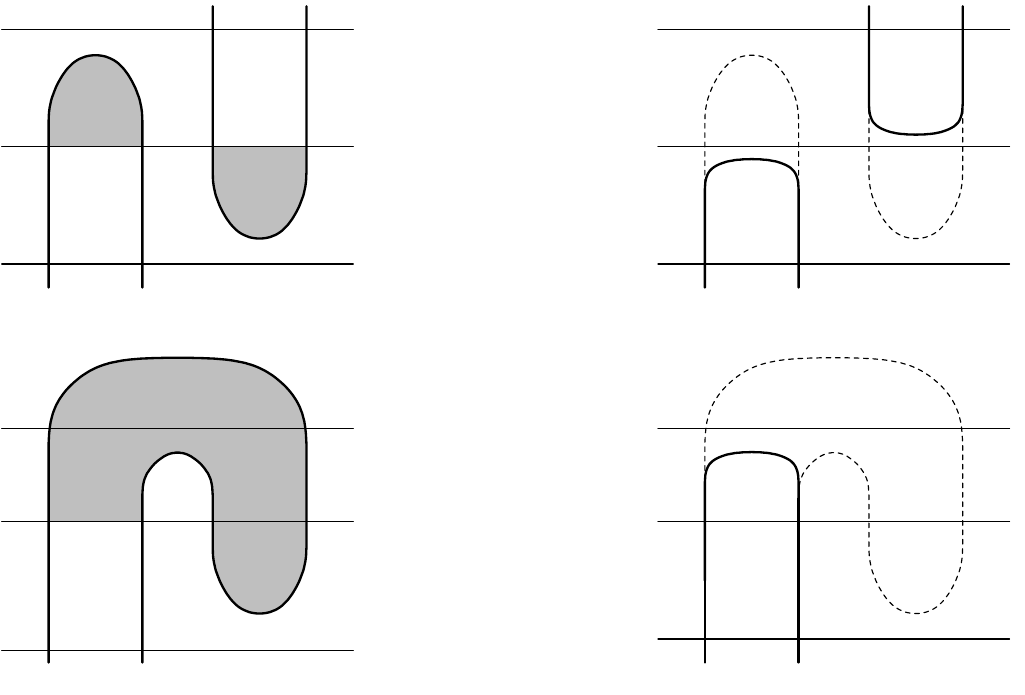}%
\end{picture}%
\setlength{\unitlength}{2960sp}%
\begingroup\makeatletter\ifx\SetFigFont\undefined%
\gdef\SetFigFont#1#2#3#4#5{%
  \reset@font\fontsize{#1}{#2pt}%
  \fontfamily{#3}\fontseries{#4}\fontshape{#5}%
  \selectfont}%
\fi\endgroup%
\begin{picture}(6537,4426)(1189,-5065)
\put(6526,-2611){\makebox(0,0)[b]{\smash{{\SetFigFont{9}{10.8}{\rmdefault}{\mddefault}{\updefault}{\color[rgb]{0,0,0}(b)}%
}}}}
\put(2326,-2611){\makebox(0,0)[b]{\smash{{\SetFigFont{9}{10.8}{\rmdefault}{\mddefault}{\updefault}{\color[rgb]{0,0,0}(a)}%
}}}}
\put(7726,-3361){\makebox(0,0)[lb]{\smash{{\SetFigFont{9}{10.8}{\rmdefault}{\mddefault}{\updefault}{\color[rgb]{0,0,0}$\hatT_1$}%
}}}}
\put(6526,-5011){\makebox(0,0)[b]{\smash{{\SetFigFont{9}{10.8}{\rmdefault}{\mddefault}{\updefault}{\color[rgb]{0,0,0}(d)}%
}}}}
\put(2326,-5011){\makebox(0,0)[b]{\smash{{\SetFigFont{9}{10.8}{\rmdefault}{\mddefault}{\updefault}{\color[rgb]{0,0,0}(c)}%
}}}}
\put(3526,-3361){\makebox(0,0)[lb]{\smash{{\SetFigFont{9}{10.8}{\rmdefault}{\mddefault}{\updefault}{\color[rgb]{0,0,0}$\hatT_1$}%
}}}}
\put(1801,-1336){\makebox(0,0)[b]{\smash{{\SetFigFont{9}{10.8}{\rmdefault}{\mddefault}{\updefault}{\color[rgb]{0,0,0}$\Delta_+$}%
}}}}
\put(2851,-1936){\makebox(0,0)[b]{\smash{{\SetFigFont{9}{10.8}{\rmdefault}{\mddefault}{\updefault}{\color[rgb]{0,0,0}$\Delta_-$}%
}}}}
\put(3526,-811){\makebox(0,0)[lb]{\smash{{\SetFigFont{9}{10.8}{\rmdefault}{\mddefault}{\updefault}{\color[rgb]{0,0,0}$\hatT_1$}%
}}}}
\put(3526,-1561){\makebox(0,0)[lb]{\smash{{\SetFigFont{9}{10.8}{\rmdefault}{\mddefault}{\updefault}{\color[rgb]{0,0,0}$\hatT_2 = \hatT$}%
}}}}
\put(3526,-2311){\makebox(0,0)[lb]{\smash{{\SetFigFont{9}{10.8}{\rmdefault}{\mddefault}{\updefault}{\color[rgb]{0,0,0}$\hatT_3$}%
}}}}
\put(1801,-3736){\makebox(0,0)[b]{\smash{{\SetFigFont{9}{10.8}{\rmdefault}{\mddefault}{\updefault}{\color[rgb]{0,0,0}$\Delta$}%
}}}}
\put(3526,-3961){\makebox(0,0)[lb]{\smash{{\SetFigFont{9}{10.8}{\rmdefault}{\mddefault}{\updefault}{\color[rgb]{0,0,0}$\hatT_2 = \hatT$}%
}}}}
\put(7726,-3961){\makebox(0,0)[lb]{\smash{{\SetFigFont{9}{10.8}{\rmdefault}{\mddefault}{\updefault}{\color[rgb]{0,0,0}$\hatT_2 = \hatT$}%
}}}}
\put(7726,-4711){\makebox(0,0)[lb]{\smash{{\SetFigFont{9}{10.8}{\rmdefault}{\mddefault}{\updefault}{\color[rgb]{0,0,0}$\hatT_3$}%
}}}}
\put(3526,-4711){\makebox(0,0)[lb]{\smash{{\SetFigFont{9}{10.8}{\rmdefault}{\mddefault}{\updefault}{\color[rgb]{0,0,0}$\hatT_3$}%
}}}}
\put(7726,-2311){\makebox(0,0)[lb]{\smash{{\SetFigFont{9}{10.8}{\rmdefault}{\mddefault}{\updefault}{\color[rgb]{0,0,0}$\hatT_3$}%
}}}}
\put(7726,-1561){\makebox(0,0)[lb]{\smash{{\SetFigFont{9}{10.8}{\rmdefault}{\mddefault}{\updefault}{\color[rgb]{0,0,0}$\hatT_2 = \hatT$}%
}}}}
\put(7726,-811){\makebox(0,0)[lb]{\smash{{\SetFigFont{9}{10.8}{\rmdefault}{\mddefault}{\updefault}{\color[rgb]{0,0,0}$\hatT_1$}%
}}}}
\end{picture}%
\caption{(a) A high disk $\Delta_+$ and a low disk $\Delta_-$\qua  (b) The result of a thinning isotopy associated to the disks $\Delta_+$ and $\Delta_-$ \qua (c) A long disk $\Delta$ \qquad (d) The result of a thinning isotopy associated to the long disk $\Delta$}
\label{fig:thinningdisks}
\end{figure}

\begin{proof}
This is standard in the theory of thin position (see eg\ \cite{gabai:fatto3mIII} and \cite{scharlemann:tpittock}).  Assume $\hatT$ admits a high disk $\Delta^+$ and a low disk $\Delta^-$ such that $\bdry \Delta^+ \cap \bdry \Delta^- - K = \emptyset$; see \fullref{fig:thinningdisks}(a).  Let $\alpha^+ = \bdry \Delta^+ \cap \hatT$ and $\alpha^- = \bdry \Delta^- \cap \hatT$.  Then $|\alpha^+ \cap \alpha^-| \leq 2$.  

Use $\Delta^+$ and $\Delta^-$ to isotop $\Delta^+ \cap K$ and $\Delta^- \cap K$ onto $\alpha^+$ and $\alpha^-$ respectively.
If $|\alpha^+ \cap \alpha^-| = 2$ then $K$ is $0$--bridge.  If $|\alpha^+ \cap \alpha^-| \leq 1$ then a further slight isotopy of $K$ will produce a Morse presentation of $K$ with reduced width; see \fullref{fig:thinningdisks}(b).  Both situations contradict our hypothesis.
\end{proof}
\newpage

Let $K$ be in Morse position, and let $\hatT$ be a noncritical level torus.  Suppose that $\hatT$ contains an arc $\alpha$ with interior disjoint from $K$ that together with an arc $\beta$ of $K$ bounds an embedded disk $\Delta$ with interior disjoint from $K$.  If $\Delta$ is not a high or low disk (and hence the interior of $\beta$ intersects $\hatT$) then we say that $\Delta$ is a {\em long disk\/} for $\hatT$ or for the arc $K \cap \Delta$.

\begin{Lemma}\label{longdisk}
Assume $K$ is not $0$--bridge.  If $\hatT$ is a noncritical level torus in a thin presentation of $K$, then $\hatT$ cannot have a long disk.
\end{Lemma}

Note that the existence of a long disk does not necessarily imply the existence of a pair of high and low disks satisfying the hypotheses of \fullref{highdisklowdisk}.

\begin{proof}
Let $\mathbb{T} = \{\hatT_{z_1}, \dots, \hatT_{z_m}\}$ be a collection of level tori used to calculate the width of $K$.

Assume $\Delta$ is a long disk for $\hatT$; see \fullref{fig:thinningdisks}(c).  Furthermore, without loss of generality, assume $N(\alpha) \cap \Delta$ is above $\hatT$.  We may use $\Delta$ to isotop the arc $\beta = \Delta \cap K$ onto $\alpha = \bdry \Delta - \Int \beta \subseteq \hatT$.  Perform a slight isotopy of the interior of this arc upwards so that it has just one critical point, a maximum.   We may then pull this critical point up to the height of the lowest critical point above $\hatT$ of the former arc $\beta$; see \fullref{fig:thinningdisks}(d).  Hence $K$ will again be in a Morse position.

Such an isotopy necessarily reduces the number of critical points of $K$.  Moreover, all the critical points are at the heights of former critical points.  Therefore there exists a proper subset of $\mathbb{T}$ that forms a suitable collection of level tori with which to calculate the width of $K$ after the isotopy.  Furthermore, the isotopy does not increase $|\hatT_{z_i} \cap K|$ for any $i$.  Therefore the isotopy decreases the width of $K$, contradicting that $K$ is in thin position.
\end{proof}

\subsection{Graphs of intersection}\label{sec:graphsofintersection}
Let $K$ be a knot of genus $g$ and order $s$ in the lens space $X$ of order $r$  whose Seifert surfaces have just one boundary component.  Let $S \subseteq X - N(K)$ be a Seifert surface for $K$ of genus $g$.  Note that $S$ is incompressible and $\bdry$--incompressible in the exterior of $K$.   Let $\what{S}$ be the surface $S$ with its boundary (abstractly) capped off by a disk.

Assume that $K$ is in thin position and that $K$ is not $0$--bridge.  
 Let $\hatT$ be a thick level for $K$ with $|\hatT \cap K|=t$, and set $T = \hatT - N(K)$.  By an isotopy of $S$ we may assume that $S$ and $T$ intersect transversely and each component of $\bdry T$ intersects $\bdry S$ exactly $s$ times.  Gabai \cite{gabai:fatto3mIII} shows that we may assume $S$ has been isotoped so that each arc component of $S \cap T$ is essential in $S$ and in $T$ (cf\ Gordon's comment in \cite{gordon:cmids} or \cite[Proposition~2.1]{gordon:tdsokils}).    We may further assume that $S$ has been chosen among all such Seifert surfaces for $K$ to minimize the number of intersections with $T$. 
 
Note that every closed component of $S \cap T$ is nontrivial in $T$.  If a component $\gamma$ of $S \cap T$ were trivial on $T$, then since $S$ is incompressible it must also be trivial on $S$.  Therefore a disk exchange would produce another Seifert surface for $K$ that has fewer intersections with $T$.

The arc components of $S \cap T$ define fat vertexed graphs $G_S$ and $G_T$ in $\what{S}$ and $\hatT$ respectively.  The {\em (fat) vertices\/}, {\em edges\/}, and {\em faces\/} of these graphs are defined as follows.   The fat vertex of $G_S$ is the single disk $\what{S} - \Int S$.  The fat vertices of $G_T$ are the disks of intersection $\hatT \cap \bar{N}(K) = \hatT - \Int T$.  The edges of $G_S$ are the arc components of $S \cap T$ as they lie on $S$, and the edges of $G_T$ are the arc components of $S \cap T$ as they lie on $T$.  The faces of $G_S$ are the connected components of $S$ cut along the edges of $G_S$, ie\ the path metric closure of the complement of the edges of $G_S$ in $S$.  Similarly, the faces of $G_T$ are the connected components of $T$ cut along the edges of $G_T$.  

Observe that simple closed curves of $S \cap T$ are not recorded in the graphs $G_S$ and $G_T$.  For example, the interior of a face of $G_S$ may intersect $T$.  Therefore each component of $S \cut T$ is contained in a face of $G_S$, though often a component of $S \cut T$ actually is a face of $G_S$.  To make the distinction from faces, we refer to the components of $S \cut T$ and $T \cut S$ as {\em regions\/}.

Orient $K$ and sequentially label the intersections of $K \cap \hatT$ from $1$ to $t$.  We may refer to the $i$--th intersection as $K_i$.  This induces a numbering of the components of $\bdry T$ as $1, 2, \dots, t$ in the order in which they appear on $\bdry \bar{N}(K)$, and hence numbers the vertices of $G_T$.  Denote these vertices $U_1, \dots, U_t$.  As there is only one vertex of $G_S$, we leave it unnumbered. 
  Since $\hatT$ is separating and $|\bdry S| = 1$, the parity rule implies that the edges of $G_T$ connect vertices of opposite parity.  An edge of $G_T$ with end points on $U_i$ and $U_j$ is an arc of $S \cap T$ with endpoints on the intersection of $\bdry S$ with the $i$--th and $j$--th components of $\bdry T$ and hence it is also an edge of $G_S$ with end points on the single vertex of $G_S$.  On the boundary of the fat vertex of $G_S$ we label the end points of this arc with $i$ and $j$.  Around the vertex of $G_S$, the labels $1, 2, \dots, t$ appear sequentially and repeat $s$ times.  Since each arc component of $S \cap T$ is essential in $S$ and in $T$, neither graph $G_S$ nor $G_T$ contains a trivial loop.

Denote by $\mathbf{t}$ the set of edge-endpoint labels $\{1, 2, \dots, t\}$ of $G_S$ for which we have the associated $\mathbf{t}$--{\em intervals\/} $(1,2), (2,3), \dots, (t{-}1, t), (t,1)$.  We may index the arc of $K \cut \hatT$ running from $K_i$ to $K_{i+1}$ by the $\mathbf{t}$--interval $(i,i{+}1)$ as $K_{(i,\,i+1)}$.  Similarly $H_{(i,\,i+1)}$ denotes the $1$--handle of $\bar{N}(K) \cut \hatT$ running from vertex $U_i$ to vertex $U_{i+1}$.  Concatenations of consecutive $\mathbf{t}$--intervals such as $(i-1,i+1)$ may be used to index longer arcs of $K$ and longer $1$--handles, eg\ $K_{(i-1,\,i+1)} = K_{(i-1,\,i)} \cup K_{(i,\,i+1)}$ and $H_{(i-1,\,i+1)} = H_{(i-1,\,i)} \cup H_{(i,\,i+1)}$.  

For each label $x$ in $\mathbf{t}$, the subgraph $\smash{G_S^x}$ of $G_S$ is the graph in $\what{S}$ consisting of the single vertex of $G_S$ and every edge of $G_S$ that has an endpoint labeled $x$.  Due to the parity rule, each graph $\smash{G_S^x}$ has exactly $s$ edges.  The faces of $\smash{G_S^x}$ are the connected components of $S$ cut along the edges of $\smash{G_S^x}$.

If $f$ is a face of $G_S$ (resp.\ $\smash{G_S^x}$ for some $x \in$ $\mathbf{t}$), then each component of $\bdry f$ consists of an alternating sequence of {\em edges\/} and {\em corners\/} where the edges are edges of $G_S$ (resp.\ $\smash{G_S^x}$) and the corners are identified with $\mathbf{t}$--intervals $(i, i{+}1)$ (resp.\ concatenated $\mathbf{t}$--intervals $(i,j)$) as they are arcs of $\bdry S$ between the labeled components $i$ and $i+1$ (resp.\ $i$ and $j$) of $\bdry T$.  The edges and corners of a region $R$ of $S \cut T$ contained in a face $f$ of $G_S$ are the edges and corners of $f$ contained in $R$.  If $R$ is a proper subset of $f$, then may have edges and corners not contained in $R$ and $R$ may have a boundary component that is not comprised of edges and corners.  We may similarly define the edges and corners of faces of $G_T$ and regions of $T \cut S$; here corners are arcs of $\bdry T$ between edges of $G_S$.

A disk face of a graph $G_S$, $\smash{G_S^x}$, or $G_T$ with $n$ edges in its boundary is an {\em $n$--gon\/}.  We commonly refer to a $2$--gon as a {\em bigon\/}, a $3$--gon as a {\em trigon\/}, a $4$--gon as a {\em tetragon\/}, and a $5$--gon as a {\em pentagon\/}.  Given an $n$--gon $f$ of $\smash{G_S^x}$, if for every $y \in \mathbf{t}$ and every $n$--gon $f'$ of $G_S^y$ such that $f' \subseteq f \subseteq S$ implies $f' = f$, then $f$ is an {\em innermost $n$--gon\/}.

To each edge of $G_S$ with endpoints labeled $i$ and $j$, we associate its {\em label pair\/} $\{i,j\}$.  A {\em Scharlemann cycle of length $n$\/} ($Sn$ cycle, for short) is a set of $n$ edges of $G_S$ each with label pair $\{i, i{+}1\}$ for some $i$ in $\mathbf{t}$ that bounds an $n$--gon of $G_S$ (with corners all $(i, i{+}1)$).  An {\em extended Scharlemann cycle of length $n$\/} is a set of $n$ edges of $\smash{G_S^x}$ each with label pair $\{x, y\}$ for some $x$ and $y$ in $\mathbf{t}$ that bounds an $n$--gon of $\smash{G_S^x}$ with corners all $(x,y)$ or all $(y,x)$ and contains a Scharlemann cycle of length $n$.  Notice that the Scharlemann cycle of length $n$ contained in the $n$--gon of an extended Scharlemann cycle of length $n$ with label pair $\{x,y\}$ has label pair $\{(x+y-1)/2, (x+y+1)/2 \}$ or $\{(x+y-1)/2 + t/2, (x+y+1)/2 +t/2\}$. If $n$ is $2$ or $3$ we will often abbreviate these terms as {\em (extended) $S2$ cycle\/} and {\em (extended) $S3$ cycle\/}.  By a {\em forked extended $S2$ cycle\/} we mean a set of three edges of some $\smash{G_S^x}$ that bounds a trigon composed of a bigon $B$ bounded by an (extended) $S2$ cycle of either $G_S^{\smash{x+1}}$ or $G_S^{\smash{x-1}}$ together with a bigon and a trigon of $G_S$ adjoined to the edges of $B$; see \fullref{fig:forkedbigona+1}.

\subsection[Outline of proof of Theorem 1.1]{Outline of proof of \fullref{main}}

\begin{proof}
Throughout the subsequent sections, unless noted otherwise, we assume that we have the following:
\begin{itemize}
\item a lens space $X$ of order $r$,
\item a height function $h$ on $X$, 
\item a knot $K$ of order $s$ in thin position with respect to $h$, 
\item a thick level $\hatT$ splitting $X$ into the two solid tori $X^+$ and $X^-$, 
\item the surface $S$ of genus $g \geq 1$ with $|\bdry S|=1$, $\what{S}$, and $T$ as defined above,
\item the graphs $G_S$, $\smash{G_S^x}$ for all $x \in \mathbf{t}$, and $G_T$ and
\item $s \geq 4g-1$.
\end{itemize}

Since $g \geq 1$ and $s \geq 4g-1$, $s \geq 3$.  Since $r \geq s$, $r \geq 3$ as well.  For technical reasons, we defer the case $r = 3$ to \fullref{sec:genus1case} at the end.  Accordingly, we will work under the assumption that $r \geq 4$ until then.  The main body of work encompasses showing that $t \leq 6$.

In \fullref{sec:bigonsandtrigons} we distill the hypothesis of our theorem into the existence of bigons and trigons in the graphs $\smash{G_S^x}$ for each $x \in \mathbf{t}$.  We adapt some fundamental lemmas from Goda and Teragaito \cite{gt:dsokwylsagok} which are then employed to understand what the bigons and trigons may look like.

In \fullref{sec:annuliandtrees} we use the existence of bigons and trigons to construct annuli that weave back and forth through $X^+$ and $X^-$ crossing $\hatT$.  Arcs of $K$ lie on these annuli.  The thinness of $K$ then gives constraints on ``how much'' of $K$ may lie on such an annulus.  With these constraints, the bigons and trigons imply the existence of a second such annulus.

In \fullref{sec:twoschcycles}, to have the requisite bigons and trigons, we reckon with two disjoint annuli containing most of the knot $K$.  The techniques of \fullref{sec:annuliandtrees} are then re-employed to obtain contradictions to the thinness of $K$.  In \fullref{doubleunfurl} we conclude that for $r \geq 4$, we have $t \leq 6$.

In \fullref{sec:tis6} we fix $K$ with $t=6$ and isotop the interior of $S$ to gain a better grasp on the faces of $G_S$.  Then we use Euler characteristic estimates and further thin position arguments to refine our understanding of the faces of $G_S$.  From this \fullref{prop:tnot6} concludes $t \neq 6$ and hence $t \leq 4$.

In \fullref{sec:bridgeposition} we consider multiple thick levels for $K$ in thin position.  Using the result that $K$ may intersect a thick level at most $4$ times, we promote the thin position of $K$ to a bridge position.  Once again using the existence of bigons and trigons, we find a thinning isotopy of $K$.  Thus we conclude in \fullref{thm:t=2} that for $r \geq 4$, we have $t = 2$.  \fullref{lem:t=2} quickly shows that $K$ is at most $1$--bridge.

Finally, in \fullref{sec:genus1case} we treat the case that $r=3$ in which case $s=3$ and $g=1$.  \fullref{thm:r=3} concludes that $K$ is at most $1$--bridge. 
\end{proof}

\section{Bigons and trigons}\label{sec:bigonsandtrigons}

Our proof of \fullref{main} is set in motion by the following lemma.
\begin{lemma}\label{musthavebigons}
For each $x \in \mathbf{t}$, $\smash{G_S^x}$ must have a bigon or trigon face.
\end{lemma}

\begin{proof}
Recall that $\smash{G_S^x}$ has $s$ edges and we are assuming $s \geq 4g-1$.  
\begin{align*}
 \chi(S) = 1- 2g &= -s + \sum_{\mbox{{\scriptsize disk faces of} } \smash{G_S^x}} \chi(\mbox{disk}) + \sum_{\mbox{{\scriptsize nondisk faces of} } \smash{G_S^x}} \chi(\mbox{nondisk}) \\
&\leq -s + \#(\mbox{disk faces})
\end{align*}
Assume $\smash{G_S^x}$ has no bigons or trigons.  So each disk face has at least four edges.  Thus 
\begin{align*} s \geq \frac{1}{2} \cdot 4 \cdot &\#(\mbox{disk faces}). \\
1-2g \leq -s + &\#(\mbox{disk faces}) \leq -s + \frac{1}{2} s,\tag*{\hbox{Hence}}\end{align*}
and so $s/2 \leq 2g-1$ or $s \leq 4g-2$, a contradiction.  
\end{proof}

We thus study how bigons and trigons of $\smash{G_S^x}$ for each $x \in \mathbf{t}$ arise within $G_S$.  Often a bigon or trigon of one graph $\smash{G_S^x}$ will contain a bigon or trigon of another graph $G_S^y$.  We say a face $f$ of $\smash{G_S^x}$ {\em accounts for the label $y$\/} if $f$ contains a bigon or trigon of $G_S^y$.

\subsection[Fundamental lemmas about G sub S]{Fundamental lemmas about $G_S$}
We adapt and build on some useful lemmas about order $2$ and order $3$ Scharlemann cycles and the faces of $G_S$ they bound from the work of Goda and Teragaito in \cite{gt:dsokwylsagok}.  They work with the case that $s=r$, but our generalization of this presents no problem here.  The main difficulty to overcome in our work that is not present in \cite{gt:dsokwylsagok} is the potential presence of simple closed curves $\gamma \in S \cap T$ that are essential in $T$ and yet trivial in both $\hatT$ and $S$.  Goda and Teragaito avoid such occurrences by choosing $\hatT$ among all Heegaard tori.  In our situation, we require $\hatT$ to be a level Heegaard torus.  Nevertheless \fullref{musthavebigons} will permit us to conclude in \fullref{circlesofintersection} that no such curve $\gamma$ exists.

We begin with some lemmas about how edges of Scharlemann cycles and extended Scharlemann cycles may lie on $\hatT$.

Let $\sigma$ be a set of edges of $G_S$.  Let $\Gamma$ be the subgraph of $G_T$ consisting of the edges of $\sigma$ and the vertices of $G_T$ to which the edges are incident.  If $\Gamma$ is contained in a disk in $\hatT$, then we say the edges of $\sigma$ {\em lie in a disk\/} in $\hatT$.  If $\Gamma$ is contained in an annulus in $\hatT$ but do not lie in a disk, then we say the edges of $\sigma$ {\em lie in an essential annulus\/} in $\hatT$.

Let $f$ be a face of $G_S$ or $\smash{G_S^x}$, and let $\sigma$ be the edges of $\bdry f$.  Let $Q$ be a two-sided surface in $X$ with product neighborhood $Q \times [-\epsilon, \epsilon]$ for small $\epsilon > 0$ so that $Q$ is identified with $Q \times \{0\}$.  If $\bdry f \cap Q = \sigma$ and $N(\sigma) \cap f$ is contained in either $Q \times [0, \epsilon]$ or $Q \times [-\epsilon, 0]$ then we say $f$ {\em lies on one side of $Q$\/} even if $\Int f$ does not.  

\begin{lemma}{\rm (cf\ \cite[Lemma 2.3]{gt:dsokwylsagok})}\qua\label{GT:L2.3}
Let $\sigma$ be a Scharlemann cycle in $G_S$ of length $p$ with label pair $\{x,x{+}1\}$ and let $f$ be the face of $G_S$ bounded by $\sigma$.  Suppose that $p \neq r$.  Then $f$ cannot lie on one side of a disk.  In particular, the edges of $\sigma$ cannot lie in a disk in $\hatT$.
\end{lemma}

\begin{proof}
Assume the edges of $\sigma$ lie in a disk $D$ and $f$ lies on one side of $D$.  (If the edges of $\sigma$ lie in disk in $\hatT$, then necessarily $f$ lies to one side of that disk.)  If $\Int f \cap D \neq \emptyset$ then by choosing $D$ smaller, we may assume $\Int f \cap D$ is a collection of simple closed curves.  Let $\xi \in \Int f \cap D$ be an innermost simple closed curve on $f$.  Alter $D$ by a disk exchange with the disk on $f$ that $\xi$ bounds.  In this manner we may produce a disk $D' \subseteq X$ of which $f$ lies to one side (containing the vertices $U_x$ and $U_{x+1}$ and the edges of $\sigma$) such that $\Int f \cap D' = \emptyset$.  Then $N(D' \cup H_{(x,\, x{+}1)} \cup f)$ is a punctured lens space of order $p$.  Since a lens space is irreducible, $X$ is a lens space of order $p$.  This is contrary to the assumption that $p \neq r$.
\end{proof}

\begin{lemma}{\rm (cf\ \cite[Lemma 2.1]{gt:dsokwylsagok})}\qua\label{GT:L2.1}  
Let $\sigma$ be an $Sp$ cycle of $G_S$, $p = 2$ or $3$, with label pair $\{x, x{+}1\}$.  Let $f$ be the face of $G_S$ bounded by $\sigma$.  Then the edges of $\sigma$ lie in an essential annulus $A$ in $\hatT$.  
Furthermore, the core of $A$ does not bound a disk in $X$.  
Indeed, the core of $A$ runs $p$ times in the longitudinal direction of the solid torus on the side of $\hatT$ to which $f$ lies.
\end{lemma}

\begin{proof}  Because we are assuming $r \geq 4$, \fullref{GT:L2.3} implies the edges of $\sigma$ do not lie in a disk in $\hatT$.  The proof of Lemma~2.1 of \cite{gt:dsokwylsagok} then applies to show $\sigma$ must lie in an essential annulus $A$ on $\hatT$.  

 Let $W = X^+$ or $X^-$ be the solid torus on the side of $\hatT$ in which $f$ lies.  If the core of $A$ bounds a disk in $X$, then it must bound a meridional disk of a solid torus on one side of $\hatT$.  By \fullref{GT:L2.3}, the core of $A$ cannot bound a disk in the solid torus $X \cut W$.  As we will show, the core of $A$ cannot bound a disk in $W$ since it runs $p$ times in the longitudinal direction of $W$.

If $\Int f \cap \hatT = \emptyset$, then the remainder of this proof follows from Lemma~2.1 of \cite{gt:dsokwylsagok}.  Because the core of $A$ runs $p$ times in the longitudinal direction of the solid torus $M = \bar{N}(A \cup H_{(x,\, x{+}1)} \cup f) \subseteq W$, the space $W \cut M$ must be a solid torus whose meridional disk crosses $\bdry M \cut A$ once.  Thus the core of $A$ must also run $p$ times in the longitudinal direction on $W$. 

If $\Int f \cap \hatT \neq \emptyset$, then by choosing $A$ smaller we may assume that $\Int f \cap  \bdry A = \emptyset$.  Let $\xi \subseteq \Int f \cap \hatT$ be an innermost simple closed curve on $f$.  Alter $\hatT$ and $A$ by a disk exchange with the disk on $f$ that $\xi$ bounds.  (Note that $\xi$ might be essential in $T$ and yet contained in $A$.)  Continuing in this manner we may produce a Heegaard torus $\hatT'$ (isotopic to $\hatT$ in $X$) and an essential annulus $A'$ on $\hatT'$ such that $\sigma$ lies in the essential annulus $A'$, $f$ lies on the same sides of $\hatT'$ and $\hatT$, and $\Int f \cap \hatT' = \emptyset$.  Let $W'$ be the solid torus on the side of $\hatT'$ containing $f$.  As in \cite{gt:dsokwylsagok}, $M = \bar{N}(A' \cup H_{\smash{(x,\,x{+}1)}} \cup f)$ is a solid torus such that $A'$ runs $p$ times in the longitudinal direction of $M$.  Thus the core of $A'$ and each component of $\bdry A'$ runs $p$ times in the longitudinal direction of $W'$.  Since $\bdry A' = \bdry A$ and $W'$ is isotopic to $W$, it follows that the core of $A$ runs $p$ times in the longitudinal direction of $W$.
\end{proof}

\begin{lemma}\label{lieinanannulus}
Let $\sigma'$ be an extended $Sp$ cycle for $p=2$ or $3$.  Assume $\sigma$ is the $Sp$ cycle contained in the face of $S$ bounded by $\sigma'$.  Then the edges of $\sigma$ and $\sigma'$ each lie in essential annuli $A$ and $A'$ respectively in $\hatT$ so that $A \cap A' = \emptyset$.
\end{lemma}

\begin{proof}
Assume $\sigma$ has label pair $\{x, x{+}1\}$.  Then $\sigma'$ has label pair $\{x-n, x+1+n\}$ (taken mod $t$) for some positive integer $n$.  We proceed by induction. 

By \fullref{GT:L2.1}, the edges of $\sigma$ lie in an essential annulus $A = A_0$ in $\hatT$.  

Let $\sigma_{i-1}$ be the extended $Sp$ cycle with label pair $\{x-(i-1), x+1+(i-1)\}$ contained in the extended $Sp$ cycle $\sigma_i$ with label pair $\{x-i, x+1+i\}$.  Assume the edges of $\sigma_{i-1}$ lie in an essential annulus $A_{i-1}$ in $\hatT$.
 There are $p$ bigons of $G_S$ whose edges give a one to one correspondence between the edges of $\sigma_{i-1}$ and the edges of $\sigma_i$.  Each bigon has the two corners $(x-i, x-i+1)$ and $(x+i, x+i+1)$. 

Of these bigons, take two whose edges in $\sigma_{i-1}$ lie in the essential annulus $A_{i-1}$ (and not in a disk).  Extend their corners radially into $H_{(x-i,\, x-i+1)}$ and $H_{(x+i,\, x+i+1)}$, and join them together to form an annulus $B$ with $\bdry B \subseteq \hatT$.  Since one component of $\bdry B$ is an essential curve contained in $A_{i-1}$, the other is either parallel on $\hatT$ to the core of $A_{i-1}$ or bounds a disk on $\hatT$.  

If one component of $\bdry B$ bounds a disk $D$ on $\hatT$, then the other component bounds the disk $D \cup B$.  Since this other component is parallel to the core of $A_{i-1}$ and hence to the core of $A$, the core of $A$ bounds a disk in $X$ contrary to \fullref{GT:L2.1}. 

 Thus the component of $\bdry B$ not in $A_{i-1}$ is an essential curve on $\hatT$ parallel to the core of $A$.  Hence the edges of $\sigma_i$ lie in an essential annulus $A_i$ in $\hatT$ and $A_{i-1} \cap A_{i} = \emptyset$.   

When $i=n$ we obtain that the edges of $\sigma'$ lie in the annulus $A_n = A'$.  Furthermore, since the label pairs of $\sigma_i$ and $\sigma_j$ for $0\leq i < j \leq n$ are distinct, $A_i \cap A_j = \emptyset$.
\end{proof}

\begin{lemma}\label{lem:twoS2cycleswithsamelabelpair}
If $\sigma_1$ and $\sigma_2$ are two $S2$ cycles with the same label pairs then the edges of $\sigma_1 \cup \sigma_2$ lie in an essential annulus.  Furthermore, within the annulus, the edges of $\sigma_1$ separate the edges of $\sigma_2$.
\end{lemma}

This lemma implies that each edge of $\sigma_1$ bounds a bigon on $\hatT$ with an edge of $\sigma_2$.

\begin{proof}
Assume the edges of $\sigma_1 \cup \sigma_2$ do not lie in an essential annulus.  Nevertheless, by \fullref{GT:L2.3} the edges $\sigma_i$ lie in an essential annulus $A_i$ for each $i=1,2$.  Let $f$ be the bigon of $G_S$ bounded by $\sigma_1$.  Let $g$ be the bigon of $G_S$ bounded by $\sigma_2$.  Let $\{x, x{+}1\}$ be the label pair of $\sigma_1$ and $\sigma_2$.  Note that the corners of $f$ and $g$ are all on $H_{(x,\, x{+}1)}$.

If $\Int(f \cup g) \cap \hatT \neq \emptyset$, then alter $\hatT$ by disk exchanges to obtain the Heegaard torus $\hatT'$ so that $\Int(f \cup g) \cap \hatT' = \emptyset$ while $\bdry(f \cup g) \cap \hatT' = \bdry(f \cup g) \cap \hatT$.  Similarly, we obtain the essential annuli $A_i'$ on $\hatT'$ in which the edges of $\sigma_i$ lie.  Let $V$ be the solid torus on the same side of $\hatT'$ as $f$ and $g$.  

Consider the solid torus $V_f = V \cut (f \cup H_{(x,\, x{+}1)})$.  Note that the core of the annulus $\hatT' \cut A_1'$ is a longitudinal curve on $\bdry V_f$.  Since $\sigma_2$ does not lie in $A_1'$, $\sigma_2$ must have an edge that is not parallel to an edge of $\sigma_1$.  Thus both corners of $g$ lie on the same rectangle of $\bdry H_{(x,\, x{+}1)} \cut (U_x \cup U_{x+1} \cup f)$ which is contained in $\bdry V_f$.  Hence the other edge of $\sigma_2$ is not parallel to an edge of $\sigma_2$.  Then the simple closed curve $\bdry g$ on $\bdry V_f$ must intersect the core of $\hatT' \cut A'$ twice.  Therefore $\bdry g$ cannot be nullhomologous in $H_1(V_f)$ contradicting that it bounds a disk in $V_f$.

Since the corners of $g$ cannot be on the same rectangle of $\bdry H_{(x,\, x{+}1)} \cap \bdry V_f$, the edges of $g$ must be incident to $U_x$ and $U_{x+1}$ on opposite sides of the edges of $\sigma_1$.  Since the edges of $\sigma_1$ and $\sigma_2$ lie in an essential annulus, the edges of $\sigma_1$ separate the edges of $\sigma_2$ within this annulus.
\end{proof}

\begin{lemma}{\rm (cf\ \cite[Lemma~2.4]{gt:dsokwylsagok})}\qua\label{oppositesides}
Let $\sigma_1$ and $\sigma_2$ be $S2$ or $S3$ cycles of $G_S$ with disjoint label pairs, and let $f_1$ and $f_2$ be the faces of $G_S$ bounded by $\sigma_1$ and $\sigma_2$ respectively.  Then the faces $f_1$ and $f_2$ lie on opposite sides of $\hatT$.
\end{lemma}

\begin{proof}
By \fullref{lieinanannulus}, the edges of $\sigma_i$ lie in an annulus $A_i$ for $i = 1, 2$.  Since the label pairs of $\sigma_1$ and $\sigma_2$ are disjoint, $A_1 \cap A_2 = \emptyset$.

We may assume $A_1$ and $A_2$ have been chosen so that $(A_1 \cup A_2) \cap (\Int f_1 \cup \Int f_2)$ is only a collection of simple closed curves.
By the minimality assumption on $|S \cap T|$ and \fullref{GT:L2.1}, $\hatT \cap (\Int f_1 \cup \Int f_2)$ may contain only simple closed curves that are essential on $T$ and yet trivial on $\hatT$.  Let $\xi \in \hatT \cap (\Int f_1 \cup \Int f_2)$ be an innermost simple closed curve on $f_1 \cup f_2$.  Alter $\hatT$ (and both $A_1$ and $A_2$) by a disk exchange with the disk on $f_1 \cup f_2$ bounded by $\xi$.  In this manner we may produce a Heegaard torus $\hatT'$ (isotopic to $\hatT$ in $X$) with annuli $A_1'$ and $A_2$' such that $\sigma_i \subseteq A_i'$ for $i=1,2$, $f_i$ lies on the same side of $\hatT'$ as $\hatT$ for $i=1,2$, and $\hatT' \cap (\Int f_1 \cup \Int f_2) = \emptyset$.

 Let $\{x_i, x_i{+}1\}$ be the label pair of $\sigma_i$.  Construct the solid torus $$V_i = \bar{N}(A_i' \cup H_{(x_i,\, x_i+1)} \cup f_i).$$  The meridian of $V_i$ intersects the core of the annulus $A_i$ algebraically $2$ or $3$ times depending on the order of the cycle $\sigma_i$.  If $f_1$ and $f_2$ lie on the same side of $\hatT$, then $f_1$ and $f_2$ lie on the same side of $\hatT'$.  Then the solid tori $V_1$ and $V_2$ are both contained in a single solid torus on the same side of $\hatT'$.  Let $A_3$ be an annulus of $\hatT' \cut (A_1' \cup A_2')$.  The manifold $V_3 = V_1 \cup \bar{N}(A_3) \cup V_2$ has toroidal boundary which intersects $\hatT'$ in the annulus $A_1' \cup A_2' \cup A_3$.  Thus $V_3$ must be a solid torus.  Due to neither meridian of $V_1$ or $V_2$ intersecting a component of $\bdry A_3$ algebraically once, $V_3$ cannot be a solid torus.  Therefore $f_1$ and $f_2$ cannot lie on the same side of $\hatT$.  
\end{proof}

\begin{figure}[ht!]
\centering
\begin{picture}(0,0)%
\includegraphics{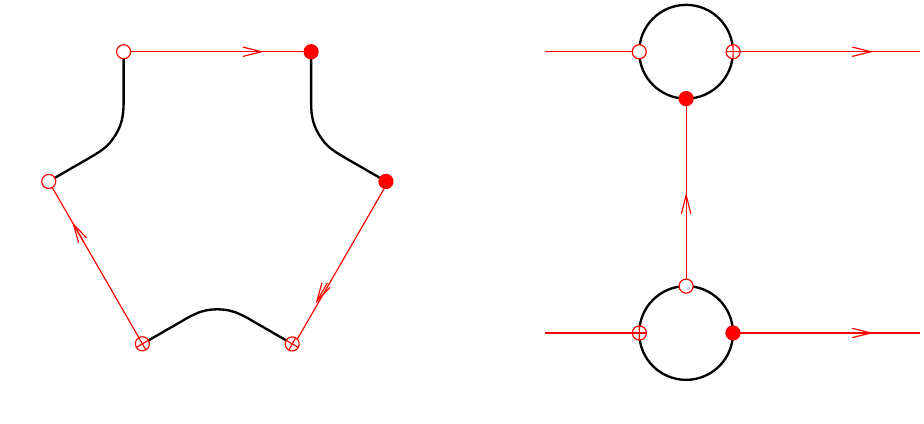}%
\end{picture}%
\setlength{\unitlength}{2960sp}%
\begingroup\makeatletter\ifx\SetFigFont\undefined%
\gdef\SetFigFont#1#2#3#4#5{%
  \reset@font\fontsize{#1}{#2pt}%
  \fontfamily{#3}\fontseries{#4}\fontshape{#5}%
  \selectfont}%
\fi\endgroup%
\begin{picture}(5904,2693)(2509,-3640)
\put(3376,-3286){\makebox(0,0)[rb]{\smash{{\SetFigFont{9}{10.8}{\rmdefault}{\mddefault}{\updefault}{\color[rgb]{0,0,0}$x$}%
}}}}
\put(4501,-1186){\makebox(0,0)[rb]{\smash{{\SetFigFont{9}{10.8}{\rmdefault}{\mddefault}{\updefault}{\color[rgb]{0,0,0}$x$}%
}}}}
\put(3226,-1261){\makebox(0,0)[rb]{\smash{{\SetFigFont{9}{10.8}{\rmdefault}{\mddefault}{\updefault}{\color[rgb]{0,0,0}$x+1$}%
}}}}
\put(7201,-3586){\makebox(0,0)[b]{\smash{{\SetFigFont{9}{10.8}{\rmdefault}{\mddefault}{\updefault}{\color[rgb]{0,0,0}(b)}%
}}}}
\put(5026,-2311){\makebox(0,0)[lb]{\smash{{\SetFigFont{9}{10.8}{\rmdefault}{\mddefault}{\updefault}{\color[rgb]{0,0,0}$x+1$}%
}}}}
\put(3901,-3586){\makebox(0,0)[b]{\smash{{\SetFigFont{9}{10.8}{\rmdefault}{\mddefault}{\updefault}{\color[rgb]{0,0,0}(a)}%
}}}}
\put(4426,-3286){\makebox(0,0)[lb]{\smash{{\SetFigFont{9}{10.8}{\rmdefault}{\mddefault}{\updefault}{\color[rgb]{0,0,0}$x+1$}%
}}}}
\put(2786,-2086){\makebox(0,0)[rb]{\smash{{\SetFigFont{9}{10.8}{\rmdefault}{\mddefault}{\updefault}{\color[rgb]{0,0,0}$x$}%
}}}}
\put(6901,-3136){\makebox(0,0)[b]{\smash{{\SetFigFont{9}{10.8}{\rmdefault}{\mddefault}{\updefault}{\color[rgb]{0,0,0}$x+1$}%
}}}}
\put(6901,-1336){\makebox(0,0)[b]{\smash{{\SetFigFont{9}{10.8}{\rmdefault}{\mddefault}{\updefault}{\color[rgb]{0,0,0}$x$}%
}}}}
\end{picture}%
\caption{(a) The trigon $F$\qua (b) The edges of $F$ on $\hatT'$}
\label{fig:funnytrigon}
\end{figure}
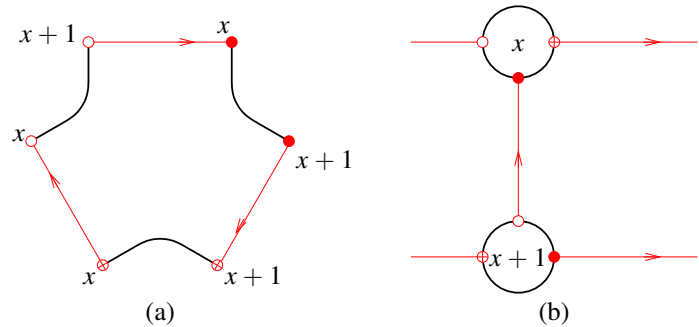

\begin{lemma} \label{lem:funnytrigon}
Let $F$ be a trigon with corners $(x, x{+}1)$ and edges as in \fullref{fig:funnytrigon}(a) such that its edges lie in an essential annulus $A$ on $\hatT$ as in \fullref{fig:funnytrigon}(b).  Let $V$ be the solid torus on the same side of $\hatT$ as $F$.  Then the core of $A$ is a meridional curve for $V$.
\end{lemma}

\begin{proof}
Note that $F$ itself is not a face of $G_S$ or any of its subgraphs since it has edges with both endpoints on the same vertex.  Nevertheless we will continue to use the language of edges and corners when talking about $F$.  

Consider $M = \bar{N}(A \cup H_{(x,\, x{+}1)} \cup F)$ formed by attaching the $2$--handle $\bar{N}(F)$ to the genus $2$ handlebody $\bar{N}(A \cup H_{(x,\, x{+}1)})$.  Let $D_z$ be a meridional disk of $\bar{N}(A)$ whose boundary intersects the edges of $F$ only twice and misses the vertices.  Let $D_H$ be the cocore of the $1$--handle $H_{(x,\, x{+}1)}$.  The fundamental group of $\bar{N}(A \cup H_{(x,\, x{+}1)})$ is then generated by curves $\zeta$ and $\eta$ dual to $D_z$ and $D_H$ respectively.  Thus we have the presentation 
\begin{align*}
\pi_1(M) &= \langle \zeta, \eta | \eta \zeta \eta^{-1} \zeta \eta = 1 \rangle \\
    &= \langle \zeta, \eta | (\zeta \eta^2)^2 = \eta^3 \rangle \\
    &= \langle a, b | a^2 = b^3 \rangle
\end{align*}    
which is the trefoil group.  Therefore $\pi_1(M)$ is not $\Z$, and $M$ is not a solid torus.

We may consider $M$ as contained in $V$ such that $\bdry M \cap \bdry V = A$.  Then $A' = \bdry M \cut A$ is properly embedded annulus in $V$.  Assume the core of $A$ is not a meridional curve for $V$.  Then the components of $\bdry A' = \bdry A$ are not meridional curves for $V$.  Thus the two components of $V \cut A'$ must be solid tori.  This is contrary to $M$ not being a solid torus. 
\end{proof}

\subsection{Bigons}

\begin{lemma}\label{bigonS2}  
If $\sigma = \{e', e''\} \subseteq \smash{G_S^x}$ bounds a bigon $f$ of $\smash{G_S^x}$, then $\sigma$ is either an $S2$ cycle or an extended $S2$ cycle.
\end{lemma}

\begin{proof}
The edges of $G_S$ on $f$ must be mutually parallel.
By relabeling if necessary, we may assume $x=1$ and the edges of $G_S$ on $f$ are $\{e'=e_1, e_2, \dots, e_n = e'' \}$ labeled successively so that, along a chosen corner of $f$, $e_i$ has label $i$ taken mod $t$.  We claim that $n \leq t$.

If $n > t+1$, then the edge $e_{t+1}$ is contained in $\smash{G_S^1}$ and yet is not on the boundary of $f$.  This contradicts that $f$ is a bigon of $\smash{G_S^1}$.

If $n=t+1$, then the label $1$ for the endpoints of each $e'$ and $e''$ occur on the same corner of $f$.  Since there are $t+1$ edges, the label $1$ must occur for the end point of some edge $e_i$ on the other corner of $f$.  Hence there is an edge of $\smash{G_S^1}$ contained in the interior of $f$ contradicting that $f$ is a bigon of $\smash{G_S^1}$.  

Since $1 < n \leq t$, $\sigma$ has label pair $\{1,n\}$.  The label $1$ of the edges $e'$ and $e''$ occur on opposite corners of $f$.  If $n=2$, then $\sigma$ is an $S2$ cycle.  If $n > 2$, then $\sigma$ is an extended $S2$ cycle.
\end{proof}

\subsection{Trigons}
The structure of trigons of $\smash{G_S^x}$ takes a bit more work to determine. 
We will first classify innermost trigons of $\smash{G_S^x}$ and then determine how trigons of another graph $G_S^y$ may contain them.  By relabeling, we may assume a given innermost trigon is a trigon of $\smash{G_S^1}$.

\begin{figure}[ht!]
\centering
\begin{picture}(0,0)%
\includegraphics{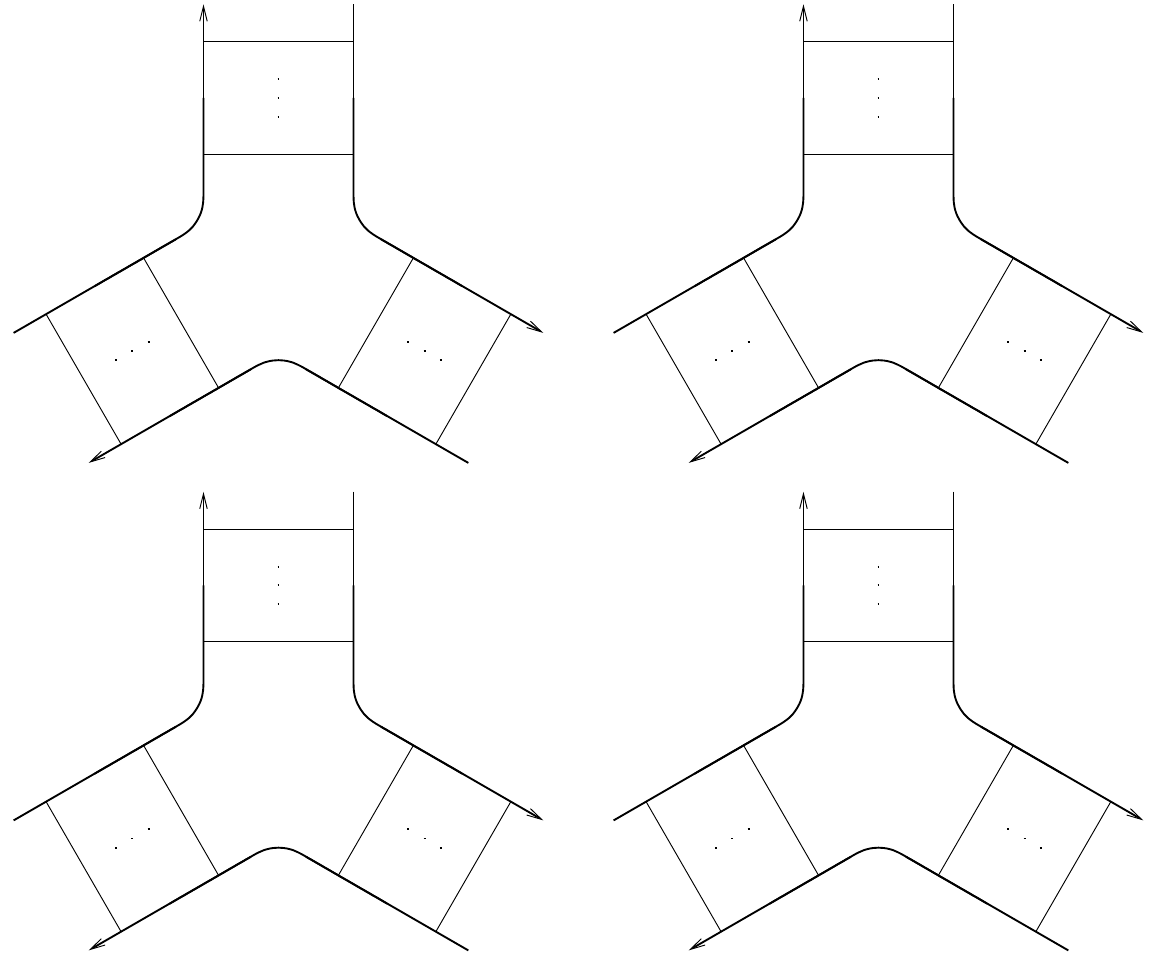}%
\end{picture}%
\setlength{\unitlength}{2368sp}%
\begingroup\makeatletter\ifx\SetFigFont\undefined%
\gdef\SetFigFont#1#2#3#4#5{%
  \reset@font\fontsize{#1}{#2pt}%
  \fontfamily{#3}\fontseries{#4}\fontshape{#5}%
  \selectfont}%
\fi\endgroup%
\begin{picture}(9253,7655)(3174,-8894)
\put(3601,-5761){\makebox(0,0)[lb]{\smash{{\SetFigFont{7}{8.4}{\rmdefault}{\mddefault}{\updefault}{\color[rgb]{0,0,0}I$'$}%
}}}}
\put(4651,-2461){\makebox(0,0)[rb]{\smash{{\SetFigFont{7}{8.4}{\rmdefault}{\mddefault}{\updefault}{\color[rgb]{0,0,0}$a+1$}%
}}}}
\put(4276,-3211){\makebox(0,0)[rb]{\smash{{\SetFigFont{7}{8.4}{\rmdefault}{\mddefault}{\updefault}{\color[rgb]{0,0,0}$a$}%
}}}}
\put(3451,-3661){\makebox(0,0)[rb]{\smash{{\SetFigFont{7}{8.4}{\rmdefault}{\mddefault}{\updefault}{\color[rgb]{0,0,0}$1$}%
}}}}
\put(4651,-1561){\makebox(0,0)[rb]{\smash{{\SetFigFont{7}{8.4}{\rmdefault}{\mddefault}{\updefault}{\color[rgb]{0,0,0}$x$}%
}}}}
\put(5851,-4486){\makebox(0,0)[rb]{\smash{{\SetFigFont{7}{8.4}{\rmdefault}{\mddefault}{\updefault}{\color[rgb]{0,0,0}$c$}%
}}}}
\put(6601,-4936){\makebox(0,0)[rb]{\smash{{\SetFigFont{7}{8.4}{\rmdefault}{\mddefault}{\updefault}{\color[rgb]{0,0,0}$1$}%
}}}}
\put(6151,-1561){\makebox(0,0)[lb]{\smash{{\SetFigFont{7}{8.4}{\rmdefault}{\mddefault}{\updefault}{\color[rgb]{0,0,0}$1$}%
}}}}
\put(6151,-2461){\makebox(0,0)[lb]{\smash{{\SetFigFont{7}{8.4}{\rmdefault}{\mddefault}{\updefault}{\color[rgb]{0,0,0}$b$}%
}}}}
\put(6526,-3211){\makebox(0,0)[lb]{\smash{{\SetFigFont{7}{8.4}{\rmdefault}{\mddefault}{\updefault}{\color[rgb]{0,0,0}$b+1$}%
}}}}
\put(7351,-3661){\makebox(0,0)[lb]{\smash{{\SetFigFont{7}{8.4}{\rmdefault}{\mddefault}{\updefault}{\color[rgb]{0,0,0}$y$}%
}}}}
\put(4951,-4486){\makebox(0,0)[lb]{\smash{{\SetFigFont{7}{8.4}{\rmdefault}{\mddefault}{\updefault}{\color[rgb]{0,0,0}$c+1$}%
}}}}
\put(4201,-4936){\makebox(0,0)[lb]{\smash{{\SetFigFont{7}{8.4}{\rmdefault}{\mddefault}{\updefault}{\color[rgb]{0,0,0}$z$}%
}}}}
\put(3601,-1861){\makebox(0,0)[lb]{\smash{{\SetFigFont{7}{8.4}{\rmdefault}{\mddefault}{\updefault}{\color[rgb]{0,0,0}I}%
}}}}
\put(9451,-2461){\makebox(0,0)[rb]{\smash{{\SetFigFont{7}{8.4}{\rmdefault}{\mddefault}{\updefault}{\color[rgb]{0,0,0}$a+1$}%
}}}}
\put(9076,-3211){\makebox(0,0)[rb]{\smash{{\SetFigFont{7}{8.4}{\rmdefault}{\mddefault}{\updefault}{\color[rgb]{0,0,0}$a$}%
}}}}
\put(8251,-3661){\makebox(0,0)[rb]{\smash{{\SetFigFont{7}{8.4}{\rmdefault}{\mddefault}{\updefault}{\color[rgb]{0,0,0}$1$}%
}}}}
\put(10651,-4486){\makebox(0,0)[rb]{\smash{{\SetFigFont{7}{8.4}{\rmdefault}{\mddefault}{\updefault}{\color[rgb]{0,0,0}$c$}%
}}}}
\put(11401,-4936){\makebox(0,0)[rb]{\smash{{\SetFigFont{7}{8.4}{\rmdefault}{\mddefault}{\updefault}{\color[rgb]{0,0,0}$y$}%
}}}}
\put(10951,-1561){\makebox(0,0)[lb]{\smash{{\SetFigFont{7}{8.4}{\rmdefault}{\mddefault}{\updefault}{\color[rgb]{0,0,0}$1$}%
}}}}
\put(10951,-2461){\makebox(0,0)[lb]{\smash{{\SetFigFont{7}{8.4}{\rmdefault}{\mddefault}{\updefault}{\color[rgb]{0,0,0}$b$}%
}}}}
\put(11326,-3211){\makebox(0,0)[lb]{\smash{{\SetFigFont{7}{8.4}{\rmdefault}{\mddefault}{\updefault}{\color[rgb]{0,0,0}$b+1$}%
}}}}
\put(12151,-3661){\makebox(0,0)[lb]{\smash{{\SetFigFont{7}{8.4}{\rmdefault}{\mddefault}{\updefault}{\color[rgb]{0,0,0}$1$}%
}}}}
\put(9751,-4486){\makebox(0,0)[lb]{\smash{{\SetFigFont{7}{8.4}{\rmdefault}{\mddefault}{\updefault}{\color[rgb]{0,0,0}$c+1$}%
}}}}
\put(9001,-4936){\makebox(0,0)[lb]{\smash{{\SetFigFont{7}{8.4}{\rmdefault}{\mddefault}{\updefault}{\color[rgb]{0,0,0}$z$}%
}}}}
\put(8401,-1861){\makebox(0,0)[lb]{\smash{{\SetFigFont{7}{8.4}{\rmdefault}{\mddefault}{\updefault}{\color[rgb]{0,0,0}II}%
}}}}
\put(9451,-1561){\makebox(0,0)[rb]{\smash{{\SetFigFont{7}{8.4}{\rmdefault}{\mddefault}{\updefault}{\color[rgb]{0,0,0}$x$}%
}}}}
\put(9451,-6361){\makebox(0,0)[rb]{\smash{{\SetFigFont{7}{8.4}{\rmdefault}{\mddefault}{\updefault}{\color[rgb]{0,0,0}$a+1$}%
}}}}
\put(9076,-7111){\makebox(0,0)[rb]{\smash{{\SetFigFont{7}{8.4}{\rmdefault}{\mddefault}{\updefault}{\color[rgb]{0,0,0}$a$}%
}}}}
\put(8251,-7561){\makebox(0,0)[rb]{\smash{{\SetFigFont{7}{8.4}{\rmdefault}{\mddefault}{\updefault}{\color[rgb]{0,0,0}$z$}%
}}}}
\put(9451,-5461){\makebox(0,0)[rb]{\smash{{\SetFigFont{7}{8.4}{\rmdefault}{\mddefault}{\updefault}{\color[rgb]{0,0,0}$x$}%
}}}}
\put(10651,-8386){\makebox(0,0)[rb]{\smash{{\SetFigFont{7}{8.4}{\rmdefault}{\mddefault}{\updefault}{\color[rgb]{0,0,0}$c$}%
}}}}
\put(11401,-8836){\makebox(0,0)[rb]{\smash{{\SetFigFont{7}{8.4}{\rmdefault}{\mddefault}{\updefault}{\color[rgb]{0,0,0}$y$}%
}}}}
\put(10951,-5461){\makebox(0,0)[lb]{\smash{{\SetFigFont{7}{8.4}{\rmdefault}{\mddefault}{\updefault}{\color[rgb]{0,0,0}$1$}%
}}}}
\put(10951,-6361){\makebox(0,0)[lb]{\smash{{\SetFigFont{7}{8.4}{\rmdefault}{\mddefault}{\updefault}{\color[rgb]{0,0,0}$b$}%
}}}}
\put(11326,-7111){\makebox(0,0)[lb]{\smash{{\SetFigFont{7}{8.4}{\rmdefault}{\mddefault}{\updefault}{\color[rgb]{0,0,0}$b+1$}%
}}}}
\put(12151,-7561){\makebox(0,0)[lb]{\smash{{\SetFigFont{7}{8.4}{\rmdefault}{\mddefault}{\updefault}{\color[rgb]{0,0,0}$1$}%
}}}}
\put(9751,-8386){\makebox(0,0)[lb]{\smash{{\SetFigFont{7}{8.4}{\rmdefault}{\mddefault}{\updefault}{\color[rgb]{0,0,0}$c+1$}%
}}}}
\put(9001,-8836){\makebox(0,0)[lb]{\smash{{\SetFigFont{7}{8.4}{\rmdefault}{\mddefault}{\updefault}{\color[rgb]{0,0,0}$1$}%
}}}}
\put(8401,-5761){\makebox(0,0)[lb]{\smash{{\SetFigFont{7}{8.4}{\rmdefault}{\mddefault}{\updefault}{\color[rgb]{0,0,0}II$'$}%
}}}}
\put(4651,-6361){\makebox(0,0)[rb]{\smash{{\SetFigFont{7}{8.4}{\rmdefault}{\mddefault}{\updefault}{\color[rgb]{0,0,0}$a+1$}%
}}}}
\put(4276,-7111){\makebox(0,0)[rb]{\smash{{\SetFigFont{7}{8.4}{\rmdefault}{\mddefault}{\updefault}{\color[rgb]{0,0,0}$a$}%
}}}}
\put(3451,-7561){\makebox(0,0)[rb]{\smash{{\SetFigFont{7}{8.4}{\rmdefault}{\mddefault}{\updefault}{\color[rgb]{0,0,0}$z$}%
}}}}
\put(4651,-5461){\makebox(0,0)[rb]{\smash{{\SetFigFont{7}{8.4}{\rmdefault}{\mddefault}{\updefault}{\color[rgb]{0,0,0}$1$}%
}}}}
\put(5851,-8386){\makebox(0,0)[rb]{\smash{{\SetFigFont{7}{8.4}{\rmdefault}{\mddefault}{\updefault}{\color[rgb]{0,0,0}$c$}%
}}}}
\put(6601,-8836){\makebox(0,0)[rb]{\smash{{\SetFigFont{7}{8.4}{\rmdefault}{\mddefault}{\updefault}{\color[rgb]{0,0,0}$y$}%
}}}}
\put(6151,-5461){\makebox(0,0)[lb]{\smash{{\SetFigFont{7}{8.4}{\rmdefault}{\mddefault}{\updefault}{\color[rgb]{0,0,0}$x$}%
}}}}
\put(6151,-6361){\makebox(0,0)[lb]{\smash{{\SetFigFont{7}{8.4}{\rmdefault}{\mddefault}{\updefault}{\color[rgb]{0,0,0}$b$}%
}}}}
\put(6526,-7111){\makebox(0,0)[lb]{\smash{{\SetFigFont{7}{8.4}{\rmdefault}{\mddefault}{\updefault}{\color[rgb]{0,0,0}$b+1$}%
}}}}
\put(7351,-7561){\makebox(0,0)[lb]{\smash{{\SetFigFont{7}{8.4}{\rmdefault}{\mddefault}{\updefault}{\color[rgb]{0,0,0}$1$}%
}}}}
\put(4951,-8386){\makebox(0,0)[lb]{\smash{{\SetFigFont{7}{8.4}{\rmdefault}{\mddefault}{\updefault}{\color[rgb]{0,0,0}$c+1$}%
}}}}
\put(4201,-8836){\makebox(0,0)[lb]{\smash{{\SetFigFont{7}{8.4}{\rmdefault}{\mddefault}{\updefault}{\color[rgb]{0,0,0}$1$}%
}}}}
\end{picture}%
\caption{The four possible types of trigons}
\label{fig:trigontypes}
\end{figure}
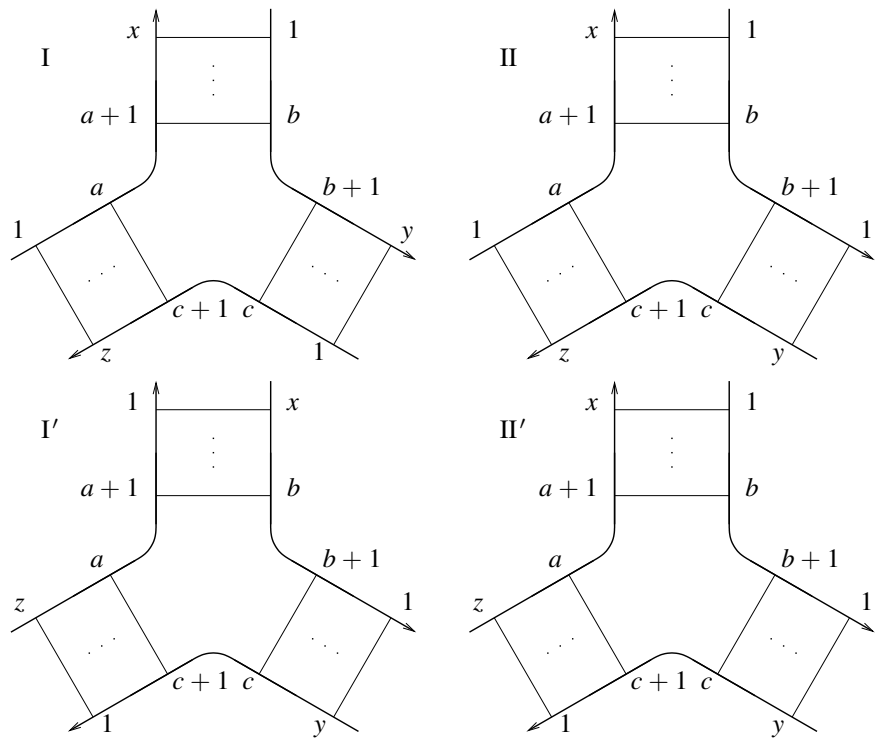

If $f$ is a trigon of $\smash{G_S^1}$, then $G_S$ on $f$ appears as one of the four types shown in \fullref{fig:trigontypes}.  The label $1$ appears only where shown since otherwise $f$ would not be a face of $\smash{G_S^1}$.  Types I$'$ and II$'$ may be obtained from types I and II respectively by suitable changes of orientations and labeling.  Therefore we may further assume an innermost trigon is a trigon of $\smash{G_S^1}$ of type I or type II.

\begin{prop} \label{innermosttrigons}
An innermost trigon of $\smash{G_S^1}$ is bounded by either an $S3$ cycle or a forked extended $S2$ cycle.
\end{prop}

\begin{proof}
To highlight the overall structure of this proof we relegate two of the subcases to \fullref{lem:ruleoutcaseIa<b} and \fullref{lem:ruleoutcaseIa>b}.

Let $f$ be an innermost trigon of $\smash{G_S^1}$.  Note that we are not concerned with whether or not $f \cap T$ contains any circle components.

{\bf Case I}\qua  The trigon $f$ is of type I.

To be innermost, we must have $a$, $b$, or $c$ be $1$.  Otherwise there would be a trigon of $G_S^2$ contained in $f$.  Say $c = 1$.  The trigon appears as in \fullref{fig:trigonIa}.

\begin{figure}[ht!]
\centering
\begin{picture}(0,0)%
\includegraphics{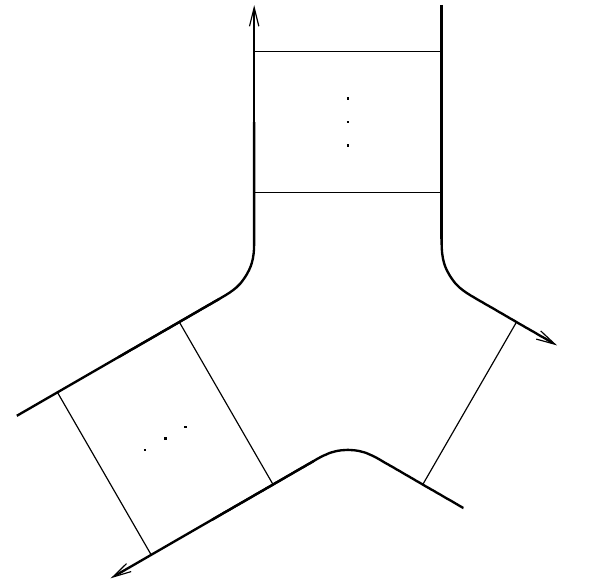}%
\end{picture}%
\setlength{\unitlength}{2960sp}%
\begingroup\makeatletter\ifx\SetFigFont\undefined%
\gdef\SetFigFont#1#2#3#4#5{%
  \reset@font\fontsize{#1}{#2pt}%
  \fontfamily{#3}\fontseries{#4}\fontshape{#5}%
  \selectfont}%
\fi\endgroup%
\begin{picture}(3828,3746)(3174,-4985)
\put(4201,-4936){\makebox(0,0)[lb]{\smash{{\SetFigFont{9}{10.8}{\rmdefault}{\mddefault}{\updefault}{\color[rgb]{0,0,0}$a+1$}%
}}}}
\put(4651,-2461){\makebox(0,0)[rb]{\smash{{\SetFigFont{9}{10.8}{\rmdefault}{\mddefault}{\updefault}{\color[rgb]{0,0,0}$a+1$}%
}}}}
\put(4276,-3211){\makebox(0,0)[rb]{\smash{{\SetFigFont{9}{10.8}{\rmdefault}{\mddefault}{\updefault}{\color[rgb]{0,0,0}$a$}%
}}}}
\put(3451,-3661){\makebox(0,0)[rb]{\smash{{\SetFigFont{9}{10.8}{\rmdefault}{\mddefault}{\updefault}{\color[rgb]{0,0,0}$1$}%
}}}}
\put(4651,-1561){\makebox(0,0)[rb]{\smash{{\SetFigFont{9}{10.8}{\rmdefault}{\mddefault}{\updefault}{\color[rgb]{0,0,0}$a+b$}%
}}}}
\put(5851,-4486){\makebox(0,0)[rb]{\smash{{\SetFigFont{9}{10.8}{\rmdefault}{\mddefault}{\updefault}{\color[rgb]{0,0,0}$1$}%
}}}}
\put(6151,-1561){\makebox(0,0)[lb]{\smash{{\SetFigFont{9}{10.8}{\rmdefault}{\mddefault}{\updefault}{\color[rgb]{0,0,0}$1$}%
}}}}
\put(6151,-2461){\makebox(0,0)[lb]{\smash{{\SetFigFont{9}{10.8}{\rmdefault}{\mddefault}{\updefault}{\color[rgb]{0,0,0}$b$}%
}}}}
\put(6526,-3211){\makebox(0,0)[lb]{\smash{{\SetFigFont{9}{10.8}{\rmdefault}{\mddefault}{\updefault}{\color[rgb]{0,0,0}$b+1$}%
}}}}
\put(4951,-4486){\makebox(0,0)[lb]{\smash{{\SetFigFont{9}{10.8}{\rmdefault}{\mddefault}{\updefault}{\color[rgb]{0,0,0}$2$}%
}}}}
\end{picture}%
\caption{The trigon of type I with $c=1$}
\label{fig:trigonIa}
\end{figure}
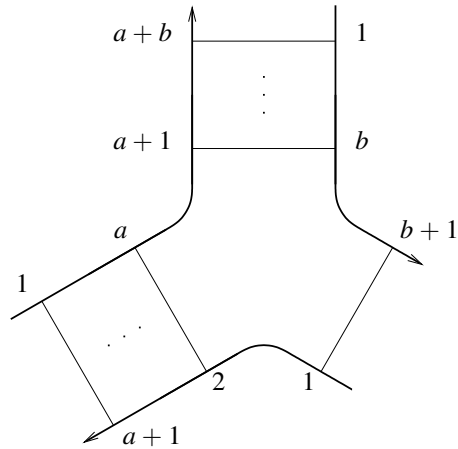

Since $f$ is a face of $\smash{G_S^1}$, we have $a+b \leq t$.  Notice that $a$ and $b$ are both odd.  One may readily check that both the cases $b=c=1\neq a$ and $a=b=c=1$ are innermost trigons corresponding to a forked extended $S2$ cycle and an $S3$ cycle respectively.  Thus assume $a \neq 1 \neq b$.

{\bf Subcase $a<b$}\qua  
This is ruled out by \fullref{lem:ruleoutcaseIa<b}.

{\bf Subcase $a=b$}\qua
In this situation there is a trigon of $G_S^{a+1}$ within $f$ contradicting that $f$ is innermost.

{\bf Subcase $a>b$}\qua
This is ruled out by \fullref{lem:ruleoutcaseIa>b}.

{\bf Case II}\qua The trigon $f$ is of type II.

For $f$ to be a face of $\smash{G_S^1}$, we must have $c+b >t$, $c+a \leq t$, and $a+b \leq t$.   These conditions imply $a < b$ and $a < c$.   The trigon appears as in \fullref{fig:trigonIIa}.  We have three main cases.

\begin{figure}[ht!]
\centering
\begin{picture}(0,0)%
\includegraphics{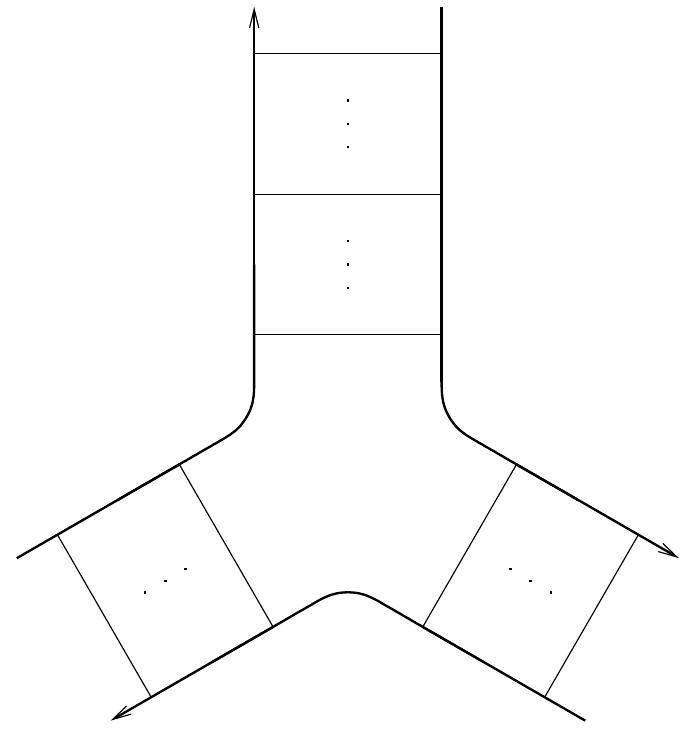}%
\end{picture}%
\setlength{\unitlength}{2960sp}%
\begingroup\makeatletter\ifx\SetFigFont\undefined%
\gdef\SetFigFont#1#2#3#4#5{%
  \reset@font\fontsize{#1}{#2pt}%
  \fontfamily{#3}\fontseries{#4}\fontshape{#5}%
  \selectfont}%
\fi\endgroup%
\begin{picture}(4453,4646)(3174,-4985)
\put(6526,-4936){\makebox(0,0)[rb]{\smash{{\SetFigFont{9}{10.8}{\rmdefault}{\mddefault}{\updefault}{\color[rgb]{0,0,0}$c+b-t$}%
}}}}
\put(6151,-736){\makebox(0,0)[lb]{\smash{{\SetFigFont{9}{10.8}{\rmdefault}{\mddefault}{\updefault}{\color[rgb]{0,0,0}$1$}%
}}}}
\put(5026,-4486){\makebox(0,0)[lb]{\smash{{\SetFigFont{9}{10.8}{\rmdefault}{\mddefault}{\updefault}{\color[rgb]{0,0,0}$c+1$}%
}}}}
\put(4276,-4936){\makebox(0,0)[lb]{\smash{{\SetFigFont{9}{10.8}{\rmdefault}{\mddefault}{\updefault}{\color[rgb]{0,0,0}$c+a$}%
}}}}
\put(6526,-3211){\makebox(0,0)[lb]{\smash{{\SetFigFont{9}{10.8}{\rmdefault}{\mddefault}{\updefault}{\color[rgb]{0,0,0}$b+1$}%
}}}}
\put(5776,-4486){\makebox(0,0)[rb]{\smash{{\SetFigFont{9}{10.8}{\rmdefault}{\mddefault}{\updefault}{\color[rgb]{0,0,0}$c$}%
}}}}
\put(3451,-3661){\makebox(0,0)[rb]{\smash{{\SetFigFont{9}{10.8}{\rmdefault}{\mddefault}{\updefault}{\color[rgb]{0,0,0}$1$}%
}}}}
\put(4651,-2536){\makebox(0,0)[rb]{\smash{{\SetFigFont{9}{10.8}{\rmdefault}{\mddefault}{\updefault}{\color[rgb]{0,0,0}$a+1$}%
}}}}
\put(4651,-1636){\makebox(0,0)[rb]{\smash{{\SetFigFont{9}{10.8}{\rmdefault}{\mddefault}{\updefault}{\color[rgb]{0,0,0}$b$}%
}}}}
\put(4651,-736){\makebox(0,0)[rb]{\smash{{\SetFigFont{9}{10.8}{\rmdefault}{\mddefault}{\updefault}{\color[rgb]{0,0,0}$a+b$}%
}}}}
\put(7351,-3661){\makebox(0,0)[lb]{\smash{{\SetFigFont{9}{10.8}{\rmdefault}{\mddefault}{\updefault}{\color[rgb]{0,0,0}$1$}%
}}}}
\put(4276,-3211){\makebox(0,0)[rb]{\smash{{\SetFigFont{9}{10.8}{\rmdefault}{\mddefault}{\updefault}{\color[rgb]{0,0,0}$a$}%
}}}}
\put(6151,-2536){\makebox(0,0)[lb]{\smash{{\SetFigFont{9}{10.8}{\rmdefault}{\mddefault}{\updefault}{\color[rgb]{0,0,0}$b$}%
}}}}
\put(6151,-1636){\makebox(0,0)[lb]{\smash{{\SetFigFont{9}{10.8}{\rmdefault}{\mddefault}{\updefault}{\color[rgb]{0,0,0}$a+1$}%
}}}}
\end{picture}%
\caption{The trigon of type II}
\label{fig:trigonIIa}
\end{figure}
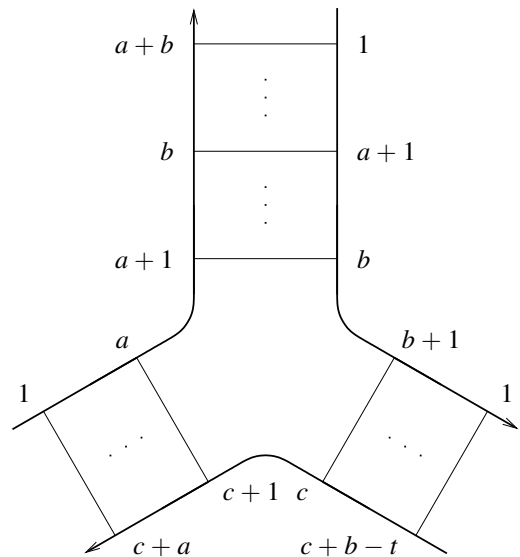

{\bf Subcase $c>b$}\qua
With some relabeling, part of the trigon appears as in Case I, $a<b$.  This gives a contradiction to \fullref{lem:ruleoutcaseIa<b}. 

{\bf Subcase $c=b$}\qua
In this situation there is a trigon of $G_S^{b+1}$ inside $f$.

{\bf Subcase $c<b$}\qua
Here, $a<c<b$.  With some relabeling, the arguments of Case I, $a>b$ apply.

Thus there are no innermost trigons of type II.

\begin{figure}[ht!]
\centering
\begin{picture}(0,0)%
\includegraphics{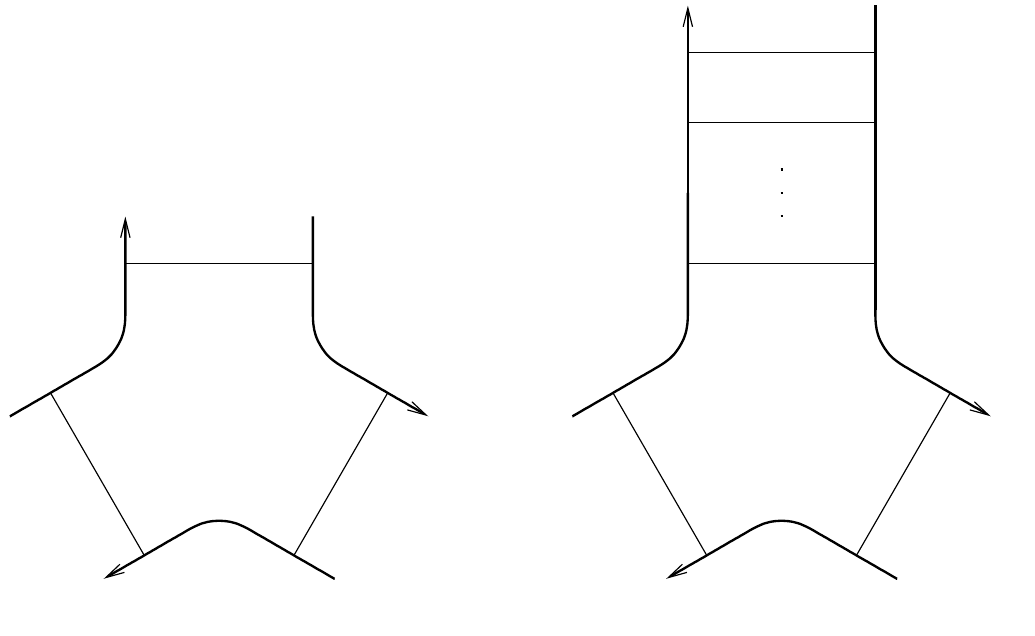}%
\end{picture}%
\setlength{\unitlength}{2960sp}%
\begingroup\makeatletter\ifx\SetFigFont\undefined%
\gdef\SetFigFont#1#2#3#4#5{%
  \reset@font\fontsize{#1}{#2pt}%
  \fontfamily{#3}\fontseries{#4}\fontshape{#5}%
  \selectfont}%
\fi\endgroup%
\begin{picture}(6603,4055)(3999,-4844)
\put(9001,-4786){\makebox(0,0)[b]{\smash{{\SetFigFont{9}{10.8}{\rmdefault}{\mddefault}{\updefault}{\color[rgb]{0,0,0}$a \neq 1$}%
}}}}
\put(7876,-3211){\makebox(0,0)[rb]{\smash{{\SetFigFont{9}{10.8}{\rmdefault}{\mddefault}{\updefault}{\color[rgb]{0,0,0}$1$}%
}}}}
\put(9451,-4486){\makebox(0,0)[rb]{\smash{{\SetFigFont{9}{10.8}{\rmdefault}{\mddefault}{\updefault}{\color[rgb]{0,0,0}$a$}%
}}}}
\put(10126,-3211){\makebox(0,0)[lb]{\smash{{\SetFigFont{9}{10.8}{\rmdefault}{\mddefault}{\updefault}{\color[rgb]{0,0,0}$a+1$}%
}}}}
\put(8551,-4486){\makebox(0,0)[lb]{\smash{{\SetFigFont{9}{10.8}{\rmdefault}{\mddefault}{\updefault}{\color[rgb]{0,0,0}$a+1$}%
}}}}
\put(8251,-2536){\makebox(0,0)[rb]{\smash{{\SetFigFont{9}{10.8}{\rmdefault}{\mddefault}{\updefault}{\color[rgb]{0,0,0}$2$}%
}}}}
\put(8251,-1636){\makebox(0,0)[rb]{\smash{{\SetFigFont{9}{10.8}{\rmdefault}{\mddefault}{\updefault}{\color[rgb]{0,0,0}$a$}%
}}}}
\put(9751,-1186){\makebox(0,0)[lb]{\smash{{\SetFigFont{9}{10.8}{\rmdefault}{\mddefault}{\updefault}{\color[rgb]{0,0,0}$1$}%
}}}}
\put(9751,-2536){\makebox(0,0)[lb]{\smash{{\SetFigFont{9}{10.8}{\rmdefault}{\mddefault}{\updefault}{\color[rgb]{0,0,0}$a$}%
}}}}
\put(8251,-1186){\makebox(0,0)[rb]{\smash{{\SetFigFont{9}{10.8}{\rmdefault}{\mddefault}{\updefault}{\color[rgb]{0,0,0}$a+1$}%
}}}}
\put(9751,-1636){\makebox(0,0)[lb]{\smash{{\SetFigFont{9}{10.8}{\rmdefault}{\mddefault}{\updefault}{\color[rgb]{0,0,0}$2$}%
}}}}
\put(4651,-2461){\makebox(0,0)[rb]{\smash{{\SetFigFont{9}{10.8}{\rmdefault}{\mddefault}{\updefault}{\color[rgb]{0,0,0}$2$}%
}}}}
\put(4276,-3211){\makebox(0,0)[rb]{\smash{{\SetFigFont{9}{10.8}{\rmdefault}{\mddefault}{\updefault}{\color[rgb]{0,0,0}$1$}%
}}}}
\put(5851,-4486){\makebox(0,0)[rb]{\smash{{\SetFigFont{9}{10.8}{\rmdefault}{\mddefault}{\updefault}{\color[rgb]{0,0,0}$1$}%
}}}}
\put(6151,-2461){\makebox(0,0)[lb]{\smash{{\SetFigFont{9}{10.8}{\rmdefault}{\mddefault}{\updefault}{\color[rgb]{0,0,0}$1$}%
}}}}
\put(6526,-3211){\makebox(0,0)[lb]{\smash{{\SetFigFont{9}{10.8}{\rmdefault}{\mddefault}{\updefault}{\color[rgb]{0,0,0}$2$}%
}}}}
\put(4951,-4486){\makebox(0,0)[lb]{\smash{{\SetFigFont{9}{10.8}{\rmdefault}{\mddefault}{\updefault}{\color[rgb]{0,0,0}$2$}%
}}}}
\end{picture}%
\caption{The two forms of innermost trigons of $\smash{G_S^1}$}
\label{fig:trigonIe}
\end{figure}
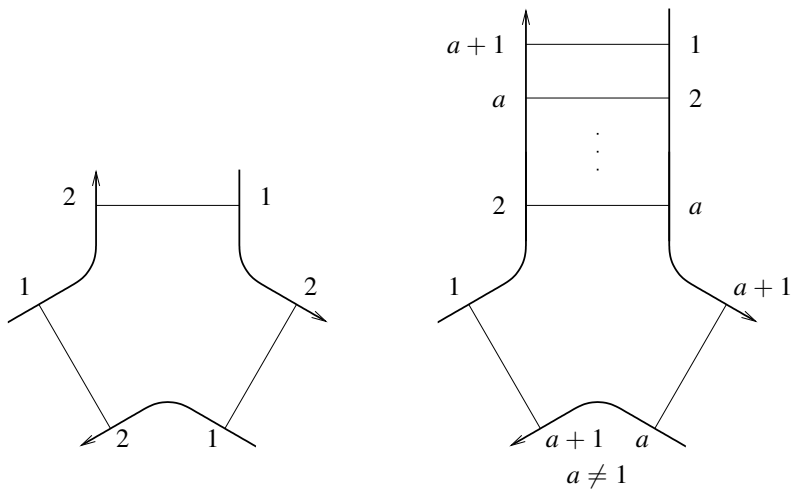

Hence, up to relabeling and changes in orientation, any innermost trigon of $\smash{G_S^1}$ appears as in \fullref{fig:trigonIe}.  Such trigons are bounded by either an $S3$ cycle or a forked extended $S2$ cycle.
\end{proof}

\begin{remark} Note that with \fullref{innermosttrigons} in hand, it follows from the definitions of $S3$ cycles and forked extended $S2$ cycles that all innermost trigons have one or two sets of label pairs for their edges.  None have three distinct label pairs.  \end{remark}

\begin{lemma}\label{lem:ruleoutcaseIa<b}
There cannot exist an innermost trigon of $\smash{G_S^1}$ of type I with $c=1$ and $a < b$.
\end{lemma}

\begin{proof}
Assume a trigon $f$ of $\smash{G_S^1}$ of type I with $c=1$ and $a < b$ does exist.

 On $f$, $G_S$ has a bigon $\kreis{A}$ with edges $e_1$ and $e_2$ having label pairs $\{2,a\}$ and $\{1,a{+}1\}$ respectively and corners $(1,2)$ and $(a,a{+}1)$, another bigon $\kreis{B}$ with edges $e_3$ and $e_4$ having label pairs $\{a, b{+}1\}$ and $\{a{+}1, b\}$ respectively and corners $(a,a{+}1)$ and $(b,b{+}1)$, and a trigon $\kreis{C}$ with edges $e_5$, $e_6$, and $e_7$ having label pairs $\{1, b{+}1\}$, $\{2,a\}$, and $\{a{+}1,b\}$ respectively and corners $(1,2)$, $(a,a{+}1)$, $(b,b{+}1)$.  See \fullref{fig:trigonalessthanb}(a).  Note that $\kreis{A}$, $\kreis{B}$, and $\kreis{C}$ all lie on the same side of $\hatT$.

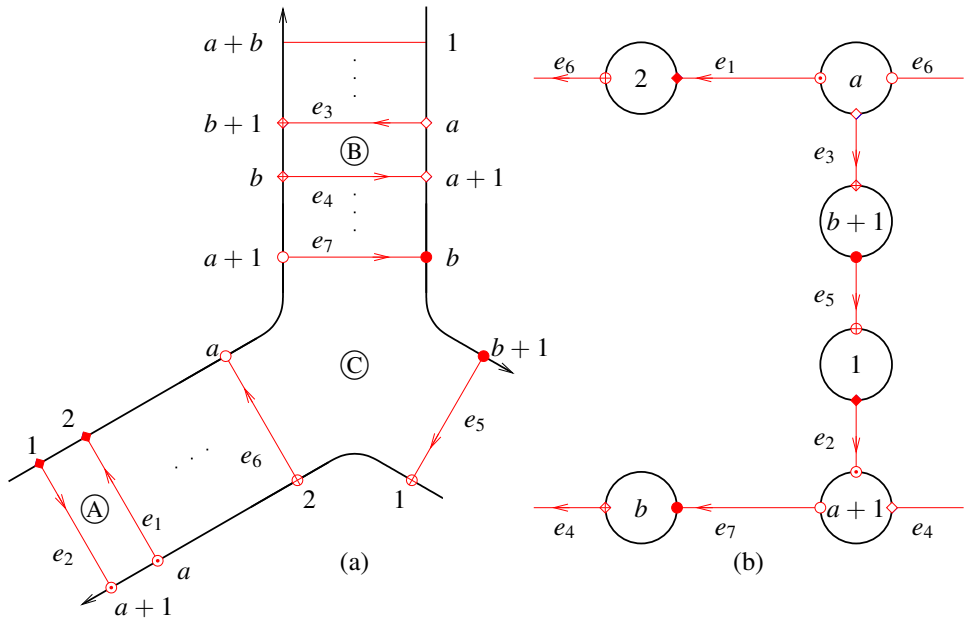
\begin{figure}[ht!]
\centering
\begin{picture}(0,0)%
\includegraphics{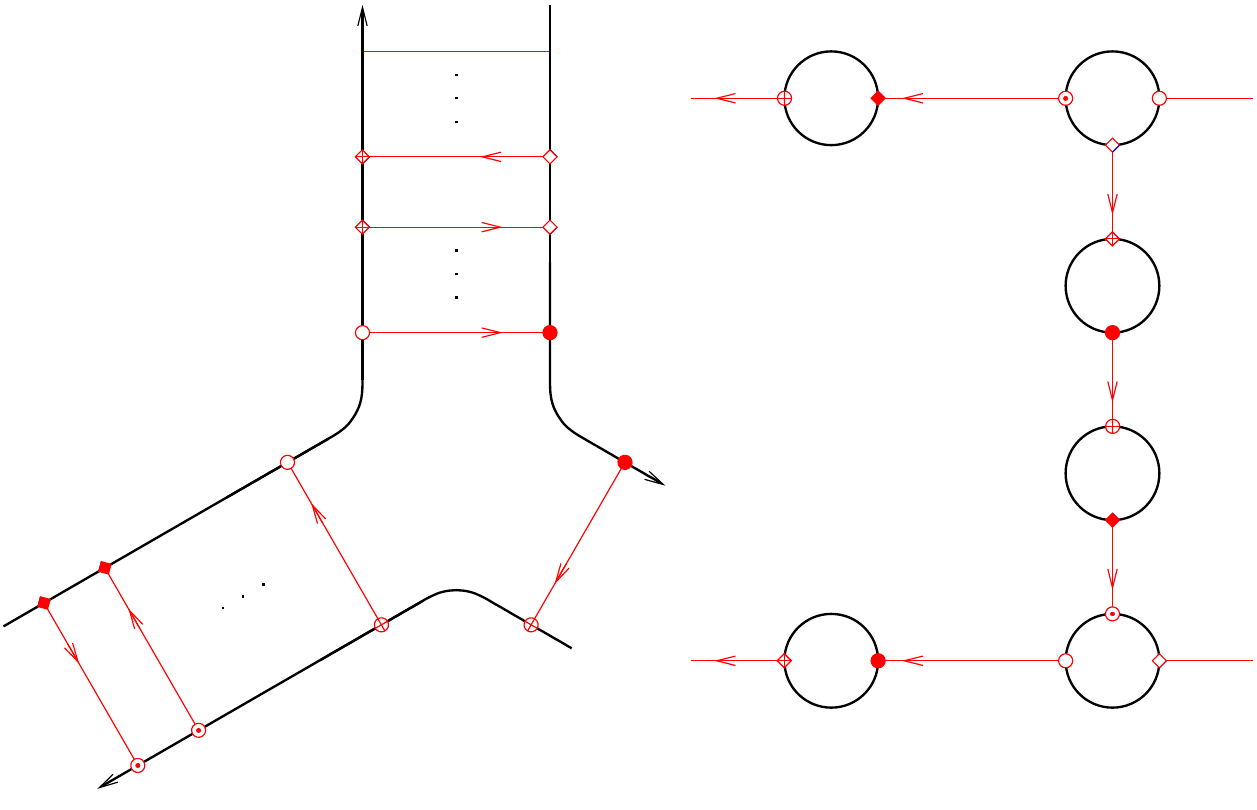}%
\end{picture}%
\setlength{\unitlength}{2960sp}%
\begingroup\makeatletter\ifx\SetFigFont\undefined%
\gdef\SetFigFont#1#2#3#4#5{%
  \reset@font\fontsize{#1}{#2pt}%
  \fontfamily{#3}\fontseries{#4}\fontshape{#5}%
  \selectfont}%
\fi\endgroup%
\begin{picture}(8032,5171)(2481,-5510)
\put(7801,-1036){\makebox(0,0)[b]{\smash{{\SetFigFont{9}{10.8}{\rmdefault}{\mddefault}{\updefault}{\color[rgb]{0,0,0}$2$}%
}}}}
\put(3376,-5461){\makebox(0,0)[lb]{\smash{{\SetFigFont{9}{10.8}{\rmdefault}{\mddefault}{\updefault}{\color[rgb]{0,0,0}$a+1$}%
}}}}
\put(2776,-4111){\makebox(0,0)[rb]{\smash{{\SetFigFont{9}{10.8}{\rmdefault}{\mddefault}{\updefault}{\color[rgb]{0,0,0}$1$}%
}}}}
\put(4951,-4561){\makebox(0,0)[lb]{\smash{{\SetFigFont{9}{10.8}{\rmdefault}{\mddefault}{\updefault}{\color[rgb]{0,0,0}$2$}%
}}}}
\put(3901,-5161){\makebox(0,0)[lb]{\smash{{\SetFigFont{9}{10.8}{\rmdefault}{\mddefault}{\updefault}{\color[rgb]{0,0,0}$a$}%
}}}}
\put(6526,-3286){\makebox(0,0)[lb]{\smash{{\SetFigFont{9}{10.8}{\rmdefault}{\mddefault}{\updefault}{\color[rgb]{0,0,0}$b+1$}%
}}}}
\put(6151,-2536){\makebox(0,0)[lb]{\smash{{\SetFigFont{9}{10.8}{\rmdefault}{\mddefault}{\updefault}{\color[rgb]{0,0,0}$b$}%
}}}}
\put(6151,-1861){\makebox(0,0)[lb]{\smash{{\SetFigFont{9}{10.8}{\rmdefault}{\mddefault}{\updefault}{\color[rgb]{0,0,0}$a+1$}%
}}}}
\put(6151,-1411){\makebox(0,0)[lb]{\smash{{\SetFigFont{9}{10.8}{\rmdefault}{\mddefault}{\updefault}{\color[rgb]{0,0,0}$a$}%
}}}}
\put(6151,-736){\makebox(0,0)[lb]{\smash{{\SetFigFont{9}{10.8}{\rmdefault}{\mddefault}{\updefault}{\color[rgb]{0,0,0}$1$}%
}}}}
\put(3076,-3886){\makebox(0,0)[rb]{\smash{{\SetFigFont{9}{10.8}{\rmdefault}{\mddefault}{\updefault}{\color[rgb]{0,0,0}$2$}%
}}}}
\put(4286,-3286){\makebox(0,0)[rb]{\smash{{\SetFigFont{9}{10.8}{\rmdefault}{\mddefault}{\updefault}{\color[rgb]{0,0,0}$a$}%
}}}}
\put(4651,-2536){\makebox(0,0)[rb]{\smash{{\SetFigFont{9}{10.8}{\rmdefault}{\mddefault}{\updefault}{\color[rgb]{0,0,0}$a+1$}%
}}}}
\put(4651,-1861){\makebox(0,0)[rb]{\smash{{\SetFigFont{9}{10.8}{\rmdefault}{\mddefault}{\updefault}{\color[rgb]{0,0,0}$b$}%
}}}}
\put(4651,-1411){\makebox(0,0)[rb]{\smash{{\SetFigFont{9}{10.8}{\rmdefault}{\mddefault}{\updefault}{\color[rgb]{0,0,0}$b+1$}%
}}}}
\put(4651,-736){\makebox(0,0)[rb]{\smash{{\SetFigFont{9}{10.8}{\rmdefault}{\mddefault}{\updefault}{\color[rgb]{0,0,0}$a+b$}%
}}}}
\put(5851,-4561){\makebox(0,0)[rb]{\smash{{\SetFigFont{9}{10.8}{\rmdefault}{\mddefault}{\updefault}{\color[rgb]{0,0,0}$1$}%
}}}}
\put(5026,-2011){\makebox(0,0)[lb]{\smash{{\SetFigFont{9}{10.8}{\rmdefault}{\mddefault}{\updefault}{\color[rgb]{0,0,0}$e_4$}%
}}}}
\put(5026,-1261){\makebox(0,0)[lb]{\smash{{\SetFigFont{9}{10.8}{\rmdefault}{\mddefault}{\updefault}{\color[rgb]{0,0,0}$e_3$}%
}}}}
\put(5026,-2386){\makebox(0,0)[lb]{\smash{{\SetFigFont{9}{10.8}{\rmdefault}{\mddefault}{\updefault}{\color[rgb]{0,0,0}$e_7$}%
}}}}
\put(6301,-3886){\makebox(0,0)[lb]{\smash{{\SetFigFont{9}{10.8}{\rmdefault}{\mddefault}{\updefault}{\color[rgb]{0,0,0}$e_5$}%
}}}}
\put(5401,-1636){\makebox(0,0)[b]{\smash{{\SetFigFont{9}{10.8}{\rmdefault}{\mddefault}{\updefault}{\color[rgb]{0,0,0}$\bigkreis{B}$}%
}}}}
\put(5401,-3436){\makebox(0,0)[b]{\smash{{\SetFigFont{9}{10.8}{\rmdefault}{\mddefault}{\updefault}{\color[rgb]{0,0,0}$\bigkreis{C}$}%
}}}}
\put(4651,-4186){\makebox(0,0)[rb]{\smash{{\SetFigFont{9}{10.8}{\rmdefault}{\mddefault}{\updefault}{\color[rgb]{0,0,0}$e_6$}%
}}}}
\put(3226,-4636){\makebox(0,0)[b]{\smash{{\SetFigFont{9}{10.8}{\rmdefault}{\mddefault}{\updefault}{\color[rgb]{0,0,0}$\bigkreis{A}$}%
}}}}
\put(3076,-5011){\makebox(0,0)[rb]{\smash{{\SetFigFont{9}{10.8}{\rmdefault}{\mddefault}{\updefault}{\color[rgb]{0,0,0}$e_2$}%
}}}}
\put(3601,-4711){\makebox(0,0)[lb]{\smash{{\SetFigFont{9}{10.8}{\rmdefault}{\mddefault}{\updefault}{\color[rgb]{0,0,0}$e_1$}%
}}}}
\put(5401,-5086){\makebox(0,0)[b]{\smash{{\SetFigFont{9}{10.8}{\rmdefault}{\mddefault}{\updefault}{\color[rgb]{0,0,0}(a)}%
}}}}
\put(9601,-1036){\makebox(0,0)[b]{\smash{{\SetFigFont{9}{10.8}{\rmdefault}{\mddefault}{\updefault}{\color[rgb]{0,0,0}$a$}%
}}}}
\put(9601,-2236){\makebox(0,0)[b]{\smash{{\SetFigFont{9}{10.8}{\rmdefault}{\mddefault}{\updefault}{\color[rgb]{0,0,0}$b+1$}%
}}}}
\put(9601,-3436){\makebox(0,0)[b]{\smash{{\SetFigFont{9}{10.8}{\rmdefault}{\mddefault}{\updefault}{\color[rgb]{0,0,0}$1$}%
}}}}
\put(9601,-4636){\makebox(0,0)[b]{\smash{{\SetFigFont{9}{10.8}{\rmdefault}{\mddefault}{\updefault}{\color[rgb]{0,0,0}$a+1$}%
}}}}
\put(7801,-4636){\makebox(0,0)[b]{\smash{{\SetFigFont{9}{10.8}{\rmdefault}{\mddefault}{\updefault}{\color[rgb]{0,0,0}$b$}%
}}}}
\put(9451,-4036){\makebox(0,0)[rb]{\smash{{\SetFigFont{9}{10.8}{\rmdefault}{\mddefault}{\updefault}{\color[rgb]{0,0,0}$e_2$}%
}}}}
\put(7051,-886){\makebox(0,0)[lb]{\smash{{\SetFigFont{9}{10.8}{\rmdefault}{\mddefault}{\updefault}{\color[rgb]{0,0,0}$e_6$}%
}}}}
\put(8401,-886){\makebox(0,0)[lb]{\smash{{\SetFigFont{9}{10.8}{\rmdefault}{\mddefault}{\updefault}{\color[rgb]{0,0,0}$e_1$}%
}}}}
\put(7051,-4786){\makebox(0,0)[lb]{\smash{{\SetFigFont{9}{10.8}{\rmdefault}{\mddefault}{\updefault}{\color[rgb]{0,0,0}$e_4$}%
}}}}
\put(8401,-4786){\makebox(0,0)[lb]{\smash{{\SetFigFont{9}{10.8}{\rmdefault}{\mddefault}{\updefault}{\color[rgb]{0,0,0}$e_7$}%
}}}}
\put(10051,-4786){\makebox(0,0)[lb]{\smash{{\SetFigFont{9}{10.8}{\rmdefault}{\mddefault}{\updefault}{\color[rgb]{0,0,0}$e_4$}%
}}}}
\put(10051,-886){\makebox(0,0)[lb]{\smash{{\SetFigFont{9}{10.8}{\rmdefault}{\mddefault}{\updefault}{\color[rgb]{0,0,0}$e_6$}%
}}}}
\put(9451,-1636){\makebox(0,0)[rb]{\smash{{\SetFigFont{9}{10.8}{\rmdefault}{\mddefault}{\updefault}{\color[rgb]{0,0,0}$e_3$}%
}}}}
\put(9451,-2836){\makebox(0,0)[rb]{\smash{{\SetFigFont{9}{10.8}{\rmdefault}{\mddefault}{\updefault}{\color[rgb]{0,0,0}$e_5$}%
}}}}
\put(8701,-5086){\makebox(0,0)[b]{\smash{{\SetFigFont{9}{10.8}{\rmdefault}{\mddefault}{\updefault}{\color[rgb]{0,0,0}(b)}%
}}}}
\end{picture}%
\caption{(a) The trigon of type I with $c=1$ and $a<b$\qua (b) The seven edges on $\hatT$}
\label{fig:trigonalessthanb}
\end{figure}

There are two extended $S2$ cycles contained in $f$:  $\sigma_1 = \{e_1, e_6\}$ with label pair $\{2,a\}$ and $\sigma_2 = \{e_3,e_7\}$ with label pair $\{a{+}1, b\}$.

Since $\sigma_1$ is an extended $S2$ cycle, by \fullref{lieinanannulus} the edges $e_1$ and $e_6$ lie in an essential annulus $A_1$.  Similarly, since $\sigma_2$ is an extended $S2$ cycle, the edges $e_3$ and $e_7$ lie in an essential annulus $A_2$.  The three edges $e_2$, $e_3$, and $e_5$ connect the vertices $U_a$ and $U_{a+1}$ (via the vertices $U_1$ and $U_{b+1}$) and thus lie in an annulus of $\hatT$ between the cores of $A_1$ and $A_2$.  These seven edges must lie on $\hatT$ as in \fullref{fig:trigonalessthanb}(b).

 Choose rectangles $\rho_{(1, 2)}$ and $\rho_{(b,\,b+1)}$ on $H_{(1,2)}$ and $H_{(b,\,b+1)}$ respectively that connect the edges of the bigons and trigon on them.  Form a larger trigon $F$ by connecting the two bigons $\kreis{A}$ and $\kreis{B}$ to the trigon $\kreis{C}$ with the rectangles $\rho_{(1,2)}$ and $\rho_{(b,\,b+1)}$.  Let $\bdry_1 F$ be the edge of $F$ with label pair $\{a,a\}$, $\bdry_2 F$ be the edge with label pair $\{a{+}1,a{+}1\}$, an $\bdry_3 F$ be the edge with label pair $\{a, a{+}1\}$.  Note that $\bdry_1 F$ and $\bdry_2 F$ each lie in an essential annulus.  See \fullref{fig:trigonF}(a) for the labeling of $\bdry F$ and \fullref{fig:trigonF}(b) for the placement of $\bdry F$ on $\hatT$.

\begin{figure}[ht!]
\centering
\begin{picture}(0,0)%
\includegraphics[scale=.95]{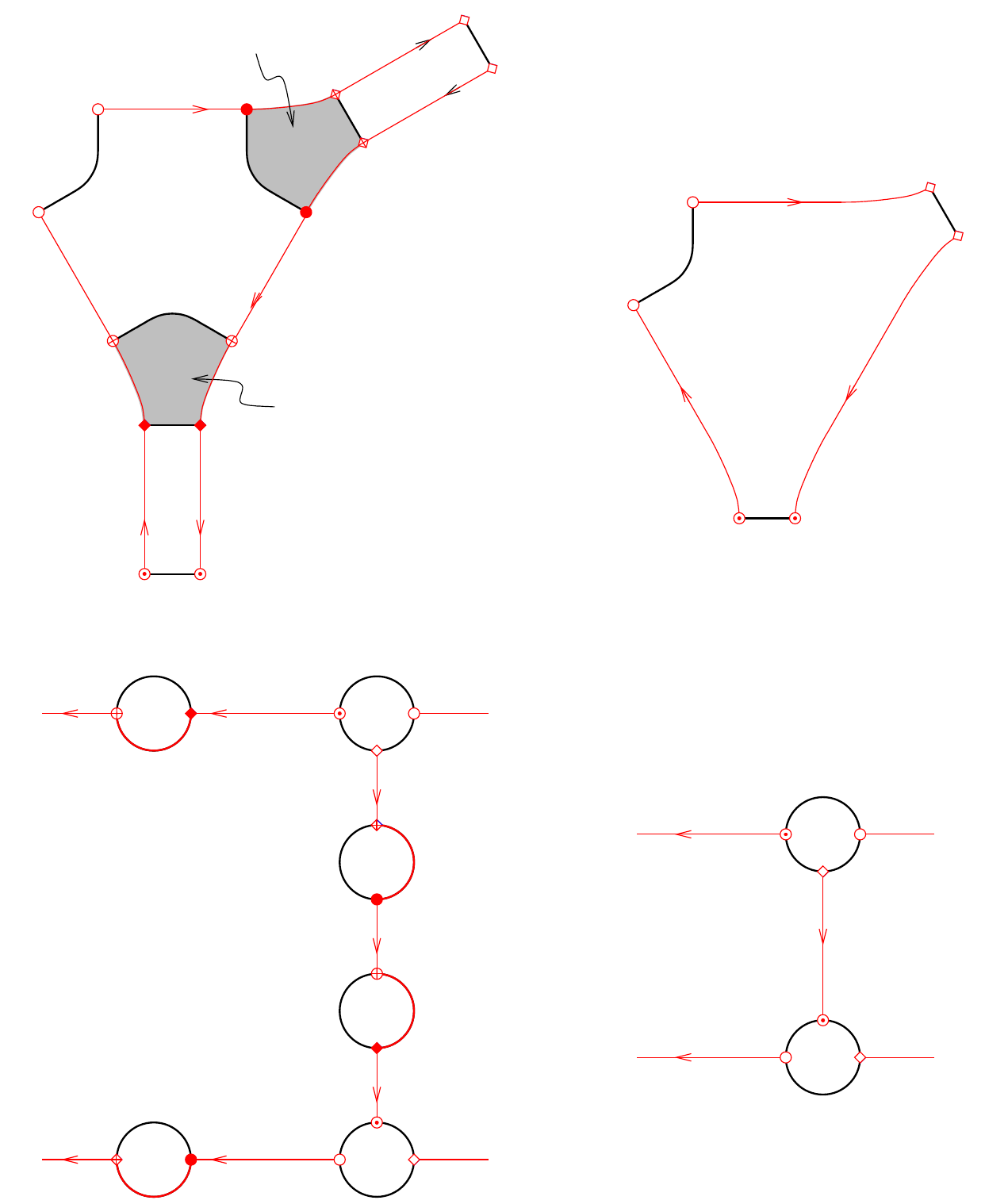}%
\end{picture}%
\setlength{\unitlength}{2812sp}%
\begingroup\makeatletter\ifx\SetFigFont\undefined%
\gdef\SetFigFont#1#2#3#4#5{%
  \reset@font\fontsize{#1}{#2pt}%
  \fontfamily{#3}\fontseries{#4}\fontshape{#5}%
  \selectfont}%
\fi\endgroup%
\begin{picture}(8118,9711)(2509,-10090)
\put(2786,-2086){\makebox(0,0)[rb]{\smash{{\SetFigFont{9}{10.8}{\rmdefault}{\mddefault}{\updefault}{\color[rgb]{0,0,0}$a$}%
}}}}
\put(3526,-1186){\makebox(0,0)[lb]{\smash{{\SetFigFont{9}{10.8}{\rmdefault}{\mddefault}{\updefault}{\color[rgb]{0,0,0}$e_7$}%
}}}}
\put(4801,-2686){\makebox(0,0)[lb]{\smash{{\SetFigFont{9}{10.8}{\rmdefault}{\mddefault}{\updefault}{\color[rgb]{0,0,0}$e_5$}%
}}}}
\put(3901,-2236){\makebox(0,0)[b]{\smash{{\SetFigFont{9}{10.8}{\rmdefault}{\mddefault}{\updefault}{\color[rgb]{0,0,0}$\bigkreis{C}$}%
}}}}
\put(3151,-2986){\makebox(0,0)[rb]{\smash{{\SetFigFont{9}{10.8}{\rmdefault}{\mddefault}{\updefault}{\color[rgb]{0,0,0}$e_6$}%
}}}}
\put(3901,-4336){\makebox(0,0)[b]{\smash{{\SetFigFont{9}{10.8}{\rmdefault}{\mddefault}{\updefault}{\color[rgb]{0,0,0}$\bigkreis{A}$}%
}}}}
\put(3601,-4636){\makebox(0,0)[rb]{\smash{{\SetFigFont{9}{10.8}{\rmdefault}{\mddefault}{\updefault}{\color[rgb]{0,0,0}$e_1$}%
}}}}
\put(4201,-4636){\makebox(0,0)[lb]{\smash{{\SetFigFont{9}{10.8}{\rmdefault}{\mddefault}{\updefault}{\color[rgb]{0,0,0}$e_2$}%
}}}}
\put(5701,-1186){\makebox(0,0)[b]{\smash{{\SetFigFont{9}{10.8}{\rmdefault}{\mddefault}{\updefault}{\color[rgb]{0,0,0}$\bigkreis{B}$}%
}}}}
\put(5476,-661){\makebox(0,0)[lb]{\smash{{\SetFigFont{9}{10.8}{\rmdefault}{\mddefault}{\updefault}{\color[rgb]{0,0,0}$e_4$}%
}}}}
\put(6076,-1336){\makebox(0,0)[lb]{\smash{{\SetFigFont{9}{10.8}{\rmdefault}{\mddefault}{\updefault}{\color[rgb]{0,0,0}$e_3$}%
}}}}
\put(3226,-1261){\makebox(0,0)[rb]{\smash{{\SetFigFont{9}{10.8}{\rmdefault}{\mddefault}{\updefault}{\color[rgb]{0,0,0}$a+1$}%
}}}}
\put(6226,-511){\makebox(0,0)[lb]{\smash{{\SetFigFont{9}{10.8}{\rmdefault}{\mddefault}{\updefault}{\color[rgb]{0,0,0}$a+1$}%
}}}}
\put(6526,-1036){\makebox(0,0)[lb]{\smash{{\SetFigFont{9}{10.8}{\rmdefault}{\mddefault}{\updefault}{\color[rgb]{0,0,0}$a$}%
}}}}
\put(4201,-5161){\makebox(0,0)[lb]{\smash{{\SetFigFont{9}{10.8}{\rmdefault}{\mddefault}{\updefault}{\color[rgb]{0,0,0}$a+1$}%
}}}}
\put(4501,-1186){\makebox(0,0)[rb]{\smash{{\SetFigFont{9}{10.8}{\rmdefault}{\mddefault}{\updefault}{\color[rgb]{0,0,0}$b$}%
}}}}
\put(3601,-3886){\makebox(0,0)[rb]{\smash{{\SetFigFont{9}{10.8}{\rmdefault}{\mddefault}{\updefault}{\color[rgb]{0,0,0}$1$}%
}}}}
\put(4201,-3886){\makebox(0,0)[lb]{\smash{{\SetFigFont{9}{10.8}{\rmdefault}{\mddefault}{\updefault}{\color[rgb]{0,0,0}$2$}%
}}}}
\put(3601,-5161){\makebox(0,0)[rb]{\smash{{\SetFigFont{9}{10.8}{\rmdefault}{\mddefault}{\updefault}{\color[rgb]{0,0,0}$a$}%
}}}}
\put(3376,-3286){\makebox(0,0)[rb]{\smash{{\SetFigFont{9}{10.8}{\rmdefault}{\mddefault}{\updefault}{\color[rgb]{0,0,0}$1$}%
}}}}
\put(4426,-3286){\makebox(0,0)[lb]{\smash{{\SetFigFont{9}{10.8}{\rmdefault}{\mddefault}{\updefault}{\color[rgb]{0,0,0}$2$}%
}}}}
\put(5176,-1111){\makebox(0,0)[rb]{\smash{{\SetFigFont{9}{10.8}{\rmdefault}{\mddefault}{\updefault}{\color[rgb]{0,0,0}$b$}%
}}}}
\put(5026,-2311){\makebox(0,0)[lb]{\smash{{\SetFigFont{9}{10.8}{\rmdefault}{\mddefault}{\updefault}{\color[rgb]{0,0,0}$b+1$}%
}}}}
\put(5476,-1711){\makebox(0,0)[lb]{\smash{{\SetFigFont{9}{10.8}{\rmdefault}{\mddefault}{\updefault}{\color[rgb]{0,0,0}$b+1$}%
}}}}
\put(5551,-6211){\makebox(0,0)[b]{\smash{{\SetFigFont{9}{10.8}{\rmdefault}{\mddefault}{\updefault}{\color[rgb]{0,0,0}$a$}%
}}}}
\put(5551,-7411){\makebox(0,0)[b]{\smash{{\SetFigFont{9}{10.8}{\rmdefault}{\mddefault}{\updefault}{\color[rgb]{0,0,0}$b+1$}%
}}}}
\put(5551,-8611){\makebox(0,0)[b]{\smash{{\SetFigFont{9}{10.8}{\rmdefault}{\mddefault}{\updefault}{\color[rgb]{0,0,0}$1$}%
}}}}
\put(5551,-9811){\makebox(0,0)[b]{\smash{{\SetFigFont{9}{10.8}{\rmdefault}{\mddefault}{\updefault}{\color[rgb]{0,0,0}$a+1$}%
}}}}
\put(3751,-9811){\makebox(0,0)[b]{\smash{{\SetFigFont{9}{10.8}{\rmdefault}{\mddefault}{\updefault}{\color[rgb]{0,0,0}$b$}%
}}}}
\put(5401,-9211){\makebox(0,0)[rb]{\smash{{\SetFigFont{9}{10.8}{\rmdefault}{\mddefault}{\updefault}{\color[rgb]{0,0,0}$e_2$}%
}}}}
\put(3001,-6061){\makebox(0,0)[lb]{\smash{{\SetFigFont{9}{10.8}{\rmdefault}{\mddefault}{\updefault}{\color[rgb]{0,0,0}$e_6$}%
}}}}
\put(4351,-6061){\makebox(0,0)[lb]{\smash{{\SetFigFont{9}{10.8}{\rmdefault}{\mddefault}{\updefault}{\color[rgb]{0,0,0}$e_1$}%
}}}}
\put(3001,-9961){\makebox(0,0)[lb]{\smash{{\SetFigFont{9}{10.8}{\rmdefault}{\mddefault}{\updefault}{\color[rgb]{0,0,0}$e_4$}%
}}}}
\put(4351,-9961){\makebox(0,0)[lb]{\smash{{\SetFigFont{9}{10.8}{\rmdefault}{\mddefault}{\updefault}{\color[rgb]{0,0,0}$e_7$}%
}}}}
\put(6001,-9961){\makebox(0,0)[lb]{\smash{{\SetFigFont{9}{10.8}{\rmdefault}{\mddefault}{\updefault}{\color[rgb]{0,0,0}$e_4$}%
}}}}
\put(6001,-6061){\makebox(0,0)[lb]{\smash{{\SetFigFont{9}{10.8}{\rmdefault}{\mddefault}{\updefault}{\color[rgb]{0,0,0}$e_6$}%
}}}}
\put(5401,-6811){\makebox(0,0)[rb]{\smash{{\SetFigFont{9}{10.8}{\rmdefault}{\mddefault}{\updefault}{\color[rgb]{0,0,0}$e_3$}%
}}}}
\put(5401,-8011){\makebox(0,0)[rb]{\smash{{\SetFigFont{9}{10.8}{\rmdefault}{\mddefault}{\updefault}{\color[rgb]{0,0,0}$e_5$}%
}}}}
\put(3751,-6211){\makebox(0,0)[b]{\smash{{\SetFigFont{9}{10.8}{\rmdefault}{\mddefault}{\updefault}{\color[rgb]{0,0,0}$2$}%
}}}}
\put(8701,-2911){\makebox(0,0)[b]{\smash{{\SetFigFont{9}{10.8}{\rmdefault}{\mddefault}{\updefault}{\color[rgb]{0,0,0}$F$}%
}}}}
\put(8851,-1936){\makebox(0,0)[b]{\smash{{\SetFigFont{9}{10.8}{\rmdefault}{\mddefault}{\updefault}{\color[rgb]{0,0,0}$\bdry_2 F$}%
}}}}
\put(9601,-3436){\makebox(0,0)[lb]{\smash{{\SetFigFont{9}{10.8}{\rmdefault}{\mddefault}{\updefault}{\color[rgb]{0,0,0}$\bdry_3 F$}%
}}}}
\put(8476,-8836){\makebox(0,0)[rb]{\smash{{\SetFigFont{9}{10.8}{\rmdefault}{\mddefault}{\updefault}{\color[rgb]{0,0,0}$\bdry_2 F$}%
}}}}
\put(8101,-3886){\makebox(0,0)[rb]{\smash{{\SetFigFont{9}{10.8}{\rmdefault}{\mddefault}{\updefault}{\color[rgb]{0,0,0}$\bdry_1 F$}%
}}}}
\put(9901,-1786){\makebox(0,0)[lb]{\smash{{\SetFigFont{9}{10.8}{\rmdefault}{\mddefault}{\updefault}{\color[rgb]{0,0,0}$a+1$}%
}}}}
\put(10351,-2386){\makebox(0,0)[lb]{\smash{{\SetFigFont{9}{10.8}{\rmdefault}{\mddefault}{\updefault}{\color[rgb]{0,0,0}$a$}%
}}}}
\put(8401,-4636){\makebox(0,0)[rb]{\smash{{\SetFigFont{9}{10.8}{\rmdefault}{\mddefault}{\updefault}{\color[rgb]{0,0,0}$a$}%
}}}}
\put(9001,-4636){\makebox(0,0)[lb]{\smash{{\SetFigFont{9}{10.8}{\rmdefault}{\mddefault}{\updefault}{\color[rgb]{0,0,0}$a+1$}%
}}}}
\put(7576,-2836){\makebox(0,0)[rb]{\smash{{\SetFigFont{9}{10.8}{\rmdefault}{\mddefault}{\updefault}{\color[rgb]{0,0,0}$a$}%
}}}}
\put(8026,-2011){\makebox(0,0)[rb]{\smash{{\SetFigFont{9}{10.8}{\rmdefault}{\mddefault}{\updefault}{\color[rgb]{0,0,0}$a+1$}%
}}}}
\put(8476,-7336){\makebox(0,0)[rb]{\smash{{\SetFigFont{9}{10.8}{\rmdefault}{\mddefault}{\updefault}{\color[rgb]{0,0,0}$\bdry_1 F$}%
}}}}
\put(9301,-7861){\makebox(0,0)[lb]{\smash{{\SetFigFont{9}{10.8}{\rmdefault}{\mddefault}{\updefault}{\color[rgb]{0,0,0}$\bdry_3 F$}%
}}}}
\put(9151,-7186){\makebox(0,0)[b]{\smash{{\SetFigFont{9}{10.8}{\rmdefault}{\mddefault}{\updefault}{\color[rgb]{0,0,0}$a$}%
}}}}
\put(6826,-8086){\makebox(0,0)[b]{\smash{{\SetFigFont{17}{20.4}{\rmdefault}{\mddefault}{\updefault}{\color[rgb]{0,0,0}$=$}%
}}}}
\put(6826,-3211){\makebox(0,0)[b]{\smash{{\SetFigFont{17}{20.4}{\rmdefault}{\mddefault}{\updefault}{\color[rgb]{0,0,0}$=$}%
}}}}
\put(6826,-10036){\makebox(0,0)[b]{\smash{{\SetFigFont{9}{10.8}{\rmdefault}{\mddefault}{\updefault}{\color[rgb]{0,0,0}(b)}%
}}}}
\put(6826,-5386){\makebox(0,0)[b]{\smash{{\SetFigFont{9}{10.8}{\rmdefault}{\mddefault}{\updefault}{\color[rgb]{0,0,0}(a)}%
}}}}
\put(9151,-8986){\makebox(0,0)[b]{\smash{{\SetFigFont{9}{10.8}{\rmdefault}{\mddefault}{\updefault}{\color[rgb]{0,0,0}$a+1$}%
}}}}
\put(4051,-736){\makebox(0,0)[lb]{\smash{{\SetFigFont{9}{10.8}{\rmdefault}{\mddefault}{\updefault}{\color[rgb]{0,0,0}$\rho_{(b,\,b+1)}$}%
}}}}
\put(4801,-3886){\makebox(0,0)[lb]{\smash{{\SetFigFont{9}{10.8}{\rmdefault}{\mddefault}{\updefault}{\color[rgb]{0,0,0}$\rho_{(1, 2)}$}%
}}}}
\end{picture}%
\caption{(a) The trigon $F$\qua (b) $\bdry F$ on $\hatT$}
\label{fig:trigonF}
\end{figure}

\fullref{lem:funnytrigon} implies that the core of the annulus in which the edges of $F$ lie bounds a meridional disk of the solid torus on the same side of $\hatT$ as $F$.  Hence the cores of $A_1$ and $A_2$ each bound meridional disks.  Let $\sigma_i'$ be the $S2$ cycle of $G_S$ in the face of $G_S$ bounded by $\sigma_i$ for $i=1,2$.  By \fullref{lieinanannulus}, the core of the annuli in which the edges of $\sigma_i'$ bound meridional disks for $i=1,2$.  This contradicts \fullref{GT:L2.1}.
\end{proof}

\begin{lemma}\label{lem:ruleoutcaseIa>b}
There cannot exist a trigon of $\smash{G_S^1}$ of type I with $c=1$ and $a > b$.
\end{lemma}

\begin{proof}
The lower left branch of \fullref{fig:trigonIa} appears now (rotated upwards) as in \fullref{fig:trigonIb}.

\begin{figure}[ht!]
\centering
\begin{picture}(0,0)%
\includegraphics{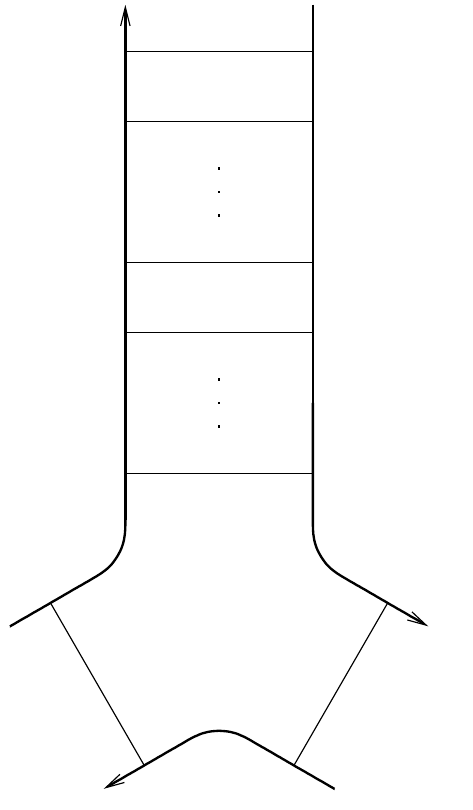}%
\end{picture}%
\setlength{\unitlength}{2960sp}%
\begingroup\makeatletter\ifx\SetFigFont\undefined%
\gdef\SetFigFont#1#2#3#4#5{%
  \reset@font\fontsize{#1}{#2pt}%
  \fontfamily{#3}\fontseries{#4}\fontshape{#5}%
  \selectfont}%
\fi\endgroup%
\begin{picture}(3003,5171)(3999,-4610)
\put(6301,-3886){\makebox(0,0)[lb]{\smash{{\SetFigFont{9}{10.8}{\rmdefault}{\mddefault}{\updefault}{\color[rgb]{0,0,0}$e_1$}%
}}}}
\put(4276,-3211){\makebox(0,0)[rb]{\smash{{\SetFigFont{9}{10.8}{\rmdefault}{\mddefault}{\updefault}{\color[rgb]{0,0,0}$1$}%
}}}}
\put(6151,-1186){\makebox(0,0)[lb]{\smash{{\SetFigFont{9}{10.8}{\rmdefault}{\mddefault}{\updefault}{\color[rgb]{0,0,0}$b$}%
}}}}
\put(6151,-1561){\makebox(0,0)[lb]{\smash{{\SetFigFont{9}{10.8}{\rmdefault}{\mddefault}{\updefault}{\color[rgb]{0,0,0}$b+1$}%
}}}}
\put(6076,164){\makebox(0,0)[lb]{\smash{{\SetFigFont{9}{10.8}{\rmdefault}{\mddefault}{\updefault}{\color[rgb]{0,0,0}$1$}%
}}}}
\put(6076,-286){\makebox(0,0)[lb]{\smash{{\SetFigFont{9}{10.8}{\rmdefault}{\mddefault}{\updefault}{\color[rgb]{0,0,0}$2$}%
}}}}
\put(6076,-2536){\makebox(0,0)[lb]{\smash{{\SetFigFont{9}{10.8}{\rmdefault}{\mddefault}{\updefault}{\color[rgb]{0,0,0}$a$}%
}}}}
\put(6526,-3211){\makebox(0,0)[lb]{\smash{{\SetFigFont{9}{10.8}{\rmdefault}{\mddefault}{\updefault}{\color[rgb]{0,0,0}$a+1$}%
}}}}
\put(4726,164){\makebox(0,0)[rb]{\smash{{\SetFigFont{9}{10.8}{\rmdefault}{\mddefault}{\updefault}{\color[rgb]{0,0,0}$a+1$}%
}}}}
\put(4726,-286){\makebox(0,0)[rb]{\smash{{\SetFigFont{9}{10.8}{\rmdefault}{\mddefault}{\updefault}{\color[rgb]{0,0,0}$a$}%
}}}}
\put(4726,-1186){\makebox(0,0)[rb]{\smash{{\SetFigFont{9}{10.8}{\rmdefault}{\mddefault}{\updefault}{\color[rgb]{0,0,0}$a-b+2$}%
}}}}
\put(4726,-1636){\makebox(0,0)[rb]{\smash{{\SetFigFont{9}{10.8}{\rmdefault}{\mddefault}{\updefault}{\color[rgb]{0,0,0}$a-b+1$}%
}}}}
\put(4726,-2536){\makebox(0,0)[rb]{\smash{{\SetFigFont{9}{10.8}{\rmdefault}{\mddefault}{\updefault}{\color[rgb]{0,0,0}$2$}%
}}}}
\put(5851,-4561){\makebox(0,0)[rb]{\smash{{\SetFigFont{9}{10.8}{\rmdefault}{\mddefault}{\updefault}{\color[rgb]{0,0,0}$b$}%
}}}}
\put(4951,-4561){\makebox(0,0)[lb]{\smash{{\SetFigFont{9}{10.8}{\rmdefault}{\mddefault}{\updefault}{\color[rgb]{0,0,0}$b+1$}%
}}}}
\put(5026,-1786){\makebox(0,0)[lb]{\smash{{\SetFigFont{9}{10.8}{\rmdefault}{\mddefault}{\updefault}{\color[rgb]{0,0,0}$e_3$}%
}}}}
\put(5026,-2386){\makebox(0,0)[lb]{\smash{{\SetFigFont{9}{10.8}{\rmdefault}{\mddefault}{\updefault}{\color[rgb]{0,0,0}$e_2$}%
}}}}
\put(5026,-1036){\makebox(0,0)[lb]{\smash{{\SetFigFont{9}{10.8}{\rmdefault}{\mddefault}{\updefault}{\color[rgb]{0,0,0}$e_4$}%
}}}}
\put(5026,-436){\makebox(0,0)[lb]{\smash{{\SetFigFont{9}{10.8}{\rmdefault}{\mddefault}{\updefault}{\color[rgb]{0,0,0}$e_5$}%
}}}}
\put(5026,314){\makebox(0,0)[lb]{\smash{{\SetFigFont{9}{10.8}{\rmdefault}{\mddefault}{\updefault}{\color[rgb]{0,0,0}$e_6$}%
}}}}
\put(4501,-3886){\makebox(0,0)[rb]{\smash{{\SetFigFont{9}{10.8}{\rmdefault}{\mddefault}{\updefault}{\color[rgb]{0,0,0}$e_0$}%
}}}}
\end{picture}%
\caption{Part of the trigon of type I with $c=1$ and $a > b$}
\label{fig:trigonIb}
\end{figure}
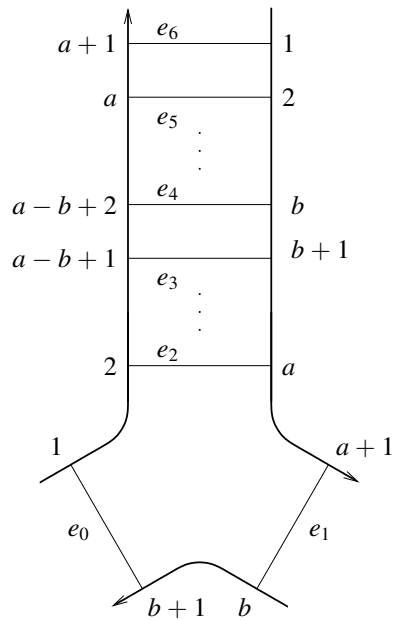

Also, since there is an extended $S2$ cycle with label pair $\{2, a\}$ contained in $f$, by \fullref{lieinanannulus} the two edges with label pair $\{2+i, a-i\}$  for each $0\leq i \leq (a+1)/2$ lie in disjoint essential annuli.  The $S2$ cycle has label pair $\{(a+1)/2, (a+3)/2\}$.

If $a-b+1 = b$ so that the $S2$ cycle has label pair $\{b, b{+}1\}$ then we may view the edges $e_0, e_1, \dots, e_5$ on $T$ as in \fullref{fig:trigonIc}.  The endpoints labeled $b+1$ of edges $e_4$, $e_0$, and $e_3$ are immediately preceded by the endpoints labeled $b$ of edges $e_3$, $e_1$, and $e_4$ respectively.  Due to orientations, given that the edges $e_4$, $e_0$, and $e_3$ appear in that order clockwise around the vertex $U_{b+1}$, the edges $e_3$, $e_1$, and $e_4$ must appear in the order counterclockwise around the vertex $U_b$ as shown.  Hence the edges $e_0$ and $e_1$ must emanate from opposite sides of the annulus in which the edges $e_3$ and $e_4$ lie. This however leaves no possible position for $e_6$. 
\begin{figure}[ht!]
\centering
\begin{picture}(0,0)%
\includegraphics{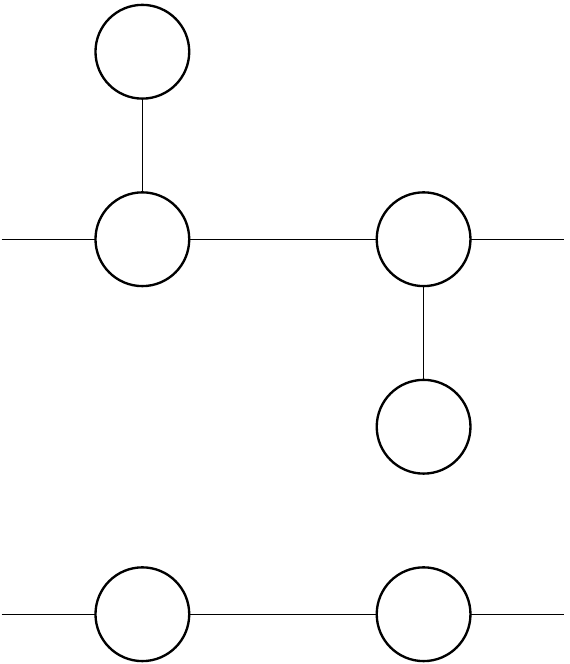}%
\end{picture}%
\setlength{\unitlength}{2960sp}%
\begingroup\makeatletter\ifx\SetFigFont\undefined%
\gdef\SetFigFont#1#2#3#4#5{%
  \reset@font\fontsize{#1}{#2pt}%
  \fontfamily{#3}\fontseries{#4}\fontshape{#5}%
  \selectfont}%
\fi\endgroup%
\begin{picture}(3624,4228)(3589,-5475)
\put(4351,-2236){\makebox(0,0)[rb]{\smash{{\SetFigFont{9}{10.8}{\rmdefault}{\mddefault}{\updefault}{\color[rgb]{0,0,0}$e_1$}%
}}}}
\put(6301,-2836){\makebox(0,0)[b]{\smash{{\SetFigFont{9}{10.8}{\rmdefault}{\mddefault}{\updefault}{\color[rgb]{0,0,0}$b+1$}%
}}}}
\put(6301,-4036){\makebox(0,0)[b]{\smash{{\SetFigFont{9}{10.8}{\rmdefault}{\mddefault}{\updefault}{\color[rgb]{0,0,0}$1$}%
}}}}
\put(3751,-5386){\makebox(0,0)[lb]{\smash{{\SetFigFont{9}{10.8}{\rmdefault}{\mddefault}{\updefault}{\color[rgb]{0,0,0}$e_5$}%
}}}}
\put(5101,-5386){\makebox(0,0)[lb]{\smash{{\SetFigFont{9}{10.8}{\rmdefault}{\mddefault}{\updefault}{\color[rgb]{0,0,0}$e_2$}%
}}}}
\put(6751,-5386){\makebox(0,0)[lb]{\smash{{\SetFigFont{9}{10.8}{\rmdefault}{\mddefault}{\updefault}{\color[rgb]{0,0,0}$e_5$}%
}}}}
\put(6151,-3436){\makebox(0,0)[rb]{\smash{{\SetFigFont{9}{10.8}{\rmdefault}{\mddefault}{\updefault}{\color[rgb]{0,0,0}$e_0$}%
}}}}
\put(4501,-2836){\makebox(0,0)[b]{\smash{{\SetFigFont{9}{10.8}{\rmdefault}{\mddefault}{\updefault}{\color[rgb]{0,0,0}$b$}%
}}}}
\put(4501,-5236){\makebox(0,0)[b]{\smash{{\SetFigFont{9}{10.8}{\rmdefault}{\mddefault}{\updefault}{\color[rgb]{0,0,0}$2$}%
}}}}
\put(4501,-1636){\makebox(0,0)[b]{\smash{{\SetFigFont{9}{10.8}{\rmdefault}{\mddefault}{\updefault}{\color[rgb]{0,0,0}$a+1$}%
}}}}
\put(6301,-5236){\makebox(0,0)[b]{\smash{{\SetFigFont{9}{10.8}{\rmdefault}{\mddefault}{\updefault}{\color[rgb]{0,0,0}$a$}%
}}}}
\put(3751,-2686){\makebox(0,0)[lb]{\smash{{\SetFigFont{9}{10.8}{\rmdefault}{\mddefault}{\updefault}{\color[rgb]{0,0,0}$e_4$}%
}}}}
\put(5101,-2686){\makebox(0,0)[lb]{\smash{{\SetFigFont{9}{10.8}{\rmdefault}{\mddefault}{\updefault}{\color[rgb]{0,0,0}$e_3$}%
}}}}
\put(6751,-2686){\makebox(0,0)[lb]{\smash{{\SetFigFont{9}{10.8}{\rmdefault}{\mddefault}{\updefault}{\color[rgb]{0,0,0}$e_4$}%
}}}}
\end{picture}%
\caption{Edges of the trigon of type I with $c=1$ and $a = 2b-1$ on $T$}
\label{fig:trigonIc}
\end{figure}
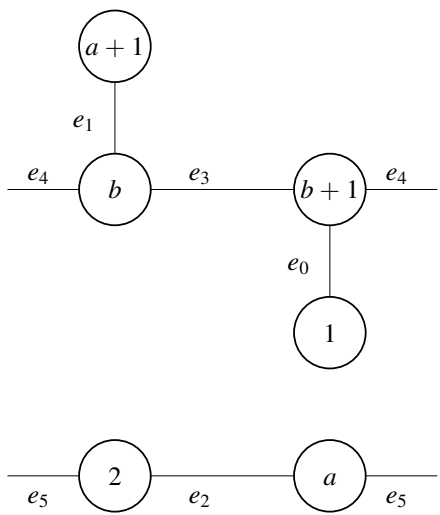

If the $S2$ cycle does not have label pair $\{b, b{+}1\}$, then let $\smash{e_3'}$ and $\smash{e_4'}$ be the other edges on $f$ that have the same labels as $e_3$ and $e_4$ respectively.  Since each of $\{e_3, \smash{e_3'}\}$ and $\{e_4, \smash{e_4'}\}$ are extended $S2$ cycles, they each lie in an essential annulus.    

Furthermore each pair of edges $\{e_3, e_4\}$ and $\{\smash{e_3'}, \smash{e_4'}\}$ bounds a bigon $B$ and $B'$ respectively on $G_S$.  We may thus form an ``annulus'' $$A = B \cup H_{(b,\,b+1)} \cup B' \cup H_{(a-b+1,\,a-b+2)}.$$  Then the positions of the edges $e_3$ and $\smash{e_3'}$ in $G_T$ dictate the placement of the edges $e_4$ and $\smash{e_4'}$.  If $e_3$ and $\smash{e_3'}$ are swapped, then $e_4$ and $\smash{e_4'}$ must also be swapped.  Otherwise $A$ would be a properly embedded nonseparating annulus in one of the Heegaard solid tori which cannot occur.  We may thus view the edges $e_0, e_1, e_3, \smash{e_3'}, e_4, \smash{e_4'}$ and $e_6$ on $T$ as in \fullref{fig:trigonId}.

Since the edges $e_2$ and $e_5$ also form an extended $S2$ cycle, they lie in an essential annulus disjoint from the previous edges.  These two edges appear either as shown in \fullref{fig:trigonId} or perhaps swapped with one another (depending on whether $e_3$ or $\smash{e_3'}$ is closer to $e_2$ on $f$).
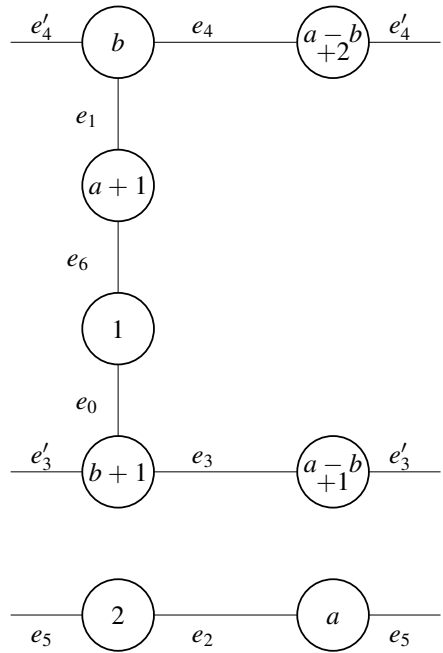
\begin{figure}[ht!]
\centering
\begin{picture}(0,0)%
\includegraphics{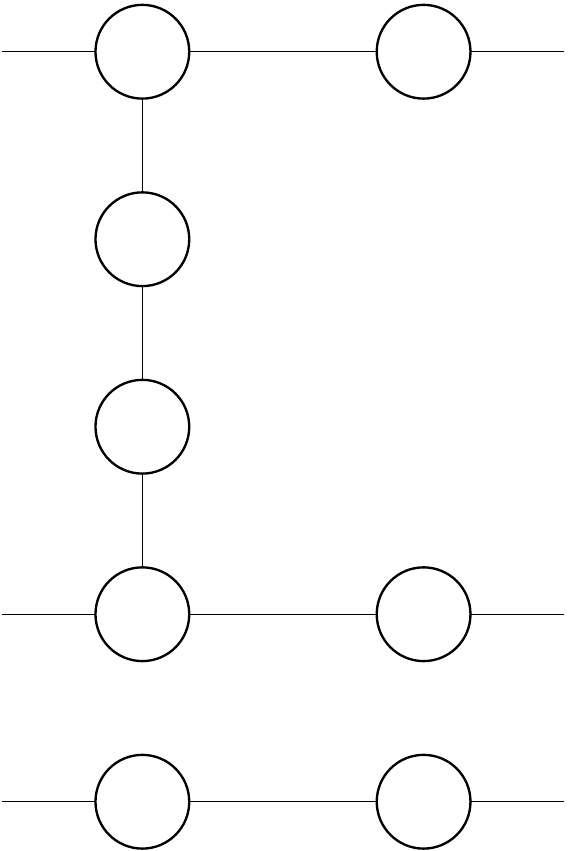}%
\end{picture}%
\setlength{\unitlength}{2960sp}%
\begingroup\makeatletter\ifx\SetFigFont\undefined%
\gdef\SetFigFont#1#2#3#4#5{%
  \reset@font\fontsize{#1}{#2pt}%
  \fontfamily{#3}\fontseries{#4}\fontshape{#5}%
  \selectfont}%
\fi\endgroup%
\begin{picture}(3654,5428)(3589,-5475)
\put(6301,-4111){\makebox(0,0)[b]{\smash{{\SetFigFont{9}{10.8}{\rmdefault}{\mddefault}{\updefault}{\color[rgb]{0,0,0}$+1$}%
}}}}
\put(3751,-5386){\makebox(0,0)[lb]{\smash{{\SetFigFont{9}{10.8}{\rmdefault}{\mddefault}{\updefault}{\color[rgb]{0,0,0}$e_5$}%
}}}}
\put(5101,-5386){\makebox(0,0)[lb]{\smash{{\SetFigFont{9}{10.8}{\rmdefault}{\mddefault}{\updefault}{\color[rgb]{0,0,0}$e_2$}%
}}}}
\put(6751,-5386){\makebox(0,0)[lb]{\smash{{\SetFigFont{9}{10.8}{\rmdefault}{\mddefault}{\updefault}{\color[rgb]{0,0,0}$e_5$}%
}}}}
\put(4501,-5236){\makebox(0,0)[b]{\smash{{\SetFigFont{9}{10.8}{\rmdefault}{\mddefault}{\updefault}{\color[rgb]{0,0,0}$2$}%
}}}}
\put(4501,-1636){\makebox(0,0)[b]{\smash{{\SetFigFont{9}{10.8}{\rmdefault}{\mddefault}{\updefault}{\color[rgb]{0,0,0}$a+1$}%
}}}}
\put(6301,-5236){\makebox(0,0)[b]{\smash{{\SetFigFont{9}{10.8}{\rmdefault}{\mddefault}{\updefault}{\color[rgb]{0,0,0}$a$}%
}}}}
\put(4501,-4036){\makebox(0,0)[b]{\smash{{\SetFigFont{9}{10.8}{\rmdefault}{\mddefault}{\updefault}{\color[rgb]{0,0,0}$b+1$}%
}}}}
\put(4501,-2836){\makebox(0,0)[b]{\smash{{\SetFigFont{9}{10.8}{\rmdefault}{\mddefault}{\updefault}{\color[rgb]{0,0,0}$1$}%
}}}}
\put(4501,-436){\makebox(0,0)[b]{\smash{{\SetFigFont{9}{10.8}{\rmdefault}{\mddefault}{\updefault}{\color[rgb]{0,0,0}$b$}%
}}}}
\put(6751,-286){\makebox(0,0)[lb]{\smash{{\SetFigFont{9}{10.8}{\rmdefault}{\mddefault}{\updefault}{\color[rgb]{0,0,0}$\smash{e_4'}$}%
}}}}
\put(3751,-286){\makebox(0,0)[lb]{\smash{{\SetFigFont{9}{10.8}{\rmdefault}{\mddefault}{\updefault}{\color[rgb]{0,0,0}$\smash{e_4'}$}%
}}}}
\put(5101,-286){\makebox(0,0)[lb]{\smash{{\SetFigFont{9}{10.8}{\rmdefault}{\mddefault}{\updefault}{\color[rgb]{0,0,0}$e_4$}%
}}}}
\put(4351,-1036){\makebox(0,0)[rb]{\smash{{\SetFigFont{9}{10.8}{\rmdefault}{\mddefault}{\updefault}{\color[rgb]{0,0,0}$e_1$}%
}}}}
\put(4351,-3436){\makebox(0,0)[rb]{\smash{{\SetFigFont{9}{10.8}{\rmdefault}{\mddefault}{\updefault}{\color[rgb]{0,0,0}$e_0$}%
}}}}
\put(3751,-3886){\makebox(0,0)[lb]{\smash{{\SetFigFont{9}{10.8}{\rmdefault}{\mddefault}{\updefault}{\color[rgb]{0,0,0}$\smash{e_3'}$}%
}}}}
\put(5101,-3886){\makebox(0,0)[lb]{\smash{{\SetFigFont{9}{10.8}{\rmdefault}{\mddefault}{\updefault}{\color[rgb]{0,0,0}$e_3$}%
}}}}
\put(6751,-3886){\makebox(0,0)[lb]{\smash{{\SetFigFont{9}{10.8}{\rmdefault}{\mddefault}{\updefault}{\color[rgb]{0,0,0}$\smash{e_3'}$}%
}}}}
\put(4276,-2236){\makebox(0,0)[rb]{\smash{{\SetFigFont{9}{10.8}{\rmdefault}{\mddefault}{\updefault}{\color[rgb]{0,0,0}$e_6$}%
}}}}
\put(6301,-361){\makebox(0,0)[b]{\smash{{\SetFigFont{9}{10.8}{\rmdefault}{\mddefault}{\updefault}{\color[rgb]{0,0,0}$a-b$}%
}}}}
\put(6301,-511){\makebox(0,0)[b]{\smash{{\SetFigFont{9}{10.8}{\rmdefault}{\mddefault}{\updefault}{\color[rgb]{0,0,0}$+2$}%
}}}}
\put(6301,-3961){\makebox(0,0)[b]{\smash{{\SetFigFont{9}{10.8}{\rmdefault}{\mddefault}{\updefault}{\color[rgb]{0,0,0}$a-b$}%
}}}}
\end{picture}%
\caption{Edges of the trigon of type I with $c=1$, $a > b$, and $a \neq 2b-1$ on $T$}
\label{fig:trigonId}
\end{figure}

Observe that the annulus $A$ separates the vertices $U_1$ and $U_2$ on $\hatT$.  Since $1$ and $b$ are both odd, $H_{(1, 2)}$ lies on the same side of $\hatT$ as $A$.  Though the interiors of $B$ and $B'$ may intersect $\hatT$ in simple closed curves, $H_{(1, 2)}$ must intersect $A$.  This cannot occur.
\end{proof}

\begin{figure}[ht!]
\centering
\begin{picture}(0,0)%
\includegraphics{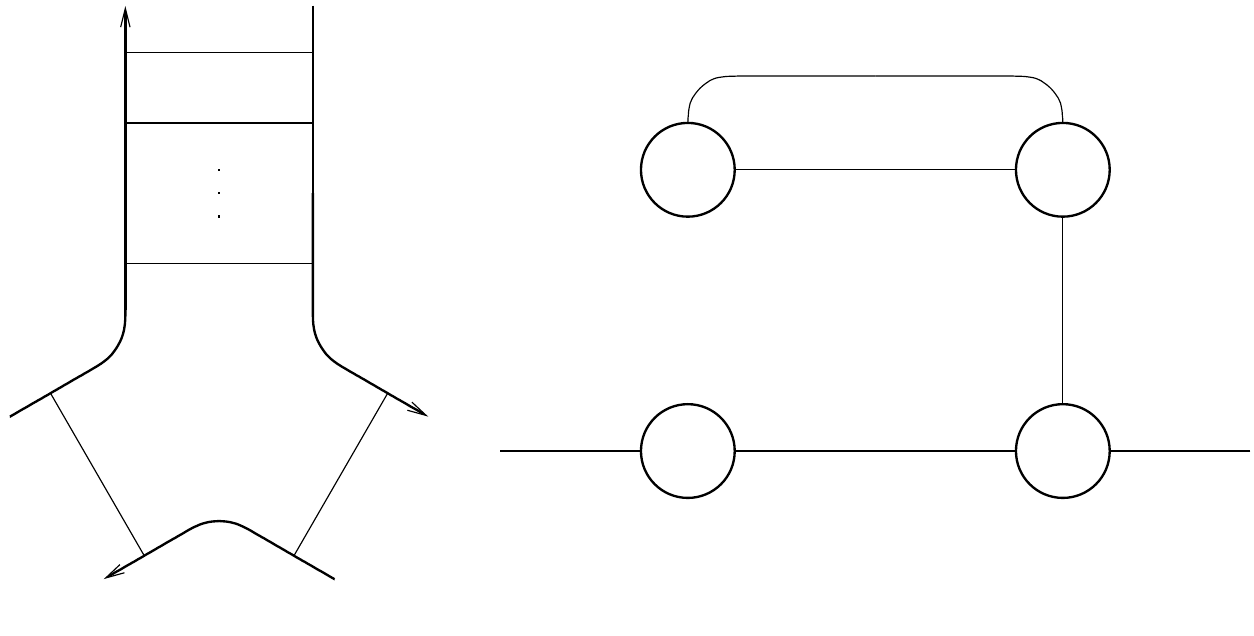}%
\end{picture}%
\setlength{\unitlength}{2960sp}%
\begingroup\makeatletter\ifx\SetFigFont\undefined%
\gdef\SetFigFont#1#2#3#4#5{%
  \reset@font\fontsize{#1}{#2pt}%
  \fontfamily{#3}\fontseries{#4}\fontshape{#5}%
  \selectfont}%
\fi\endgroup%
\begin{picture}(8014,4051)(3999,-4840)
\put(9601,-4786){\makebox(0,0)[b]{\smash{{\SetFigFont{9}{10.8}{\rmdefault}{\mddefault}{\updefault}{\color[rgb]{0,0,0}(b)}%
}}}}
\put(4276,-3211){\makebox(0,0)[rb]{\smash{{\SetFigFont{9}{10.8}{\rmdefault}{\mddefault}{\updefault}{\color[rgb]{0,0,0}$1$}%
}}}}
\put(6076,-2536){\makebox(0,0)[lb]{\smash{{\SetFigFont{9}{10.8}{\rmdefault}{\mddefault}{\updefault}{\color[rgb]{0,0,0}$a$}%
}}}}
\put(6526,-3211){\makebox(0,0)[lb]{\smash{{\SetFigFont{9}{10.8}{\rmdefault}{\mddefault}{\updefault}{\color[rgb]{0,0,0}$a+1$}%
}}}}
\put(4726,-2536){\makebox(0,0)[rb]{\smash{{\SetFigFont{9}{10.8}{\rmdefault}{\mddefault}{\updefault}{\color[rgb]{0,0,0}$2$}%
}}}}
\put(6076,-1636){\makebox(0,0)[lb]{\smash{{\SetFigFont{9}{10.8}{\rmdefault}{\mddefault}{\updefault}{\color[rgb]{0,0,0}$2$}%
}}}}
\put(6076,-1186){\makebox(0,0)[lb]{\smash{{\SetFigFont{9}{10.8}{\rmdefault}{\mddefault}{\updefault}{\color[rgb]{0,0,0}$1$}%
}}}}
\put(4726,-1636){\makebox(0,0)[rb]{\smash{{\SetFigFont{9}{10.8}{\rmdefault}{\mddefault}{\updefault}{\color[rgb]{0,0,0}$a$}%
}}}}
\put(4726,-1186){\makebox(0,0)[rb]{\smash{{\SetFigFont{9}{10.8}{\rmdefault}{\mddefault}{\updefault}{\color[rgb]{0,0,0}$a+1$}%
}}}}
\put(5851,-4561){\makebox(0,0)[rb]{\smash{{\SetFigFont{9}{10.8}{\rmdefault}{\mddefault}{\updefault}{\color[rgb]{0,0,0}$1$}%
}}}}
\put(4951,-4561){\makebox(0,0)[lb]{\smash{{\SetFigFont{9}{10.8}{\rmdefault}{\mddefault}{\updefault}{\color[rgb]{0,0,0}$2$}%
}}}}
\put(4501,-4036){\makebox(0,0)[rb]{\smash{{\SetFigFont{9}{10.8}{\rmdefault}{\mddefault}{\updefault}{\color[rgb]{0,0,0}$e_4$}%
}}}}
\put(6226,-4036){\makebox(0,0)[lb]{\smash{{\SetFigFont{9}{10.8}{\rmdefault}{\mddefault}{\updefault}{\color[rgb]{0,0,0}$e_5$}%
}}}}
\put(4951,-1786){\makebox(0,0)[lb]{\smash{{\SetFigFont{9}{10.8}{\rmdefault}{\mddefault}{\updefault}{\color[rgb]{0,0,0}$e_2$}%
}}}}
\put(4951,-2386){\makebox(0,0)[lb]{\smash{{\SetFigFont{9}{10.8}{\rmdefault}{\mddefault}{\updefault}{\color[rgb]{0,0,0}$e_3$}%
}}}}
\put(4951,-1036){\makebox(0,0)[lb]{\smash{{\SetFigFont{9}{10.8}{\rmdefault}{\mddefault}{\updefault}{\color[rgb]{0,0,0}$e_1$}%
}}}}
\put(5401,-1411){\makebox(0,0)[b]{\smash{{\SetFigFont{9}{10.8}{\rmdefault}{\mddefault}{\updefault}{\color[rgb]{0,0,0}$\bigkreis{A}$}%
}}}}
\put(5401,-3436){\makebox(0,0)[b]{\smash{{\SetFigFont{9}{10.8}{\rmdefault}{\mddefault}{\updefault}{\color[rgb]{0,0,0}$\bigkreis{D}$}%
}}}}
\put(5401,-4786){\makebox(0,0)[b]{\smash{{\SetFigFont{9}{10.8}{\rmdefault}{\mddefault}{\updefault}{\color[rgb]{0,0,0}(a)}%
}}}}
\put(10801,-1936){\makebox(0,0)[b]{\smash{{\SetFigFont{9}{10.8}{\rmdefault}{\mddefault}{\updefault}{\color[rgb]{0,0,0}$1$}%
}}}}
\put(10801,-3736){\makebox(0,0)[b]{\smash{{\SetFigFont{9}{10.8}{\rmdefault}{\mddefault}{\updefault}{\color[rgb]{0,0,0}$2$}%
}}}}
\put(8401,-3736){\makebox(0,0)[b]{\smash{{\SetFigFont{9}{10.8}{\rmdefault}{\mddefault}{\updefault}{\color[rgb]{0,0,0}$a$}%
}}}}
\put(8401,-1936){\makebox(0,0)[b]{\smash{{\SetFigFont{9}{10.8}{\rmdefault}{\mddefault}{\updefault}{\color[rgb]{0,0,0}$a+1$}%
}}}}
\put(9301,-3886){\makebox(0,0)[lb]{\smash{{\SetFigFont{9}{10.8}{\rmdefault}{\mddefault}{\updefault}{\color[rgb]{0,0,0}$e_3$}%
}}}}
\put(11401,-3886){\makebox(0,0)[lb]{\smash{{\SetFigFont{9}{10.8}{\rmdefault}{\mddefault}{\updefault}{\color[rgb]{0,0,0}$e_2$}%
}}}}
\put(7501,-3886){\makebox(0,0)[lb]{\smash{{\SetFigFont{9}{10.8}{\rmdefault}{\mddefault}{\updefault}{\color[rgb]{0,0,0}$e_2$}%
}}}}
\put(9301,-1186){\makebox(0,0)[lb]{\smash{{\SetFigFont{9}{10.8}{\rmdefault}{\mddefault}{\updefault}{\color[rgb]{0,0,0}$e_5$}%
}}}}
\put(10876,-2761){\makebox(0,0)[lb]{\smash{{\SetFigFont{9}{10.8}{\rmdefault}{\mddefault}{\updefault}{\color[rgb]{0,0,0}$e_4$}%
}}}}
\put(9301,-2086){\makebox(0,0)[lb]{\smash{{\SetFigFont{9}{10.8}{\rmdefault}{\mddefault}{\updefault}{\color[rgb]{0,0,0}$e_1$}%
}}}}
\end{picture}%
\caption{(a) A trigon bounded by a forked extended $S2$ cycle of $\smash{G_S^1}$\qua  (b) Edges of $G_S$ on a trigon bounded by a forked extended $S2$ cycle of $\smash{G_S^1}$ as they lie on $G_T$}
\label{fig:forkedbigon}
\end{figure}
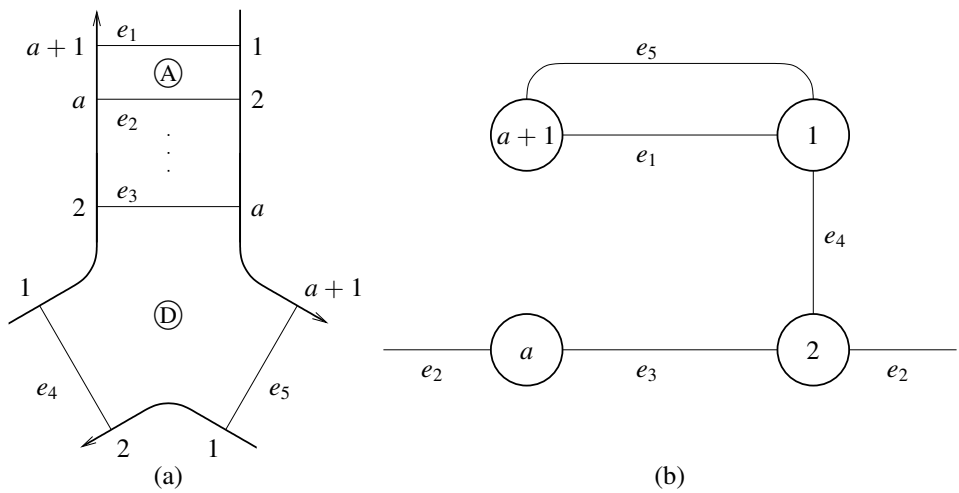

\begin{lemma} \label{forkededges}
If an innermost trigon $g$ of $\smash{G_S^1}$ is bounded by a forked extended $S2$ cycle with edges as in \fullref{fig:forkedbigon}(a), then the edges $e_1, \dots, e_5$ appear on $G_T$ as in \fullref{fig:forkedbigon}(b) (up to symmetries).
\end{lemma}

\begin{proof}
To prove, we examine the possible configurations of these five edges on $G_T$ and see what the bigon $\kreis{A}$ and the trigon $\kreis{D}$ bounded by them on $g$ imply.

The edges $e_2$ and $e_3$ form an (extended) $S2$ cycle $\sigma$.  By \fullref{lieinanannulus} they lie in an essential annulus $A$.  Thus we may begin by assuming these two edges are as shown in \fullref{fig:forkedbigon}(b).   The edge $e_4$ connects vertices $U_1$ and $U_2$, and the edge $e_1$ connects the vertices $U_1$ and $U_{a+1}$.  Without loss of generality, we may assume $e_4$ and $e_1$ are also as shown in \fullref{fig:forkedbigon}(b).  It remains for us to determine the position of $e_5$ (relative to the first four edges).

The endpoints labeled $2$ of edges $e_3$, $e_4$, and $e_2$  are immediately preceded by the endpoints labeled $1$ of edges $e_4$, $e_5$, and $e_1$ respectively.  Due to orientations, given that the edges $e_3$, $e_4$, and $e_2$ appear in that order clockwise around the vertex $U_2$, the edges $e_4$, $e_5$, and $e_1$ must appear in that order counterclockwise around the vertex $U_1$.  Therefore the edge $e_5$ may appear either as in \fullref{fig:forkedbigon}(b) or as in \fullref{fig:forkededges2}.  
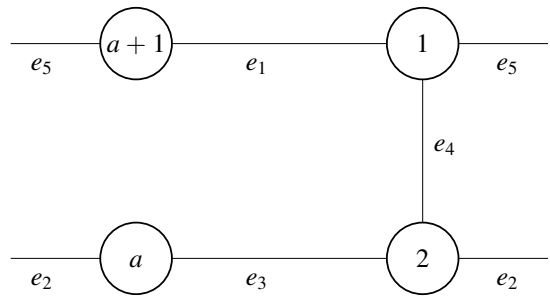
\begin{figure}[ht!]
\centering
\begin{picture}(0,0)%
\includegraphics{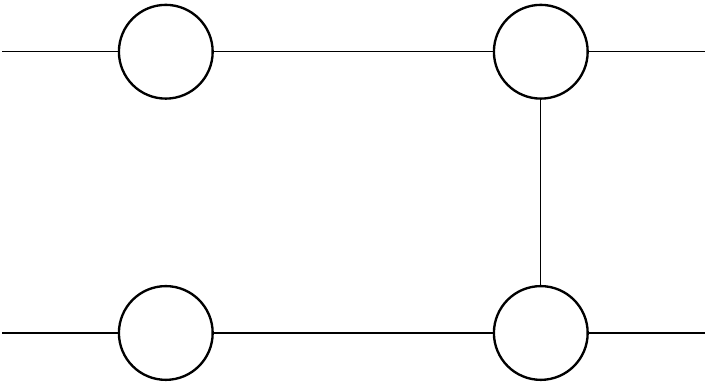}%
\end{picture}%
\setlength{\unitlength}{2960sp}%
\begingroup\makeatletter\ifx\SetFigFont\undefined%
\gdef\SetFigFont#1#2#3#4#5{%
  \reset@font\fontsize{#1}{#2pt}%
  \fontfamily{#3}\fontseries{#4}\fontshape{#5}%
  \selectfont}%
\fi\endgroup%
\begin{picture}(4524,2428)(2389,-3075)
\put(2551,-1186){\makebox(0,0)[lb]{\smash{{\SetFigFont{9}{10.8}{\rmdefault}{\mddefault}{\updefault}{\color[rgb]{0,0,0}$e_5$}%
}}}}
\put(5851,-1036){\makebox(0,0)[b]{\smash{{\SetFigFont{9}{10.8}{\rmdefault}{\mddefault}{\updefault}{\color[rgb]{0,0,0}$1$}%
}}}}
\put(5851,-2836){\makebox(0,0)[b]{\smash{{\SetFigFont{9}{10.8}{\rmdefault}{\mddefault}{\updefault}{\color[rgb]{0,0,0}$2$}%
}}}}
\put(3451,-2836){\makebox(0,0)[b]{\smash{{\SetFigFont{9}{10.8}{\rmdefault}{\mddefault}{\updefault}{\color[rgb]{0,0,0}$a$}%
}}}}
\put(3451,-1036){\makebox(0,0)[b]{\smash{{\SetFigFont{9}{10.8}{\rmdefault}{\mddefault}{\updefault}{\color[rgb]{0,0,0}$a+1$}%
}}}}
\put(4351,-2986){\makebox(0,0)[lb]{\smash{{\SetFigFont{9}{10.8}{\rmdefault}{\mddefault}{\updefault}{\color[rgb]{0,0,0}$e_3$}%
}}}}
\put(6451,-2986){\makebox(0,0)[lb]{\smash{{\SetFigFont{9}{10.8}{\rmdefault}{\mddefault}{\updefault}{\color[rgb]{0,0,0}$e_2$}%
}}}}
\put(2551,-2986){\makebox(0,0)[lb]{\smash{{\SetFigFont{9}{10.8}{\rmdefault}{\mddefault}{\updefault}{\color[rgb]{0,0,0}$e_2$}%
}}}}
\put(5926,-1861){\makebox(0,0)[lb]{\smash{{\SetFigFont{9}{10.8}{\rmdefault}{\mddefault}{\updefault}{\color[rgb]{0,0,0}$e_4$}%
}}}}
\put(6451,-1186){\makebox(0,0)[lb]{\smash{{\SetFigFont{9}{10.8}{\rmdefault}{\mddefault}{\updefault}{\color[rgb]{0,0,0}$e_5$}%
}}}}
\put(4351,-1186){\makebox(0,0)[lb]{\smash{{\SetFigFont{9}{10.8}{\rmdefault}{\mddefault}{\updefault}{\color[rgb]{0,0,0}$e_1$}%
}}}}
\end{picture}%
\caption{Another possible placement of the edge $e_5$}
\label{fig:forkededges2}
\end{figure}

Assume the five edges appear as in \fullref{fig:forkededges2}.  Extend the $(a, a{+}1)$ corner of $\kreis{A}$ and $\kreis{D}$ radially inward through $H_{\smash{(a,\,a{+}1)}}$ to its core $K_{\smash{(a,\,a{+}1)}}$.  Join $\kreis{A}$ to $\kreis{D}$ along this common corner to form the trigon $F$.  

By \fullref{lem:funnytrigon}, the core of the essential annulus $A \subseteq \hatT$ in which the edges of $F$ lie bounds a meridional disk of the solid torus on the same side of $\hatT$ as $F$.  Thus the core of $A$ bounds a meridional disk.  By \fullref{lieinanannulus} the annulus in which the edges of the $S2$ cycle that lies in between the edges of $\sigma$ lie must also bound a meridional disk.  This however contradicts \fullref{GT:L2.1}.

Thus the edges of $\kreis{A}$ and $\kreis{D}$ must appear as in \fullref{fig:forkedbigon}(b).
\end{proof}

\begin{lemma}\label{thinningforkplusbigon}
Let $g$ be a trigon of $\smash{G_S^1}$ with edges as in \fullref{fig:forkedbigon}(a).  If the interior of the bigon on $\hatT$ bounded by edges $e_1$ and $e_5$ (and the vertices $U_1$ and $U_{a+1}$) is disjoint from $K$, then there cannot be a bigon of $G_S$ attached to $e_4$, the edge with label pair $\{1,2\}$.
\end{lemma}

\begin{proof}
Note that since $g$ is bounded by a forked extended $S2$ cycle, $t \geq 4$.
Assume there is a bigon of $G_S$ attached to $e_4$.  Let $g'$ be $g$ with this bigon attached.  Let $\Delta$ be the bigon on $\hatT$ bounded by the edges $e_1$ and $e_5$.  Assume $\Int \Delta \cap K = \emptyset$.  Then, by the minimality assumption on $|S \cap T|$, we have $\Int \Delta \cap g' = \emptyset$.  The corners of $\Delta$ are arcs $u_1$ of $\bdry U_1$ and $u_{a+1}$ of $\bdry U_{a+1}$.  We first consider the case that $a+1 \geq 6$ (and hence $t \geq 6$).  

\begin{figure}[ht!]
\centering
\begin{picture}(0,0)%
\includegraphics{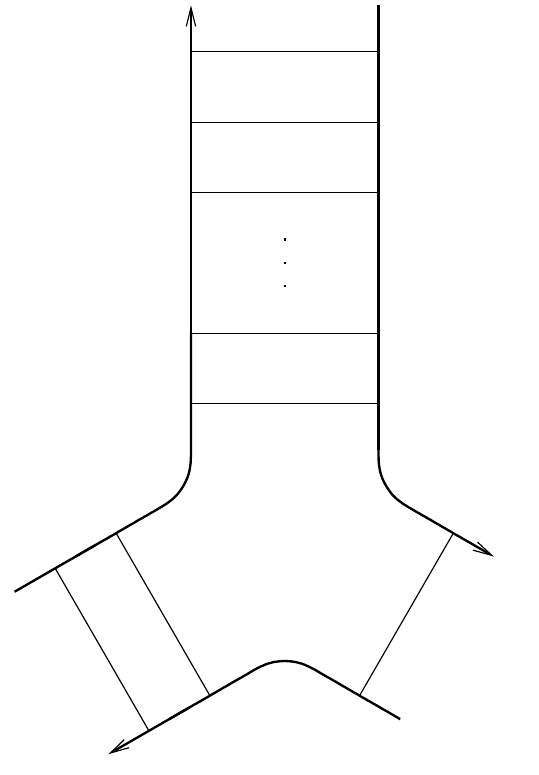}%
\end{picture}%
\setlength{\unitlength}{2960sp}%
\begingroup\makeatletter\ifx\SetFigFont\undefined%
\gdef\SetFigFont#1#2#3#4#5{%
  \reset@font\fontsize{#1}{#2pt}%
  \fontfamily{#3}\fontseries{#4}\fontshape{#5}%
  \selectfont}%
\fi\endgroup%
\begin{picture}(3423,4946)(3579,-4835)
\put(4426,-4036){\makebox(0,0)[b]{\smash{{\SetFigFont{9}{10.8}{\rmdefault}{\mddefault}{\updefault}{\color[rgb]{0,0,0}$\bigkreis{E}$}%
}}}}
\put(4276,-3211){\makebox(0,0)[rb]{\smash{{\SetFigFont{9}{10.8}{\rmdefault}{\mddefault}{\updefault}{\color[rgb]{0,0,0}$1$}%
}}}}
\put(6076,-2536){\makebox(0,0)[lb]{\smash{{\SetFigFont{9}{10.8}{\rmdefault}{\mddefault}{\updefault}{\color[rgb]{0,0,0}$a$}%
}}}}
\put(6526,-3211){\makebox(0,0)[lb]{\smash{{\SetFigFont{9}{10.8}{\rmdefault}{\mddefault}{\updefault}{\color[rgb]{0,0,0}$a+1$}%
}}}}
\put(4951,-4561){\makebox(0,0)[lb]{\smash{{\SetFigFont{9}{10.8}{\rmdefault}{\mddefault}{\updefault}{\color[rgb]{0,0,0}$2$}%
}}}}
\put(5851,-4561){\makebox(0,0)[rb]{\smash{{\SetFigFont{9}{10.8}{\rmdefault}{\mddefault}{\updefault}{\color[rgb]{0,0,0}$1$}%
}}}}
\put(4576,-4786){\makebox(0,0)[lb]{\smash{{\SetFigFont{9}{10.8}{\rmdefault}{\mddefault}{\updefault}{\color[rgb]{0,0,0}$3$}%
}}}}
\put(6076,-736){\makebox(0,0)[lb]{\smash{{\SetFigFont{9}{10.8}{\rmdefault}{\mddefault}{\updefault}{\color[rgb]{0,0,0}$2$}%
}}}}
\put(6076,-286){\makebox(0,0)[lb]{\smash{{\SetFigFont{9}{10.8}{\rmdefault}{\mddefault}{\updefault}{\color[rgb]{0,0,0}$1$}%
}}}}
\put(6076,-1186){\makebox(0,0)[lb]{\smash{{\SetFigFont{9}{10.8}{\rmdefault}{\mddefault}{\updefault}{\color[rgb]{0,0,0}$3$}%
}}}}
\put(4726,-2536){\makebox(0,0)[rb]{\smash{{\SetFigFont{9}{10.8}{\rmdefault}{\mddefault}{\updefault}{\color[rgb]{0,0,0}$2$}%
}}}}
\put(6076,-2086){\makebox(0,0)[lb]{\smash{{\SetFigFont{9}{10.8}{\rmdefault}{\mddefault}{\updefault}{\color[rgb]{0,0,0}$a-1$}%
}}}}
\put(4726,-2086){\makebox(0,0)[rb]{\smash{{\SetFigFont{9}{10.8}{\rmdefault}{\mddefault}{\updefault}{\color[rgb]{0,0,0}$3$}%
}}}}
\put(4726,-736){\makebox(0,0)[rb]{\smash{{\SetFigFont{9}{10.8}{\rmdefault}{\mddefault}{\updefault}{\color[rgb]{0,0,0}$a$}%
}}}}
\put(4726,-286){\makebox(0,0)[rb]{\smash{{\SetFigFont{9}{10.8}{\rmdefault}{\mddefault}{\updefault}{\color[rgb]{0,0,0}$a+1$}%
}}}}
\put(4726,-1186){\makebox(0,0)[rb]{\smash{{\SetFigFont{9}{10.8}{\rmdefault}{\mddefault}{\updefault}{\color[rgb]{0,0,0}$a-1$}%
}}}}
\put(3826,-3436){\makebox(0,0)[rb]{\smash{{\SetFigFont{9}{10.8}{\rmdefault}{\mddefault}{\updefault}{\color[rgb]{0,0,0}$t$}%
}}}}
\put(4951,-136){\makebox(0,0)[lb]{\smash{{\SetFigFont{9}{10.8}{\rmdefault}{\mddefault}{\updefault}{\color[rgb]{0,0,0}$e_1$}%
}}}}
\put(4951,-586){\makebox(0,0)[lb]{\smash{{\SetFigFont{9}{10.8}{\rmdefault}{\mddefault}{\updefault}{\color[rgb]{0,0,0}$e_2$}%
}}}}
\put(4951,-2686){\makebox(0,0)[lb]{\smash{{\SetFigFont{9}{10.8}{\rmdefault}{\mddefault}{\updefault}{\color[rgb]{0,0,0}$e_3$}%
}}}}
\put(6001,-3961){\makebox(0,0)[rb]{\smash{{\SetFigFont{9}{10.8}{\rmdefault}{\mddefault}{\updefault}{\color[rgb]{0,0,0}$e_5$}%
}}}}
\put(4801,-3961){\makebox(0,0)[lb]{\smash{{\SetFigFont{9}{10.8}{\rmdefault}{\mddefault}{\updefault}{\color[rgb]{0,0,0}$e_4$}%
}}}}
\put(4951,-1336){\makebox(0,0)[lb]{\smash{{\SetFigFont{9}{10.8}{\rmdefault}{\mddefault}{\updefault}{\color[rgb]{0,0,0}$e_6$}%
}}}}
\put(4951,-1936){\makebox(0,0)[lb]{\smash{{\SetFigFont{9}{10.8}{\rmdefault}{\mddefault}{\updefault}{\color[rgb]{0,0,0}$e_7$}%
}}}}
\put(4276,-4411){\makebox(0,0)[rb]{\smash{{\SetFigFont{9}{10.8}{\rmdefault}{\mddefault}{\updefault}{\color[rgb]{0,0,0}$e_8$}%
}}}}
\put(5401,-961){\makebox(0,0)[b]{\smash{{\SetFigFont{9}{10.8}{\rmdefault}{\mddefault}{\updefault}{\color[rgb]{0,0,0}$\bigkreis{B}$}%
}}}}
\put(5401,-511){\makebox(0,0)[b]{\smash{{\SetFigFont{9}{10.8}{\rmdefault}{\mddefault}{\updefault}{\color[rgb]{0,0,0}$\bigkreis{A}$}%
}}}}
\put(5401,-2311){\makebox(0,0)[b]{\smash{{\SetFigFont{9}{10.8}{\rmdefault}{\mddefault}{\updefault}{\color[rgb]{0,0,0}$\bigkreis{C}$}%
}}}}
\put(5401,-3436){\makebox(0,0)[b]{\smash{{\SetFigFont{9}{10.8}{\rmdefault}{\mddefault}{\updefault}{\color[rgb]{0,0,0}$\bigkreis{D}$}%
}}}}
\end{picture}%
\caption{The bigons and trigon of $g'$}
\label{fig:forkplusbigon}
\end{figure}
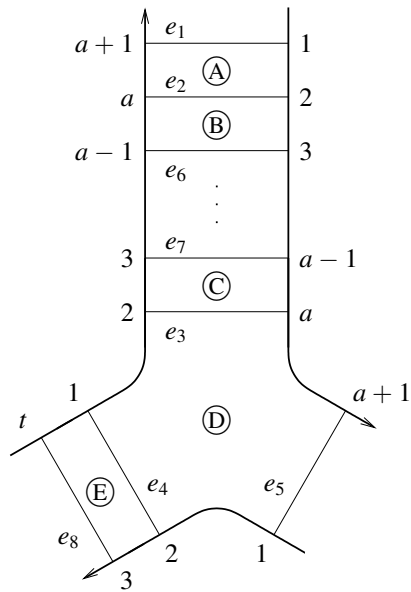
Label the bigons and the trigon of $G_S$ on $g'$ as $\kreis{A}, \kreis{B}, \dots, \kreis{E}$ as shown in \fullref{fig:forkplusbigon}. 
Let $\rho_{(1,2)}$ be the rectangle on $\bdry H_{(1,2)} \cut (U_1 \cup U_2)$ and $u_2$ be the arc of $\bdry U_2$ so that $\rho_{(1,2)}$ is bounded by $u_1$, the $(1,2)$ corner of $\kreis{A}$, $u_2$, and the corner $(1,2)$ of $\kreis{D}$.
 Let $\rho_{\smash{(a,\,a{+}1)}}$ be the rectangle on $\bdry H_{\smash{(a,\,a{+}1)}} \cut (U_a \cup U_{a+1})$ and $u_a$ be the arc of $\bdry U_a$ so that $\rho_{\smash{(a,\,a{+}1)}}$ is bounded by $u_{a+1}$, the $(a, a{+}1)$ corner of $\kreis{A}$, $u_{a}$, and a $(a, a{+}1)$ corner of $\kreis{D}$.
Let $\rho_{(2,3)}$ be the rectangle on $\bdry H_{(2,3)} \cut (U_2 \cup U_3)$ and $u_3$ be the arc of $\bdry U_3$ so that $\rho_{(2,3)}$ is bounded by $u_2$, the $(2,3)$ corner of $\kreis{B}$, $u_3$, and the $(2,3)$ corner of $\kreis{E}$.
Let $\rho_{(a-1,\,a)}$ be the rectangle on $\bdry H_{(a-1,\,a)} \cut (U_{a-1} \cup U_a)$ and $u_{a-1}$ be the arc of $\bdry U_{a-1}$ so that $\rho_{(a-1,\,a)}$ is bounded by $u_{a}$, the $(a{-}1, a)$ corner of $\kreis{B}$, $u_{a-1}$, and the $(a{-}1, a)$ corner of $\kreis{C}$.

\begin{figure}[ht!]
\centering
\begin{picture}(0,0)%
\includegraphics{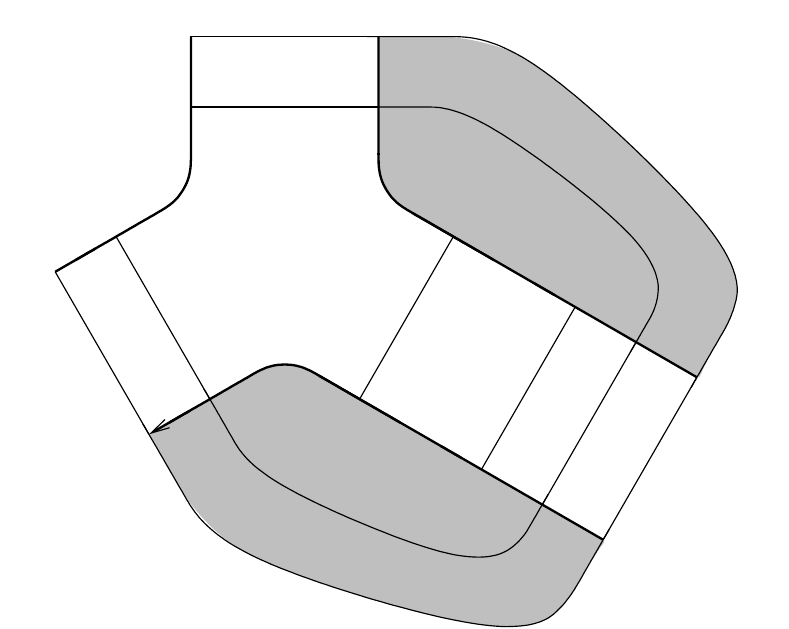}%
\end{picture}%
\setlength{\unitlength}{2960sp}%
\begingroup\makeatletter\ifx\SetFigFont\undefined%
\gdef\SetFigFont#1#2#3#4#5{%
  \reset@font\fontsize{#1}{#2pt}%
  \fontfamily{#3}\fontseries{#4}\fontshape{#5}%
  \selectfont}%
\fi\endgroup%
\begin{picture}(5027,4019)(3579,-5799)
\put(7201,-4411){\makebox(0,0)[b]{\smash{{\SetFigFont{9}{10.8}{\rmdefault}{\mddefault}{\updefault}{\color[rgb]{0,0,0}$\bigkreis{A}$}%
}}}}
\put(6076,-2386){\makebox(0,0)[lb]{\smash{{\SetFigFont{9}{10.8}{\rmdefault}{\mddefault}{\updefault}{\color[rgb]{0,0,0}$a$}%
}}}}
\put(8176,-4186){\makebox(0,0)[lb]{\smash{{\SetFigFont{9}{10.8}{\rmdefault}{\mddefault}{\updefault}{\color[rgb]{0,0,0}$a-1$}%
}}}}
\put(7801,-3961){\makebox(0,0)[lb]{\smash{{\SetFigFont{9}{10.8}{\rmdefault}{\mddefault}{\updefault}{\color[rgb]{0,0,0}$a$}%
}}}}
\put(6901,-5161){\makebox(0,0)[rb]{\smash{{\SetFigFont{9}{10.8}{\rmdefault}{\mddefault}{\updefault}{\color[rgb]{0,0,0}$2$}%
}}}}
\put(4501,-4786){\makebox(0,0)[rb]{\smash{{\SetFigFont{9}{10.8}{\rmdefault}{\mddefault}{\updefault}{\color[rgb]{0,0,0}$3$}%
}}}}
\put(7426,-5461){\makebox(0,0)[lb]{\smash{{\SetFigFont{9}{10.8}{\rmdefault}{\mddefault}{\updefault}{\color[rgb]{0,0,0}$3$}%
}}}}
\put(5101,-4486){\makebox(0,0)[lb]{\smash{{\SetFigFont{9}{10.8}{\rmdefault}{\mddefault}{\updefault}{\color[rgb]{0,0,0}$2$}%
}}}}
\put(4276,-3211){\makebox(0,0)[rb]{\smash{{\SetFigFont{9}{10.8}{\rmdefault}{\mddefault}{\updefault}{\color[rgb]{0,0,0}$1$}%
}}}}
\put(6526,-3211){\makebox(0,0)[lb]{\smash{{\SetFigFont{9}{10.8}{\rmdefault}{\mddefault}{\updefault}{\color[rgb]{0,0,0}$a+1$}%
}}}}
\put(5851,-4561){\makebox(0,0)[rb]{\smash{{\SetFigFont{9}{10.8}{\rmdefault}{\mddefault}{\updefault}{\color[rgb]{0,0,0}$1$}%
}}}}
\put(4726,-2086){\makebox(0,0)[rb]{\smash{{\SetFigFont{9}{10.8}{\rmdefault}{\mddefault}{\updefault}{\color[rgb]{0,0,0}$3$}%
}}}}
\put(3826,-3436){\makebox(0,0)[rb]{\smash{{\SetFigFont{9}{10.8}{\rmdefault}{\mddefault}{\updefault}{\color[rgb]{0,0,0}$t$}%
}}}}
\put(4951,-2686){\makebox(0,0)[lb]{\smash{{\SetFigFont{9}{10.8}{\rmdefault}{\mddefault}{\updefault}{\color[rgb]{0,0,0}$e_3$}%
}}}}
\put(6001,-3961){\makebox(0,0)[rb]{\smash{{\SetFigFont{9}{10.8}{\rmdefault}{\mddefault}{\updefault}{\color[rgb]{0,0,0}$e_5$}%
}}}}
\put(4801,-3961){\makebox(0,0)[lb]{\smash{{\SetFigFont{9}{10.8}{\rmdefault}{\mddefault}{\updefault}{\color[rgb]{0,0,0}$e_4$}%
}}}}
\put(4951,-1936){\makebox(0,0)[lb]{\smash{{\SetFigFont{9}{10.8}{\rmdefault}{\mddefault}{\updefault}{\color[rgb]{0,0,0}$e_7$}%
}}}}
\put(4276,-4411){\makebox(0,0)[rb]{\smash{{\SetFigFont{9}{10.8}{\rmdefault}{\mddefault}{\updefault}{\color[rgb]{0,0,0}$e_8$}%
}}}}
\put(5401,-2311){\makebox(0,0)[b]{\smash{{\SetFigFont{9}{10.8}{\rmdefault}{\mddefault}{\updefault}{\color[rgb]{0,0,0}$\bigkreis{C}$}%
}}}}
\put(5401,-3436){\makebox(0,0)[b]{\smash{{\SetFigFont{9}{10.8}{\rmdefault}{\mddefault}{\updefault}{\color[rgb]{0,0,0}$\bigkreis{D}$}%
}}}}
\put(4426,-4036){\makebox(0,0)[b]{\smash{{\SetFigFont{9}{10.8}{\rmdefault}{\mddefault}{\updefault}{\color[rgb]{0,0,0}$\bigkreis{E}$}%
}}}}
\put(6601,-5011){\makebox(0,0)[rb]{\smash{{\SetFigFont{9}{10.8}{\rmdefault}{\mddefault}{\updefault}{\color[rgb]{0,0,0}$1$}%
}}}}
\put(7276,-3586){\makebox(0,0)[lb]{\smash{{\SetFigFont{9}{10.8}{\rmdefault}{\mddefault}{\updefault}{\color[rgb]{0,0,0}$a+1$}%
}}}}
\put(4726,-2536){\makebox(0,0)[rb]{\smash{{\SetFigFont{9}{10.8}{\rmdefault}{\mddefault}{\updefault}{\color[rgb]{0,0,0}$2$}%
}}}}
\put(7576,-5236){\makebox(0,0)[lb]{\smash{{\SetFigFont{9}{10.8}{\rmdefault}{\mddefault}{\updefault}{\color[rgb]{0,0,0}$e_6$}%
}}}}
\put(7201,-5011){\makebox(0,0)[lb]{\smash{{\SetFigFont{9}{10.8}{\rmdefault}{\mddefault}{\updefault}{\color[rgb]{0,0,0}$e_2$}%
}}}}
\put(6751,-4786){\makebox(0,0)[lb]{\smash{{\SetFigFont{9}{10.8}{\rmdefault}{\mddefault}{\updefault}{\color[rgb]{0,0,0}$e_1$}%
}}}}
\put(7501,-4711){\makebox(0,0)[b]{\smash{{\SetFigFont{9}{10.8}{\rmdefault}{\mddefault}{\updefault}{\color[rgb]{0,0,0}$\bigkreis{B}$}%
}}}}
\put(6076,-1936){\makebox(0,0)[lb]{\smash{{\SetFigFont{9}{10.8}{\rmdefault}{\mddefault}{\updefault}{\color[rgb]{0,0,0}$a-1$}%
}}}}
\put(6526,-4111){\makebox(0,0)[b]{\smash{{\SetFigFont{9}{10.8}{\rmdefault}{\mddefault}{\updefault}{\color[rgb]{0,0,0}$\Delta$}%
}}}}
\put(5926,-4861){\makebox(0,0)[b]{\smash{{\SetFigFont{9}{10.8}{\rmdefault}{\mddefault}{\updefault}{\color[rgb]{0,0,0}$\rho_{(1, 2)}$}%
}}}}
\put(6526,-2836){\makebox(0,0)[b]{\smash{{\SetFigFont{9}{10.8}{\rmdefault}{\mddefault}{\updefault}{\color[rgb]{0,0,0}$\rho_{(a,\,a+1)}$}%
}}}}
\put(6826,-2386){\makebox(0,0)[b]{\smash{{\SetFigFont{9}{10.8}{\rmdefault}{\mddefault}{\updefault}{\color[rgb]{0,0,0}$\rho_{(a-1,\,a)}$}%
}}}}
\put(5326,-5161){\makebox(0,0)[b]{\smash{{\SetFigFont{9}{10.8}{\rmdefault}{\mddefault}{\updefault}{\color[rgb]{0,0,0}$\rho_{(2, 3)}$}%
}}}}
\end{picture}%
\caption{The bigons and trigon of $g'$ assembled to form a long disk $D$}
\label{fig:thinforkplusbigon}
\end{figure}
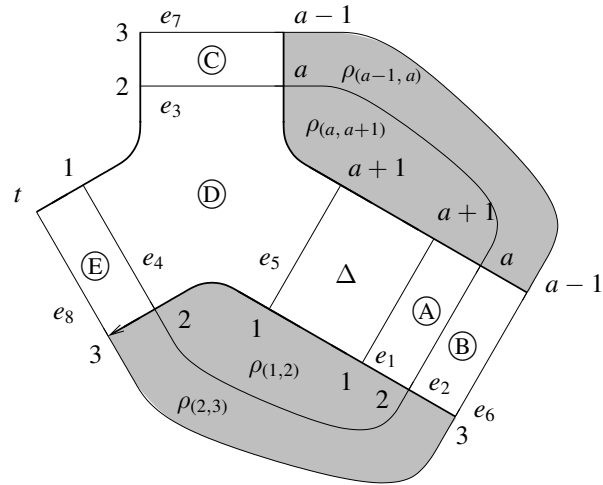
Assemble $\kreis{A}, \kreis{B}, \dots, \kreis{E}$, $\Delta$, $\rho_{(1,2)}$, $\rho_{(2,3)}$, $\rho_{(a-1,\,a)}$, and $\rho_{\smash{(a,\,a{+}1)}}$ to form the embedded disk $D$ as shown in \fullref{fig:thinforkplusbigon}.  The boundary of $D$ is composed of an arc $\alpha = e_7 \cup u_{a-1} \cup e_6 \cup u_3 \cup e_8$ and an arc $\beta = (H_{(t,1)} \cap \kreis{E}) \cup (H_{(1,2)} \cap \kreis{D}) \cup (H_{(2,3)} \cap \kreis{C})$.  The arc $\alpha$ lies on $\hatT$.  The arc $\beta$ lies on $\bdry H_{(t,3)}$ (recall $H_{(t,3)} = H_{(t,1)} \cup H_{(1,2)} \cup H_{(2,3)}$) and can be radially extended into $H_{(t,3)}$ to lie on the arc $K_{(t,3)}$.  Hence we may take $\beta$ to lie on $K$.  Therefore $D$ is a long disk.  By \fullref{longdisk}, such a disk cannot exist.

If $a+1 =4$, then $\kreis{B}$ and $\kreis{C}$ are the same bigon of $G_S$ on $g$.  The two corners of $\kreis{B} = \kreis{C}$ however are distinct arcs on $H_{(2,3)}$.   Also notationally $\rho_{(a-1,\,a)}$ is confused with $\rho_{(2,3)}$.  In this case, let us refer to the former as $\rho_{\smash{(2,3)}}'$.  Further notice that the interiors of $\rho_{(2,3)}$ and $\rho_{\smash{(2,3)}}'$ do not intersect.  We may assemble the disk $D$ using two copies of $\kreis{B}$.  A slight isotopy which fixes $\beta$ and keeps $\alpha$ on $\hatT$ makes $D$ embedded.  (Alternatively, in this construction of $D$, one may regard $\kreis{B}$ and $\kreis{C}$ as slight pushoffs of the same bigon to opposite sides.)  Again, $D$ is a long disk.
\end{proof}

\begin{figure}[ht!]
\centering
\begin{picture}(0,0)%
\includegraphics{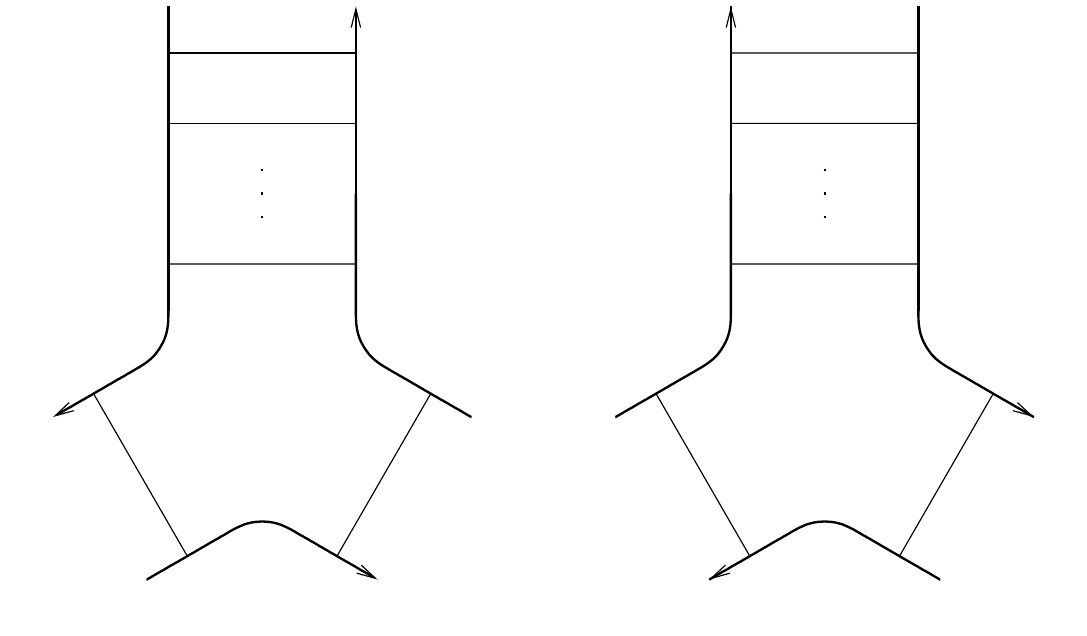}%
\end{picture}%
\setlength{\unitlength}{2960sp}%
\begingroup\makeatletter\ifx\SetFigFont\undefined%
\gdef\SetFigFont#1#2#3#4#5{%
  \reset@font\fontsize{#1}{#2pt}%
  \fontfamily{#3}\fontseries{#4}\fontshape{#5}%
  \selectfont}%
\fi\endgroup%
\begin{picture}(6953,4051)(3724,-4840)
\put(9526,-4561){\makebox(0,0)[rb]{\smash{{\SetFigFont{9}{10.8}{\rmdefault}{\mddefault}{\updefault}{\color[rgb]{0,0,0}$a+1$}%
}}}}
\put(5401,-3436){\makebox(0,0)[b]{\smash{{\SetFigFont{9}{10.8}{\rmdefault}{\mddefault}{\updefault}{\color[rgb]{0,0,0}$\bigkreis{D}'$}%
}}}}
    \put(5401,-1411){\makebox(0,0)[b]{\smash{{\SetFigFont{9}{10.8}{\rmdefault}{\mddefault}{\updefault}{\color[rgb]{0,0,0}$\bigkreis{A}'$}%
}}}}
\put(4951,-1786){\makebox(0,0)[lb]{\smash{{\SetFigFont{9}{10.8}{\rmdefault}{\mddefault}{\updefault}{\color[rgb]{0,0,0}$\smash{e_2'}$}%
}}}}
\put(4951,-1036){\makebox(0,0)[lb]{\smash{{\SetFigFont{9}{10.8}{\rmdefault}{\mddefault}{\updefault}{\color[rgb]{0,0,0}$\smash{e_1'}$}%
}}}}
\put(4951,-2386){\makebox(0,0)[lb]{\smash{{\SetFigFont{9}{10.8}{\rmdefault}{\mddefault}{\updefault}{\color[rgb]{0,0,0}$\smash{e_3'}$}%
}}}}
\put(6226,-4036){\makebox(0,0)[lb]{\smash{{\SetFigFont{9}{10.8}{\rmdefault}{\mddefault}{\updefault}{\color[rgb]{0,0,0}$e_5'$}%
}}}}
\put(4576,-4036){\makebox(0,0)[rb]{\smash{{\SetFigFont{9}{10.8}{\rmdefault}{\mddefault}{\updefault}{\color[rgb]{0,0,0}$\smash{e_4'}$}%
}}}}
\put(4726,-1186){\makebox(0,0)[rb]{\smash{{\SetFigFont{9}{10.8}{\rmdefault}{\mddefault}{\updefault}{\color[rgb]{0,0,0}$b$}%
}}}}
\put(4726,-1636){\makebox(0,0)[rb]{\smash{{\SetFigFont{9}{10.8}{\rmdefault}{\mddefault}{\updefault}{\color[rgb]{0,0,0}$b+1$}%
}}}}
\put(4726,-2536){\makebox(0,0)[rb]{\smash{{\SetFigFont{9}{10.8}{\rmdefault}{\mddefault}{\updefault}{\color[rgb]{0,0,0}$a$}%
}}}}
\put(4201,-3286){\makebox(0,0)[rb]{\smash{{\SetFigFont{9}{10.8}{\rmdefault}{\mddefault}{\updefault}{\color[rgb]{0,0,0}$a+1$}%
}}}}
\put(6076,-1186){\makebox(0,0)[lb]{\smash{{\SetFigFont{9}{10.8}{\rmdefault}{\mddefault}{\updefault}{\color[rgb]{0,0,0}$a+1$}%
}}}}
\put(6076,-1636){\makebox(0,0)[lb]{\smash{{\SetFigFont{9}{10.8}{\rmdefault}{\mddefault}{\updefault}{\color[rgb]{0,0,0}$a$}%
}}}}
\put(6076,-2536){\makebox(0,0)[lb]{\smash{{\SetFigFont{9}{10.8}{\rmdefault}{\mddefault}{\updefault}{\color[rgb]{0,0,0}$b+1$}%
}}}}
\put(6601,-3286){\makebox(0,0)[lb]{\smash{{\SetFigFont{9}{10.8}{\rmdefault}{\mddefault}{\updefault}{\color[rgb]{0,0,0}$b$}%
}}}}
\put(4951,-4561){\makebox(0,0)[lb]{\smash{{\SetFigFont{9}{10.8}{\rmdefault}{\mddefault}{\updefault}{\color[rgb]{0,0,0}$a$}%
}}}}
\put(5851,-4561){\makebox(0,0)[rb]{\smash{{\SetFigFont{9}{10.8}{\rmdefault}{\mddefault}{\updefault}{\color[rgb]{0,0,0}$a+1$}%
}}}}
\put(9001,-3436){\makebox(0,0)[b]{\smash{{\SetFigFont{9}{10.8}{\rmdefault}{\mddefault}{\updefault}{\color[rgb]{0,0,0}$\bigkreis{D}'$}%
}}}}
\put(9001,-1411){\makebox(0,0)[b]{\smash{{\SetFigFont{9}{10.8}{\rmdefault}{\mddefault}{\updefault}{\color[rgb]{0,0,0}$\bigkreis{A}'$}%
}}}}
\put(8551,-1786){\makebox(0,0)[lb]{\smash{{\SetFigFont{9}{10.8}{\rmdefault}{\mddefault}{\updefault}{\color[rgb]{0,0,0}$\smash{e_2'}$}%
}}}}
\put(8551,-1036){\makebox(0,0)[lb]{\smash{{\SetFigFont{9}{10.8}{\rmdefault}{\mddefault}{\updefault}{\color[rgb]{0,0,0}$\smash{e_1'}$}%
}}}}
\put(8551,-2386){\makebox(0,0)[lb]{\smash{{\SetFigFont{9}{10.8}{\rmdefault}{\mddefault}{\updefault}{\color[rgb]{0,0,0}$\smash{e_3'}$}%
}}}}
\put(9826,-4036){\makebox(0,0)[lb]{\smash{{\SetFigFont{9}{10.8}{\rmdefault}{\mddefault}{\updefault}{\color[rgb]{0,0,0}$e_5'$}%
}}}}
\put(8176,-4036){\makebox(0,0)[rb]{\smash{{\SetFigFont{9}{10.8}{\rmdefault}{\mddefault}{\updefault}{\color[rgb]{0,0,0}$\smash{e_4'}$}%
}}}}
\put(8326,-1186){\makebox(0,0)[rb]{\smash{{\SetFigFont{9}{10.8}{\rmdefault}{\mddefault}{\updefault}{\color[rgb]{0,0,0}$b+1$}%
}}}}
\put(8326,-1636){\makebox(0,0)[rb]{\smash{{\SetFigFont{9}{10.8}{\rmdefault}{\mddefault}{\updefault}{\color[rgb]{0,0,0}$b$}%
}}}}
\put(8326,-2536){\makebox(0,0)[rb]{\smash{{\SetFigFont{9}{10.8}{\rmdefault}{\mddefault}{\updefault}{\color[rgb]{0,0,0}$a+2$}%
}}}}
\put(7801,-3286){\makebox(0,0)[rb]{\smash{{\SetFigFont{9}{10.8}{\rmdefault}{\mddefault}{\updefault}{\color[rgb]{0,0,0}$a+1$}%
}}}}
\put(9676,-1186){\makebox(0,0)[lb]{\smash{{\SetFigFont{9}{10.8}{\rmdefault}{\mddefault}{\updefault}{\color[rgb]{0,0,0}$a+1$}%
}}}}
\put(9676,-1636){\makebox(0,0)[lb]{\smash{{\SetFigFont{9}{10.8}{\rmdefault}{\mddefault}{\updefault}{\color[rgb]{0,0,0}$a+2$}%
}}}}
\put(9676,-2536){\makebox(0,0)[lb]{\smash{{\SetFigFont{9}{10.8}{\rmdefault}{\mddefault}{\updefault}{\color[rgb]{0,0,0}$b$}%
}}}}
\put(10201,-3286){\makebox(0,0)[lb]{\smash{{\SetFigFont{9}{10.8}{\rmdefault}{\mddefault}{\updefault}{\color[rgb]{0,0,0}$b+1$}%
}}}}
\put(8476,-4561){\makebox(0,0)[lb]{\smash{{\SetFigFont{9}{10.8}{\rmdefault}{\mddefault}{\updefault}{\color[rgb]{0,0,0}$a+2$}%
}}}}
\put(5200,-4900){\makebox(0,0)[lb]{\smash{{\SetFigFont{9}{10.8}{\rmdefault}{\mddefault}{\updefault}{\color[rgb]{0,0,0}(a)}%
}}}}
\put(8800,-4900){\makebox(0,0)[lb]{\smash{{\SetFigFont{9}{10.8}{\rmdefault}{\mddefault}{\updefault}{\color[rgb]{0,0,0}(b)}%
}}}}
\end{picture}%
\caption{The two possible trigons bounded by a forked extended $S2$ cycle of $G_S^{a+1}$}
\label{fig:forkedbigona+1}
\end{figure}
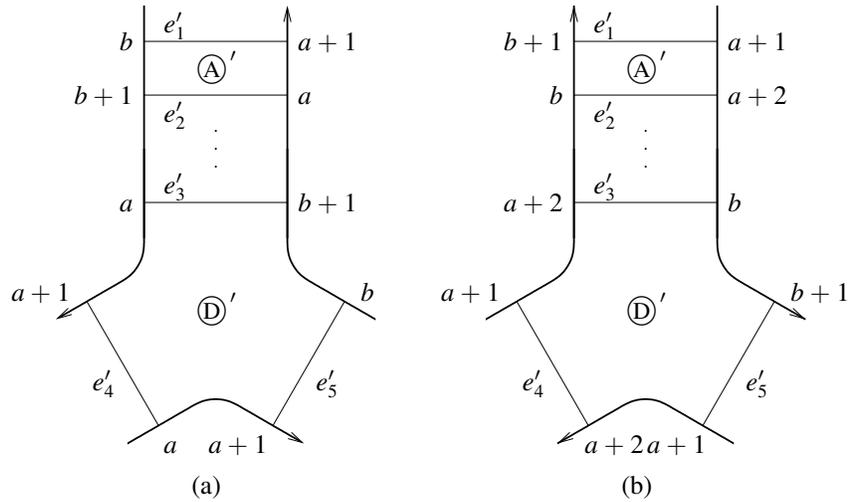

\begin{lemma}\label{thinningtwoforks}
Let $g$ be the trigon bounded by a forked extended $S2$ cycle of $\smash{G_S^1}$ with edges as in \fullref{fig:forkedbigon}(a).  Let $f$ be the trigon bounded by a forked extended $S2$ cycle of $G_S^{a+1}$.  Its edges must be as in \fullref{fig:forkedbigona+1} (a) or (b).  Let $\Delta_g$ be the bigon on $\hatT$ bounded by the edges $e_1$ and $e_5$ of $g$, and let $\Delta_f$ be the bigon on $\hatT$ bounded by the edges $\smash{e_1'}$ and $e_5'$ of $f$.  If the interiors of $\Delta_g$ and $\Delta_f$ are disjoint from $K$, then the edge $\smash{e_4'}$ has label pair $\{a, a{+}1\}$.
\end{lemma}

\begin{proof}
Assume the interiors of $\Delta_f$ and $\Delta_g$ are disjoint from $K$.  By the minimality assumption on $|S \cap T|$, we have $\Int (\Delta_f \cup \Delta_g) \cap (f \cup g) = \emptyset$.

Assume the edge $\smash{e_4'}$ of $f$ does not have label pair $\{a, a{+}1\}$.  Therefore it must have label pair $\{a{+}1, a{+}2\}$ and $f$ is as in \fullref{fig:forkedbigona+1}(b).  Assemble disks $D_g$ and $D_f$ in a manner similar to what is done in the proof of \fullref{thinningforkplusbigon}.  The disk $D_g$ is comprised of $\kreis{A}$, $\kreis{D}$, $\Delta_g$, a rectangle $\rho_{(a,\,a+1)}$ on $\bdry H_{\smash{(a,\,a{+}1)}}$, and rectangle $\rho_{(1, 2)}$ on $\bdry H_{(1,2)}$.  The disk $D_f$ is comprised of $\kreis{A}'$, $\kreis{D}'$, $\Delta_f$, a rectangle $\rho_{(a+1,\,a+2)}$ on $\bdry H_{(a+1,\,a+2)}$, and a rectangle $\rho_{(b,\,b+1)}$ on $\bdry H_{(b,\,b+1)}$.  These two disks are shown in \fullref{fig:twoforkshighlowdisks}.
\begin{figure}[ht!]
\centering
\begin{picture}(0,0)%
\includegraphics{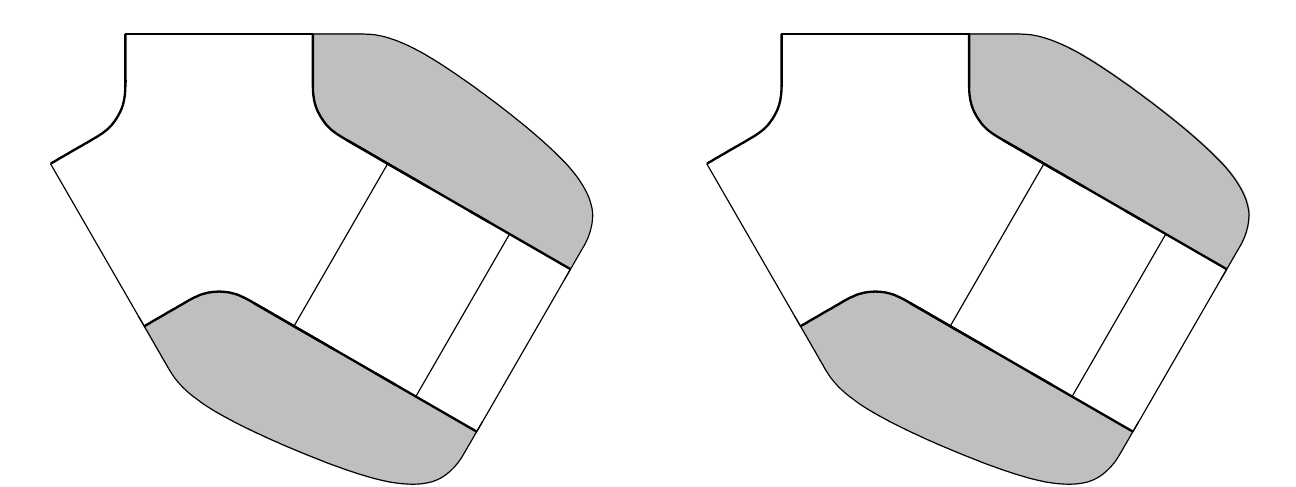}%
\end{picture}%
\setlength{\unitlength}{2960sp}%
\begingroup\makeatletter\ifx\SetFigFont\undefined%
\gdef\SetFigFont#1#2#3#4#5{%
  \reset@font\fontsize{#1}{#2pt}%
  \fontfamily{#3}\fontseries{#4}\fontshape{#5}%
  \selectfont}%
\fi\endgroup%
\begin{picture}(8278,3100)(3999,-5354)
\put(4276,-3211){\makebox(0,0)[rb]{\smash{{\SetFigFont{9}{10.8}{\rmdefault}{\mddefault}{\updefault}{\color[rgb]{0,0,0}$1$}%
}}}}
\put(6526,-3211){\makebox(0,0)[lb]{\smash{{\SetFigFont{9}{10.8}{\rmdefault}{\mddefault}{\updefault}{\color[rgb]{0,0,0}$a+1$}%
}}}}
\put(5851,-4561){\makebox(0,0)[rb]{\smash{{\SetFigFont{9}{10.8}{\rmdefault}{\mddefault}{\updefault}{\color[rgb]{0,0,0}$1$}%
}}}}
\put(4951,-2686){\makebox(0,0)[lb]{\smash{{\SetFigFont{9}{10.8}{\rmdefault}{\mddefault}{\updefault}{\color[rgb]{0,0,0}$e_3$}%
}}}}
\put(6001,-3961){\makebox(0,0)[rb]{\smash{{\SetFigFont{9}{10.8}{\rmdefault}{\mddefault}{\updefault}{\color[rgb]{0,0,0}$e_5$}%
}}}}
\put(4801,-3961){\makebox(0,0)[lb]{\smash{{\SetFigFont{9}{10.8}{\rmdefault}{\mddefault}{\updefault}{\color[rgb]{0,0,0}$e_4$}%
}}}}
\put(5401,-3436){\makebox(0,0)[b]{\smash{{\SetFigFont{9}{10.8}{\rmdefault}{\mddefault}{\updefault}{\color[rgb]{0,0,0}$\bigkreis{D}$}%
}}}}
\put(6601,-5011){\makebox(0,0)[rb]{\smash{{\SetFigFont{9}{10.8}{\rmdefault}{\mddefault}{\updefault}{\color[rgb]{0,0,0}$1$}%
}}}}
\put(7276,-3586){\makebox(0,0)[lb]{\smash{{\SetFigFont{9}{10.8}{\rmdefault}{\mddefault}{\updefault}{\color[rgb]{0,0,0}$a+1$}%
}}}}
\put(4726,-2536){\makebox(0,0)[rb]{\smash{{\SetFigFont{9}{10.8}{\rmdefault}{\mddefault}{\updefault}{\color[rgb]{0,0,0}$2$}%
}}}}
\put(7201,-5011){\makebox(0,0)[lb]{\smash{{\SetFigFont{9}{10.8}{\rmdefault}{\mddefault}{\updefault}{\color[rgb]{0,0,0}$e_2$}%
}}}}
\put(6751,-4786){\makebox(0,0)[lb]{\smash{{\SetFigFont{9}{10.8}{\rmdefault}{\mddefault}{\updefault}{\color[rgb]{0,0,0}$e_1$}%
}}}}
\put(7201,-4411){\makebox(0,0)[b]{\smash{{\SetFigFont{9}{10.8}{\rmdefault}{\mddefault}{\updefault}{\color[rgb]{0,0,0}$\bigkreis{A}$}%
}}}}
\put(6076,-2386){\makebox(0,0)[lb]{\smash{{\SetFigFont{9}{10.8}{\rmdefault}{\mddefault}{\updefault}{\color[rgb]{0,0,0}$a$}%
}}}}
\put(7801,-3961){\makebox(0,0)[lb]{\smash{{\SetFigFont{9}{10.8}{\rmdefault}{\mddefault}{\updefault}{\color[rgb]{0,0,0}$a$}%
}}}}
\put(6901,-5161){\makebox(0,0)[rb]{\smash{{\SetFigFont{9}{10.8}{\rmdefault}{\mddefault}{\updefault}{\color[rgb]{0,0,0}$2$}%
}}}}
\put(5101,-4486){\makebox(0,0)[lb]{\smash{{\SetFigFont{9}{10.8}{\rmdefault}{\mddefault}{\updefault}{\color[rgb]{0,0,0}$2$}%
}}}}
\put(6601,-4036){\makebox(0,0)[b]{\smash{{\SetFigFont{9}{10.8}{\rmdefault}{\mddefault}{\updefault}{\color[rgb]{0,0,0}$\Delta_g$}%
}}}}
\put(8476,-3211){\makebox(0,0)[rb]{\smash{{\SetFigFont{9}{10.8}{\rmdefault}{\mddefault}{\updefault}{\color[rgb]{0,0,0}$a+1$}%
}}}}
\put(10726,-3211){\makebox(0,0)[lb]{\smash{{\SetFigFont{9}{10.8}{\rmdefault}{\mddefault}{\updefault}{\color[rgb]{0,0,0}$b+1$}%
}}}}
\put(9151,-2686){\makebox(0,0)[lb]{\smash{{\SetFigFont{9}{10.8}{\rmdefault}{\mddefault}{\updefault}{\color[rgb]{0,0,0}$\smash{e_3'}$}%
}}}}
\put(10201,-3961){\makebox(0,0)[rb]{\smash{{\SetFigFont{9}{10.8}{\rmdefault}{\mddefault}{\updefault}{\color[rgb]{0,0,0}$e_5'$}%
}}}}
\put(9001,-3961){\makebox(0,0)[lb]{\smash{{\SetFigFont{9}{10.8}{\rmdefault}{\mddefault}{\updefault}{\color[rgb]{0,0,0}$\smash{e_4'}$}%
}}}}
\put(9601,-3436){\makebox(0,0)[b]{\smash{{\SetFigFont{9}{10.8}{\rmdefault}{\mddefault}{\updefault}{\color[rgb]{0,0,0}$\bigkreis{D}'$}%
}}}}
\put(10801,-5011){\makebox(0,0)[rb]{\smash{{\SetFigFont{9}{10.8}{\rmdefault}{\mddefault}{\updefault}{\color[rgb]{0,0,0}$a+1$}%
}}}}
\put(11476,-3586){\makebox(0,0)[lb]{\smash{{\SetFigFont{9}{10.8}{\rmdefault}{\mddefault}{\updefault}{\color[rgb]{0,0,0}$b+1$}%
}}}}
\put(8926,-2536){\makebox(0,0)[rb]{\smash{{\SetFigFont{9}{10.8}{\rmdefault}{\mddefault}{\updefault}{\color[rgb]{0,0,0}$a+2$}%
}}}}
\put(11401,-5011){\makebox(0,0)[lb]{\smash{{\SetFigFont{9}{10.8}{\rmdefault}{\mddefault}{\updefault}{\color[rgb]{0,0,0}$\smash{e_2'}$}%
}}}}
\put(10951,-4786){\makebox(0,0)[lb]{\smash{{\SetFigFont{9}{10.8}{\rmdefault}{\mddefault}{\updefault}{\color[rgb]{0,0,0}$\smash{e_1'}$}%
}}}}
\put(11401,-4411){\makebox(0,0)[b]{\smash{{\SetFigFont{9}{10.8}{\rmdefault}{\mddefault}{\updefault}{\color[rgb]{0,0,0}$\bigkreis{A}'$}%
}}}}
\put(10276,-2386){\makebox(0,0)[lb]{\smash{{\SetFigFont{9}{10.8}{\rmdefault}{\mddefault}{\updefault}{\color[rgb]{0,0,0}$b$}%
}}}}
\put(12001,-3961){\makebox(0,0)[lb]{\smash{{\SetFigFont{9}{10.8}{\rmdefault}{\mddefault}{\updefault}{\color[rgb]{0,0,0}$b$}%
}}}}
\put(11101,-5161){\makebox(0,0)[rb]{\smash{{\SetFigFont{9}{10.8}{\rmdefault}{\mddefault}{\updefault}{\color[rgb]{0,0,0}$a+2$}%
}}}}
\put(10801,-4036){\makebox(0,0)[b]{\smash{{\SetFigFont{9}{10.8}{\rmdefault}{\mddefault}{\updefault}{\color[rgb]{0,0,0}$\Delta_f$}%
}}}}
\put(9001,-5161){\makebox(0,0)[b]{\smash{{\SetFigFont{9}{10.8}{\rmdefault}{\mddefault}{\updefault}{\color[rgb]{0,0,0}$D_f$}%
}}}}
\put(4801,-5161){\makebox(0,0)[b]{\smash{{\SetFigFont{9}{10.8}{\rmdefault}{\mddefault}{\updefault}{\color[rgb]{0,0,0}$D_g$}%
}}}}
\put(10126,-4561){\makebox(0,0)[rb]{\smash{{\SetFigFont{9}{10.8}{\rmdefault}{\mddefault}{\updefault}{\color[rgb]{0,0,0}$a+1$}%
}}}}
\put(9301,-4411){\makebox(0,0)[lb]{\smash{{\SetFigFont{9}{10.8}{\rmdefault}{\mddefault}{\updefault}{\color[rgb]{0,0,0}$a+2$}%
}}}}
\end{picture}%
\caption{The construction of a high disk and a low disk}
\label{fig:twoforkshighlowdisks}
\end{figure}
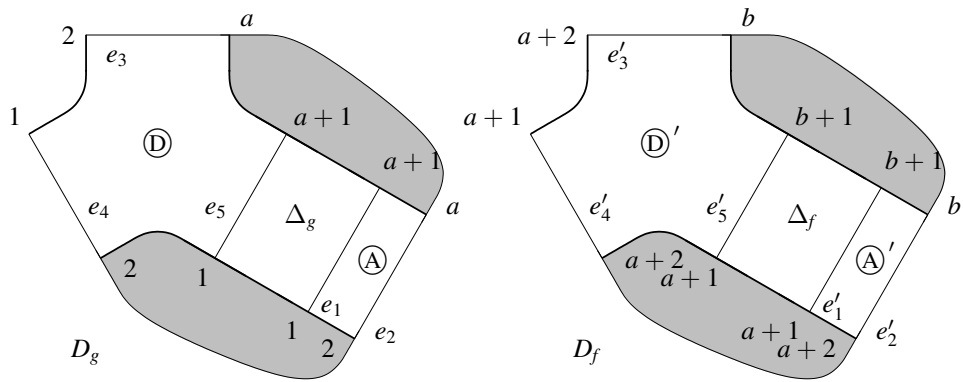
Note that the arcs $K_{(1,2)}$ and $K_{(a+1,\,a+2)}$ lie on opposite sides of $\hatT$.  We may extend arcs of the boundaries of the two disks $D_g$ and $D_f$  radially into $H_{(1,2)}$ and $H_{(a+1,\,a+2)}$ so that $\bdry D_g = K_{(1,2)} \cup (D_g \cap \hatT)$ and $\bdry D_f = K_{(a+1,\,a+2)} \cup (D_f \cap \hatT)$.  Even if $b=2$, one may check that the arcs $\bdry D_g \cap \hatT$ and $\bdry D_f \cap \hatT$ are disjoint.  Therefore the pair of disks $D_g$ and $D_f$ are a pair of disjoint high and low disks for $K$.  This is a contradiction to \fullref{highdisklowdisk}.

Therefore $f$ is as in \fullref{fig:forkedbigona+1}(a), and the edge $\smash{e_4'}$ has label pair $\{a, a{+}1\}$.
\end{proof}


\begin{lemma} \label{lem:shortpairsofforks}
Assume for $i \neq j$ that $\sigma_i \subseteq G_S^i$ and $\sigma_j \subseteq G_S^j$ are two forked extended $S2$ cycles such that for each (extended) $S2$ cycle contained in the face bounded by one of $\sigma_i$ or $\sigma_j$ there is an (extended) $S2$ cycle with the same label pair contained in the face bounded by the other forked extended $S2$ cycle.  Then the faces bounded by $\sigma_i$ and $\sigma_j$ each contain only an $S2$ cycle and no extended $S2$ cycles.
\end{lemma}

\begin{proof}
Assume $\smash{G_S^1}$ has a forked extended $S2$ cycle $\sigma_1$ that bounds a trigon $g$ on $\smash{G_S^1}$ as in \fullref{fig:forkedbigon}(a).  By \fullref{forkededges} the five edges $e_1$, $e_2$, $e_3$, $e_4$, and $e_5$ lie on $\hatT$ as in \fullref{fig:forkedbigon}(b).  For each $n = 0, \dots, (a-3)/2$ there is an (extended) $S2$ cycle contained in $g$ with label pair $\{2+n,a-n\}$.  

Assume there is a forked extended $S2$ cycle $\sigma_j$ of $G_S^j$ bounding a trigon $f$ such that for each $n = 0, \dots, (a-3)/2$ there is an (extended) $S2$ cycle contained in $f$ with label pair $\{2+n,a-n\}$.  By hypothesis $j \neq 1$.  Hence $j = a+1$.  Therefore $f$ appears as in \fullref{fig:forkedbigona+1}(a) with $b=1$.

We now consider how the edges of $\kreis{A}'$ and $\kreis{D}'$ lie with respect to the edges of $g$ as shown in \fullref{fig:forkedbigon}(b).  Note that since $\kreis{A} \neq \kreis{D}'$ and $\kreis{D} \neq \kreis{A}'$ the edges $\smash{e_1'}, \dots, e_5'$ each must be distinct from the edges $e_1, \dots, e_5$.  Also, by \fullref{lem:twoS2cycleswithsamelabelpair} and \fullref{lieinanannulus}, the edges $e_2$, $e_3$, $\smash{e_2'}$, and $\smash{e_3'}$ together lie in an essential annulus.

The edges $e_1$ and $\smash{e_1'}$ have the same label pair $\{1, a{+}1\}$.  Either they lie in a disk or they lie in an essential annulus in $\hatT$.  

{\bf Case 1}\qua If $e_1$ and $\smash{e_1'}$ lie in an essential annulus, then either (a)  $\smash{e_1'}$ is incident to $U_1$ in the arc of $\bdry U_1$ between $e_4$ and $e_5$ that does not intersect $e_1$ or (b) $\smash{e_1'}$ is incident to $U_1$ in the arc of $\bdry U_1$ between $e_4$ and $e_1$ that does not intersect $e_5$. 

For (a), the corners of $\kreis{A}'$ force $\smash{e_2'}$ to be incident to the vertices $U_2$ and $U_a$ on opposite sides of $e_2 \cup e_3 \cup U_2 \cup U_a$.  Thus $\smash{e_2'}$ cannot not lie in the essential annulus which contains $e_2$ and $e_3$ contrary to \fullref{lem:twoS2cycleswithsamelabelpair}.  

For (b), the corners of $\kreis{A}'$ force $\smash{e_2'}$ to be incident to the vertices $U_2$ and $U_a$ on the side of $e_2 \cup e_3 \cup U_2 \cup U_a$ to which $e_4$ is not incident.  By \fullref{lem:twoS2cycleswithsamelabelpair} the edge $\smash{e_3'}$ must lie on the other side of $e_2 \cup e_3 \cup U_2 \cup U_a$.  Following the $(a, a{+}1)$ corner of $\kreis{D}'$ from $\smash{e_3'}$ to $\smash{e_4'}$ implies that $\smash{e_4'}$ lies in the disk of $\hatT$ between $e_1$ and $e_5$.  This however contradicts that $\smash{e_4'}$ has label pair $\{a, a{+}1\}$.

Therefore the edges $e_1$ and $\smash{e_1'}$ cannot lie in an essential annulus.

{\bf Case 2}\qua  The edges $e_1$ and $\smash{e_1'}$ lie in a disk.  There are three possibilities for the placement of $\smash{e_1'}$ with respect to the edges of $g$.  Either (a) $e_5$ lies in the disk between $e_1$ and $\smash{e_1'}$, (b) $e_1$ lies in the disk between $\smash{e_1'}$ and $e_5$ or (c) $\smash{e_1'}$ lies in the disk between $e_1$ and $e_5$.

For (a), edge $\smash{e_1'}$ is incident to $U_1$ and $U_{a+1}$ in the same arcs as in Case 1(a).  Thus again edge $\smash{e_2'}$ cannot lie in the essential annulus which contains $e_2$ and $e_3$ contrary to \fullref{lem:twoS2cycleswithsamelabelpair}.

For (b), edge $\smash{e_1'}$ is incident to $U_1$ and $U_{a+1}$ in the same arcs as in case 2(a).  Thus again edge $\smash{e_4'}$ is forced to have the wrong label pair.

For (c), the corners of $\kreis{A}'$ force $\smash{e_2'}$ to be incident to the vertices $U_2$ and $U_a$ on the side of $e_2 \cup e_3 \cup U_2 \cup U_a$ to which $e_4$ is incident.  

The edges $e_2$ and $\smash{e_2'}$ must lie in a disk.  Otherwise $e_2$ and $\smash{e_2'}$ lie in an essential annulus with core parallel to the core of the essential annulus in which the edges $e_2$ and $e_3$ lie.  Then on $\hatT$, $e_1$ and $\smash{e_1'}$ bound a bigon $\Delta$.  With the appropriate rectangles $\rho_{(1,2)}$ and $\rho_{(a,\,a+1)}$ on $\bdry H_{(1,2)}$ and $\bdry H_{(a,\,a+1)}$ respectively we may form the disk $D = \kreis{A} \cup \kreis{A}' \cup \rho_{(1, 2)} \cup \rho_{(a,\,a+1)} \cup \Delta$ whose boundary is parallel to the core of the essential annulus in which the edges $e_2$ and $\smash{e_2'}$ lie.  Then by \fullref{lieinanannulus} $\bdry D$ is parallel to the core of the annulus in which an $S2$ cycle lies.  This contradicts \fullref{GT:L2.1}.

Since the edges $e_2$ and $\smash{e_2'}$ lie in a disk, $\smash{e_2'}$ is incident to $U_2$ in the arc of $\bdry U_2$ between $e_2$ and $e_4$ that does not contain $e_3$.  Therefore $\smash{e_3'}$ is incident to $U_2$ and $U_a$ in the arcs of $\bdry U_2$ and $\bdry U_a$ between $e_2$ and $e_3$ in which $\smash{e_2'}$ is not incident.  Following the corners of $\kreis{D}'$ puts edge $\smash{e_4'}$ incident to $U_a$ on the same side of $e_2$ and $e_3$ as $\smash{e_3'}$ and edge $e_5'$ so that $e_1$ lies in the disk between $\smash{e_1'}$ and $e_5'$.  This configuration is shown in \fullref{fig:oneconfiguration}.  Let $\Gamma$ be this subgraph of $G_T$ consisting of the edges $e_1, \dots, e_5$ and $e_1, \dots, e_5$ and the vertices $U_1$, $U_2$, $U_a$, and $U_{a+1}$.
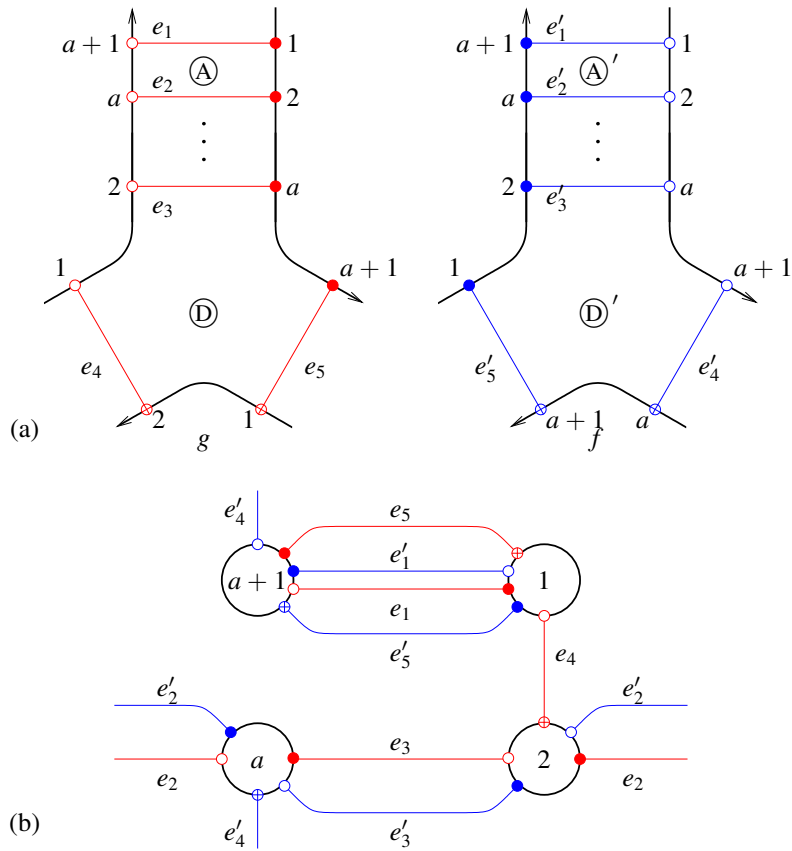
\begin{figure}[ht!]
\centering
\begin{picture}(0,0)%
\includegraphics{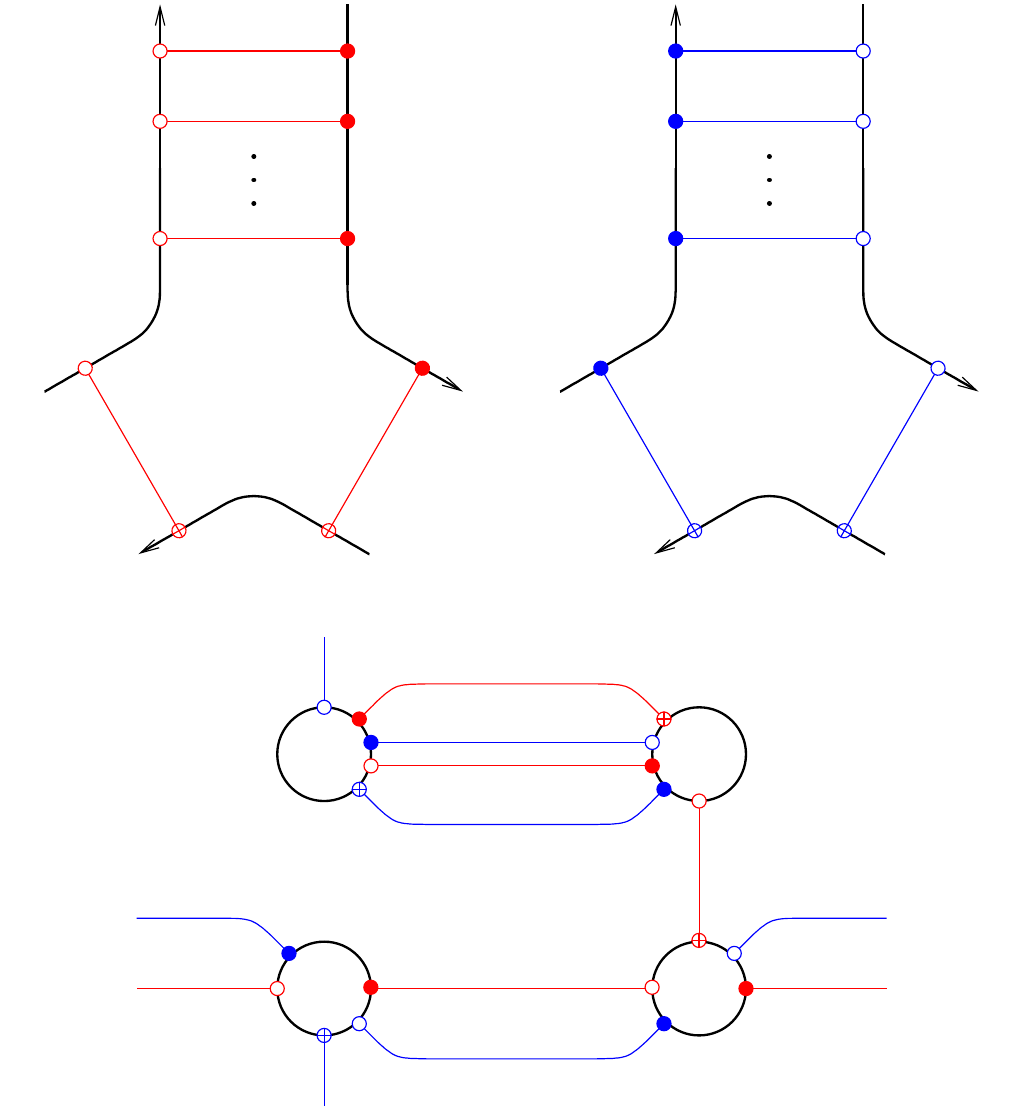}%
\end{picture}%
\setlength{\unitlength}{2960sp}%
\begingroup\makeatletter\ifx\SetFigFont\undefined%
\gdef\SetFigFont#1#2#3#4#5{%
  \reset@font\fontsize{#1}{#2pt}%
  \fontfamily{#3}\fontseries{#4}\fontshape{#5}%
  \selectfont}%
\fi\endgroup%
\begin{picture}(6525,7084)(3777,-8023)
\put(8701,-1561){\makebox(0,0)[b]{\smash{{\SetFigFont{9}{10.8}{\rmdefault}{\mddefault}{\updefault}{\color[rgb]{0,0,0}$\bigkreis{A}'$}%
}}}}
\put(9526,-4036){\makebox(0,0)[lb]{\smash{{\SetFigFont{9}{10.8}{\rmdefault}{\mddefault}{\updefault}{\color[rgb]{0,0,0}$\smash{e_4'}$}%
}}}}
\put(7876,-4036){\makebox(0,0)[rb]{\smash{{\SetFigFont{9}{10.8}{\rmdefault}{\mddefault}{\updefault}{\color[rgb]{0,0,0}$e_5'$}%
}}}}
\put(4576,-4036){\makebox(0,0)[rb]{\smash{{\SetFigFont{9}{10.8}{\rmdefault}{\mddefault}{\updefault}{\color[rgb]{0,0,0}$e_4$}%
}}}}
\put(8251,-2611){\makebox(0,0)[lb]{\smash{{\SetFigFont{9}{10.8}{\rmdefault}{\mddefault}{\updefault}{\color[rgb]{0,0,0}$\smash{e_3'}$}%
}}}}
\put(5851,-4486){\makebox(0,0)[rb]{\smash{{\SetFigFont{9}{10.8}{\rmdefault}{\mddefault}{\updefault}{\color[rgb]{0,0,0}$1$}%
}}}}
\put(4276,-3211){\makebox(0,0)[rb]{\smash{{\SetFigFont{9}{10.8}{\rmdefault}{\mddefault}{\updefault}{\color[rgb]{0,0,0}$1$}%
}}}}
\put(4726,-2536){\makebox(0,0)[rb]{\smash{{\SetFigFont{9}{10.8}{\rmdefault}{\mddefault}{\updefault}{\color[rgb]{0,0,0}$2$}%
}}}}
\put(6076,-2536){\makebox(0,0)[lb]{\smash{{\SetFigFont{9}{10.8}{\rmdefault}{\mddefault}{\updefault}{\color[rgb]{0,0,0}$a$}%
}}}}
\put(6526,-3211){\makebox(0,0)[lb]{\smash{{\SetFigFont{9}{10.8}{\rmdefault}{\mddefault}{\updefault}{\color[rgb]{0,0,0}$a+1$}%
}}}}
\put(7576,-3211){\makebox(0,0)[rb]{\smash{{\SetFigFont{9}{10.8}{\rmdefault}{\mddefault}{\updefault}{\color[rgb]{0,0,0}$1$}%
}}}}
\put(8026,-2536){\makebox(0,0)[rb]{\smash{{\SetFigFont{9}{10.8}{\rmdefault}{\mddefault}{\updefault}{\color[rgb]{0,0,0}$2$}%
}}}}
\put(9376,-2536){\makebox(0,0)[lb]{\smash{{\SetFigFont{9}{10.8}{\rmdefault}{\mddefault}{\updefault}{\color[rgb]{0,0,0}$a$}%
}}}}
\put(9826,-3211){\makebox(0,0)[lb]{\smash{{\SetFigFont{9}{10.8}{\rmdefault}{\mddefault}{\updefault}{\color[rgb]{0,0,0}$a+1$}%
}}}}
\put(9151,-4486){\makebox(0,0)[rb]{\smash{{\SetFigFont{9}{10.8}{\rmdefault}{\mddefault}{\updefault}{\color[rgb]{0,0,0}$a$}%
}}}}
\put(6226,-4036){\makebox(0,0)[lb]{\smash{{\SetFigFont{9}{10.8}{\rmdefault}{\mddefault}{\updefault}{\color[rgb]{0,0,0}$e_5$}%
}}}}
\put(4951,-4486){\makebox(0,0)[lb]{\smash{{\SetFigFont{9}{10.8}{\rmdefault}{\mddefault}{\updefault}{\color[rgb]{0,0,0}$2$}%
}}}}
\put(8251,-4486){\makebox(0,0)[lb]{\smash{{\SetFigFont{9}{10.8}{\rmdefault}{\mddefault}{\updefault}{\color[rgb]{0,0,0}$a+1$}%
}}}}
\put(8251,-7336){\makebox(0,0)[b]{\smash{{\SetFigFont{9}{10.8}{\rmdefault}{\mddefault}{\updefault}{\color[rgb]{0,0,0}$2$}%
}}}}
\put(5851,-7336){\makebox(0,0)[b]{\smash{{\SetFigFont{9}{10.8}{\rmdefault}{\mddefault}{\updefault}{\color[rgb]{0,0,0}$a$}%
}}}}
\put(5851,-5836){\makebox(0,0)[b]{\smash{{\SetFigFont{9}{10.8}{\rmdefault}{\mddefault}{\updefault}{\color[rgb]{0,0,0}$a+1$}%
}}}}
\put(8251,-5836){\makebox(0,0)[b]{\smash{{\SetFigFont{9}{10.8}{\rmdefault}{\mddefault}{\updefault}{\color[rgb]{0,0,0}$1$}%
}}}}
\put(7051,-7186){\makebox(0,0)[b]{\smash{{\SetFigFont{9}{10.8}{\rmdefault}{\mddefault}{\updefault}{\color[rgb]{0,0,0}$e_3$}%
}}}}
\put(9001,-7486){\makebox(0,0)[b]{\smash{{\SetFigFont{9}{10.8}{\rmdefault}{\mddefault}{\updefault}{\color[rgb]{0,0,0}$e_2$}%
}}}}
\put(7051,-7936){\makebox(0,0)[b]{\smash{{\SetFigFont{9}{10.8}{\rmdefault}{\mddefault}{\updefault}{\color[rgb]{0,0,0}$\smash{e_3'}$}%
}}}}
\put(5101,-7486){\makebox(0,0)[b]{\smash{{\SetFigFont{9}{10.8}{\rmdefault}{\mddefault}{\updefault}{\color[rgb]{0,0,0}$e_2$}%
}}}}
\put(5776,-5236){\makebox(0,0)[rb]{\smash{{\SetFigFont{9}{10.8}{\rmdefault}{\mddefault}{\updefault}{\color[rgb]{0,0,0}$\smash{e_4'}$}%
}}}}
\put(5776,-7936){\makebox(0,0)[rb]{\smash{{\SetFigFont{9}{10.8}{\rmdefault}{\mddefault}{\updefault}{\color[rgb]{0,0,0}$\smash{e_4'}$}%
}}}}
\put(7051,-5236){\makebox(0,0)[b]{\smash{{\SetFigFont{9}{10.8}{\rmdefault}{\mddefault}{\updefault}{\color[rgb]{0,0,0}$e_5$}%
}}}}
\put(7051,-6436){\makebox(0,0)[b]{\smash{{\SetFigFont{9}{10.8}{\rmdefault}{\mddefault}{\updefault}{\color[rgb]{0,0,0}$e_5'$}%
}}}}
\put(7051,-6061){\makebox(0,0)[b]{\smash{{\SetFigFont{9}{10.8}{\rmdefault}{\mddefault}{\updefault}{\color[rgb]{0,0,0}$e_1$}%
}}}}
\put(7051,-5611){\makebox(0,0)[b]{\smash{{\SetFigFont{9}{10.8}{\rmdefault}{\mddefault}{\updefault}{\color[rgb]{0,0,0}$\smash{e_1'}$}%
}}}}
\put(5101,-6736){\makebox(0,0)[b]{\smash{{\SetFigFont{9}{10.8}{\rmdefault}{\mddefault}{\updefault}{\color[rgb]{0,0,0}$\smash{e_2'}$}%
}}}}
\put(9001,-6736){\makebox(0,0)[b]{\smash{{\SetFigFont{9}{10.8}{\rmdefault}{\mddefault}{\updefault}{\color[rgb]{0,0,0}$\smash{e_2'}$}%
}}}}
\put(8326,-6436){\makebox(0,0)[lb]{\smash{{\SetFigFont{9}{10.8}{\rmdefault}{\mddefault}{\updefault}{\color[rgb]{0,0,0}$e_4$}%
}}}}
\put(4951,-2686){\makebox(0,0)[lb]{\smash{{\SetFigFont{9}{10.8}{\rmdefault}{\mddefault}{\updefault}{\color[rgb]{0,0,0}$e_3$}%
}}}}
\put(5401,-4636){\makebox(0,0)[b]{\smash{{\SetFigFont{9}{10.8}{\rmdefault}{\mddefault}{\updefault}{\color[rgb]{0,0,0}$g$}%
}}}}
\put(8701,-4636){\makebox(0,0)[b]{\smash{{\SetFigFont{9}{10.8}{\rmdefault}{\mddefault}{\updefault}{\color[rgb]{0,0,0}$f$}%
}}}}
\put(3901,-4561){\makebox(0,0)[b]{\smash{{\SetFigFont{9}{10.8}{\rmdefault}{\mddefault}{\updefault}{\color[rgb]{0,0,0}(a)}%
}}}}
\put(3901,-7861){\makebox(0,0)[b]{\smash{{\SetFigFont{9}{10.8}{\rmdefault}{\mddefault}{\updefault}{\color[rgb]{0,0,0}(b)}%
}}}}
\put(5401,-1561){\makebox(0,0)[b]{\smash{{\SetFigFont{9}{10.8}{\rmdefault}{\mddefault}{\updefault}{\color[rgb]{0,0,0}$\bigkreis{A}$}%
}}}}
\put(5401,-3586){\makebox(0,0)[b]{\smash{{\SetFigFont{9}{10.8}{\rmdefault}{\mddefault}{\updefault}{\color[rgb]{0,0,0}$\bigkreis{D}$}%
}}}}
\put(8701,-3586){\makebox(0,0)[b]{\smash{{\SetFigFont{9}{10.8}{\rmdefault}{\mddefault}{\updefault}{\color[rgb]{0,0,0}$\bigkreis{D}'$}%
}}}}
\put(6076,-1336){\makebox(0,0)[lb]{\smash{{\SetFigFont{9}{10.8}{\rmdefault}{\mddefault}{\updefault}{\color[rgb]{0,0,0}$1$}%
}}}}
\put(6076,-1786){\makebox(0,0)[lb]{\smash{{\SetFigFont{9}{10.8}{\rmdefault}{\mddefault}{\updefault}{\color[rgb]{0,0,0}$2$}%
}}}}
\put(4726,-1786){\makebox(0,0)[rb]{\smash{{\SetFigFont{9}{10.8}{\rmdefault}{\mddefault}{\updefault}{\color[rgb]{0,0,0}$a$}%
}}}}
\put(4726,-1336){\makebox(0,0)[rb]{\smash{{\SetFigFont{9}{10.8}{\rmdefault}{\mddefault}{\updefault}{\color[rgb]{0,0,0}$a+1$}%
}}}}
\put(4951,-1636){\makebox(0,0)[lb]{\smash{{\SetFigFont{9}{10.8}{\rmdefault}{\mddefault}{\updefault}{\color[rgb]{0,0,0}$e_2$}%
}}}}
\put(4951,-1186){\makebox(0,0)[lb]{\smash{{\SetFigFont{9}{10.8}{\rmdefault}{\mddefault}{\updefault}{\color[rgb]{0,0,0}$e_1$}%
}}}}
\put(8251,-1186){\makebox(0,0)[lb]{\smash{{\SetFigFont{9}{10.8}{\rmdefault}{\mddefault}{\updefault}{\color[rgb]{0,0,0}$\smash{e_1'}$}%
}}}}
\put(9376,-1336){\makebox(0,0)[lb]{\smash{{\SetFigFont{9}{10.8}{\rmdefault}{\mddefault}{\updefault}{\color[rgb]{0,0,0}$1$}%
}}}}
\put(8026,-1786){\makebox(0,0)[rb]{\smash{{\SetFigFont{9}{10.8}{\rmdefault}{\mddefault}{\updefault}{\color[rgb]{0,0,0}$a$}%
}}}}
\put(8026,-1336){\makebox(0,0)[rb]{\smash{{\SetFigFont{9}{10.8}{\rmdefault}{\mddefault}{\updefault}{\color[rgb]{0,0,0}$a+1$}%
}}}}
\put(9376,-1786){\makebox(0,0)[lb]{\smash{{\SetFigFont{9}{10.8}{\rmdefault}{\mddefault}{\updefault}{\color[rgb]{0,0,0}$2$}%
}}}}
\put(8251,-1636){\makebox(0,0)[lb]{\smash{{\SetFigFont{9}{10.8}{\rmdefault}{\mddefault}{\updefault}{\color[rgb]{0,0,0}$\smash{e_2'}$}%
}}}}
\end{picture}%
\caption{(a) The trigons $g$ and $f$\qua  (b) The subgraph $\Gamma$ consisting of the edges $e_1, \dots, e_5, \smash{e_1'}, \dots, e_5'$ and the vertices $U_1$, $U_2$, $U_a$, and $U_{a+1}$ on $\hatT$}
\label{fig:oneconfiguration}
\end{figure}

  If $a > 3$ then there is an (extended) $S2$ cycle with label pair $\{3, a{-}1\}$ which by \fullref{lieinanannulus} lies in an essential annulus.  Thus $\Gamma$  must lie in an essential annulus.  This is a contradiction.  Thus $a=3$.  This proves the lemma.
\end{proof}


\begin{prop}\label{lieinannulusorboundbigon}
If a trigon $f$ of $\smash{G_S^x}$ is not innermost, then there exists two edges of $f$ that have the same label pair.  These edges either lie in an essential annulus on $\hatT$ or bound a bigon on $\hatT$ whose interior intersects $K$.  Moreover, if the edges of $f$ are not an extended $S3$ cycle, then $f$ contains a forked extended $S2$ cycle whose edge with label pair distinct from the other two is not an edge of $f$.
\end{prop}

\begin{proof}
Assume $f$ is a trigon of $\smash{G_S^x}$ that is not an innermost trigon.  Let $g$ be the trigon innermost on $f$.  Therefore $f$ is obtained from $g$ by attaching bigon faces of $G_S$.  By \fullref{innermosttrigons}, the edges of $g$ are an $S3$ cycle or a forked extended $S2$ cycle.   

If the edges of $g$ are an $S3$ cycle, then to obtain $f$ an equal number of bigons of $G_S$ must be attached to all three edges of $g$.  Hence the edges of $f$ are an extended $S3$ cycle and lie in an essential annulus on $\hatT$.

If the edges of $g$ are a forked extended $S2$ cycle, then relabel so that $g$ is a trigon of $\smash{G_S^1}$ as in \fullref{fig:forkedbigon}(a) with one edge $e_4$ having label pair $\{1, 2\}$ and the other two edges $e_1$ and $e_5$ having label pair $\{1, a{+}1\}$.  Note $a \neq 1$.  Then $f$ is obtained by attaching bigon faces of $G_S$ to one, two, or three edges of $g$.  Let $\Delta_g$ be the bigon on $\hatT$ bounded by $e_1$ and $e_5$ (and the vertices $U_1$ and $U_a$).   \fullref{thinningforkplusbigon} implies that the interior of $\Delta_g$ must intersect $K$ if a bigon of $G_S$ is incident to $e_4$.

{\bf Case I}\qua
If $f$ is obtained by attaching bigon faces of $G_S$ to just one edge of $g$, then the bigons are attached to $e_4$.  Otherwise edges with label pairs $\{1,2\}$ and $\{1, a{+}1\}$ must be on the boundary of $f$.  Since $a \neq 1$, this implies that $f$ is a trigon of $\smash{G_S^1}$, contradicting that $f$ properly contains the trigon $g$ of $\smash{G_S^1}$.

Attaching bigons to the edge $e_4$ of $g$ with label pair $\{1, 2\}$ implies that $f$ must be a trigon of $G_S^{a+1}$.  Furthermore, $e_1$ and $e_5$ are still edges of $f$.  These edges bound $\Delta_g$.  According to \fullref{thinningforkplusbigon} the interior of $\Delta_g$ must intersect $K$.  This satisfies the conclusion of the proposition.

{\bf Case II}\qua
If $f$ is obtained from $g$ by attaching bigons of $G_S$ to the two edges with label pair $\{1, a{+}1\}$, then $f$ must be a trigon of either $\smash{G_S^1}$ or $G_S^2$.  This cannot occur since $g$ is a trigon of $\smash{G_S^1}$ and its interior contains edges of $G_S^2$.

Thus $f$ must be obtained from $g$ by attaching bigons of $G_S$ to the edge with label pair $\{1,2\}$ and one of the edges with label pair $\{1, a{+}1\}$.  Because the third edge has label pair $\{1, a{+}1\}$,  $f$ must be a trigon of either $\smash{G_S^1}$ or $\smash{G_S^{a+1}}$.  It can be neither as the interior of $f$ will contain an edge of $g$ with label pair $\{1, a{+}1\}$.  Thus this case cannot occur.

{\bf Case III}\qua
Assume $f$ is obtained from $g$ by attaching bigons of $G_S$ to all three edges of $g$.  The number of bigons attached to each edge is not necessarily uniform.
\begin{figure}[ht!]
\centering
\begin{picture}(0,0)%
\includegraphics{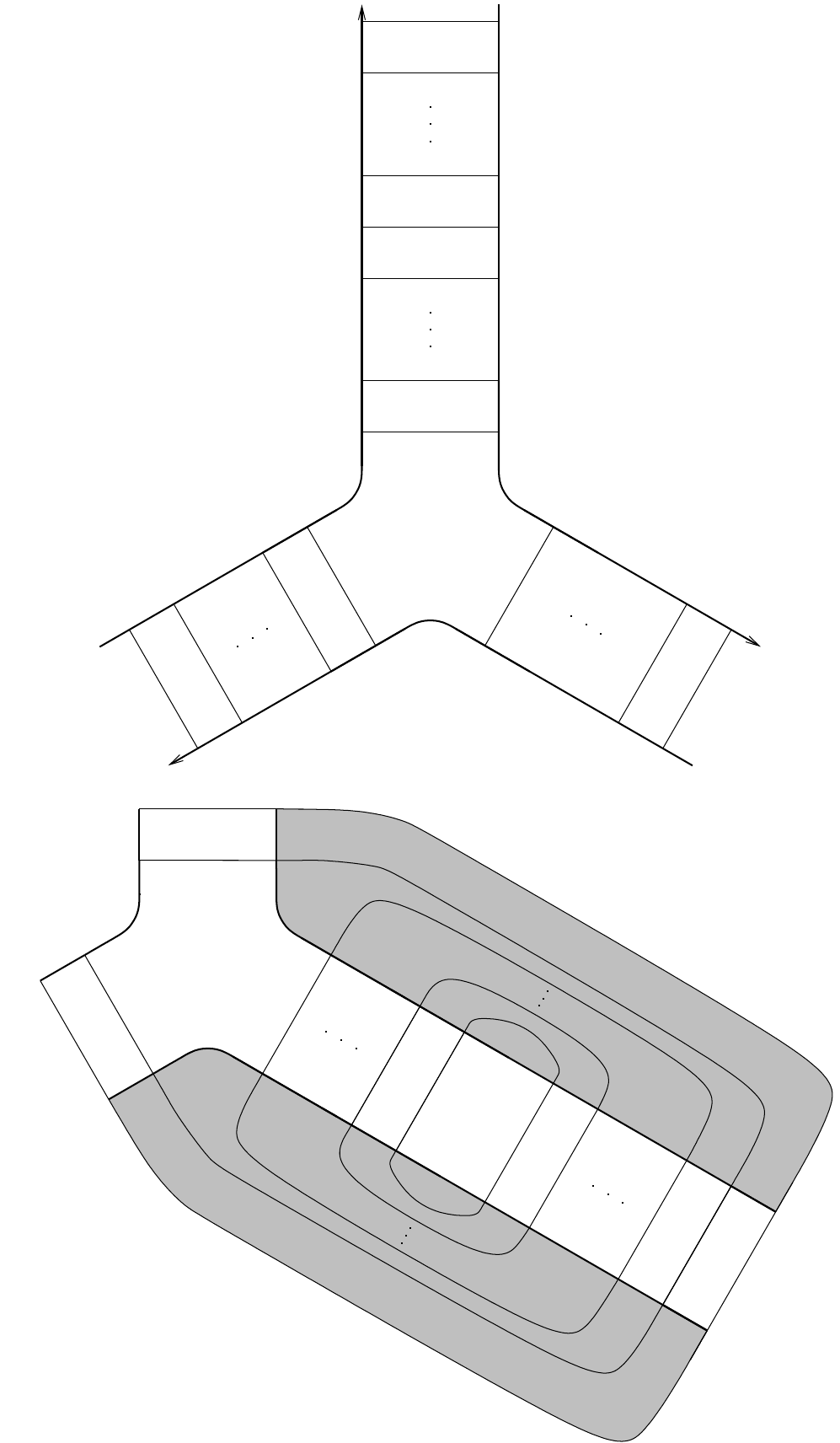}%
\end{picture}%
\setlength{\unitlength}{2565sp}%
\begingroup\makeatletter\ifx\SetFigFont\undefined%
\gdef\SetFigFont#1#2#3#4#5{%
  \reset@font\fontsize{#1}{#2pt}%
  \fontfamily{#3}\fontseries{#4}\fontshape{#5}%
  \selectfont}%
\fi\endgroup%
\begin{picture}(7352,12625)(1629,-11314)
\put(2326,-6961){\makebox(0,0)[rb]{\smash{{\SetFigFont{8}{9.6}{\rmdefault}{\mddefault}{\updefault}{\color[rgb]{0,0,0}$1$}%
}}}}
\put(2776,-6286){\makebox(0,0)[rb]{\smash{{\SetFigFont{8}{9.6}{\rmdefault}{\mddefault}{\updefault}{\color[rgb]{0,0,0}$2$}%
}}}}
\put(5851,-8536){\makebox(0,0)[b]{\smash{{\SetFigFont{8}{9.6}{\rmdefault}{\mddefault}{\updefault}{\color[rgb]{0,0,0}$\Delta_f$}%
}}}}
\put(5776,-7636){\makebox(0,0)[lb]{\smash{{\SetFigFont{8}{9.6}{\rmdefault}{\mddefault}{\updefault}{\color[rgb]{0,0,0}$a+b+1$}%
}}}}
\put(6526,-8086){\makebox(0,0)[lb]{\smash{{\SetFigFont{8}{9.6}{\rmdefault}{\mddefault}{\updefault}{\color[rgb]{0,0,0}$a+b+1$}%
}}}}
\put(5326,-7411){\makebox(0,0)[lb]{\smash{{\SetFigFont{8}{9.6}{\rmdefault}{\mddefault}{\updefault}{\color[rgb]{0,0,0}$a+b$}%
}}}}
\put(2551,-8536){\makebox(0,0)[rb]{\smash{{\SetFigFont{8}{9.6}{\rmdefault}{\mddefault}{\updefault}{\color[rgb]{0,0,0}$3$}%
}}}}
\put(3076,-8236){\makebox(0,0)[lb]{\smash{{\SetFigFont{8}{9.6}{\rmdefault}{\mddefault}{\updefault}{\color[rgb]{0,0,0}$2$}%
}}}}
\put(3826,-8236){\makebox(0,0)[rb]{\smash{{\SetFigFont{8}{9.6}{\rmdefault}{\mddefault}{\updefault}{\color[rgb]{0,0,0}$1$}%
}}}}
\put(4576,-8686){\makebox(0,0)[rb]{\smash{{\SetFigFont{8}{9.6}{\rmdefault}{\mddefault}{\updefault}{\color[rgb]{0,0,0}$t-b+2$}%
}}}}
\put(5026,-8911){\makebox(0,0)[rb]{\smash{{\SetFigFont{8}{9.6}{\rmdefault}{\mddefault}{\updefault}{\color[rgb]{0,0,0}$t-b+1$}%
}}}}
\put(6901,-10036){\makebox(0,0)[rb]{\smash{{\SetFigFont{8}{9.6}{\rmdefault}{\mddefault}{\updefault}{\color[rgb]{0,0,0}$1$}%
}}}}
\put(7351,-10261){\makebox(0,0)[rb]{\smash{{\SetFigFont{8}{9.6}{\rmdefault}{\mddefault}{\updefault}{\color[rgb]{0,0,0}$2$}%
}}}}
\put(7876,-10561){\makebox(0,0)[lb]{\smash{{\SetFigFont{8}{9.6}{\rmdefault}{\mddefault}{\updefault}{\color[rgb]{0,0,0}$3$}%
}}}}
\put(8176,-9061){\makebox(0,0)[lb]{\smash{{\SetFigFont{8}{9.6}{\rmdefault}{\mddefault}{\updefault}{\color[rgb]{0,0,0}$a$}%
}}}}
\put(7501,-8686){\makebox(0,0)[lb]{\smash{{\SetFigFont{8}{9.6}{\rmdefault}{\mddefault}{\updefault}{\color[rgb]{0,0,0}$a+1$}%
}}}}
\put(8551,-9361){\makebox(0,0)[lb]{\smash{{\SetFigFont{8}{9.6}{\rmdefault}{\mddefault}{\updefault}{\color[rgb]{0,0,0}$a-1$}%
}}}}
\put(6976,-8386){\makebox(0,0)[lb]{\smash{{\SetFigFont{8}{9.6}{\rmdefault}{\mddefault}{\updefault}{\color[rgb]{0,0,0}$a+b$}%
}}}}
\put(6376,-9736){\makebox(0,0)[rb]{\smash{{\SetFigFont{8}{9.6}{\rmdefault}{\mddefault}{\updefault}{\color[rgb]{0,0,0}$t-b+2$}%
}}}}
\put(6001,-9436){\makebox(0,0)[rb]{\smash{{\SetFigFont{8}{9.6}{\rmdefault}{\mddefault}{\updefault}{\color[rgb]{0,0,0}$t-b+1$}%
}}}}
\put(4651,-7036){\makebox(0,0)[lb]{\smash{{\SetFigFont{8}{9.6}{\rmdefault}{\mddefault}{\updefault}{\color[rgb]{0,0,0}$a+1$}%
}}}}
\put(4126,-6136){\makebox(0,0)[lb]{\smash{{\SetFigFont{8}{9.6}{\rmdefault}{\mddefault}{\updefault}{\color[rgb]{0,0,0}$a$}%
}}}}
\put(4051,-5686){\makebox(0,0)[lb]{\smash{{\SetFigFont{8}{9.6}{\rmdefault}{\mddefault}{\updefault}{\color[rgb]{0,0,0}$a-1$}%
}}}}
\put(2776,-5761){\makebox(0,0)[rb]{\smash{{\SetFigFont{8}{9.6}{\rmdefault}{\mddefault}{\updefault}{\color[rgb]{0,0,0}$3$}%
}}}}
\put(1876,-7261){\makebox(0,0)[rb]{\smash{{\SetFigFont{8}{9.6}{\rmdefault}{\mddefault}{\updefault}{\color[rgb]{0,0,0}$t$}%
}}}}
\put(7351,-5461){\makebox(0,0)[rb]{\smash{{\SetFigFont{8}{9.6}{\rmdefault}{\mddefault}{\updefault}{\color[rgb]{0,0,0}$t-b+1$}%
}}}}
\put(5401,-5386){\makebox(0,0)[b]{\smash{{\SetFigFont{8}{9.6}{\rmdefault}{\mddefault}{\updefault}{\color[rgb]{0,0,0}(a)}%
}}}}
\put(5401,-11161){\makebox(0,0)[b]{\smash{{\SetFigFont{8}{9.6}{\rmdefault}{\mddefault}{\updefault}{\color[rgb]{0,0,0}(b)}%
}}}}
\put(8101,-4111){\makebox(0,0)[lb]{\smash{{\SetFigFont{8}{9.6}{\rmdefault}{\mddefault}{\updefault}{\color[rgb]{0,0,0}$a+b+1$}%
}}}}
\put(7651,-3886){\makebox(0,0)[lb]{\smash{{\SetFigFont{8}{9.6}{\rmdefault}{\mddefault}{\updefault}{\color[rgb]{0,0,0}$a+b$}%
}}}}
\put(6601,-3211){\makebox(0,0)[lb]{\smash{{\SetFigFont{8}{9.6}{\rmdefault}{\mddefault}{\updefault}{\color[rgb]{0,0,0}$a+1$}%
}}}}
\put(6151,-2536){\makebox(0,0)[lb]{\smash{{\SetFigFont{8}{9.6}{\rmdefault}{\mddefault}{\updefault}{\color[rgb]{0,0,0}$a$}%
}}}}
\put(6151,-2086){\makebox(0,0)[lb]{\smash{{\SetFigFont{8}{9.6}{\rmdefault}{\mddefault}{\updefault}{\color[rgb]{0,0,0}$a-1$}%
}}}}
\put(6151,-1186){\makebox(0,0)[lb]{\smash{{\SetFigFont{8}{9.6}{\rmdefault}{\mddefault}{\updefault}{\color[rgb]{0,0,0}$3$}%
}}}}
\put(6151,-736){\makebox(0,0)[lb]{\smash{{\SetFigFont{8}{9.6}{\rmdefault}{\mddefault}{\updefault}{\color[rgb]{0,0,0}$2$}%
}}}}
\put(6151,-286){\makebox(0,0)[lb]{\smash{{\SetFigFont{8}{9.6}{\rmdefault}{\mddefault}{\updefault}{\color[rgb]{0,0,0}$1$}%
}}}}
\put(6151,614){\makebox(0,0)[lb]{\smash{{\SetFigFont{8}{9.6}{\rmdefault}{\mddefault}{\updefault}{\color[rgb]{0,0,0}$t-b+2$}%
}}}}
\put(6151,1064){\makebox(0,0)[lb]{\smash{{\SetFigFont{8}{9.6}{\rmdefault}{\mddefault}{\updefault}{\color[rgb]{0,0,0}$t-b+1$}%
}}}}
\put(4951,-4561){\makebox(0,0)[lb]{\smash{{\SetFigFont{8}{9.6}{\rmdefault}{\mddefault}{\updefault}{\color[rgb]{0,0,0}$2$}%
}}}}
\put(4576,-4786){\makebox(0,0)[lb]{\smash{{\SetFigFont{8}{9.6}{\rmdefault}{\mddefault}{\updefault}{\color[rgb]{0,0,0}$3$}%
}}}}
\put(3901,-5236){\makebox(0,0)[lb]{\smash{{\SetFigFont{8}{9.6}{\rmdefault}{\mddefault}{\updefault}{\color[rgb]{0,0,0}$b+1$}%
}}}}
\put(3451,-5461){\makebox(0,0)[lb]{\smash{{\SetFigFont{8}{9.6}{\rmdefault}{\mddefault}{\updefault}{\color[rgb]{0,0,0}$b+2$}%
}}}}
\put(6976,-5236){\makebox(0,0)[rb]{\smash{{\SetFigFont{8}{9.6}{\rmdefault}{\mddefault}{\updefault}{\color[rgb]{0,0,0}$t-b+2$}%
}}}}
\put(5851,-4561){\makebox(0,0)[rb]{\smash{{\SetFigFont{8}{9.6}{\rmdefault}{\mddefault}{\updefault}{\color[rgb]{0,0,0}$1$}%
}}}}
\put(2701,-4111){\makebox(0,0)[rb]{\smash{{\SetFigFont{8}{9.6}{\rmdefault}{\mddefault}{\updefault}{\color[rgb]{0,0,0}$t-b+1$}%
}}}}
\put(3151,-3886){\makebox(0,0)[rb]{\smash{{\SetFigFont{8}{9.6}{\rmdefault}{\mddefault}{\updefault}{\color[rgb]{0,0,0}$t-b+2$}%
}}}}
\put(3901,-3436){\makebox(0,0)[rb]{\smash{{\SetFigFont{8}{9.6}{\rmdefault}{\mddefault}{\updefault}{\color[rgb]{0,0,0}$t$}%
}}}}
\put(4201,-3211){\makebox(0,0)[rb]{\smash{{\SetFigFont{8}{9.6}{\rmdefault}{\mddefault}{\updefault}{\color[rgb]{0,0,0}$1$}%
}}}}
\put(4651,-2536){\makebox(0,0)[rb]{\smash{{\SetFigFont{8}{9.6}{\rmdefault}{\mddefault}{\updefault}{\color[rgb]{0,0,0}$2$}%
}}}}
\put(4651,-2086){\makebox(0,0)[rb]{\smash{{\SetFigFont{8}{9.6}{\rmdefault}{\mddefault}{\updefault}{\color[rgb]{0,0,0}$3$}%
}}}}
\put(4651,-1186){\makebox(0,0)[rb]{\smash{{\SetFigFont{8}{9.6}{\rmdefault}{\mddefault}{\updefault}{\color[rgb]{0,0,0}$a-1$}%
}}}}
\put(4651,-736){\makebox(0,0)[rb]{\smash{{\SetFigFont{8}{9.6}{\rmdefault}{\mddefault}{\updefault}{\color[rgb]{0,0,0}$a$}%
}}}}
\put(4651,1064){\makebox(0,0)[rb]{\smash{{\SetFigFont{8}{9.6}{\rmdefault}{\mddefault}{\updefault}{\color[rgb]{0,0,0}$a+b+1$}%
}}}}
\put(4651,614){\makebox(0,0)[rb]{\smash{{\SetFigFont{8}{9.6}{\rmdefault}{\mddefault}{\updefault}{\color[rgb]{0,0,0}$a+b$}%
}}}}
\put(4651,-286){\makebox(0,0)[rb]{\smash{{\SetFigFont{8}{9.6}{\rmdefault}{\mddefault}{\updefault}{\color[rgb]{0,0,0}$a+1$}%
}}}}
\end{picture}%
\caption{(a) A forked extended $S2$ cycle with $b$ bigons attached to each edge of the face $g$ it bounds \qua (b) The construction of a long disk}
\label{fig:3extforkedbigon}
\end{figure}
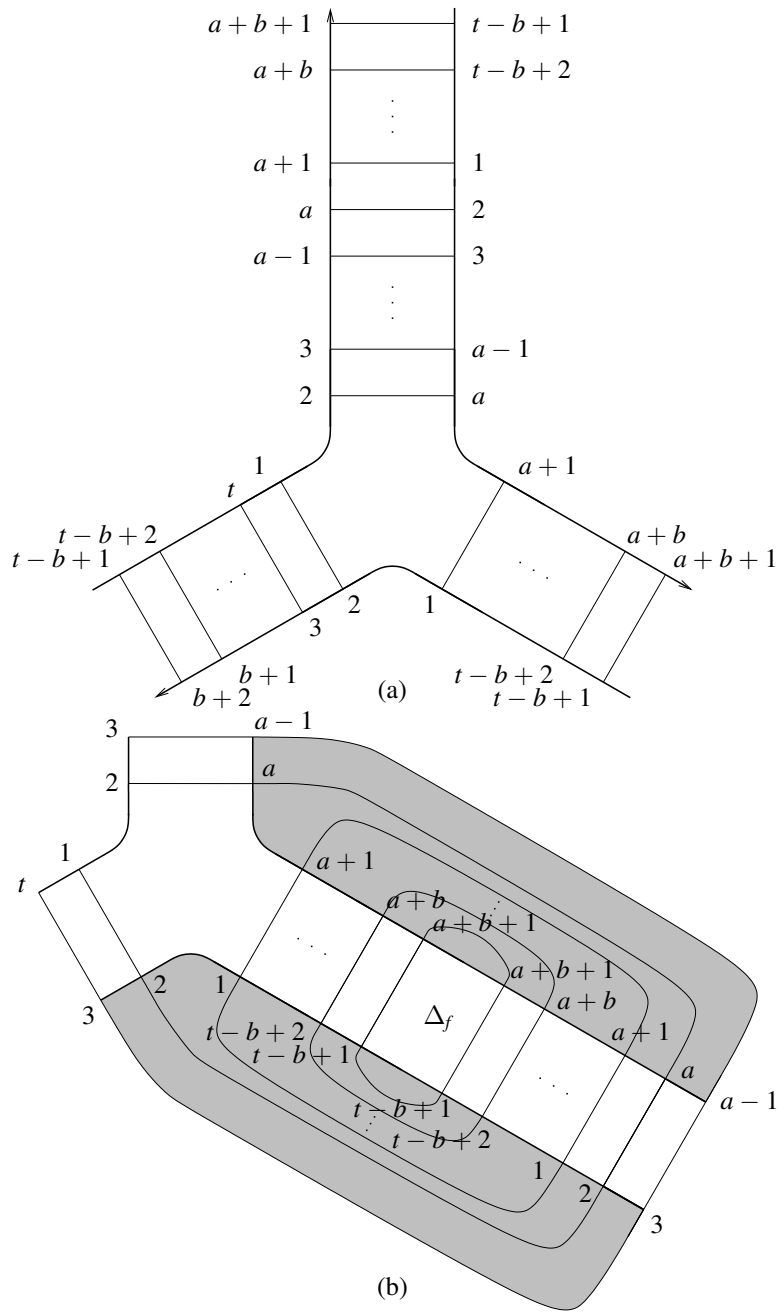

If $b$ bigons are attached to each edge of $g$, then $f$ is a trigon of $G_S^{t-b+1}$.  Then $f$ has two edges with label pair $\{a+b+1,t-b+1\}$ and one edge with label pair $\{t-b+1,b+2\}$.  The trigon $f$ appears as in \fullref{fig:3extforkedbigon}(a).  Note that $b+2 < a+b+1$ since $a > 1$, and $a+b+1 < t-b+1$ since otherwise there would be an edge of $G_S^{t-b+1}$ in the interior of $f$.

Either the two edges of $f$ with label pair $\{a+b+1, t-b+1\}$ lie in an essential annulus (satisfying the proposition) or they bound a bigon $\Delta_f$ on $\hatT$.  If $\Int \Delta_f \cap K = \emptyset$ then we may construct a long disk in a manner similar to the constructions for $a+1 \geq 6$ and $a+1=4$ in \fullref{thinningforkplusbigon}.  See \fullref{fig:3extforkedbigon}(b).  Thus $\Int \Delta_f \cap K \neq \emptyset$ satisfying the proposition.

If a different number of bigons is attached to each edge of $g$ to form $f$, then assume $b$ is the minimum number of bigons attached to an edge of $g$.  Hence all edges of $g$ \clearpage have $b$ bigons attached, and either one or two edges have more than $b$ edges attached.  Let $g'$ be the trigon formed by attaching $b$ bigons to each of the edges of $g$.  The trigon $g'$ appears as in \fullref{fig:3extforkedbigon}(a).  

If $f$ is obtained from $g'$ by attaching bigons to two of the edges of $g'$, then one may check that the arguments of Case II apply.   If $f$ is obtained from $g'$ by attaching bigons to just one of the edges, then one may check that the initial arguments of Case I also apply to imply that extra bigons are attached to the edge of $g'$ with label pair $\{b+2, t-b+1\}$.  Thus $f$ is a trigon of $G_S^{a+b+1}$. 

Either the two edges of $f$ with label pair $\{a+b+1, t-b+1\}$ lie in an essential annulus (satisfying the proposition) or they bound a bigon $\Delta_f$ on $\hatT$.  If $\Int \Delta_f \cap K = \emptyset$ then we may construct a long disk in a manner similar to the construction in \fullref{thinningforkplusbigon}.  Again, see \fullref{fig:3extforkedbigon}(b).  Thus $\Int \Delta_f \cap K \neq \emptyset$ satisfying the proposition.
\end{proof}

\begin{lemma}\label{lem:emptybigon}
If two edges of a trigon of some $\smash{G_S^x}$ bound a bigon $B$ on $\hatT$, then $\Int B \cap K = \emptyset$.
\end{lemma}

\begin{proof}
Assume otherwise.  Then for some $x$ there is a trigon $F$ of $\smash{G_S^x}$ with  two edges that bound a bigon $B$ in $\hatT$ with $\Int B \cap K \neq \emptyset$ such that (*) any other bigon in $\hatT$ bounded by two edges of a trigon of some $G_S^y$ and contained in $B$ has interior disjoint from $K$.

Let $U_1$ be a vertex of $G_T$ in $\Int B$.  By \fullref{musthavebigons} the graph $\smash{G_S^1}$ must have a bigon or trigon face $g$.  Because $U_1$ is contained in the interior of $B$, no pair of edges bounding $g$ may lie in an essential annulus on $\hatT$.  By \fullref{GT:L2.1} and \fullref{lieinanannulus}, $g$ must be a trigon that is not bounded by an $S3$ cycle or an extended $S3$ cycle.  If $g$ is not an innermost trigon then by \fullref{lieinannulusorboundbigon} a pair of its edges must bound a bigon $B'$ on $\hatT$ whose interior intersects $K$.  This bigon $B'$ however must be contained in $B$, contradicting (*).  Therefore $g$ must be an innermost trigon.  Since it cannot be bounded by an $S3$ cycle, \fullref{innermosttrigons} implies that it must be bounded by a forked extended $S2$ cycle.  We may assume $g$ has edges as in \fullref{fig:forkedbigon}(a) which appear on $G_T$ as in \fullref{fig:forkedbigon}(b) as described by \fullref{forkededges}.  

Because $U_1$ is in the interior of $B$ and the edges $e_2$ and $e_3$ lie in an essential annulus on $\hatT$, the edges of $G_T$ bounding $B$ must both be incident to $U_2$.  Because edges of $G_T$ connect vertices of opposite parity, the edges of $G_T$ bounding $B$ cannot also be incident to $U_{a+1}$.  Therefore they must be incident to some other vertex, say $U_z$, which may be $U_a$.  (Note that this means $x=2$ or $x=z$.)  One such configuration with $z \neq a$ is depicted in \fullref{fig:bigonconfig}. If $z = a$ then it may be the case that either $e_2$ or $e_3$ bounds a bigon on $\hatT$ with an edge of $B$. 
\begin{figure}[ht!]
\centering
\begin{picture}(0,0)%
\includegraphics{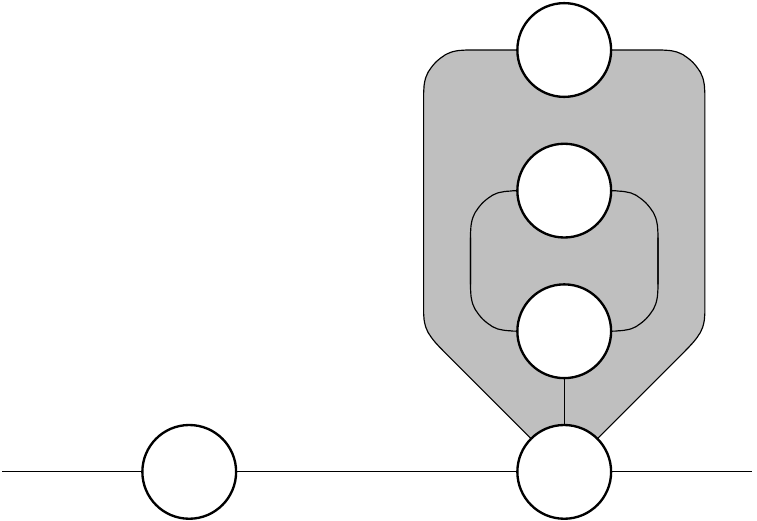}%
\end{picture}%
\setlength{\unitlength}{2960sp}%
\begingroup\makeatletter\ifx\SetFigFont\undefined%
\gdef\SetFigFont#1#2#3#4#5{%
  \reset@font\fontsize{#1}{#2pt}%
  \fontfamily{#3}\fontseries{#4}\fontshape{#5}%
  \selectfont}%
\fi\endgroup%
\begin{picture}(4824,3328)(7189,-3675)
\put(10801,-3436){\makebox(0,0)[b]{\smash{{\SetFigFont{9}{10.8}{\rmdefault}{\mddefault}{\updefault}{\color[rgb]{0,0,0}$2$}%
}}}}
\put(8401,-3436){\makebox(0,0)[b]{\smash{{\SetFigFont{9}{10.8}{\rmdefault}{\mddefault}{\updefault}{\color[rgb]{0,0,0}$a$}%
}}}}
\put(9301,-3586){\makebox(0,0)[lb]{\smash{{\SetFigFont{9}{10.8}{\rmdefault}{\mddefault}{\updefault}{\color[rgb]{0,0,0}$e_3$}%
}}}}
\put(11401,-3586){\makebox(0,0)[lb]{\smash{{\SetFigFont{9}{10.8}{\rmdefault}{\mddefault}{\updefault}{\color[rgb]{0,0,0}$e_2$}%
}}}}
\put(7501,-3586){\makebox(0,0)[lb]{\smash{{\SetFigFont{9}{10.8}{\rmdefault}{\mddefault}{\updefault}{\color[rgb]{0,0,0}$e_2$}%
}}}}
\put(10801,-2536){\makebox(0,0)[b]{\smash{{\SetFigFont{9}{10.8}{\rmdefault}{\mddefault}{\updefault}{\color[rgb]{0,0,0}$1$}%
}}}}
\put(10801,-1636){\makebox(0,0)[b]{\smash{{\SetFigFont{9}{10.8}{\rmdefault}{\mddefault}{\updefault}{\color[rgb]{0,0,0}$a+1$}%
}}}}
\put(10801,-736){\makebox(0,0)[b]{\smash{{\SetFigFont{9}{10.8}{\rmdefault}{\mddefault}{\updefault}{\color[rgb]{0,0,0}$z$}%
}}}}
\put(10276,-2086){\makebox(0,0)[lb]{\smash{{\SetFigFont{9}{10.8}{\rmdefault}{\mddefault}{\updefault}{\color[rgb]{0,0,0}$e_1$}%
}}}}
\put(11326,-2086){\makebox(0,0)[rb]{\smash{{\SetFigFont{9}{10.8}{\rmdefault}{\mddefault}{\updefault}{\color[rgb]{0,0,0}$e_5$}%
}}}}
\put(10876,-2986){\makebox(0,0)[lb]{\smash{{\SetFigFont{9}{10.8}{\rmdefault}{\mddefault}{\updefault}{\color[rgb]{0,0,0}$e_4$}%
}}}}
\end{picture}%
\caption{For $z \neq a$ the bigon $B$ is shown in grey with the edges of $g$ on $\hatT$.}
\label{fig:bigonconfig}
\end{figure}
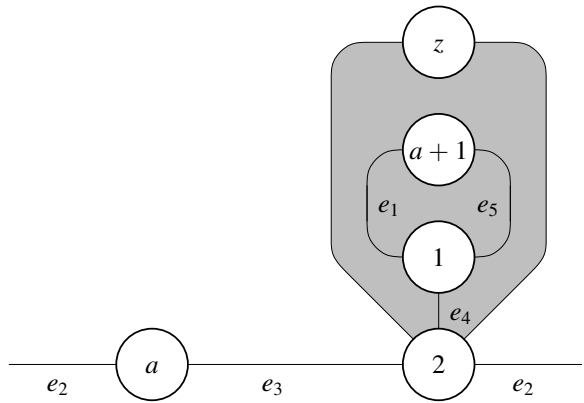

Observe that the vertex $U_{a+1}$ must be contained in the interior of $B$.  Therefore the bigon $\Delta_g$ on $\hatT$ bounded by the two edges $e_1$ and $e_5$ of $g$ is contained in $B$.  By (*) the interior of $\Delta_g$ must be disjoint from $K$.  

Furthermore the graph $G_S^{a+1}$ must have a trigon face $f$ bounded by a forked extended $S2$ cycle by the same arguments used above for the trigon $g$.  The edges of $f$ then appear as in \fullref{fig:forkedbigona+1}(a) or (b).  Because the bigon $\Delta_f$ on $\hatT$ bounded by the two edges $\smash{e_1'}$ and $e_5'$ of $f$ is contained in $B$, (*) implies its interior must be disjoint from $K$.  Since the interiors of both $\Delta_g$ and $\Delta_f$ are disjoint from $K$, \fullref{thinningtwoforks} implies that the edges of $G_S$ on $f$ must appear as in \fullref{fig:forkedbigona+1}(a).  Moreover, since the edge $\smash{e_4'}$ has label pair $\{a, a{+}1\}$, it must be the case that $z=a$.  Note that the two edges $\smash{e_2'}$ and $\smash{e_3'}$ both have label pairs $\{a, b{+}1\}$ and lie in an essential annulus on $\hatT$.  

{\bf Claim}\qua $b=1$

Observe that the vertex $U_b$ is contained in $B$.  Note that $b \neq 2, a+1$ since edges of $G_T$ connect vertices of opposite parity, and $b \neq a$ since the edges of $f$ would not form a forked extended $S2$ cycle otherwise.  If $b \neq 1$ then the above arguments apply once again to show that $G_S^b$ contains a trigon $h$ bounded by a forked extended $S2$ cycle.  As with the trigons $g$ and $f$, two edges of $h$ bound a bigon $\Delta_h$ on $\hatT$ which is contained in $B$ and thus has interior disjoint from $K$.  With $f$ playing the role of $g$, $h$ playing the role of $f$, and adjusting the labeling accordingly, \fullref{thinningtwoforks} implies that $h$ has an edge, say $\smash{e_4''}$, with label pair $\{b, b{+}1\}$.  Since $b \neq 1$, $b+1 \neq 2$.  Therefore in order for $\smash{e_2'}$ and $\smash{e_3'}$ to lie in an essential annulus on $\hatT$, the vertex $U_{b+1}$ must be disjoint from $B$.  Yet this contradicts that the edge $\smash{e_4''}$ connects the vertices $U_b$ and $U_{b+1}$.  Hence $b = 1$.

\newpage
Since $b=1$ the two forked extended $S2$ cycles bounding $g$ and $f$ fit the hypotheses of \fullref{lem:shortpairsofforks}.  Indeed the notation that we are currently using agrees with the notation used in the proof of \fullref{lem:shortpairsofforks}.  As a consequence, we have that $a=3$ and the subgraph $\Gamma$ of $G_T$ consisting of the edges $e_1, \dots, e_5$ and $\smash{e_1'}, \dots, e_5'$ and the vertices $U_1 \dots U_4$ does not lie in an essential annulus on $\hatT$.

Since $z=a=3$, the edges of $B$ have label pair $\{2, 3\}$.  Therefore $F$ is a trigon of either $G_S^2$ or $G_S^3$.  The subgraph $\Gamma$ and the two edges of $B$ are shown in \fullref{fig:twoshortforkedextendedS2}.  Because of the symmetry between the edges of $f$ and $g$, we may assume $F$ is a trigon of $G_S^2$.

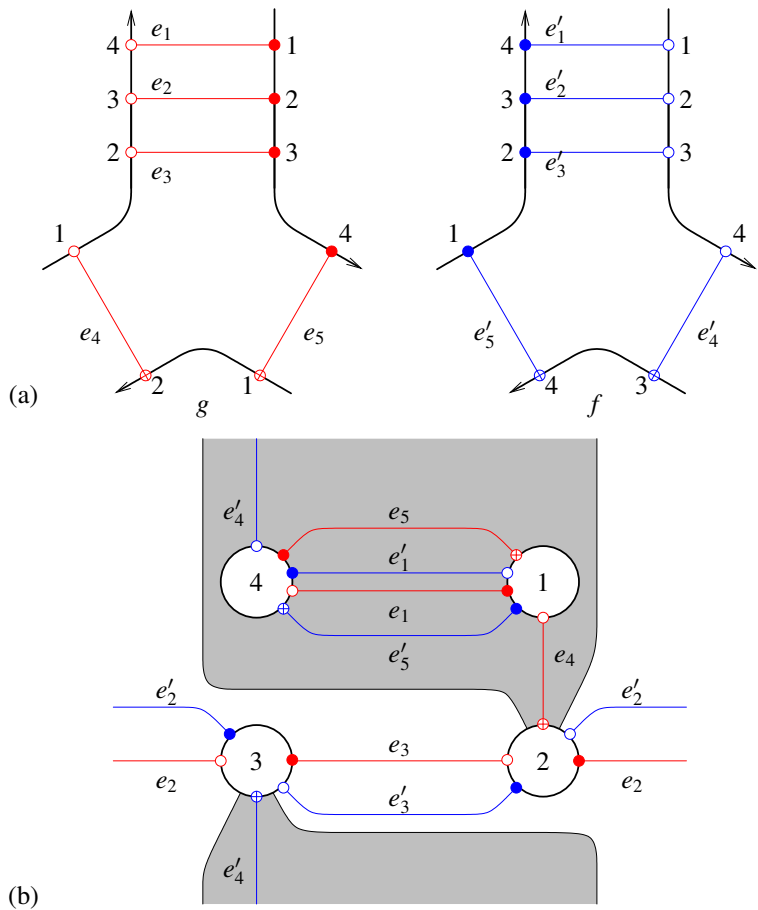
\begin{figure}[ht!]
\centering
\begin{picture}(0,0)%
\includegraphics{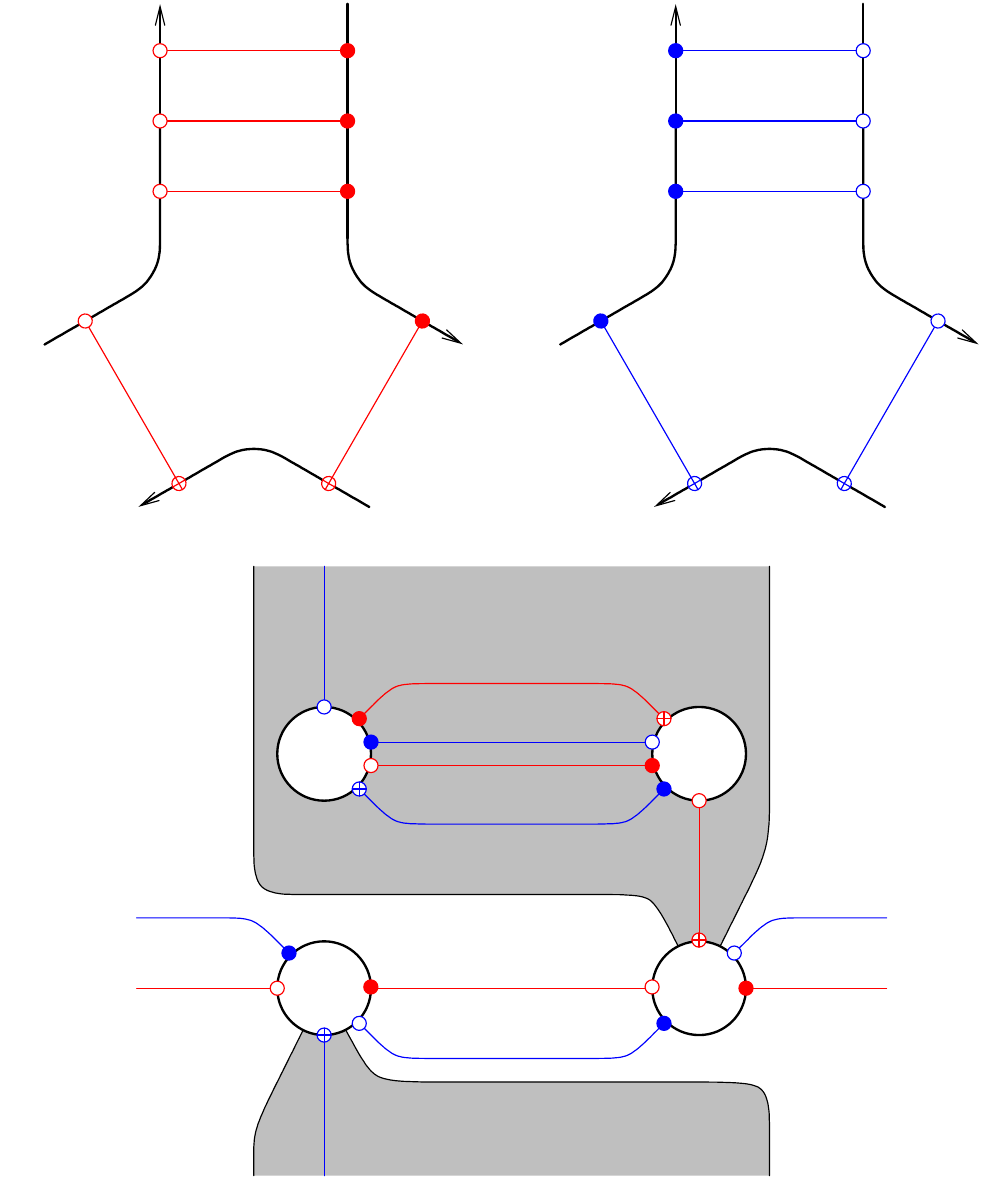}%
\end{picture}%
\setlength{\unitlength}{2960sp}%
\begingroup\makeatletter\ifx\SetFigFont\undefined%
\gdef\SetFigFont#1#2#3#4#5{%
  \reset@font\fontsize{#1}{#2pt}%
  \fontfamily{#3}\fontseries{#4}\fontshape{#5}%
  \selectfont}%
\fi\endgroup%
\begin{picture}(6325,7576)(3777,-8815)
\put(7051,-7936){\makebox(0,0)[b]{\smash{{\SetFigFont{9}{10.8}{\rmdefault}{\mddefault}{\updefault}{\color[rgb]{0,0,0}$\smash{e_3'}$}%
}}}}
\put(9526,-4036){\makebox(0,0)[lb]{\smash{{\SetFigFont{9}{10.8}{\rmdefault}{\mddefault}{\updefault}{\color[rgb]{0,0,0}$\smash{e_4'}$}%
}}}}
\put(7876,-4036){\makebox(0,0)[rb]{\smash{{\SetFigFont{9}{10.8}{\rmdefault}{\mddefault}{\updefault}{\color[rgb]{0,0,0}$e_5'$}%
}}}}
\put(4576,-4036){\makebox(0,0)[rb]{\smash{{\SetFigFont{9}{10.8}{\rmdefault}{\mddefault}{\updefault}{\color[rgb]{0,0,0}$e_4$}%
}}}}
\put(8251,-2611){\makebox(0,0)[lb]{\smash{{\SetFigFont{9}{10.8}{\rmdefault}{\mddefault}{\updefault}{\color[rgb]{0,0,0}$\smash{e_3'}$}%
}}}}
\put(5851,-4486){\makebox(0,0)[rb]{\smash{{\SetFigFont{9}{10.8}{\rmdefault}{\mddefault}{\updefault}{\color[rgb]{0,0,0}$1$}%
}}}}
\put(4276,-3211){\makebox(0,0)[rb]{\smash{{\SetFigFont{9}{10.8}{\rmdefault}{\mddefault}{\updefault}{\color[rgb]{0,0,0}$1$}%
}}}}
\put(4726,-2536){\makebox(0,0)[rb]{\smash{{\SetFigFont{9}{10.8}{\rmdefault}{\mddefault}{\updefault}{\color[rgb]{0,0,0}$2$}%
}}}}
\put(6076,-2536){\makebox(0,0)[lb]{\smash{{\SetFigFont{9}{10.8}{\rmdefault}{\mddefault}{\updefault}{\color[rgb]{0,0,0}$3$}%
}}}}
\put(6526,-3211){\makebox(0,0)[lb]{\smash{{\SetFigFont{9}{10.8}{\rmdefault}{\mddefault}{\updefault}{\color[rgb]{0,0,0}$4$}%
}}}}
\put(7576,-3211){\makebox(0,0)[rb]{\smash{{\SetFigFont{9}{10.8}{\rmdefault}{\mddefault}{\updefault}{\color[rgb]{0,0,0}$1$}%
}}}}
\put(8026,-2536){\makebox(0,0)[rb]{\smash{{\SetFigFont{9}{10.8}{\rmdefault}{\mddefault}{\updefault}{\color[rgb]{0,0,0}$2$}%
}}}}
\put(9376,-2536){\makebox(0,0)[lb]{\smash{{\SetFigFont{9}{10.8}{\rmdefault}{\mddefault}{\updefault}{\color[rgb]{0,0,0}$3$}%
}}}}
\put(9826,-3211){\makebox(0,0)[lb]{\smash{{\SetFigFont{9}{10.8}{\rmdefault}{\mddefault}{\updefault}{\color[rgb]{0,0,0}$4$}%
}}}}
\put(9151,-4486){\makebox(0,0)[rb]{\smash{{\SetFigFont{9}{10.8}{\rmdefault}{\mddefault}{\updefault}{\color[rgb]{0,0,0}$3$}%
}}}}
\put(6226,-4036){\makebox(0,0)[lb]{\smash{{\SetFigFont{9}{10.8}{\rmdefault}{\mddefault}{\updefault}{\color[rgb]{0,0,0}$e_5$}%
}}}}
\put(4951,-4486){\makebox(0,0)[lb]{\smash{{\SetFigFont{9}{10.8}{\rmdefault}{\mddefault}{\updefault}{\color[rgb]{0,0,0}$2$}%
}}}}
\put(8251,-4486){\makebox(0,0)[lb]{\smash{{\SetFigFont{9}{10.8}{\rmdefault}{\mddefault}{\updefault}{\color[rgb]{0,0,0}$4$}%
}}}}
\put(4951,-2686){\makebox(0,0)[lb]{\smash{{\SetFigFont{9}{10.8}{\rmdefault}{\mddefault}{\updefault}{\color[rgb]{0,0,0}$e_3$}%
}}}}
\put(6076,-1636){\makebox(0,0)[lb]{\smash{{\SetFigFont{9}{10.8}{\rmdefault}{\mddefault}{\updefault}{\color[rgb]{0,0,0}$1$}%
}}}}
\put(6076,-2086){\makebox(0,0)[lb]{\smash{{\SetFigFont{9}{10.8}{\rmdefault}{\mddefault}{\updefault}{\color[rgb]{0,0,0}$2$}%
}}}}
\put(4726,-2086){\makebox(0,0)[rb]{\smash{{\SetFigFont{9}{10.8}{\rmdefault}{\mddefault}{\updefault}{\color[rgb]{0,0,0}$3$}%
}}}}
\put(4726,-1636){\makebox(0,0)[rb]{\smash{{\SetFigFont{9}{10.8}{\rmdefault}{\mddefault}{\updefault}{\color[rgb]{0,0,0}$4$}%
}}}}
\put(4951,-1936){\makebox(0,0)[lb]{\smash{{\SetFigFont{9}{10.8}{\rmdefault}{\mddefault}{\updefault}{\color[rgb]{0,0,0}$e_2$}%
}}}}
\put(4951,-1486){\makebox(0,0)[lb]{\smash{{\SetFigFont{9}{10.8}{\rmdefault}{\mddefault}{\updefault}{\color[rgb]{0,0,0}$e_1$}%
}}}}
\put(8251,-1486){\makebox(0,0)[lb]{\smash{{\SetFigFont{9}{10.8}{\rmdefault}{\mddefault}{\updefault}{\color[rgb]{0,0,0}$\smash{e_1'}$}%
}}}}
\put(9376,-1636){\makebox(0,0)[lb]{\smash{{\SetFigFont{9}{10.8}{\rmdefault}{\mddefault}{\updefault}{\color[rgb]{0,0,0}$1$}%
}}}}
\put(8026,-2086){\makebox(0,0)[rb]{\smash{{\SetFigFont{9}{10.8}{\rmdefault}{\mddefault}{\updefault}{\color[rgb]{0,0,0}$3$}%
}}}}
\put(8026,-1636){\makebox(0,0)[rb]{\smash{{\SetFigFont{9}{10.8}{\rmdefault}{\mddefault}{\updefault}{\color[rgb]{0,0,0}$4$}%
}}}}
\put(9376,-2086){\makebox(0,0)[lb]{\smash{{\SetFigFont{9}{10.8}{\rmdefault}{\mddefault}{\updefault}{\color[rgb]{0,0,0}$2$}%
}}}}
\put(8251,-1936){\makebox(0,0)[lb]{\smash{{\SetFigFont{9}{10.8}{\rmdefault}{\mddefault}{\updefault}{\color[rgb]{0,0,0}$\smash{e_2'}$}%
}}}}
\put(5401,-4636){\makebox(0,0)[b]{\smash{{\SetFigFont{9}{10.8}{\rmdefault}{\mddefault}{\updefault}{\color[rgb]{0,0,0}$g$}%
}}}}
\put(8701,-4636){\makebox(0,0)[b]{\smash{{\SetFigFont{9}{10.8}{\rmdefault}{\mddefault}{\updefault}{\color[rgb]{0,0,0}$f$}%
}}}}
\put(8251,-7636){\makebox(0,0)[b]{\smash{{\SetFigFont{9}{10.8}{\rmdefault}{\mddefault}{\updefault}{\color[rgb]{0,0,0}$2$}%
}}}}
\put(5851,-7636){\makebox(0,0)[b]{\smash{{\SetFigFont{9}{10.8}{\rmdefault}{\mddefault}{\updefault}{\color[rgb]{0,0,0}$3$}%
}}}}
\put(5851,-6136){\makebox(0,0)[b]{\smash{{\SetFigFont{9}{10.8}{\rmdefault}{\mddefault}{\updefault}{\color[rgb]{0,0,0}$4$}%
}}}}
\put(8251,-6136){\makebox(0,0)[b]{\smash{{\SetFigFont{9}{10.8}{\rmdefault}{\mddefault}{\updefault}{\color[rgb]{0,0,0}$1$}%
}}}}
\put(7051,-7486){\makebox(0,0)[b]{\smash{{\SetFigFont{9}{10.8}{\rmdefault}{\mddefault}{\updefault}{\color[rgb]{0,0,0}$e_3$}%
}}}}
\put(9001,-7786){\makebox(0,0)[b]{\smash{{\SetFigFont{9}{10.8}{\rmdefault}{\mddefault}{\updefault}{\color[rgb]{0,0,0}$e_2$}%
}}}}
\put(5101,-7786){\makebox(0,0)[b]{\smash{{\SetFigFont{9}{10.8}{\rmdefault}{\mddefault}{\updefault}{\color[rgb]{0,0,0}$e_2$}%
}}}}
\put(5776,-5536){\makebox(0,0)[rb]{\smash{{\SetFigFont{9}{10.8}{\rmdefault}{\mddefault}{\updefault}{\color[rgb]{0,0,0}$\smash{e_4'}$}%
}}}}
\put(7051,-5536){\makebox(0,0)[b]{\smash{{\SetFigFont{9}{10.8}{\rmdefault}{\mddefault}{\updefault}{\color[rgb]{0,0,0}$e_5$}%
}}}}
\put(7051,-6736){\makebox(0,0)[b]{\smash{{\SetFigFont{9}{10.8}{\rmdefault}{\mddefault}{\updefault}{\color[rgb]{0,0,0}$e_5'$}%
}}}}
\put(7051,-6361){\makebox(0,0)[b]{\smash{{\SetFigFont{9}{10.8}{\rmdefault}{\mddefault}{\updefault}{\color[rgb]{0,0,0}$e_1$}%
}}}}
\put(7051,-5911){\makebox(0,0)[b]{\smash{{\SetFigFont{9}{10.8}{\rmdefault}{\mddefault}{\updefault}{\color[rgb]{0,0,0}$\smash{e_1'}$}%
}}}}
\put(5101,-7036){\makebox(0,0)[b]{\smash{{\SetFigFont{9}{10.8}{\rmdefault}{\mddefault}{\updefault}{\color[rgb]{0,0,0}$\smash{e_2'}$}%
}}}}
\put(9001,-7036){\makebox(0,0)[b]{\smash{{\SetFigFont{9}{10.8}{\rmdefault}{\mddefault}{\updefault}{\color[rgb]{0,0,0}$\smash{e_2'}$}%
}}}}
\put(8326,-6736){\makebox(0,0)[lb]{\smash{{\SetFigFont{9}{10.8}{\rmdefault}{\mddefault}{\updefault}{\color[rgb]{0,0,0}$e_4$}%
}}}}
\put(5776,-8536){\makebox(0,0)[rb]{\smash{{\SetFigFont{9}{10.8}{\rmdefault}{\mddefault}{\updefault}{\color[rgb]{0,0,0}$\smash{e_4'}$}%
}}}}
\put(3901,-4561){\makebox(0,0)[b]{\smash{{\SetFigFont{9}{10.8}{\rmdefault}{\mddefault}{\updefault}{\color[rgb]{0,0,0}(a)}%
}}}}
\put(3901,-8761){\makebox(0,0)[b]{\smash{{\SetFigFont{9}{10.8}{\rmdefault}{\mddefault}{\updefault}{\color[rgb]{0,0,0}(b)}%
}}}}
\end{picture}%
\caption{(a) The two forked extended $S2$ cycles $g$ and $f$ \qua (b) The bigon $B$ in grey and the edges of $g$ and $f$ on $\hatT$}
\label{fig:twoshortforkedextendedS2}
\end{figure}
\newpage

If $F$ is not an innermost trigon, then by \fullref{lieinannulusorboundbigon} either $F$ is bounded by an extended $S3$ cycle or contains a forked extended $S2$ cycle.  In either case, $F$ must contain an (extended) $S3$ cycle or an (extended) $S2$ cycle $\sigma$ for which neither $2$ nor $3$ belong to its label pair.  By \fullref{lieinanannulus}, the edges of $\sigma$ must lie in an essential annulus on $\hatT$.   But this essential annulus must be disjoint from $\Gamma$ which contradicts that $\Gamma$ does not itself lie in an essential annulus.

Since $F$ must be an innermost trigon, by \fullref{innermosttrigons} it is bounded by an $S3$ cycle or a forked extended $S2$ cycle.  

If $F$ is bounded by an $S3$ cycle, then it has label pair $\{2,3\}$.  By \fullref{GT:L2.1} its third edge must lie in an essential annulus $A$ with the two edges that bound $B$.  Beginning with the two edges of $B$, by following the endpoints of the edges along the corners of $F$ (along $\bdry H_{2,3}$),  we find that all three edges of $F$ must encounter the vertex $U_2$ from the same side of the essential annulus $A'$ in which edges $e_2$ and $e_3$ lie.  Following along the corners of $F$ to the vertex $U_3$, we find that the edges of $F$ must encounter the vertex $U_3$ from the other side of the annulus $A'$.  Therefore the core curves of each $A$ and $A'$ are isotopic.  But since $e_2$ and $e_3$ bound an $S2$ cycle in the solid torus on the same side of $\hatT$ as $F$, \fullref{GT:L2.1} implies that the cores of $A$ and $A'$ run $3$ and $2$ times respectively in the longitudinal direction of this solid torus.  Because these cores are isotopic on $\hatT$, this is a contradiction.

If $F$ is bounded by a forked extended $S2$ cycle of $G_S^2$, then by \fullref{forkededges} the edge with label pair distinct from the other two has label pair $\{2, 1\}$ or label pair $\{2, 3\}$.  But since the two of its edges that bound $B$ have the label pair $\{2,3\}$, the third edge must have label pair $\{2, 1\}$.  Therefore there must be an (extended) $S2$ cycle with label pair $\{1, 4\}$ on $F$.  By \fullref{lieinanannulus} the two edges of this $S2$ cycle must lie in an essential annulus on $\hatT$.  Yet since the vertices $U_1$ and $U_4$ are both contained in the interior of the bigon $B$, this cannot occur.
\end{proof}

\begin{prop}\label{prop:threetrigontypes}
A trigon of $\smash{G_S^x}$ is bounded by either an $S3$ cycle, an extended $S3$ cycle, or a forked extended $S2$ cycle.
\end{prop}

\begin{proof}
Let $f$ be a trigon of $\smash{G_S^x}$. 
By \fullref{innermosttrigons}, $f$ is innermost if and only if it is bounded by either an $S3$ cycle or a forked extended $S2$ cycle.  If $f$ is not innermost, then by \fullref{lieinannulusorboundbigon} either it is bounded by an extended $S3$ cycle or $f$ contains a forked extended $S2$ cycle whose edge with label pair distinct from the other two is not an edge of $f$.  

Let us assume $f$ is as in this last situation since it is the only one we must rule out.  Let $g$ be the forked extended $S2$ cycle in $f$ and assume it is labeled as in \fullref{fig:forkedbigon}(a).  By \fullref{lem:emptybigon} the bigon on $\hatT$ bounded by edges $e_1$ and $e_5$ of $g$ has interior disjoint from $K$.  Since edge $e_4$ is not an edge on the boundary of $f$, within $f$ there must be a bigon of $G_S$ with $e_4$ as an edge.  This contradicts \fullref{thinningforkplusbigon}.
\end{proof}


\begin{lemma}\label{lem:bigonsandtrigons}
Each graph $\smash{G_S^x}$ contains at least one of the following: an $S2$ cycle, an extended $S2$ cycle, an $S3$ cycle, an extended $S3$ cycle, a forked extended $S2$ cycle.
\end{lemma}

\begin{proof}
Since \fullref{musthavebigons} implies that each graph $\smash{G_S^x}$ must contain a bigon or a trigon, this lemma follows from \fullref{bigonS2} and \fullref{prop:threetrigontypes}.
\end{proof}

\begin{lemma}\label{circlesofintersection}
$S \cap T$ contains no simple closed curves that bound disks in $X$.
\end{lemma}

\begin{proof}
A simple closed curve of $S \cap T$ is either essential on $\hatT$ or bounds a disk on $\hatT$.

Assume there is a simple closed curve $\gamma \in S \cap T$ that is an essential curve on $\hatT$ and bounds a disk in $X$.  Via \fullref{lem:bigonsandtrigons} $G_S$ must have an $S2$ or $S3$ cycle $\sigma$.  By \fullref{GT:L2.1}, the edges of $\sigma$ lie in an essential annulus $A$ on $\hatT$.  Since $\gamma$ must be disjoint from $A$, it must be parallel to the core of $A$.  Thus the core of $A$ bounds a disk in $X$.  This contradicts \fullref{GT:L2.1}.

Assume there is a simple closed curve $\gamma \in S \cap T$ that bounds a disk $D \subseteq \hatT$.  

Assume there is a vertex, say $U_x$, of $G_T$ in $D$.  By \fullref{lem:bigonsandtrigons}, $\smash{G_S^x}$ contains either an $S2$ cycle, an extended $S2$ cycle, an $S3$ cycle, an extended $S3$ cycle, or a forked extended $S2$ cycle.  By \fullref{GT:L2.1} and \fullref{lieinanannulus} only the edges of a forked extended $S2$ cycle may not lie in an essential annulus.  Hence $\smash{G_S^x}$ must contain a forked extended $S2$ cycle $\sigma$ that bounds a face of $\smash{G_S^x}$ with edges as in \fullref{fig:forkedbigon}.  Since an edge of $\sigma$ must be incident to a vertex that has an (extended) $S2$ cycle also incident to it, $\sigma$ cannot lie in $D$ by \fullref{lieinanannulus}.  Thus there are no vertices of $G_T$ in $D$.  

Hence $D \cap K = \emptyset$.  Since lens spaces are irreducible, $\gamma$ must also bound a disk $D' \subseteq S$ that is isotopic rel--$\bdry$ in $X - N(K)$ to $D$.  This contradicts the minimality assumption on $|S \cap T|$.
\end{proof}

\subsection[Similar forked extended S2 cycles]{Similar forked extended $S2$ cycles} 

With \fullref{lem:bigonsandtrigons} in hand, we may refine our understanding of \fullref{lem:shortpairsofforks}.  In particular, as we shall soon see, the hypotheses of \fullref{lem:shortpairsofforks} hold true only if $t=4$.  The following proposition will be relevant in \fullref{sec:twoschcycles} for the proof of \fullref{labelaccount}.

\begin{prop}  \label{prop:nosaladforks}
Assume $\sigma_i \subseteq G_S^i$ and $\sigma_j \subseteq G_S^j$ are two forked extended $S2$ cycles such that for each (extended) $S2$ cycle contained in the face bounded by one of $\sigma_i$ or $\sigma_j$ there is an (extended) $S2$ cycle with the same label pair contained in the face bounded by the other forked extended $S2$ cycle.  If $t \geq 6$ then $i=j$.
\end{prop}

\begin{proof}
Let us return to the notation set up in \fullref{lem:shortpairsofforks}.  Specifically, note that for $f$, $g$, and $\Gamma$ we set $a = 3$.  One may care to refer to \fullref{fig:twoshortforkedextendedS2} (disregarding the bigon $B$) instead of \fullref{fig:oneconfiguration} for a depiction of the arrangement of their edges.

\fullref{lem:bigonsandtrigons} implies that for each $x \in \mathbf{t}$ $\smash{G_S^x}$ contains a face $F$ bounded by an (extended) $S2$ cycle, an (extended) $S3$ cycle, or a forked extended $S2$ cycle.  In any of these cases, there exists an $S2$ cycle or an $S3$ cycle $\sigma$ of $G_S$ on $F$.  
By \fullref{GT:L2.1} the edges of an $S2$ or $S3$ cycle lie in an essential annulus.   Since the subgraph $\Gamma$ does not lie in an essential annulus (see the proof of \fullref{lem:shortpairsofforks}), the edges of $\sigma$ must be incident to at least one vertex in $\Gamma$.  Therefore $\sigma$ must have label pair $\{t, 1\}, \{1, 2\}, \{2, 3\}, \{3, 4\}$, or $\{4, 5\}$.  By \fullref{oppositesides} the faces of $G_S$ bounded by Scharlemann cycles of order $2$ or $3$ with disjoint label pairs must lie on opposite sides of $\hatT$.  Since each $f$ and $g$ contain an $S2$ cycle with label pair $\{2, 3\}$, $\sigma$ may only have label pair $\{1, 2\}$, $\{2, 3\}$, or $\{3, 4\}$. 

{\em Claim:\/}  $\sigma$ cannot have label pair $\{1, 2\}$ or $\{3, 4\}$.

Assume $\sigma$ has label pair $\{1, 2\}$.  Note that $e_4$ cannot belong to $\sigma$.  Let $c$ be a corner of the face of $G_S$ bounded by $\sigma$.  Let $E_1$ and $E_2$ be the edges incident to $c$ at the vertices $U_1$ and $U_2$ respectively.   $E_1$ is incident to $U_1$ on the arc of $\bdry U_1$ either between $e_4$ and $e_5'$ or $e_4$ and $e_5$ in which no other edge of $\Gamma$ is incident.  

In the former case the edge $E_2$ is incident to $U_2$ in the arc of $\bdry U_2$ between $e_3$ and $\smash{e_3'}$ in which no other edge of $\Gamma$ is incident.  Thus the other end is incident to $U_3$.  This contradicts that the label pair of $E_2$ is $\{1, 2\}$.

In the latter case the edge $E_2$ is incident to $U_2$ in the arc of $\bdry U_2$ between $e_3$ and $e_4$ in which no other edge of $\Gamma$ is incident.  Thus $E_1$ and $E_2$ lie in the same disk of $\hatT \cut \Gamma$.  Therefore every edge of $\sigma$ lies in a disk.  This contradicts \fullref{GT:L2.1}.

Hence $\sigma$ cannot have label pair $\{1, 2\}$.  Due to the symmetry of $f$, $g$, and $\Gamma$, the same argument prohibits $\sigma$ from having label pair $\{3, 4\}$.  This proves the claim.

Due to the above claim, the $S2$ or $S3$ cycle of $G_S$ in any bigon or trigon of $\smash{G_S^x}$ for any $x \in \mathbf{t}$ must have label pair $\{2, 3\}$.  If for $5 \leq x \leq t$ the graph $\smash{G_S^x}$ has an extended $S2$ cycle or extended $S3$ cycle $\sigma'$, then the label pair of $\sigma'$ cannot include any of the labels $1$, $2$, $3$, or $4$.  Thus $\sigma'$ must lie in a disk in the complement of $\Gamma$.  This contradicts \fullref{lieinanannulus}.  

Thus $G_5$ and $G_t$ both have forked extended $S2$ cycles which bound trigons containing extended $S2$ cycles with label pair $\{1, 4\}$.  This contradicts \fullref{lem:shortpairsofforks}.  
\end{proof}

\section{Annuli and trees}\label{sec:annuliandtrees}

\subsection{Construction of annuli, related complexes, and trees}
Assume the labeling of the vertices of $G_T$ has been chosen so that $K_{(t, 1)}$ is contained in $X^-$.  Recall that by \fullref{circlesofintersection} the interior of each face of $G_S$ is disjoint from $\hatT$.

Assume $1 < k \leq t/2$ and $\smash{G_S^k}$ contains an extended $Sp$ cycle $\sigma$ for $p=2$ or $3$ to which the innermost Scharlemann cycle has label pair $\{1, t\}$.  The edges of $\sigma$ thus have label pair $\{k, t-k+1\}$.  

Let $B$ denote the face of $\smash{G_S^k}$ bounded by $\sigma$.  The edges of $G_S$ divide $B$ into  $p(k-1)$ bigons plus the $p$--gon bounded by the Scharlemann cycle.  Let $B_1$ be the $p$--gon bounded by the Scharlemann cycle.  If $p=2$, for $i \geq 2$ let $B_i'$ and $B_i''$ be the bigons bounded by the edges with label pairs $\{i-1, t-i+2\}$ and $\{i, t-i+1\}$.  If $p=3$, for $i \geq 2$ let $B_i'$, $B_i''$, and $B_i'''$ be the bigons bounded by the edges with label pairs $\{i-1, t-i+2\}$ and $\{i, t-i+1\}$ chosen so that 
\begin{itemize}
\item $\bdry B_i' \cap \bdry B_{i+1}' \neq \emptyset$,  $\bdry B_i'' \cap \bdry B_{i+1}'' \neq \emptyset$, $\bdry B_i''' \cap \bdry B_{i+1}''' \neq \emptyset$ and
\item the edges $\bdry B_i' \cap \bdry B_{i+1}'$ and $\bdry B_i'' \cap \bdry B_{i+1}''$ do not lie in a disk; see \fullref{lieinanannulus}.
\end{itemize}
For each $2 \leq i \leq k$, form the annulus $A_i$ from $B_i' \cup B_i'' \cup H_{(i-1,\,i)} \cup H_{(t-i+1,\, t-i+2)}$ by shrinking each $H_{(i-1,\,i)}$ and $H_{(t-i+1,\, t-i+2)}$ radially to their cores $K_{(i-1,\,i)}$ and $K_{(t-i+1,\, t-i+2)}$.  The annulus $A_i$ is contained in $X^+$ if and only if $i$ is even.  Let $a_i$ be the curve of $\bdry A_i$ that contains the point $K_i$, and set $a_1$ to be the curve of $\bdry A_2$ that contains $K_1$.

\begin{figure}[ht!]
\centering
\begin{picture}(0,0)%
\includegraphics{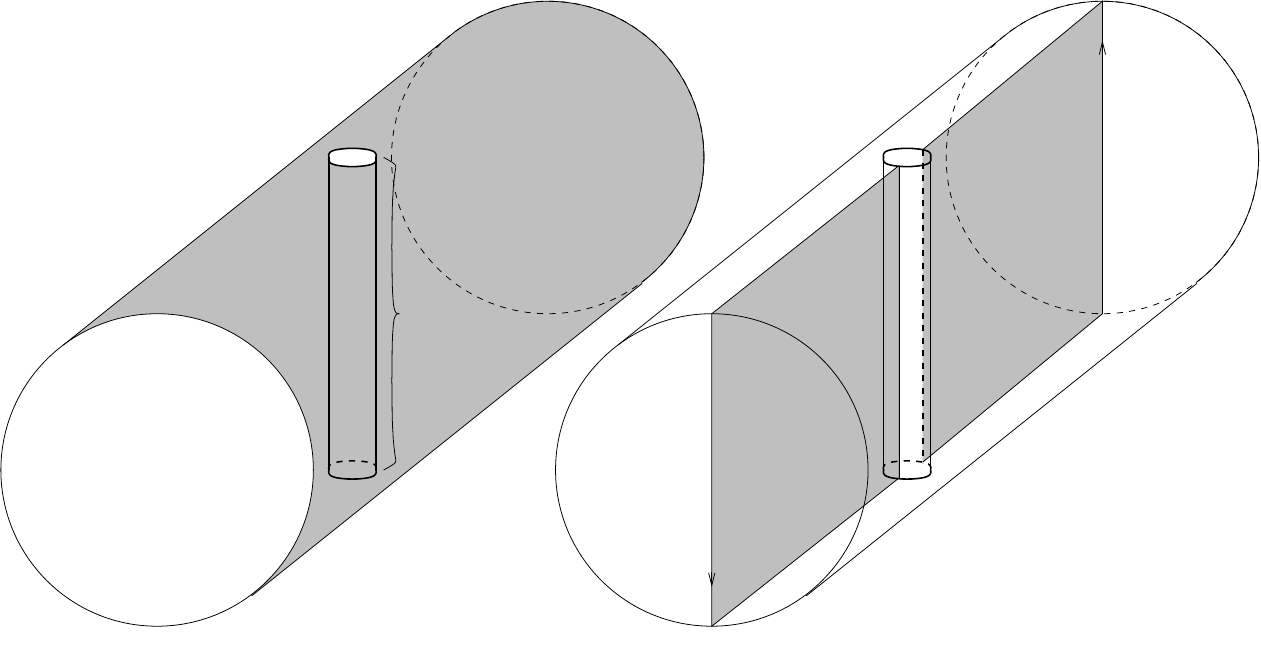}%
\end{picture}%
\setlength{\unitlength}{1973sp}%
\begingroup\makeatletter\ifx\SetFigFont\undefined%
\gdef\SetFigFont#1#2#3#4#5{%
  \reset@font\fontsize{#1}{#2pt}%
  \fontfamily{#3}\fontseries{#4}\fontshape{#5}%
  \selectfont}%
\fi\endgroup%
\begin{picture}(12091,6366)(1268,-6790)
\put(3001,-3361){\makebox(0,0)[b]{\smash{{\SetFigFont{6}{7.2}{\rmdefault}{\mddefault}{\updefault}{\color[rgb]{0,0,0}$\hatT$}%
}}}}
\put(5176,-3511){\makebox(0,0)[lb]{\smash{{\SetFigFont{6}{7.2}{\rmdefault}{\mddefault}{\updefault}{\color[rgb]{0,0,0}$H_{(t, 1)}$}%
}}}}
\put(2776,-5011){\makebox(0,0)[b]{\smash{{\SetFigFont{6}{7.2}{\rmdefault}{\mddefault}{\updefault}{\color[rgb]{0,0,0}$X^-$}%
}}}}
\put(4651,-1996){\makebox(0,0)[b]{\smash{{\SetFigFont{6}{7.2}{\rmdefault}{\mddefault}{\updefault}{\color[rgb]{0,0,0}$1$}%
}}}}
\put(4651,-4973){\makebox(0,0)[b]{\smash{{\SetFigFont{6}{7.2}{\rmdefault}{\mddefault}{\updefault}{\color[rgb]{0,0,0}$t$}%
}}}}
\put(8326,-3961){\makebox(0,0)[lb]{\smash{{\SetFigFont{6}{7.2}{\rmdefault}{\mddefault}{\updefault}{\color[rgb]{0,0,0}$B_1$}%
}}}}
\put(11176,-4861){\makebox(0,0)[lb]{\smash{{\SetFigFont{6}{7.2}{\rmdefault}{\mddefault}{\updefault}{\color[rgb]{0,0,0}$p=2$}%
}}}}
\put(2776,-6736){\makebox(0,0)[b]{\smash{{\SetFigFont{6}{7.2}{\rmdefault}{\mddefault}{\updefault}{\color[rgb]{0,0,0}(a)}%
}}}}
\put(8251,-6736){\makebox(0,0)[b]{\smash{{\SetFigFont{6}{7.2}{\rmdefault}{\mddefault}{\updefault}{\color[rgb]{0,0,0}(b)}%
}}}}
\end{picture}%
\caption{(a) The solid torus $X^-$ with $H_{(t, 1)} = \bar{N}(K_{(t, 1)})$\qua  (b) The bigon $B_1$ for $p=2$}
\label{fig:S2cycle}
\end{figure}
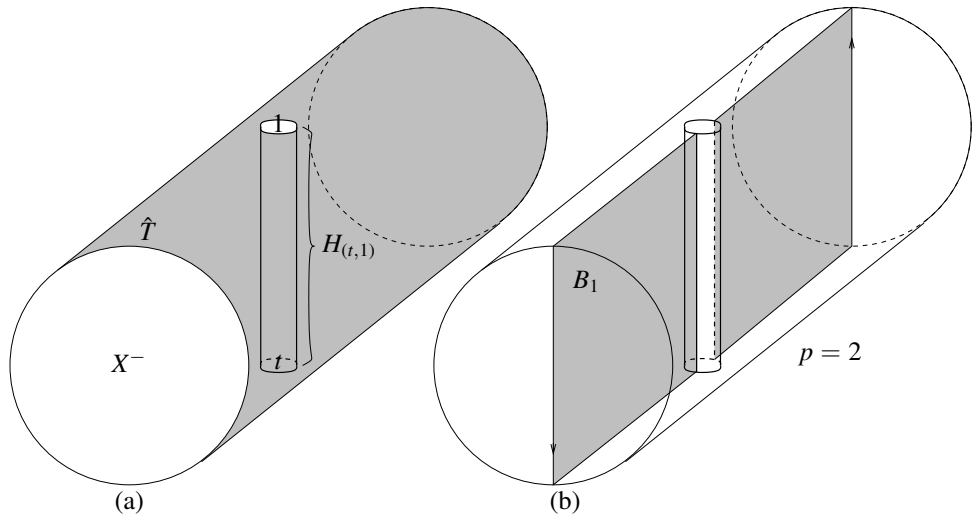

\begin{figure}[ht!]
\centering
\begin{picture}(0,0)%
\includegraphics{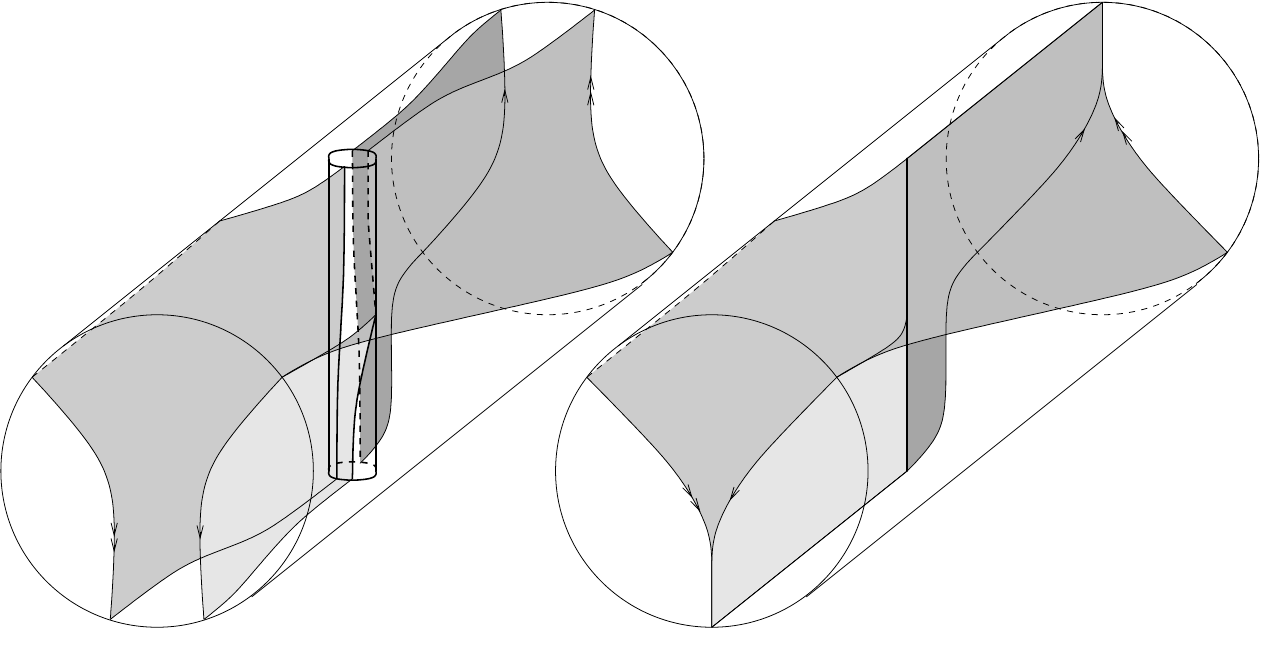}%
\end{picture}%
\setlength{\unitlength}{1973sp}%
\begingroup\makeatletter\ifx\SetFigFont\undefined%
\gdef\SetFigFont#1#2#3#4#5{%
  \reset@font\fontsize{#1}{#2pt}%
  \fontfamily{#3}\fontseries{#4}\fontshape{#5}%
  \selectfont}%
\fi\endgroup%
\begin{picture}(12091,6291)(2468,-8365)
\put(4651,-6961){\makebox(0,0)[lb]{\smash{{\SetFigFont{6}{7.2}{\rmdefault}{\mddefault}{\updefault}{\color[rgb]{0,0,0}$B_1$}%
}}}}
\put(3076,-4261){\makebox(0,0)[lb]{\smash{{\SetFigFont{6}{7.2}{\rmdefault}{\mddefault}{\updefault}{\color[rgb]{0,0,0}$p=3$}%
}}}}
\put(11701,-5686){\makebox(0,0)[lb]{\smash{{\SetFigFont{6}{7.2}{\rmdefault}{\mddefault}{\updefault}{\color[rgb]{0,0,0}$A_1$}%
}}}}
\put(3976,-8311){\makebox(0,0)[b]{\smash{{\SetFigFont{6}{7.2}{\rmdefault}{\mddefault}{\updefault}{\color[rgb]{0,0,0}(a)}%
}}}}
\put(9376,-8311){\makebox(0,0)[b]{\smash{{\SetFigFont{6}{7.2}{\rmdefault}{\mddefault}{\updefault}{\color[rgb]{0,0,0}(b)}%
}}}}
\end{picture}%
\caption{(a) The trigon $B_1$ for $p=3$ \qua (b) The complex $A_1$ for $p=3$}
\label{fig:S3cycle}
\end{figure}
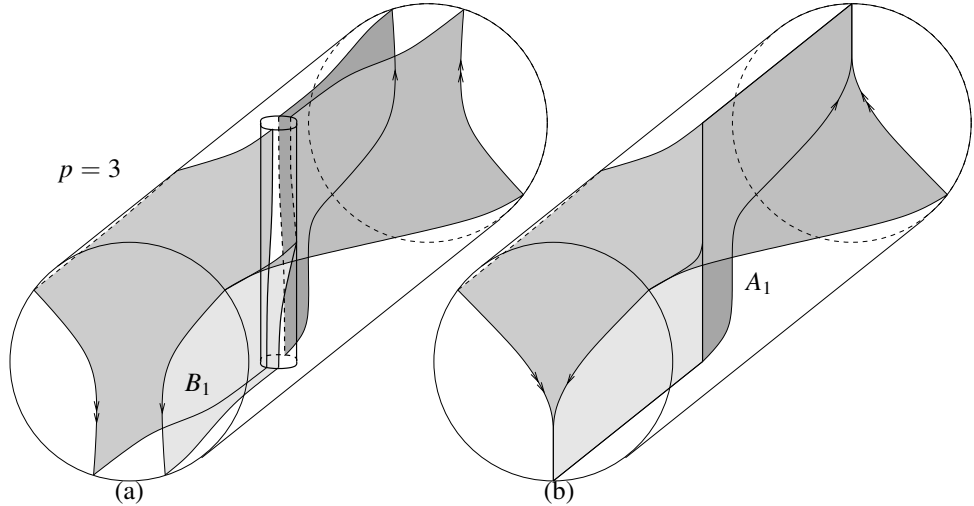

\label{$A_1$}
If $p=2$, let $A_1$ be the M\"obius band formed from $B_1 \cup H_{(t, 1)}$ by shrinking $H_{(t, 1)}$ to its core radially.  \fullref{fig:S2cycle} shows (a) the solid torus $X^-$ with $H_{(t, 1)}$ and (b) how $B_1$ sits in $X^-$.  Such a construction of a M\"obius band is done in the proof of Lemma~2.5 of \cite{gt:dsokwylsagok}.

If $p=3$, let $A_1$ be the complex formed from $B_1 \cup H_{(t, 1)}$ by first isotoping the edge $B_1 \cap B_2'''$ across the disk component of $\hatT \cut (B_1 \cup H_{(t, 1)})$ (which has interior disjoint from $K$ by \fullref{lem:emptybigon}) keeping the corners of $B_1$ on $H_{(t, 1)}$.  See \fullref{fig:S3cycle}(a) for an example of the placement of $B_1$ in $X^-$; see also \fullref{fig:S2cycle}(a).  Next identify a small collar neighborhood in $B_1$ of the two edges that have been isotoped together.  Then after shrinking $H_{(t, 1)}$ to its core radially, the resulting complex $A_1$ intersects $\hatT$ as the curve $a_1$.  See \fullref{fig:S3cycle}(b).  Note that after chopping along a suitable meridional disk $D$, $(X^- \cut D, \bar{N}(A_1) \cut D)$ is homeomorphic to $(I \times D^2, I \times \bar{N}(Y))$ where $Y$ is the complex in the standard disk $D^2$ formed by three radii; see \fullref{GT:L2.1}.
For both of these cases, $\bdry \bar{N}(A_1) \cap X^-$ is an annulus that double covers $A_1$ (except along $K_{(t, 1)}$  if $p=3$ where it is triple covered by the annulus).

Define the annulus $A = \bigcup_{i=2}^k A_i$.   By \fullref{lieinanannulus} $\bigcup_{i=2}^{k} \bdry A_i = \{a_1, \dots, a_k\}$ is a collection of $k$ essential simple closed curves on $\hatT$ that are mutually disjoint and parallel.  Furthermore, by \fullref{GT:L2.1} these curves are not meridional in either $X^+$ or $X^-$.  These curves then divide $\hatT$ into $k$ annuli $T_j$ so that $\hatT - A = \bigcup_{j=1}^k \Int T_j$.

  Recall that two $3$--manifolds with toroidal boundaries attached together along an annulus that is incompressible in each is a solid torus if and only if each of the manifolds is itself a solid torus and one of their meridians crosses the annulus exactly once.  Each annulus $A_i$ (for $i \neq 1$) thus separates $X^\pm$ (where $\pm = +$ or $-$ depending on the parity of $i$) into two solid tori where at least one of the meridians of the solid tori crosses the core curve of $A_i$ just once.   This implies that $A_i$ is isotopic in $X^\pm$ rel--$\bdry$ onto $\hatT$.  Let $V_i$ be the solid torus of $X^\pm \cut A_i$ through which this isotopy occurs.  In the event that both solid tori of $X^\pm \cut A_i$ would work, choose the $V_i$ so that if $\Int V_i \cap \Int V_j$ is nonempty, then it is either all of $\Int V_i$ or all of $\Int V_j$.  Note that none of the $V_i$ contain $A_1$.

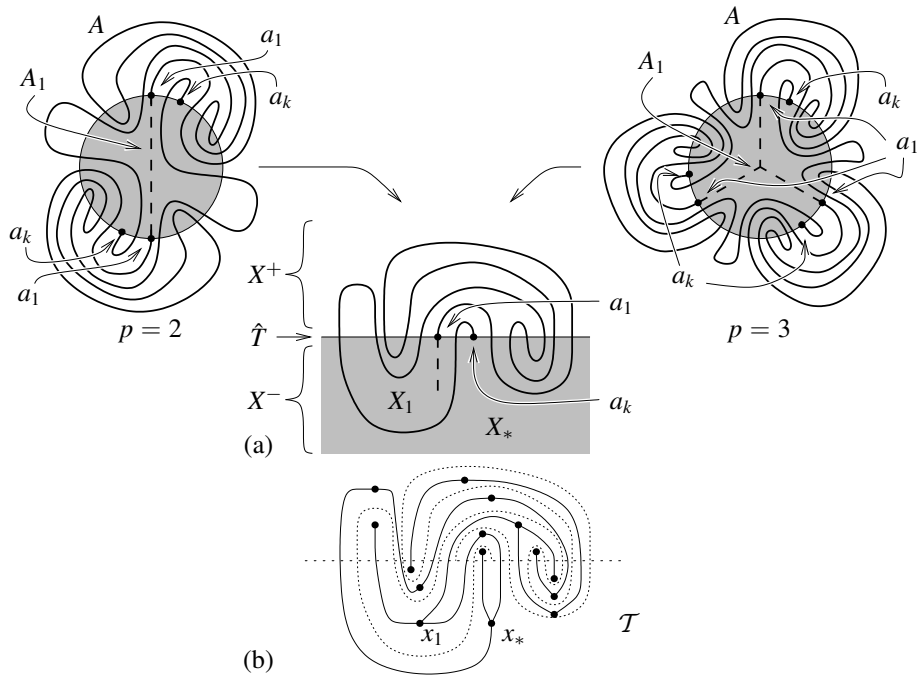
\begin{figure}[ht!]
\centering
\begin{picture}(0,0)%
\includegraphics{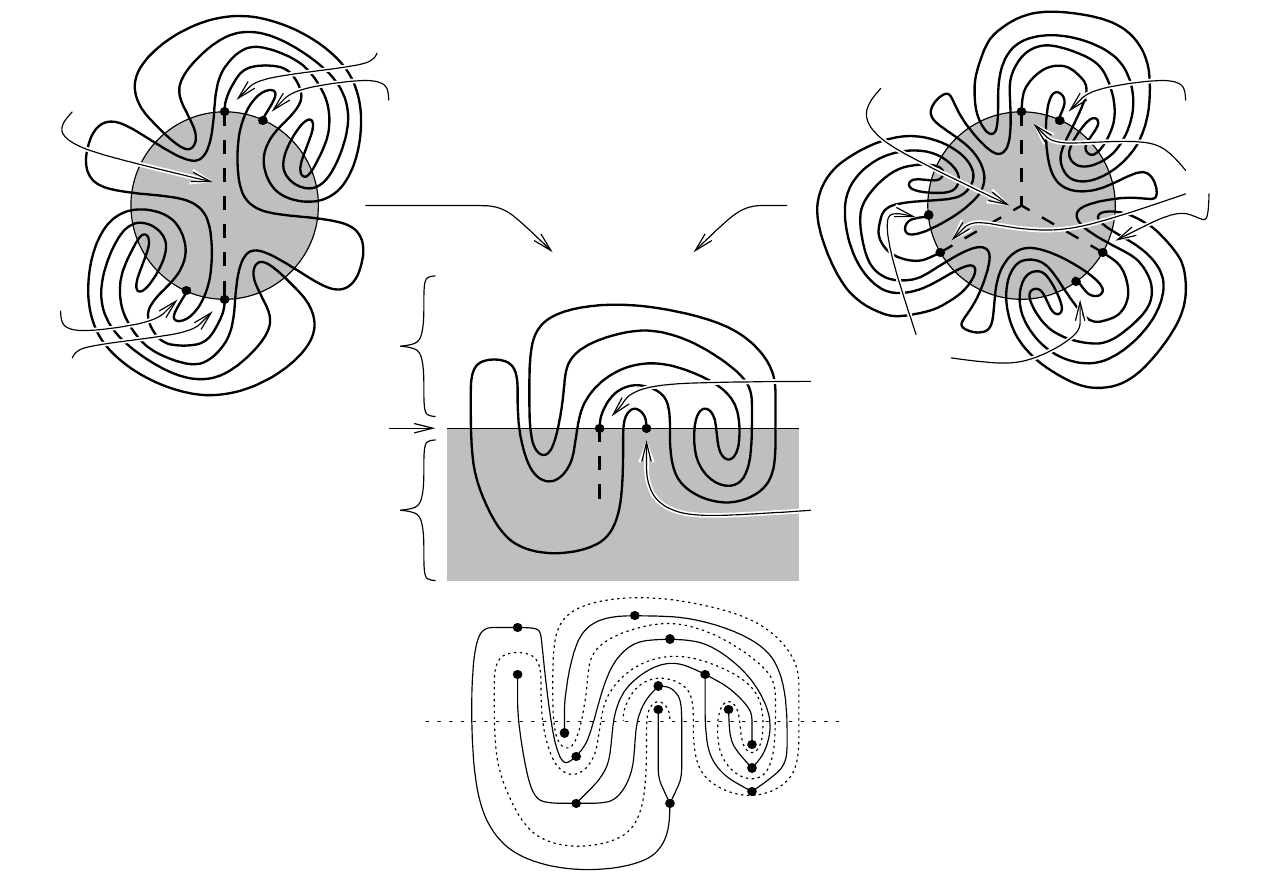}%
\end{picture}%
\setlength{\unitlength}{2960sp}%
\begingroup\makeatletter\ifx\SetFigFont\undefined%
\gdef\SetFigFont#1#2#3#4#5{%
  \reset@font\fontsize{#1}{#2pt}%
  \fontfamily{#3}\fontseries{#4}\fontshape{#5}%
  \selectfont}%
\fi\endgroup%
\begin{picture}(8123,5570)(664,-6423)
\put(3001,-6361){\makebox(0,0)[b]{\smash{{\SetFigFont{9}{10.8}{\rmdefault}{\mddefault}{\updefault}{\color[rgb]{0,0,0}(b)}%
}}}}
\put(5026,-6136){\makebox(0,0)[lb]{\smash{{\SetFigFont{9}{10.8}{\rmdefault}{\mddefault}{\updefault}{\color[rgb]{0,0,0}$x_*$}%
}}}}
\put(4351,-6136){\makebox(0,0)[lb]{\smash{{\SetFigFont{9}{10.8}{\rmdefault}{\mddefault}{\updefault}{\color[rgb]{0,0,0}$x_1$}%
}}}}
\put(6001,-6061){\makebox(0,0)[lb]{\smash{{\SetFigFont{9}{10.8}{\rmdefault}{\mddefault}{\updefault}{\color[rgb]{0,0,0}$\T$}%
}}}}
\put(3001,-4561){\makebox(0,0)[b]{\smash{{\SetFigFont{9}{10.8}{\rmdefault}{\mddefault}{\updefault}{\color[rgb]{0,0,0}(a)}%
}}}}
\put(1126,-1486){\makebox(0,0)[b]{\smash{{\SetFigFont{8}{9.6}{\rmdefault}{\mddefault}{\updefault}{\color[rgb]{0,0,0}$A_1$}%
}}}}
\put(1651,-1036){\makebox(0,0)[b]{\smash{{\SetFigFont{8}{9.6}{\rmdefault}{\mddefault}{\updefault}{\color[rgb]{0,0,0}$A$}%
}}}}
\put(3076,-1636){\makebox(0,0)[lb]{\smash{{\SetFigFont{9}{10.8}{\rmdefault}{\mddefault}{\updefault}{\color[rgb]{0,0,0}$a_k$}%
}}}}
\put(1201,-3286){\makebox(0,0)[rb]{\smash{{\SetFigFont{9}{10.8}{\rmdefault}{\mddefault}{\updefault}{\color[rgb]{0,0,0}$a_1$}%
}}}}
\put(1126,-2761){\makebox(0,0)[rb]{\smash{{\SetFigFont{9}{10.8}{\rmdefault}{\mddefault}{\updefault}{\color[rgb]{0,0,0}$a_k$}%
}}}}
\put(3001,-1111){\makebox(0,0)[lb]{\smash{{\SetFigFont{9}{10.8}{\rmdefault}{\mddefault}{\updefault}{\color[rgb]{0,0,0}$a_1$}%
}}}}
\put(2101,-3586){\makebox(0,0)[b]{\smash{{\SetFigFont{9}{10.8}{\rmdefault}{\mddefault}{\updefault}{\color[rgb]{0,0,0}$p = 2$}%
}}}}
\put(4201,-4186){\makebox(0,0)[b]{\smash{{\SetFigFont{9}{10.8}{\rmdefault}{\mddefault}{\updefault}{\color[rgb]{0,0,0}$X_1$}%
}}}}
\put(5026,-4411){\makebox(0,0)[b]{\smash{{\SetFigFont{9}{10.8}{\rmdefault}{\mddefault}{\updefault}{\color[rgb]{0,0,0}$X_*$}%
}}}}
\put(5926,-4186){\makebox(0,0)[lb]{\smash{{\SetFigFont{9}{10.8}{\rmdefault}{\mddefault}{\updefault}{\color[rgb]{0,0,0}$a_k$}%
}}}}
\put(5926,-3361){\makebox(0,0)[lb]{\smash{{\SetFigFont{9}{10.8}{\rmdefault}{\mddefault}{\updefault}{\color[rgb]{0,0,0}$a_1$}%
}}}}
\put(3076,-3661){\makebox(0,0)[rb]{\smash{{\SetFigFont{9}{10.8}{\rmdefault}{\mddefault}{\updefault}{\color[rgb]{0,0,0}$\hatT$}%
}}}}
\put(3226,-3136){\makebox(0,0)[rb]{\smash{{\SetFigFont{9}{10.8}{\rmdefault}{\mddefault}{\updefault}{\color[rgb]{0,0,0}$X^+$}%
}}}}
\put(3226,-4186){\makebox(0,0)[rb]{\smash{{\SetFigFont{9}{10.8}{\rmdefault}{\mddefault}{\updefault}{\color[rgb]{0,0,0}$X^-$}%
}}}}
\put(6301,-1336){\makebox(0,0)[b]{\smash{{\SetFigFont{8}{9.6}{\rmdefault}{\mddefault}{\updefault}{\color[rgb]{0,0,0}$A_1$}%
}}}}
\put(6976,-961){\makebox(0,0)[b]{\smash{{\SetFigFont{8}{9.6}{\rmdefault}{\mddefault}{\updefault}{\color[rgb]{0,0,0}$A$}%
}}}}
\put(6676,-3136){\makebox(0,0)[rb]{\smash{{\SetFigFont{9}{10.8}{\rmdefault}{\mddefault}{\updefault}{\color[rgb]{0,0,0}$a_k$}%
}}}}
\put(8326,-2011){\makebox(0,0)[lb]{\smash{{\SetFigFont{9}{10.8}{\rmdefault}{\mddefault}{\updefault}{\color[rgb]{0,0,0}$a_1$}%
}}}}
\put(8176,-1636){\makebox(0,0)[lb]{\smash{{\SetFigFont{9}{10.8}{\rmdefault}{\mddefault}{\updefault}{\color[rgb]{0,0,0}$a_k$}%
}}}}
\put(7201,-3586){\makebox(0,0)[b]{\smash{{\SetFigFont{9}{10.8}{\rmdefault}{\mddefault}{\updefault}{\color[rgb]{0,0,0}$p = 3$}%
}}}}
\end{picture}%
\caption{(a) A schematic of the annulus $A$\qua  (b)  The corresponding graph $\T$}
    \label{fig:crosssection}
\end{figure}

Consider the collection of solid tori $\mathcal{X} = X \cut A \cut \hatT$ formed by chopping $X$ along both $A$ and $\hatT$. 
  The boundaries of all but at most two $X_l \in \mathcal{X}$ are alternately comprised of the annuli $A_i$ and the annuli $T_j$.  Let $X_1 \subseteq X^-$ be the solid torus of $\mathcal{X}$ that contains $A_1$.  The boundary of $X_1$ intersects the curve $a_1$ but not $\Int A_2$.  Let $X_* \in \mathcal{X}$ be the solid torus whose boundary intersects $a_k$ but not $\Int A_k$.  Since $k>1$, each of $\bdry X_*$ and $\bdry X_1$ contains more than one of the annuli $T_j$.  The meridian of each $X_l \in \mathcal{X}$ intersects the core curve of an annulus $T_j$ on its boundary exactly once except for at most two, one of which is necessarily $X_1 \subseteq X^-$ and the other we shall denote $X_0 \subseteq X^+$.  Note that $X_0 \cap \Int V_i = \emptyset$.  It may be the case that $X_* = X_0$ or $X_* = X_1$.  See \fullref{fig:crosssection}(a) for a schematic example of the annulus $A$, the complex $A_1$, $X_1$, and $X_*$ for both $p = 2$ and $3$.  We will continue draw schematics of the annulus $A$ et al.\ in the quotiented manner indicated in \fullref{fig:crosssection}(a).

Form a graph $\T$ with vertices $x_l$ corresponding to solid tori $X_l \in \mathcal{X}$ and edges $t_j$ corresponding to annuli $T_j$ so that an edge $t_j$ connects vertices $x_l$ and $x_m$ if the annulus $T_j$ is contained in the boundary of each $X_l$ and $X_m$.  \fullref{fig:crosssection}(b) shows the graph $\T$ corresponding to the annulus $A$ in \fullref{fig:crosssection}(a).  The vertices $x_l$ may be marked as $+$ or $-$ according to whether the corresponding $X_l$ is contained in $X^+$ or $X^-$.  Adjacent vertices must have opposite signs.

\begin{lemma}
The graph $\T$ is a tree.
\end{lemma}

\begin{proof}
If $\T$ is not a tree, then there is a cycle of distinct edges $t_{j_1}, t_{j_2}, \dots, t_{j_w}$ such that $\bdry t_{j_{i}} \cap \bdry t_{j_{i+1}} = x_{l_i}$ (unless $s = w$ in which case $\bdry t_{j_1} = \bdry t_{j_2}$) and $j_i = j_{i'}$ only if $i = i'$.  

Let $R_{l_i}$ be an annulus properly contained in the solid torus $X_{l_i}$ corresponding to the vertex $x_{l_i}$ whose boundary components are the core curves of the annuli $T_{j_i}$ and $T_{j_{i+1}}$.   Since the annuli are connected in a cycle, $R = \bigcup_{i=1}^w R_{l_i}$ is a torus that is disjoint from $A$.  Since $R$ is contained in the lens space $X$, it must be separating.  On the boundary of $X_{l_i}$ the annuli $T_{j_i}$ and $T_{j_{i+1}}$ are necessarily separated by at least two of the curves $\{a_1, \dots, a_k\}$ since each $\bdry T_{j_i}$ and $\bdry T_{j_{i+1}}$ are pairs of curves in $\{a_1, \dots, a_k\}$ and $\Int T_{j_i} \cap \Int T_{j_{i+1}} = \emptyset$.  Specifically, $X_{l_i} - R_{l_i}$ has two components each of which has nontrivial intersection with $A$.  This however contradicts that $A$ is connected and $R$ is separating.
\end{proof}

Since $\T$ is a tree and there are $k$ annuli $T_j$, $\T$ has $k+1$ vertices.  Hence there are $k+1$ solid tori $X_l$.  Furthermore, since $T$ divides $K$ into $t$ arcs (the arcs $K_{(1, 2)}, K_{(2, 3)}, \dots, K_{(t, 1)}$) of which $2(k-1)$ (the arcs $K_{(t-k+1,\, t-k+2)}, \dots, K_{(t-1,\,t)}$ and $K_{(1, 2)}, \dots, K_{(k-1,\,k)}$) are contained in $A$, the remaining $t-2(k-1)$ arcs (the arcs $K_{(k,\,k+1)}, \dots, K_{(t-k,\, t-k+1)}$ and $K_{(t,1)}$) are disjoint from $A$ (except at the four points $K_1$, $K_t$, $K_k$, and $K_{t-k+1}$).  

Aside from $K_{(t,1)}$, these remaining arcs together form the arc $K_{(k,\,t-k+1)}$ which runs from $\bdry A$ inside $X_*$ through the interior of the solid tori $X_l$ crossing between them by passing through the $T_j$ eventually returning to $\bdry A$ in $X_*$.  Due to $\T$ being a tree, $K_{(k,\,t-k+1)}$ eventually returns through the same $T_j$ from which it entered an $X_l$.  Therefore $|K_{(k,\,t-k+1)} \cap \Int T_j|$ is even for every $j$. 

 Consider the subtree $\T_* \subseteq \T$ consisting of the vertices and edges corresponding to the $X_l$ and $T_j$ with which $K_{(k,\,t-k+1)}$ has nonempty intersection.  Since $K_{(k,\,t-k+1)}$ intersects the union $\smash{\bigcup_{j=1}^k} \Int T_j$ a total of $(t-k+1) - k - 1 = t-2k$ times and $K_{(k,\,t-k+1)}$ intersects each $\Int T_j$ an even number of times, $\T_*$ has at most $t/2 - k$ edges.  Root both trees $\T$ and $\T_*$ at the vertex $x_*$.

\begin{remark}  \label{evenintersections} 
Indeed, for each edge $t_j$ of $\T_*$ there is a positive even number of intersections of the interior of $K_{(k,\,t-k+1)}$ with the corresponding annulus $T_j \subseteq \hatT \cut A$.   Notice that if $t/2 - k < k$ (ie\ $t/4 < k$) then there exists at least one $X_l$ (and hence some $T_j$) which $K_{(k,\,t-k+1)}$ does not intersect.  To rephrase, if $K_{(k,\,t-k+1)}$ intersects every $X_l \in \mathcal{X}$, then $t/4 \geq k$ and the face $B$ contains (extended) $S2$ or $S3$ cycles for the $2k (\leq t/2)$ graphs $G_S^i$ with $t-k+1 \leq i \leq t$ or $1 \leq i \leq k$.
\end{remark}

For either $\T$ or $\T_*$ we say a vertex other than $x_*$ of valency $1$ is a {\em leaf\/} and a (nonleaf) vertex other than $x_*$ all of whose adjacent vertices except at most one are leaves is a {\em penultimate leaf\/}.

\subsection{Initial constraints on the trees due to thinness}

\begin{lemma} \label{leafisotopy} 
Let $x_m$ be a leaf of $\T$ not in $\T_*$ which is not $x_0$ and corresponds to the solid torus $X_m$.   If $X_m \subseteq X^+$ (resp.\ $X^-$), then there is a high disk (resp.\ low disk) in $X_m$ for each of the two arcs of $K \cap \bdry X_m$.  
\end{lemma}

\begin{proof}
We first remark that $x_1$ cannot be a leaf unless $k=1$.  Nevertheless, we are assuming $k > 1$.  Since $x_m \not \in \T_*$, $x_m \not \in X_*$.  Therefore the boundary of $X_m$ is formed by two annuli, say $T_m$ and $A_m$.  Since $X_m \neq X_0$, the meridian of $X_m$ crosses the core curves of $T_m$ and $A_m$ each once.  In particular, there are two disjoint meridional disks each with boundary consisting of one of the transverse arcs $K_{(m-1,\,m)}$ and $K_{(t-m+1,\, t-m+2)}$ on $A_m$ and an arc on $T_m$.  Since $K \cap \Int X_m = \emptyset$, these are high (or low) disks. 
\end{proof}

\begin{lemma}\label{leafcancel}
Any two leaves of $\T$ not in $\T_*$ neither of which is $x_0$ must have the same sign.
\end{lemma}

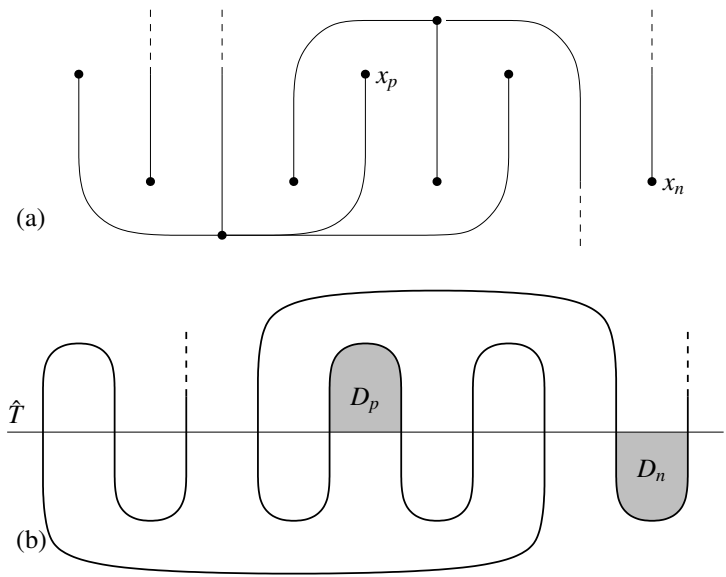
\begin{figure}[ht!]
\centering
\begin{picture}(0,0)%
\includegraphics{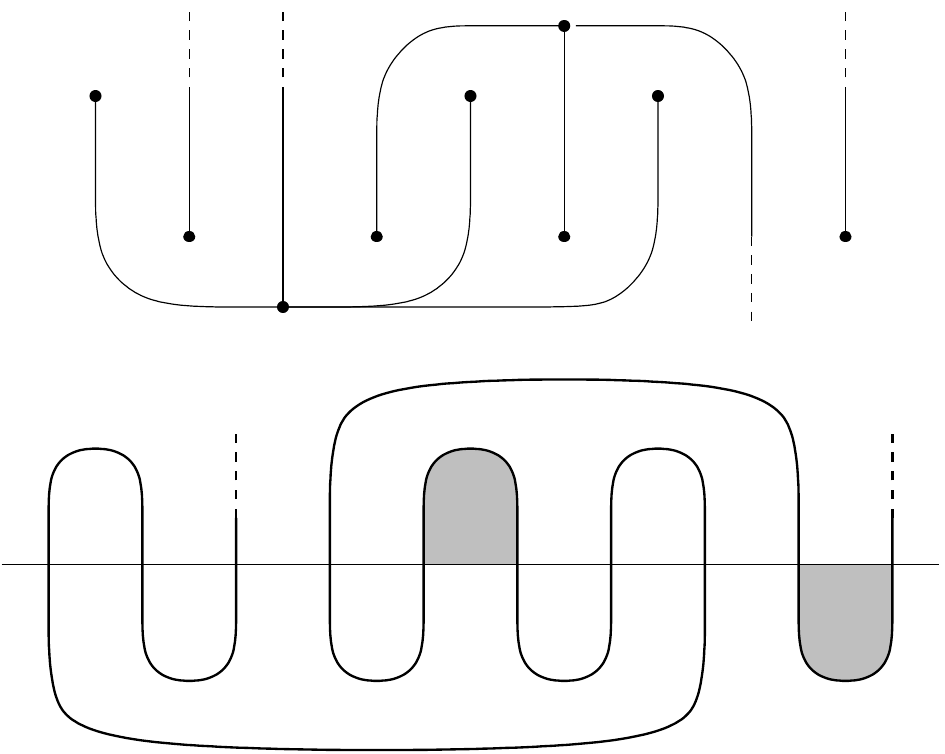}%
\end{picture}%
\setlength{\unitlength}{2960sp}%
\begingroup\makeatletter\ifx\SetFigFont\undefined%
\gdef\SetFigFont#1#2#3#4#5{%
  \reset@font\fontsize{#1}{#2pt}%
  \fontfamily{#3}\fontseries{#4}\fontshape{#5}%
  \selectfont}%
\fi\endgroup%
\begin{picture}(6024,4820)(-911,-6519)
\put(2176,-2386){\makebox(0,0)[lb]{\smash{{\SetFigFont{9}{10.8}{\rmdefault}{\mddefault}{\updefault}{\color[rgb]{0,0,0}$x_p$}%
}}}}
\put(4576,-3286){\makebox(0,0)[lb]{\smash{{\SetFigFont{9}{10.8}{\rmdefault}{\mddefault}{\updefault}{\color[rgb]{0,0,0}$x_n$}%
}}}}
\put(4501,-5686){\makebox(0,0)[b]{\smash{{\SetFigFont{9}{10.8}{\rmdefault}{\mddefault}{\updefault}{\color[rgb]{0,0,0}$D_n$}%
}}}}
\put(2101,-5086){\makebox(0,0)[b]{\smash{{\SetFigFont{9}{10.8}{\rmdefault}{\mddefault}{\updefault}{\color[rgb]{0,0,0}$D_p$}%
}}}}
\put(-824,-3586){\makebox(0,0)[lb]{\smash{{\SetFigFont{9}{10.8}{\rmdefault}{\mddefault}{\updefault}{\color[rgb]{0,0,0}(a)}%
}}}}
\put(-749,-5236){\makebox(0,0)[rb]{\smash{{\SetFigFont{9}{10.8}{\rmdefault}{\mddefault}{\updefault}{\color[rgb]{0,0,0}$\hatT$}%
}}}}
\put(-824,-6286){\makebox(0,0)[lb]{\smash{{\SetFigFont{9}{10.8}{\rmdefault}{\mddefault}{\updefault}{\color[rgb]{0,0,0}(b)}%
}}}}
\end{picture}%
\caption{(a) Parts of $\T$ with $x_n$ and $x_p$ \qua  (b) The corresponding schematic of $A$ in $X$ illustrating $D_n$ and $D_p$}
\label{fig:oppleaves}
\end{figure}

\begin{proof}
Assume otherwise.  
Let $x_p$ and $x_n$ be two leaves of $\T$ of signs $+$ and $-$ respectively that are not in $\T_*$.  Assume neither is $x_0$.  The two vertices cannot be joined by an edge since $\T$ is a tree containing the vertex $x_0$.  Therefore the edges $t_p$ and $t_n$ incident to $x_p$ and $x_n$ respectively are distinct.  Let $X_p, X_n, T_p,$ and $T_n$ be the corresponding solid tori and annuli.  Hence the annuli $T_n \subseteq X_n$ and $T_p \subseteq X_p$ are distinct.  By \fullref{leafisotopy}, there is a high disk $D_p$ such that $D_p \cap \hatT \subseteq T_p$ and a low disk $D_n$ such that $D_n \cap \hatT \subseteq T_n$.  For example, see \fullref{fig:oppleaves}.  Since $K$ is in thin position this contradicts \fullref{highdisklowdisk}.
\end{proof}

\begin{lemma} \label{penultimate}
If $\T$ has a penultimate leaf not in $\T_*$, then the vertex $x_0$ must be among the penultimate leaf and leaves to which it is adjacent.
\end{lemma}

\begin{proof}
Let $x_p$ be a penultimate leaf of $\T$ not in $\T_*$ and let $x_{l_1}, \dots, x_{l_{n}}$ be the leaves adjacent to it.  Assume $x_0$ is none of these.  Also recall that $x_1$ cannot be a leaf unless $k=1$.

Let $T_{l_i}$, for $i = 1, \dots, n$, be the annulus on $\hatT$ corresponding to the edge of $\T$ which connects $x_{l_i}$ to $x_p$.  Let $T_p$ be the remaining annulus of $\hatT \cap X_p$.  Note that $T_p$ corresponds to the edge $t_p$ which separates $\{x_p, x_{l_1}, \dots, x_{l_n} \}$ from the rest of the vertices of $\T$ and from $\T_*$ in particular.  Since meridional disks of the $X_{l_i}$ may be isotoped to intersect the cores of the $T_{l_i}$ once, the manifold $X_p'$ obtained by joining each $X_{l_i}$ to $X_p$ along $T_{l_i}$ is itself a solid torus.

{\bf Case 1}\qua  Assume $x_p \neq x_1$.

In this case $\bdry X_p'$ is composed of the annulus $T_p$ and a subannulus $A'= \bigcup_{i=m}^{m+2n} A_i$ of $A$ for some $m \geq 2$.  This annulus contains the arc $K_{(m-1,\, m+2n)}$ and the arc $K_{(t-m-2n+1,\, t-m+2)}$.  See \fullref{fig:penultimatepnot1}.

\begin{figure}[ht!]
\centering
\begin{picture}(0,0)%
\includegraphics{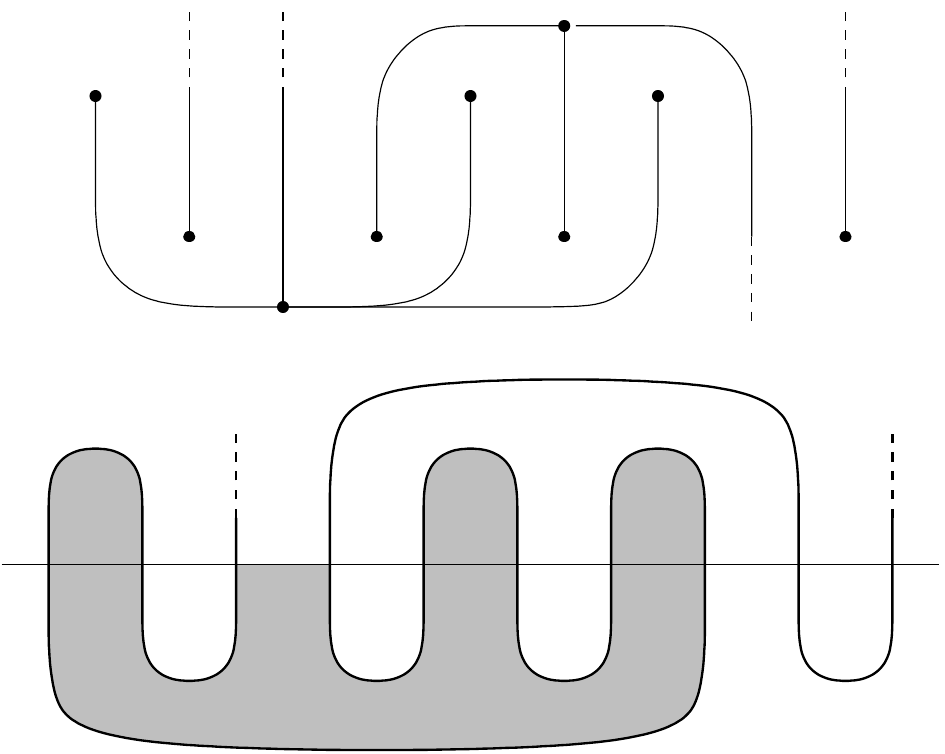}%
\end{picture}%
\setlength{\unitlength}{2960sp}%
\begingroup\makeatletter\ifx\SetFigFont\undefined%
\gdef\SetFigFont#1#2#3#4#5{%
  \reset@font\fontsize{#1}{#2pt}%
  \fontfamily{#3}\fontseries{#4}\fontshape{#5}%
  \selectfont}%
\fi\endgroup%
\begin{picture}(6024,4820)(1189,-6219)
\put(3076,-3586){\makebox(0,0)[lb]{\smash{{\SetFigFont{9}{10.8}{\rmdefault}{\mddefault}{\updefault}{\color[rgb]{0,0,0}$x_p$}%
}}}}
\put(1726,-1861){\makebox(0,0)[lb]{\smash{{\SetFigFont{9}{10.8}{\rmdefault}{\mddefault}{\updefault}{\color[rgb]{0,0,0}$x_{l_1}$}%
}}}}
\put(4126,-1861){\makebox(0,0)[lb]{\smash{{\SetFigFont{9}{10.8}{\rmdefault}{\mddefault}{\updefault}{\color[rgb]{0,0,0}$x_{l_2}$}%
}}}}
\put(5326,-1861){\makebox(0,0)[lb]{\smash{{\SetFigFont{9}{10.8}{\rmdefault}{\mddefault}{\updefault}{\color[rgb]{0,0,0}$x_{l_3}$}%
}}}}
\put(3076,-2611){\makebox(0,0)[lb]{\smash{{\SetFigFont{9}{10.8}{\rmdefault}{\mddefault}{\updefault}{\color[rgb]{0,0,0}$t_p$}%
}}}}
\put(1876,-2611){\makebox(0,0)[lb]{\smash{{\SetFigFont{9}{10.8}{\rmdefault}{\mddefault}{\updefault}{\color[rgb]{0,0,0}$t_{l_1}$}%
}}}}
\put(4276,-2611){\makebox(0,0)[lb]{\smash{{\SetFigFont{9}{10.8}{\rmdefault}{\mddefault}{\updefault}{\color[rgb]{0,0,0}$t_{l_2}$}%
}}}}
\put(5476,-2611){\makebox(0,0)[lb]{\smash{{\SetFigFont{9}{10.8}{\rmdefault}{\mddefault}{\updefault}{\color[rgb]{0,0,0}$t_{l_3}$}%
}}}}
\put(1801,-5236){\makebox(0,0)[b]{\smash{{\SetFigFont{9}{10.8}{\rmdefault}{\mddefault}{\updefault}{\color[rgb]{0,0,0}$T_{l_1}$}%
}}}}
\put(3001,-5236){\makebox(0,0)[b]{\smash{{\SetFigFont{9}{10.8}{\rmdefault}{\mddefault}{\updefault}{\color[rgb]{0,0,0}$T_p$}%
}}}}
\put(4201,-5236){\makebox(0,0)[b]{\smash{{\SetFigFont{9}{10.8}{\rmdefault}{\mddefault}{\updefault}{\color[rgb]{0,0,0}$T_{l_2}$}%
}}}}
\put(5401,-5236){\makebox(0,0)[b]{\smash{{\SetFigFont{9}{10.8}{\rmdefault}{\mddefault}{\updefault}{\color[rgb]{0,0,0}$T_{l_3}$}%
}}}}
\put(1801,-4711){\makebox(0,0)[b]{\smash{{\SetFigFont{9}{10.8}{\rmdefault}{\mddefault}{\updefault}{\color[rgb]{0,0,0}$X_{l_1}$}%
}}}}
\put(4201,-4711){\makebox(0,0)[b]{\smash{{\SetFigFont{9}{10.8}{\rmdefault}{\mddefault}{\updefault}{\color[rgb]{0,0,0}$X_{l_2}$}%
}}}}
\put(5401,-4711){\makebox(0,0)[b]{\smash{{\SetFigFont{9}{10.8}{\rmdefault}{\mddefault}{\updefault}{\color[rgb]{0,0,0}$X_{l_3}$}%
}}}}
\put(3076,-5986){\makebox(0,0)[lb]{\smash{{\SetFigFont{9}{10.8}{\rmdefault}{\mddefault}{\updefault}{\color[rgb]{0,0,0}$X_p$}%
}}}}
\put(1201,-5986){\makebox(0,0)[lb]{\smash{{\SetFigFont{9}{10.8}{\rmdefault}{\mddefault}{\updefault}{\color[rgb]{0,0,0}(b)}%
}}}}
\put(1201,-3286){\makebox(0,0)[lb]{\smash{{\SetFigFont{9}{10.8}{\rmdefault}{\mddefault}{\updefault}{\color[rgb]{0,0,0}(a)}%
}}}}
\end{picture}%
\caption{(a) Part of $\T$ where $x_p$ is a penultimate leaf not in $\T_*$ \qua  (b)  The corresponding schematic of $A$ in $X$ illustrating $X_p'$ and its meridional long disk}
\label{fig:penultimatepnot1}
\end{figure}
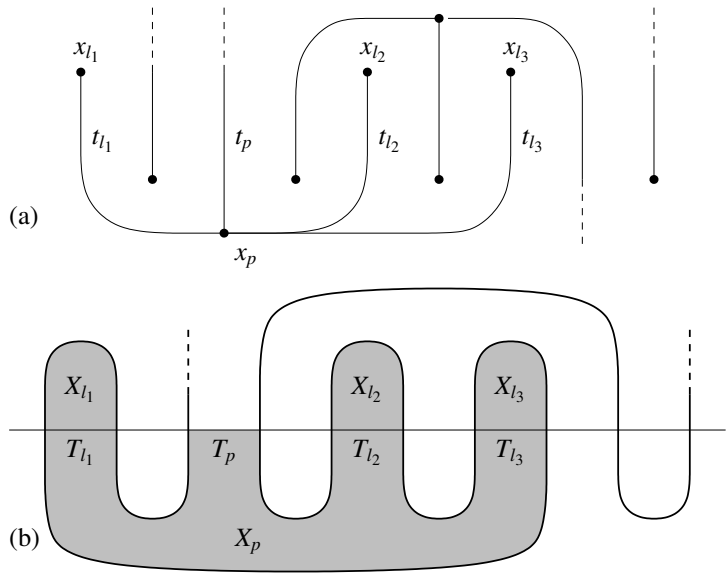

Since $X_p \neq X_0$ or $X_1$, meridional disks may be isotoped to intersect the cores of each annulus of $\hatT \cap X_p$ exactly once.  Thus a meridional disk of $X_p'$ may be isotoped to intersect the core of $T_p$ once and intersect $A'$ in the arc $K_{(m-1,\, m+2n)}$ or $K_{(t-m-2n+1,\, t-m+2)}$.  This meridional disk is a long disk.  Its existence contradicts \fullref{longdisk}.

{\bf Case 2}\qua Assume $x_p = x_1$.

A meridional disk of $X_p$ may be isotoped to intersect the core of each annulus of $\hatT \cap X_p$ two or three times depending on the order of the Scharlemann cycle at hand, so the method of Case 1 does not apply.  See \fullref{fig:penultimatepis1}(a) and (b).

Nevertheless, in this case the boundary $\bdry X_p'$ is composed of the annulus $T_p$ and a subannulus $A'= \smash{\bigcup_{i=2}^{2n+1}} A_i \subseteq A$.  This annulus contains the arcs $K_{(1, \,2n+1)}$ and $K_{(t-2n,\,t)}$.
 Because $x_1 = x_p \in \T \cut \T_*$, the only arc of $K$ in the interior of $X_p'$ is $K_{(t,1)}$.  Recall that $K_{(t,1)} \subseteq A_1$.

Let $D'$ be a meridional disk of $X_p'$.   See \fullref{fig:penultimatepis1}(b) and (c) for illustrations.  Isotop $K_{(t-2n,\, 2n+1)}$ within $A' \cup A_1$ to lie on $D'$.   This puts $K_{(t-2n,\,t)} \cup K_{(1, \,2n+1)} \subseteq \bdry D'$ and $K_{(t, 1)}$ as a properly embedded arc on $D'$.  Let $D$ be one of the two components of $D' \cut K_{(t, 1)}$ so that $\bdry D$ is the union of one arc on $\hatT$ and an arc of $K$.  See \fullref{fig:penultimatepis1}(d).  The disk $D$ is a long disk.  Since $K$ is in thin position, its existence contradicts \fullref{longdisk}.
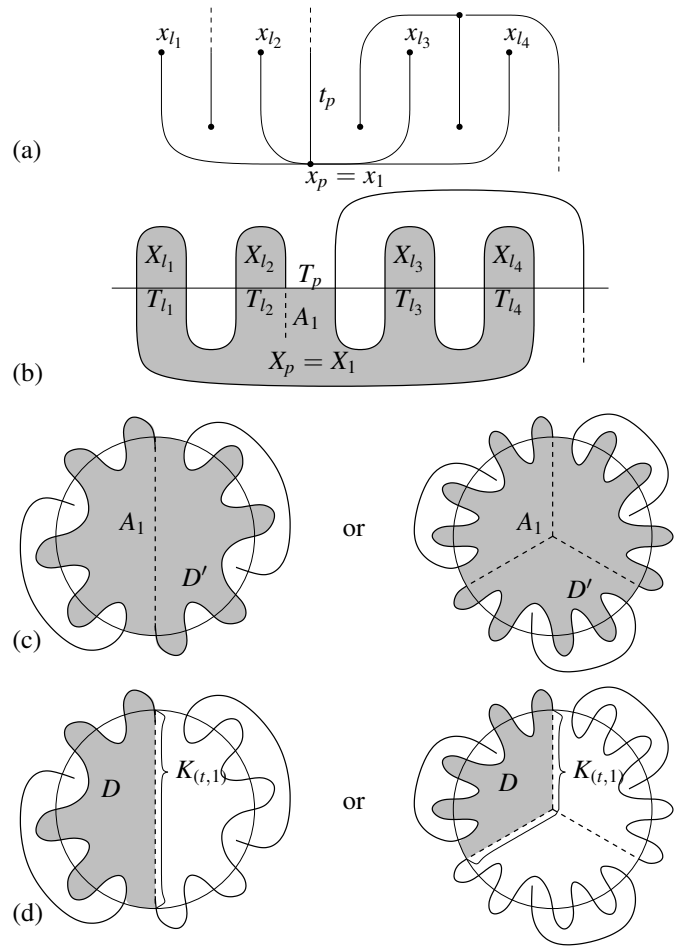
\begin{figure}[ht!]
\centering
\begin{picture}(0,0)%
\includegraphics[scale=.80]{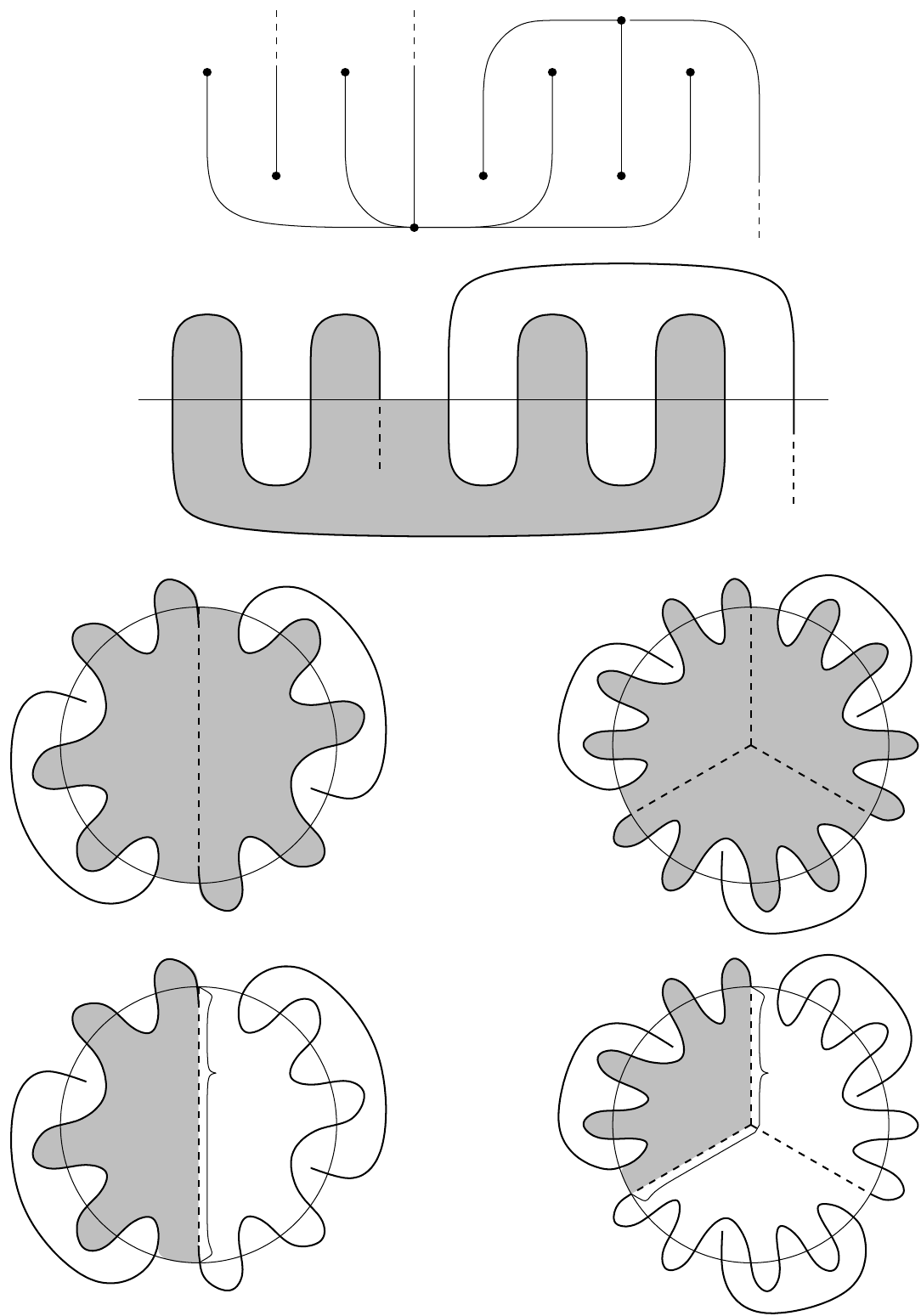}%
\end{picture}%
\setlength{\unitlength}{2052sp}%
\begingroup\makeatletter\ifx\SetFigFont\undefined%
\gdef\SetFigFont#1#2#3#4#5{%
  \reset@font\fontsize{#1}{#2pt}%
  \fontfamily{#3}\fontseries{#4}\fontshape{#5}%
  \selectfont}%
\fi\endgroup%
\begin{picture}(8000,11422)(676,-16121)
\put(4351,-5911){\makebox(0,0)[lb]{\smash{{\SetFigFont{8}{9.6}{\rmdefault}{\mddefault}{\updefault}{\color[rgb]{0,0,0}$t_p$}%
}}}}
\put(2401,-5161){\makebox(0,0)[lb]{\smash{{\SetFigFont{8}{9.6}{\rmdefault}{\mddefault}{\updefault}{\color[rgb]{0,0,0}$x_{l_1}$}%
}}}}
\put(3601,-5161){\makebox(0,0)[lb]{\smash{{\SetFigFont{8}{9.6}{\rmdefault}{\mddefault}{\updefault}{\color[rgb]{0,0,0}$x_{l_2}$}%
}}}}
\put(5401,-5161){\makebox(0,0)[lb]{\smash{{\SetFigFont{8}{9.6}{\rmdefault}{\mddefault}{\updefault}{\color[rgb]{0,0,0}$x_{l_3}$}%
}}}}
\put(6601,-5161){\makebox(0,0)[lb]{\smash{{\SetFigFont{8}{9.6}{\rmdefault}{\mddefault}{\updefault}{\color[rgb]{0,0,0}$x_{l_4}$}%
}}}}
\put(4201,-6886){\makebox(0,0)[lb]{\smash{{\SetFigFont{8}{9.6}{\rmdefault}{\mddefault}{\updefault}{\color[rgb]{0,0,0}$x_p=x_1$}%
}}}}
\put(2476,-8386){\makebox(0,0)[b]{\smash{{\SetFigFont{8}{9.6}{\rmdefault}{\mddefault}{\updefault}{\color[rgb]{0,0,0}$T_{l_1}$}%
}}}}
\put(3751,-9136){\makebox(0,0)[lb]{\smash{{\SetFigFont{8}{9.6}{\rmdefault}{\mddefault}{\updefault}{\color[rgb]{0,0,0}$X_p=X_1$}%
}}}}
\put(2476,-7861){\makebox(0,0)[b]{\smash{{\SetFigFont{8}{9.6}{\rmdefault}{\mddefault}{\updefault}{\color[rgb]{0,0,0}$X_{l_1}$}%
}}}}
\put(4051,-8611){\makebox(0,0)[lb]{\smash{{\SetFigFont{8}{9.6}{\rmdefault}{\mddefault}{\updefault}{\color[rgb]{0,0,0}$A_1$}%
}}}}
\put(6676,-8386){\makebox(0,0)[b]{\smash{{\SetFigFont{8}{9.6}{\rmdefault}{\mddefault}{\updefault}{\color[rgb]{0,0,0}$T_{l_4}$}%
}}}}
\put(6676,-7861){\makebox(0,0)[b]{\smash{{\SetFigFont{8}{9.6}{\rmdefault}{\mddefault}{\updefault}{\color[rgb]{0,0,0}$X_{l_4}$}%
}}}}
\put(4276,-8086){\makebox(0,0)[b]{\smash{{\SetFigFont{8}{9.6}{\rmdefault}{\mddefault}{\updefault}{\color[rgb]{0,0,0}$T_p$}%
}}}}
\put(3676,-8386){\makebox(0,0)[b]{\smash{{\SetFigFont{8}{9.6}{\rmdefault}{\mddefault}{\updefault}{\color[rgb]{0,0,0}$T_{l_2}$}%
}}}}
\put(3676,-7861){\makebox(0,0)[b]{\smash{{\SetFigFont{8}{9.6}{\rmdefault}{\mddefault}{\updefault}{\color[rgb]{0,0,0}$X_{l_2}$}%
}}}}
\put(5476,-8386){\makebox(0,0)[b]{\smash{{\SetFigFont{8}{9.6}{\rmdefault}{\mddefault}{\updefault}{\color[rgb]{0,0,0}$T_{l_3}$}%
}}}}
\put(5476,-7861){\makebox(0,0)[b]{\smash{{\SetFigFont{8}{9.6}{\rmdefault}{\mddefault}{\updefault}{\color[rgb]{0,0,0}$X_{l_3}$}%
}}}}
\put(7126,-11086){\makebox(0,0)[rb]{\smash{{\SetFigFont{8}{9.6}{\rmdefault}{\mddefault}{\updefault}{\color[rgb]{0,0,0}$A_1$}%
}}}}
\put(7351,-11911){\makebox(0,0)[lb]{\smash{{\SetFigFont{8}{9.6}{\rmdefault}{\mddefault}{\updefault}{\color[rgb]{0,0,0}$D'$}%
}}}}
\put(2626,-14086){\makebox(0,0)[lb]{\smash{{\SetFigFont{8}{9.6}{\rmdefault}{\mddefault}{\updefault}{\color[rgb]{0,0,0}$K_{(t, 1)}$}%
}}}}
\put(1726,-14311){\makebox(0,0)[lb]{\smash{{\SetFigFont{8}{9.6}{\rmdefault}{\mddefault}{\updefault}{\color[rgb]{0,0,0}$D$}%
}}}}
\put(4801,-11161){\makebox(0,0)[b]{\smash{{\SetFigFont{8}{9.6}{\rmdefault}{\mddefault}{\updefault}{\color[rgb]{0,0,0}or}%
}}}}
\put(4801,-14461){\makebox(0,0)[b]{\smash{{\SetFigFont{8}{9.6}{\rmdefault}{\mddefault}{\updefault}{\color[rgb]{0,0,0}or}%
}}}}
\put(2701,-11761){\makebox(0,0)[lb]{\smash{{\SetFigFont{8}{9.6}{\rmdefault}{\mddefault}{\updefault}{\color[rgb]{0,0,0}$D'$}%
}}}}
\put(7426,-14086){\makebox(0,0)[lb]{\smash{{\SetFigFont{8}{9.6}{\rmdefault}{\mddefault}{\updefault}{\color[rgb]{0,0,0}$K_{(t, 1)}$}%
}}}}
\put(6526,-14236){\makebox(0,0)[lb]{\smash{{\SetFigFont{8}{9.6}{\rmdefault}{\mddefault}{\updefault}{\color[rgb]{0,0,0}$D$}%
}}}}
\put(676,-12511){\makebox(0,0)[lb]{\smash{{\SetFigFont{8}{9.6}{\rmdefault}{\mddefault}{\updefault}{\color[rgb]{0,0,0}(c)}%
}}}}
\put(676,-15811){\makebox(0,0)[lb]{\smash{{\SetFigFont{8}{9.6}{\rmdefault}{\mddefault}{\updefault}{\color[rgb]{0,0,0}(d)}%
}}}}
\put(676,-9286){\makebox(0,0)[lb]{\smash{{\SetFigFont{8}{9.6}{\rmdefault}{\mddefault}{\updefault}{\color[rgb]{0,0,0}(b)}%
}}}}
\put(676,-6586){\makebox(0,0)[lb]{\smash{{\SetFigFont{8}{9.6}{\rmdefault}{\mddefault}{\updefault}{\color[rgb]{0,0,0}(a)}%
}}}}
\put(2326,-11086){\makebox(0,0)[rb]{\smash{{\SetFigFont{8}{9.6}{\rmdefault}{\mddefault}{\updefault}{\color[rgb]{0,0,0}$A_1$}%
}}}}
\end{picture}%
\caption{(a) Part of $\T$ where $x_1$ is a penultimate leaf not in $\T_*$ \qua (b)  The corresponding schematic of $A$ in $X$ illustrating $X_p'$ and $D'$.  (c) ``Unquotiented'' pictures of the meridional disk $D'$ of $X_p'$ when $A_1$ is constructed from an order $2$ or $3$ Scharlemann cycle.  (d)  The long disk $D$ as a subdisk of $D'$ with $K_{(t, 1)}$ as an arc of $A_1 \cap D'$}
\label{fig:penultimatepis1}
\end{figure}
\end{proof}

\subsection{The unfurling isotopy and further constraints on the trees}

Let $R$ be a torus in a $3$--manifold $Y$ that bounds a solid torus $V$.  Because $R$ bounds the solid torus $V$, any Dehn twist along $R$ is isotopic to the identity.

Assume $R$ is divided into two annuli $T_R$ and $A_R$ by the two parallel curves  $a_{\rm out}$ and $a_{\rm in}$.  Let $K$ be a knot in $Y$ so that $K \cap R$ is a collection of spanning arcs of $A_R$ with $K \cap N(a_{\rm out}) \subseteq V$ and $K \cap N(a_{\rm in}) \subseteq Y - \Int V$.

\begin{Lemma} \label{main-unfurl} 
There is an ambient isotopy $\Theta_u \co Y \times [0,2\pi] \to Y$ so that 
\begin{enumerate}
\item $\Theta_0$ is the identity,\vspace{-2pt}
\item $\Theta_{2\pi} (K-R) = \Theta_0 (K-R)$ and\vspace{-2pt}
\item $\Theta_{2\pi}(K \cap R)$ is a collection of transverse arcs of $T_R$.\vspace{-2pt}
\end{enumerate}
\end{Lemma}

\begin{figure}[ht!]
\centering
\begin{picture}(0,0)%
\includegraphics{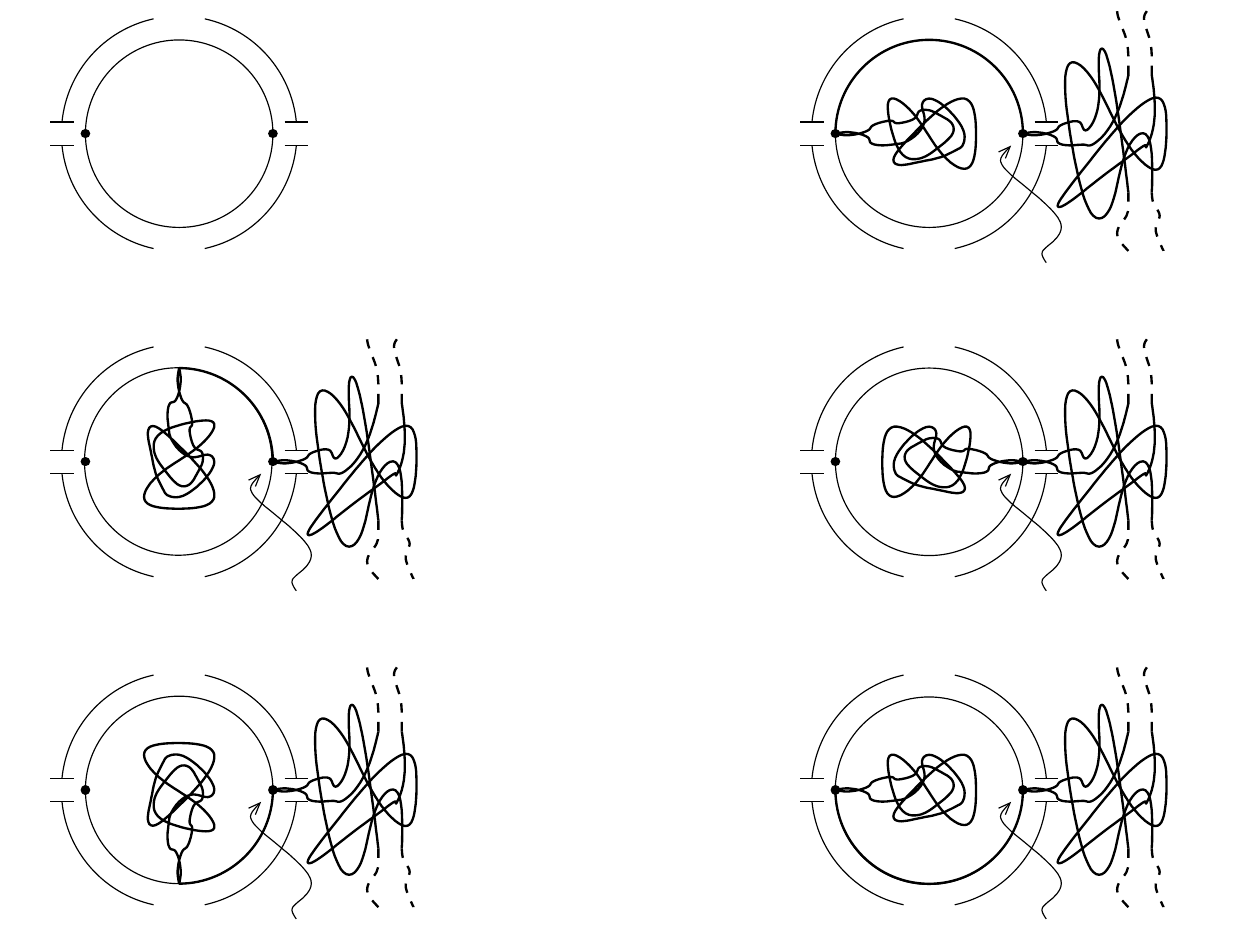}%
\end{picture}%
\setlength{\unitlength}{2960sp}%
\begingroup\makeatletter\ifx\SetFigFont\undefined%
\gdef\SetFigFont#1#2#3#4#5{%
  \reset@font\fontsize{#1}{#2pt}%
  \fontfamily{#3}\fontseries{#4}\fontshape{#5}%
  \selectfont}%
\fi\endgroup%
\begin{picture}(7920,6080)(1854,-7394)
\put(2251,-7111){\makebox(0,0)[lb]{\smash{{\SetFigFont{9}{10.8}{\rmdefault}{\mddefault}{\updefault}{\color[rgb]{0,0,0}$\Theta_{3\pi/2}$}%
}}}}
\put(3001,-2986){\makebox(0,0)[b]{\smash{{\SetFigFont{9}{10.8}{\rmdefault}{\mddefault}{\updefault}{\color[rgb]{0,0,0}$T_R$}%
}}}}
\put(3001,-1486){\makebox(0,0)[b]{\smash{{\SetFigFont{9}{10.8}{\rmdefault}{\mddefault}{\updefault}{\color[rgb]{0,0,0}$A_R$}%
}}}}
\put(3001,-2236){\makebox(0,0)[b]{\smash{{\SetFigFont{9}{10.8}{\rmdefault}{\mddefault}{\updefault}{\color[rgb]{0,0,0}$V$}%
}}}}
\put(2176,-2236){\makebox(0,0)[rb]{\smash{{\SetFigFont{9}{10.8}{\rmdefault}{\mddefault}{\updefault}{\color[rgb]{0,0,0}$a_{\mbox{in}}$}%
}}}}
\put(3826,-2236){\makebox(0,0)[lb]{\smash{{\SetFigFont{9}{10.8}{\rmdefault}{\mddefault}{\updefault}{\color[rgb]{0,0,0}$a_{\mbox{out}}$}%
}}}}
\put(3001,-5086){\makebox(0,0)[b]{\smash{{\SetFigFont{9}{10.8}{\rmdefault}{\mddefault}{\updefault}{\color[rgb]{0,0,0}$T_R$}%
}}}}
\put(3001,-3586){\makebox(0,0)[b]{\smash{{\SetFigFont{9}{10.8}{\rmdefault}{\mddefault}{\updefault}{\color[rgb]{0,0,0}$A_R$}%
}}}}
\put(2176,-4336){\makebox(0,0)[rb]{\smash{{\SetFigFont{9}{10.8}{\rmdefault}{\mddefault}{\updefault}{\color[rgb]{0,0,0}$a_{\mbox{in}}$}%
}}}}
\put(3751,-5236){\makebox(0,0)[lb]{\smash{{\SetFigFont{9}{10.8}{\rmdefault}{\mddefault}{\updefault}{\color[rgb]{0,0,0}$a_{\mbox{out}}$}%
}}}}
\put(4651,-4261){\makebox(0,0)[lb]{\smash{{\SetFigFont{9}{10.8}{\rmdefault}{\mddefault}{\updefault}{\color[rgb]{0,0,0}$K$}%
}}}}
\put(3001,-7186){\makebox(0,0)[b]{\smash{{\SetFigFont{9}{10.8}{\rmdefault}{\mddefault}{\updefault}{\color[rgb]{0,0,0}$T_R$}%
}}}}
\put(3001,-5686){\makebox(0,0)[b]{\smash{{\SetFigFont{9}{10.8}{\rmdefault}{\mddefault}{\updefault}{\color[rgb]{0,0,0}$A_R$}%
}}}}
\put(2176,-6436){\makebox(0,0)[rb]{\smash{{\SetFigFont{9}{10.8}{\rmdefault}{\mddefault}{\updefault}{\color[rgb]{0,0,0}$a_{\mbox{in}}$}%
}}}}
\put(3751,-7336){\makebox(0,0)[lb]{\smash{{\SetFigFont{9}{10.8}{\rmdefault}{\mddefault}{\updefault}{\color[rgb]{0,0,0}$a_{\mbox{out}}$}%
}}}}
\put(4651,-6361){\makebox(0,0)[lb]{\smash{{\SetFigFont{9}{10.8}{\rmdefault}{\mddefault}{\updefault}{\color[rgb]{0,0,0}$K$}%
}}}}
\put(7801,-2986){\makebox(0,0)[b]{\smash{{\SetFigFont{9}{10.8}{\rmdefault}{\mddefault}{\updefault}{\color[rgb]{0,0,0}$T_R$}%
}}}}
\put(7801,-1486){\makebox(0,0)[b]{\smash{{\SetFigFont{9}{10.8}{\rmdefault}{\mddefault}{\updefault}{\color[rgb]{0,0,0}$A_R$}%
}}}}
\put(6976,-2236){\makebox(0,0)[rb]{\smash{{\SetFigFont{9}{10.8}{\rmdefault}{\mddefault}{\updefault}{\color[rgb]{0,0,0}$a_{\mbox{in}}$}%
}}}}
\put(8551,-3136){\makebox(0,0)[lb]{\smash{{\SetFigFont{9}{10.8}{\rmdefault}{\mddefault}{\updefault}{\color[rgb]{0,0,0}$a_{\mbox{out}}$}%
}}}}
\put(9451,-2161){\makebox(0,0)[lb]{\smash{{\SetFigFont{9}{10.8}{\rmdefault}{\mddefault}{\updefault}{\color[rgb]{0,0,0}$K$}%
}}}}
\put(7801,-7186){\makebox(0,0)[b]{\smash{{\SetFigFont{9}{10.8}{\rmdefault}{\mddefault}{\updefault}{\color[rgb]{0,0,0}$T_R$}%
}}}}
\put(7801,-5686){\makebox(0,0)[b]{\smash{{\SetFigFont{9}{10.8}{\rmdefault}{\mddefault}{\updefault}{\color[rgb]{0,0,0}$A_R$}%
}}}}
\put(6976,-6436){\makebox(0,0)[rb]{\smash{{\SetFigFont{9}{10.8}{\rmdefault}{\mddefault}{\updefault}{\color[rgb]{0,0,0}$a_{\mbox{in}}$}%
}}}}
\put(8551,-7336){\makebox(0,0)[lb]{\smash{{\SetFigFont{9}{10.8}{\rmdefault}{\mddefault}{\updefault}{\color[rgb]{0,0,0}$a_{\mbox{out}}$}%
}}}}
\put(9451,-6361){\makebox(0,0)[lb]{\smash{{\SetFigFont{9}{10.8}{\rmdefault}{\mddefault}{\updefault}{\color[rgb]{0,0,0}$K$}%
}}}}
\put(8551,-5236){\makebox(0,0)[lb]{\smash{{\SetFigFont{9}{10.8}{\rmdefault}{\mddefault}{\updefault}{\color[rgb]{0,0,0}$a_{\mbox{out}}$}%
}}}}
\put(9451,-4261){\makebox(0,0)[lb]{\smash{{\SetFigFont{9}{10.8}{\rmdefault}{\mddefault}{\updefault}{\color[rgb]{0,0,0}$K$}%
}}}}
\put(7801,-5086){\makebox(0,0)[b]{\smash{{\SetFigFont{9}{10.8}{\rmdefault}{\mddefault}{\updefault}{\color[rgb]{0,0,0}$T_R$}%
}}}}
\put(7801,-3586){\makebox(0,0)[b]{\smash{{\SetFigFont{9}{10.8}{\rmdefault}{\mddefault}{\updefault}{\color[rgb]{0,0,0}$A_R$}%
}}}}
\put(6976,-4336){\makebox(0,0)[rb]{\smash{{\SetFigFont{9}{10.8}{\rmdefault}{\mddefault}{\updefault}{\color[rgb]{0,0,0}$a_{\mbox{in}}$}%
}}}}
\put(2101,-2911){\makebox(0,0)[rb]{\smash{{\SetFigFont{9}{10.8}{\rmdefault}{\mddefault}{\updefault}{\color[rgb]{0,0,0}(a)}%
}}}}
\put(6901,-2911){\makebox(0,0)[rb]{\smash{{\SetFigFont{9}{10.8}{\rmdefault}{\mddefault}{\updefault}{\color[rgb]{0,0,0}(b)}%
}}}}
\put(6901,-5011){\makebox(0,0)[rb]{\smash{{\SetFigFont{9}{10.8}{\rmdefault}{\mddefault}{\updefault}{\color[rgb]{0,0,0}(d)}%
}}}}
\put(2101,-5011){\makebox(0,0)[rb]{\smash{{\SetFigFont{9}{10.8}{\rmdefault}{\mddefault}{\updefault}{\color[rgb]{0,0,0}(c)}%
}}}}
\put(6901,-7111){\makebox(0,0)[rb]{\smash{{\SetFigFont{9}{10.8}{\rmdefault}{\mddefault}{\updefault}{\color[rgb]{0,0,0}(f)}%
}}}}
\put(2101,-7111){\makebox(0,0)[rb]{\smash{{\SetFigFont{9}{10.8}{\rmdefault}{\mddefault}{\updefault}{\color[rgb]{0,0,0}(e)}%
}}}}
\put(2251,-5011){\makebox(0,0)[lb]{\smash{{\SetFigFont{9}{10.8}{\rmdefault}{\mddefault}{\updefault}{\color[rgb]{0,0,0}$\Theta_{\pi/2}$}%
}}}}
\put(7051,-5011){\makebox(0,0)[lb]{\smash{{\SetFigFont{9}{10.8}{\rmdefault}{\mddefault}{\updefault}{\color[rgb]{0,0,0}$\Theta_{\pi}$}%
}}}}
\put(7051,-7111){\makebox(0,0)[lb]{\smash{{\SetFigFont{9}{10.8}{\rmdefault}{\mddefault}{\updefault}{\color[rgb]{0,0,0}$\Theta_{2\pi}$}%
}}}}
\end{picture}%
\caption{(a) A (not necessarily meridional) ``cross section'' of the solid torus $V$ \qua  (b) The cross section shown with $K$ at $u=0$\qua  (c)  The isotopy at $u=\pi/2$\qua  (d)  The isotopy at $u=\pi$\qua (e) The isotopy at $u=3\pi/2$ \qua(f) The finished isotopy at $u=2\pi$}
\label{fig:unfurlingschematic}
\end{figure}
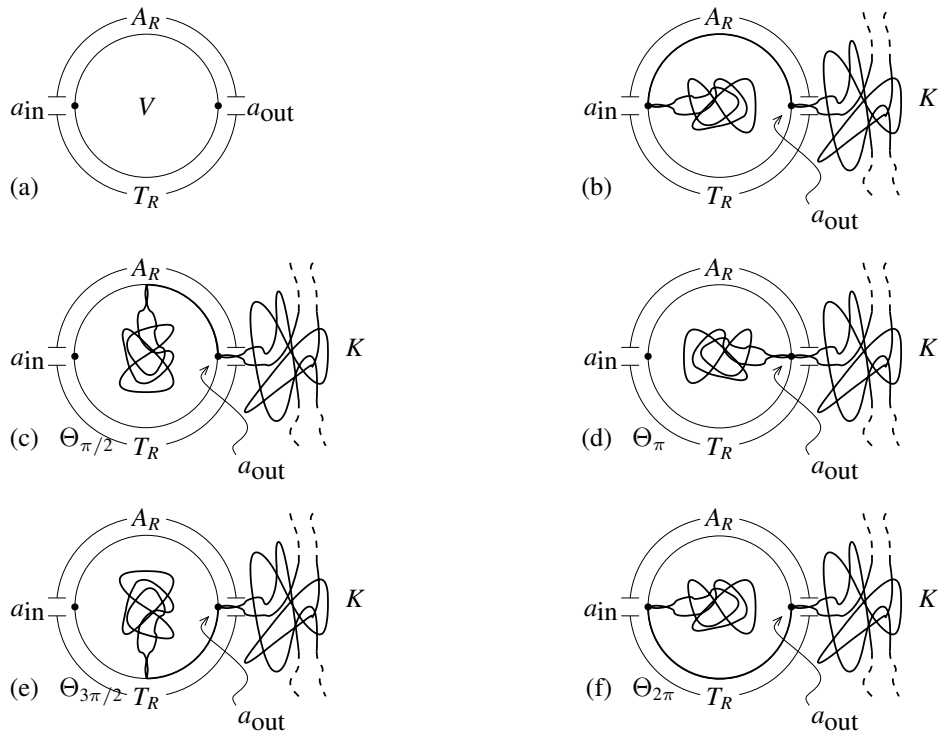

We refer to an ambient isotopy of $K$ via such a Dehn twist as an {\em unfurling\/}.  This unfurling isotopy is shown schematically in \fullref{fig:unfurlingschematic}.  Note that $a_{\rm out}$ and $a_{\rm in}$ need not be longitudinal curves on $R =\bdry V$.

\begin{proof}
Let $\gamma_A$ be an arc of $A_R \cap K$, $\gamma_T$ be an arc of $T_R$ connecting the endpoints of $\gamma_A$, and $\gamma=\gamma_A \cup \gamma_T$ be the simple closed curve oriented so that $\gamma_A$ runs from $a_{\rm out}$ to $a_{\rm in}$.  Take $\Theta_u$ to be the Dehn twist along $R$ in the direction of $\gamma$.  We may view \fullref{fig:unfurlingschematic} as a ``cross section'' of $V$  along $\gamma$.  The conclusions of the lemma are immediate following the definition of the Dehn twist.  
\end{proof}

Let $T_R$ be an annulus of $\hatT \cut A$.  Then we have $\bdry T_R = a_m
\cup a_n$ where $m<n$.  Let $A_R = \smash{\bigcup_{i=m+1}^n A_{i}}$.
\begin{lemma} \label{dividingtorus}
The torus $R = A_R \cup T_R$ bounds a solid torus in $X$ that does not contain $X_0$.
\end{lemma}

\begin{proof}
Since $R$ is a torus in a lens space, it separates $X$ into two pieces, say $V$ and $W$.  Deleting the edge $t_R$ corresponding to the annulus $T_R$ from $\T$ yields two trees $\T_V$ and $\T_W$ each consisting of vertices $x_i$ corresponding to solid tori $X_i \in \mathcal{X}$ in $V$ and $W$ respectively.

Since only one of $V$ and $W$ may contain $X_0$, only one of $\T_V$ and $\T_W$ may contain $x_0$.  Without loss of generality, assume the tree $\T_V$ does not contain the vertex $x_0$.  Though $\T_V$ may contain $x_1$, all other vertices of $\T_V$ then correspond to solid tori $X_i \neq X_0$ or $X_1$ whose meridians traverse each annulus of $X_i \cap \hatT$ (and each annulus of $X_i \cap A$) exactly once.  Therefore the union $X_V$ of the $X_i$ corresponding to $x_i \in \T_V$ along their common annuli $X_i \cap \hatT$ is a solid torus.    Note that $\bdry X_V = R$, and the solid torus $X_V$ is indeed $V$.
\end{proof}

\begin{lemma} \label{unfurling}
 $x_1  \in \T_*$
\end{lemma}

\begin{proof}
Assume $x_1 \not \in \T_*$.
We may assume $k \geq 3$ since if $k \leq 2$ then either $x_1 \in \T_*$ or $t=2$. 

Put a transverse direction on $A$ respecting the ordering of the $A_i$.  For odd (resp.\ even) $i$, the transverse direction at $a_i$ points above (resp.\ below) $\hatT$.  Consider the collection $\A$ of curves $a_i$ on the boundary of $X_1$.  If $k$ is even, then $a_k \not \in \A$ since otherwise $x_1 \in \T_*$.  For each $a_i \in \A$ with $i$ even, the annulus $A_{i+1}$ must be contained in $\bdry X_1$.  Since $a_1 \in \A$, there are more odd indexed curves than even indexed curves in $\A$.  Hence there are two curves $a_m$ and $a_n$ in $\A$ cobounding an annulus $T_R \subseteq \hatT \cut A$ in $\bdry X_1$ with $m<n$ and both $m$ and $n$ odd. 

Consider the annulus $A_R = \bigcup_{i=m+1}^n A_i$ and the torus $R = A_R \cup T_R$.  See \fullref{fig:unfurl1}(a).  By \fullref{dividingtorus}, $R$ bounds a solid torus, say $V$, that does not contain $X_0$.    

  The assumption that $x_1 \not \in \T_*$ is equivalent to the statement that the solid torus $X_1$ is disjoint from the arc $K_{(k,\,t-k+1)}$ of $K$.  Hence $T_R$ and moreover $R$ are disjoint from this arc.   Therefore $K$ only intersects $R$ as two transverse arcs of $A_R$.  At one of the curves $a_m$ and $a_n$, $K$ continues into $V$; at the other curve, $K$ continues away from $V$.  We may apply \fullref{main-unfurl} to obtain an isotopy of $K$. See \fullref{fig:unfurl1}(b). 

\begin{figure}[ht!]
\centering
\begin{picture}(0,0)%
\includegraphics{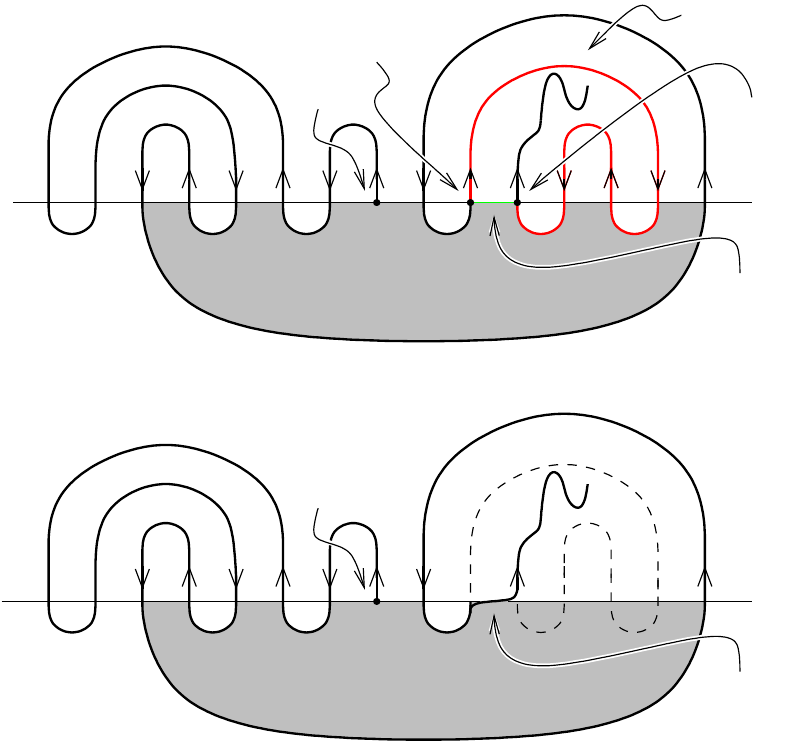}%
\end{picture}%
\setlength{\unitlength}{2960sp}%
\begingroup\makeatletter\ifx\SetFigFont\undefined%
\gdef\SetFigFont#1#2#3#4#5{%
  \reset@font\fontsize{#1}{#2pt}%
  \fontfamily{#3}\fontseries{#4}\fontshape{#5}%
  \selectfont}%
\fi\endgroup%
\begin{picture}(5042,4755)(1189,-6219)
\put(5926,-3436){\makebox(0,0)[b]{\smash{{\SetFigFont{9}{10.8}{\rmdefault}{\mddefault}{\updefault}{\color[rgb]{0,0,0}$T_R$}%
}}}}
\put(3526,-1786){\makebox(0,0)[b]{\smash{{\SetFigFont{9}{10.8}{\rmdefault}{\mddefault}{\updefault}{\color[rgb]{0,0,0}$a_m$}%
}}}}
\put(6001,-2236){\makebox(0,0)[b]{\smash{{\SetFigFont{9}{10.8}{\rmdefault}{\mddefault}{\updefault}{\color[rgb]{0,0,0}$a_n$}%
}}}}
\put(1201,-2686){\makebox(0,0)[b]{\smash{{\SetFigFont{9}{10.8}{\rmdefault}{\mddefault}{\updefault}{\color[rgb]{0,0,0}$\hatT$}%
}}}}
\put(2101,-3436){\makebox(0,0)[b]{\smash{{\SetFigFont{9}{10.8}{\rmdefault}{\mddefault}{\updefault}{\color[rgb]{0,0,0}$A$}%
}}}}
\put(3826,-3361){\makebox(0,0)[b]{\smash{{\SetFigFont{9}{10.8}{\rmdefault}{\mddefault}{\updefault}{\color[rgb]{0,0,0}$X_1$}%
}}}}
\put(5551,-1636){\makebox(0,0)[lb]{\smash{{\SetFigFont{9}{10.8}{\rmdefault}{\mddefault}{\updefault}{\color[rgb]{0,0,0}$A_R$}%
}}}}
\put(5926,-5986){\makebox(0,0)[b]{\smash{{\SetFigFont{9}{10.8}{\rmdefault}{\mddefault}{\updefault}{\color[rgb]{0,0,0}$\Theta_{2\pi}(A_R)$}%
}}}}
\put(1201,-5236){\makebox(0,0)[b]{\smash{{\SetFigFont{9}{10.8}{\rmdefault}{\mddefault}{\updefault}{\color[rgb]{0,0,0}$\hatT$}%
}}}}
\put(2101,-5986){\makebox(0,0)[b]{\smash{{\SetFigFont{9}{10.8}{\rmdefault}{\mddefault}{\updefault}{\color[rgb]{0,0,0}$A$}%
}}}}
\put(3826,-5911){\makebox(0,0)[b]{\smash{{\SetFigFont{9}{10.8}{\rmdefault}{\mddefault}{\updefault}{\color[rgb]{0,0,0}$X_1$}%
}}}}
\put(3226,-2086){\makebox(0,0)[b]{\smash{{\SetFigFont{9}{10.8}{\rmdefault}{\mddefault}{\updefault}{\color[rgb]{0,0,0}$a_1$}%
}}}}
\put(3226,-4636){\makebox(0,0)[b]{\smash{{\SetFigFont{9}{10.8}{\rmdefault}{\mddefault}{\updefault}{\color[rgb]{0,0,0}$a_1$}%
}}}}
\put(1201,-6136){\makebox(0,0)[lb]{\smash{{\SetFigFont{9}{10.8}{\rmdefault}{\mddefault}{\updefault}{\color[rgb]{0,0,0}(b)}%
}}}}
\put(1201,-3586){\makebox(0,0)[lb]{\smash{{\SetFigFont{9}{10.8}{\rmdefault}{\mddefault}{\updefault}{\color[rgb]{0,0,0}(a)}%
}}}}
\end{picture}%
\caption{(a) The construction of the unfurling torus\qua (b) After the unfurling isotopy}
\label{fig:unfurl1}
\end{figure}
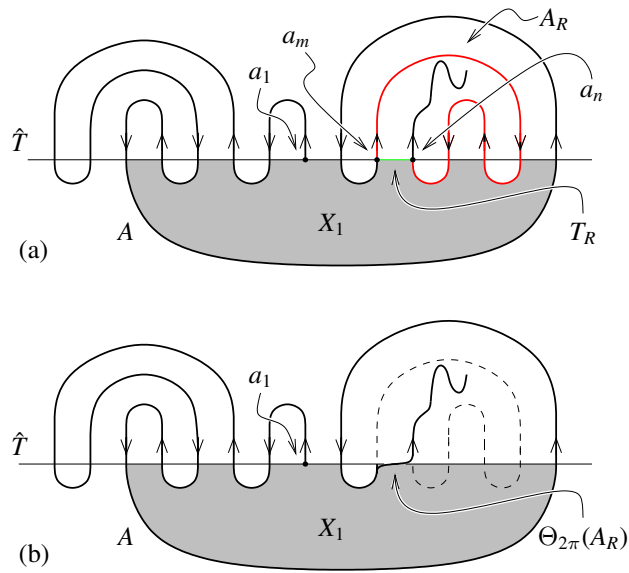

After this unfurling isotopy, the arcs once of $K \cap A_R$ are now arcs of $K \cap T_R$ and may be nudged to be transverse to the height function.  The new set of critical levels is a strict subset of the former set of critical levels.  Furthermore the intersection number of $K$ with any of the remaining critical levels has not increased.  Therefore the unfurling isotopy has decreased the width of $K$, contradicting the thinness of $K$.
\end{proof}

\begin{Lemma} \label{unfurling2}
The vertices of $\T$ not in $\T_*$ are leaves of $\T$ of the same sign.
\end{Lemma}

\begin{figure}[ht!]
\centering
\begin{picture}(0,0)%
\includegraphics[scale=.85]{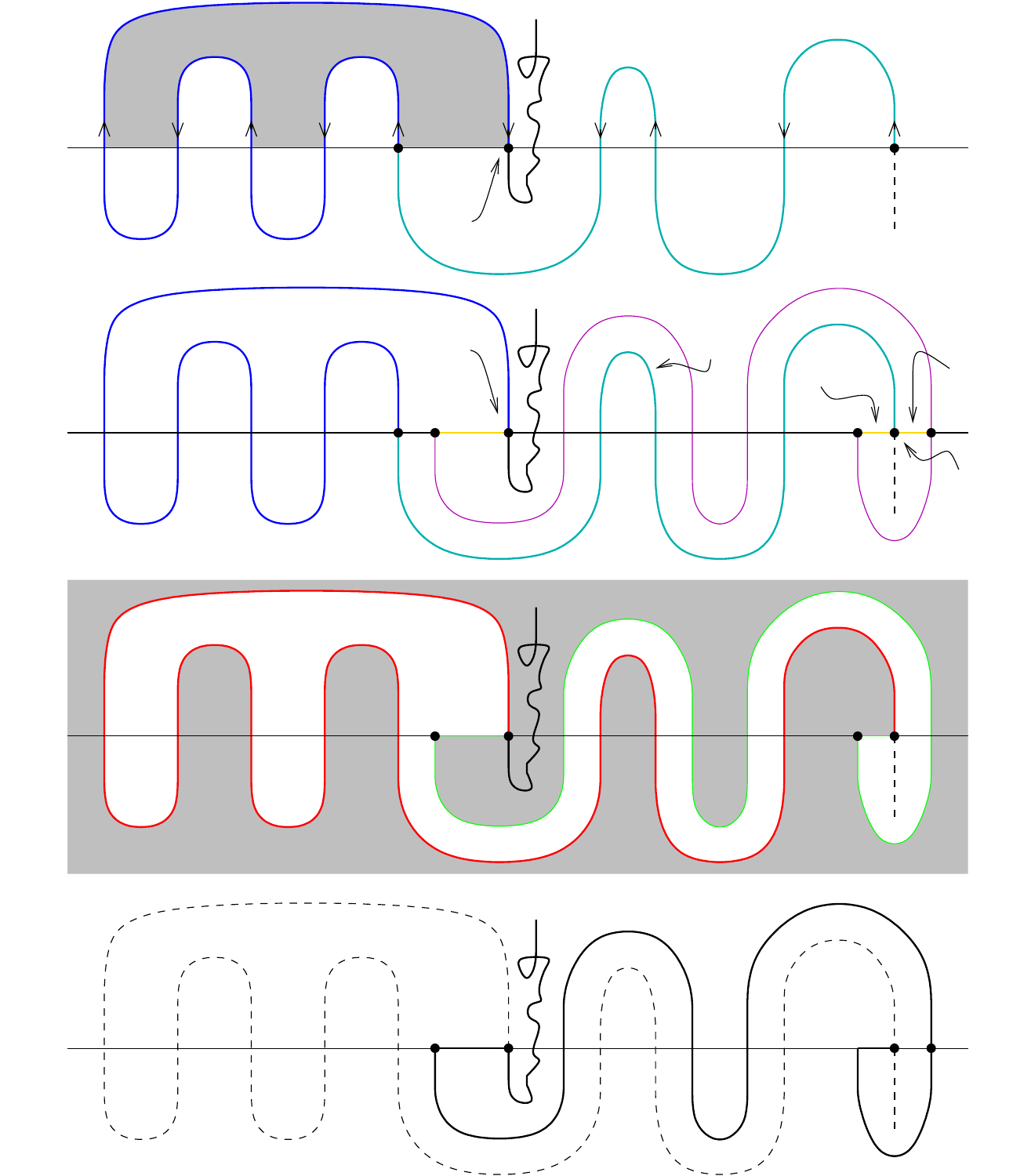}%
\end{picture}%
\setlength{\unitlength}{2516sp}%
\begingroup\makeatletter\ifx\SetFigFont\undefined%
\gdef\SetFigFont#1#2#3#4#5{%
  \reset@font\fontsize{#1}{#2pt}%
  \fontfamily{#3}\fontseries{#4}\fontshape{#5}%
  \selectfont}%
\fi\endgroup%
\begin{picture}(8258,9599)(-550,-8827)
\put(6676,-1036){\makebox(0,0)[rb]{\smash{{\SetFigFont{9}{10.8}{\rmdefault}{\mddefault}{\updefault}{\color[rgb]{0,0,0}$A_1$}%
}}}}
\put(2626,-661){\makebox(0,0)[rb]{\smash{{\SetFigFont{9}{10.8}{\rmdefault}{\mddefault}{\updefault}{\color[rgb]{0,0,0}$a_m$}%
}}}}
\put(6826,-661){\makebox(0,0)[lb]{\smash{{\SetFigFont{9}{10.8}{\rmdefault}{\mddefault}{\updefault}{\color[rgb]{0,0,0}$a_1$}%
}}}}
\put(1801,389){\makebox(0,0)[b]{\smash{{\SetFigFont{9}{10.8}{\rmdefault}{\mddefault}{\updefault}{\color[rgb]{0,0,0}$X_0$}%
}}}}
\put(4801,164){\makebox(0,0)[lb]{\smash{{\SetFigFont{9}{10.8}{\rmdefault}{\mddefault}{\updefault}{\color[rgb]{0,0,0}$A_{(1,\,m)}$}%
}}}}
\put(151,-361){\makebox(0,0)[rb]{\smash{{\SetFigFont{9}{10.8}{\rmdefault}{\mddefault}{\updefault}{\color[rgb]{0,0,0}$\hatT$}%
}}}}
\put(3151,-586){\makebox(0,0)[b]{\smash{{\SetFigFont{9}{10.8}{\rmdefault}{\mddefault}{\updefault}{\color[rgb]{0,0,0}$T_0$}%
}}}}
\put(2626,-2986){\makebox(0,0)[rb]{\smash{{\SetFigFont{9}{10.8}{\rmdefault}{\mddefault}{\updefault}{\color[rgb]{0,0,0}$a_m$}%
}}}}
\put(3301,-2911){\makebox(0,0)[b]{\smash{{\SetFigFont{9}{10.8}{\rmdefault}{\mddefault}{\updefault}{\color[rgb]{0,0,0}$T_0'$}%
}}}}
\put(3001,-2686){\makebox(0,0)[b]{\smash{{\SetFigFont{9}{10.8}{\rmdefault}{\mddefault}{\updefault}{\color[rgb]{0,0,0}$a_m'$}%
}}}}
\put(7126,-2686){\makebox(0,0)[lb]{\smash{{\SetFigFont{9}{10.8}{\rmdefault}{\mddefault}{\updefault}{\color[rgb]{0,0,0}$a_1'$}%
}}}}
\put(6376,-2986){\makebox(0,0)[rb]{\smash{{\SetFigFont{9}{10.8}{\rmdefault}{\mddefault}{\updefault}{\color[rgb]{0,0,0}$a_1''$}%
}}}}
\put(7201,-3211){\makebox(0,0)[lb]{\smash{{\SetFigFont{9}{10.8}{\rmdefault}{\mddefault}{\updefault}{\color[rgb]{0,0,0}$a_1$}%
}}}}
\put(6451,-2311){\makebox(0,0)[rb]{\smash{{\SetFigFont{9}{10.8}{\rmdefault}{\mddefault}{\updefault}{\color[rgb]{0,0,0}$T_1''$}%
}}}}
\put(7201,-2311){\makebox(0,0)[lb]{\smash{{\SetFigFont{9}{10.8}{\rmdefault}{\mddefault}{\updefault}{\color[rgb]{0,0,0}$T_1'$}%
}}}}
\put(4951,-1861){\makebox(0,0)[lb]{\smash{{\SetFigFont{9}{10.8}{\rmdefault}{\mddefault}{\updefault}{\color[rgb]{0,0,0}$A_{(1,\,m)}'$}%
}}}}
\put(5101,-2161){\makebox(0,0)[lb]{\smash{{\SetFigFont{9}{10.8}{\rmdefault}{\mddefault}{\updefault}{\color[rgb]{0,0,0}$A_{(1,\,m)}$}%
}}}}
\put(3526,-4786){\makebox(0,0)[rb]{\smash{{\SetFigFont{9}{10.8}{\rmdefault}{\mddefault}{\updefault}{\color[rgb]{0,0,0}$A_R$}%
}}}}
\put(5326,-4486){\makebox(0,0)[b]{\smash{{\SetFigFont{9}{10.8}{\rmdefault}{\mddefault}{\updefault}{\color[rgb]{0,0,0}$T_R$}%
}}}}
\put(1201,-6286){\makebox(0,0)[b]{\smash{{\SetFigFont{9}{10.8}{\rmdefault}{\mddefault}{\updefault}{\color[rgb]{0,0,0}$V$}%
}}}}
\put( 76,-6211){\makebox(0,0)[lb]{\smash{{\SetFigFont{9}{10.8}{\rmdefault}{\mddefault}{\updefault}{\color[rgb]{0,0,0}(c)}%
}}}}
\put( 76,-3736){\makebox(0,0)[lb]{\smash{{\SetFigFont{9}{10.8}{\rmdefault}{\mddefault}{\updefault}{\color[rgb]{0,0,0}(b)}%
}}}}
\put( 76,-1486){\makebox(0,0)[lb]{\smash{{\SetFigFont{9}{10.8}{\rmdefault}{\mddefault}{\updefault}{\color[rgb]{0,0,0}(a)}%
}}}}
\put( 76,-8686){\makebox(0,0)[lb]{\smash{{\SetFigFont{9}{10.8}{\rmdefault}{\mddefault}{\updefault}{\color[rgb]{0,0,0}(d)}%
}}}}
\put(3301,-2161){\makebox(0,0)[rb]{\smash{{\SetFigFont{9}{10.8}{\rmdefault}{\mddefault}{\updefault}{\color[rgb]{0,0,0}$a_n$}%
}}}}
\put(526,-1711){\makebox(0,0)[rb]{\smash{{\SetFigFont{9}{10.8}{\rmdefault}{\mddefault}{\updefault}{\color[rgb]{0,0,0}$A_{(m,\,n)}$}%
}}}}
\put(526,614){\makebox(0,0)[rb]{\smash{{\SetFigFont{9}{10.8}{\rmdefault}{\mddefault}{\updefault}{\color[rgb]{0,0,0}$A_{(m,\,n)}$}%
}}}}
\put(3301,-1111){\makebox(0,0)[rb]{\smash{{\SetFigFont{9}{10.8}{\rmdefault}{\mddefault}{\updefault}{\color[rgb]{0,0,0}$a_n$}%
}}}}
\end{picture}%
\caption{(a) The annuli $A_{(1,\,m)}$, $A_{(m,\,n)}$ and $T_0$ with the solid torus $X_0$\qua (b)  The annulus $A_{(1,\,m)}'$ and the annulus $T_0'$\qua  (c) The annuli $A_R$ and $T_R$ and the solid torus $V$\qua (d) The result of unfurling along $R = A_R \cup T_R$} 
\label{fig:unfurl2}
\end{figure}
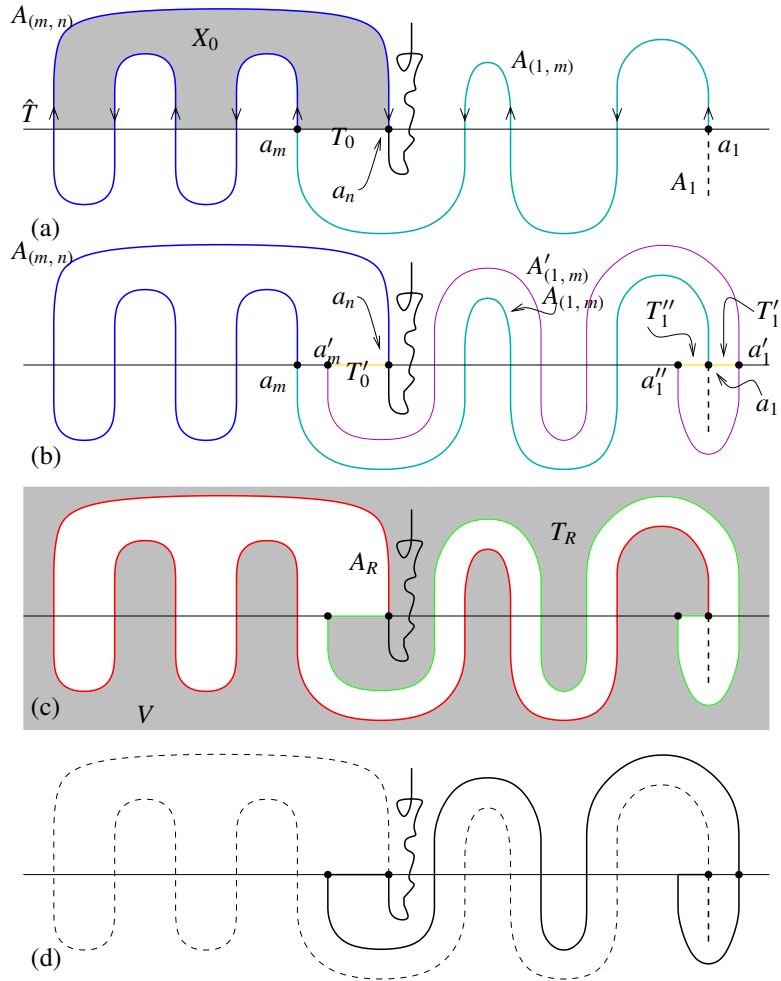

\begin{proof}
By \fullref{unfurling}, $x_1 \in \T_*$.  If $x_0 \in \T_*$, then \fullref{leafcancel} and \fullref{penultimate} imply the conclusion.  Thus we may assume $x_0 \not \in \T_*$.

Let $t_0$ be the edge of $\T$ incident to $x_0$ that separates $x_0$ from $x_*$, $T_0$ be the annulus of $\hatT \cut A$ corresponding to $t_0$, and $\bdry T_0 = a_m \cup a_n$ with $m < n$.  Form the two annuli $A_{1,\,m} = \smash{\bigcup_{i=2}^m} A_i$ and $A_{m,\,n} = \smash{\bigcup_{i=m+1}^n} A_i$.  See \fullref{fig:unfurl2}(a).  Let $\smash{A_{1,\,m}'}$ be a slight push off of $A_{1,\,m}$ with boundaries $\smash{a_1'}$ and $\smash{a_m'}$ so that $\smash{a_m'} \subseteq T_0$.  Also let $\smash{A_1'}$ be $(\bdry \smash{\bar{N}}(A_1) \cap X_1) \cut \hatT$ where one component of $\bdry \smash{A_1'}$ is $\smash{a_1'}$ and the other, say $\smash{a_1''}$, is a slight push off of $a_1$ to its other side.  Let $\smash{T_1'}$ be the annulus between $a_1$ and $\smash{a_1'}$ and $\smash{T_1''}$ be the annulus between $a_1$ and $\smash{a_1''}$ so that $\smash{T_1'} \cap \smash{T_1''} = a_1$.  Let $\smash{T_0'}$ be the annulus in $T_0$ bounded by $\smash{a_m'}$ and $a_n$.  See \fullref{fig:unfurl2}(b).  We now form the torus 
\[R = A_{1,\,m} \cup A_{m,\,n} \cup T_0' \cup \smash{A_{1,\,m}'} \cup A_1' \cup T_1''\] as shown in \fullref{fig:unfurl2}(c). 
By construction $A_1$ and $X_0$ both lie on the same side of $R$.  It follows (as in \fullref{dividingtorus}) that $R$ bounds a solid torus, say $V$, on its other side.  

Divide $R$ into the two annuli $A_R = A_{1,\,m} \cup A_{m,\,n}$ (which is a subannulus of $A$) and $T_R = \smash{T_0'} \cup \smash{A_{1,\,m}'} \cup \smash{A_1'} \cup \smash{T_1''}$.  
Notice that $K \cap R = K_{(1,\,n)} \cup K_{(t-n+1,\,t)}$ are two transverse arcs of $A_R$, $K \cap \Int T_R = \emptyset$, and $K \cut V = K_{(t, 1)}$ with $\bdry K_{(t, 1)} \subseteq a_n$.  \fullref{main-unfurl} applies.  The result of this unfurling isotopy is depicted in \fullref{fig:unfurl2}(d).

Form the annuli $A_Q = A_{(1,\,m)}' \cup A_1' \cup T_1''$ and $T_Q = A_{(1,\,m)} \cup (T_0 \cut T_0')$.    See \fullref{fig:unfurl3}(a).  The torus $Q = A_Q \cup T_Q$ bounds a solid torus that contains $A_1$.  Notice that $Q_A = T_R$.  Hence after unfurling along $R$, $Q_A$ is a subannulus of the isotoped $A_R$, $K \cap T_Q = \emptyset$, and \fullref{main-unfurl} applies again. 

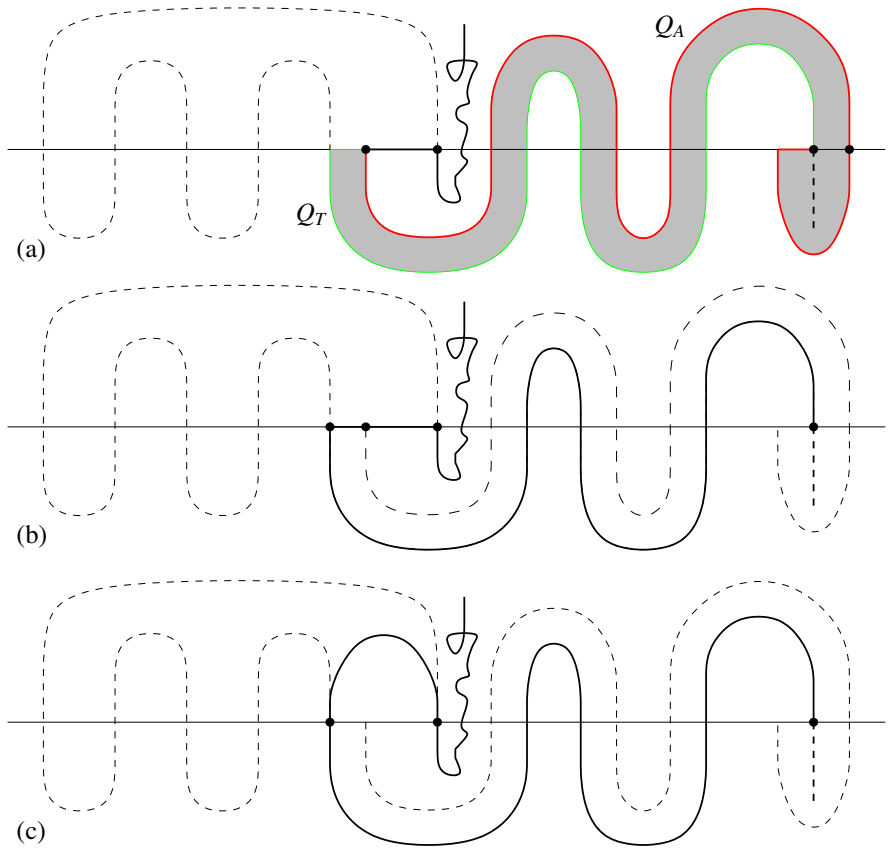
\begin{figure}[ht!]
\centering
\begin{picture}(0,0)%
\includegraphics{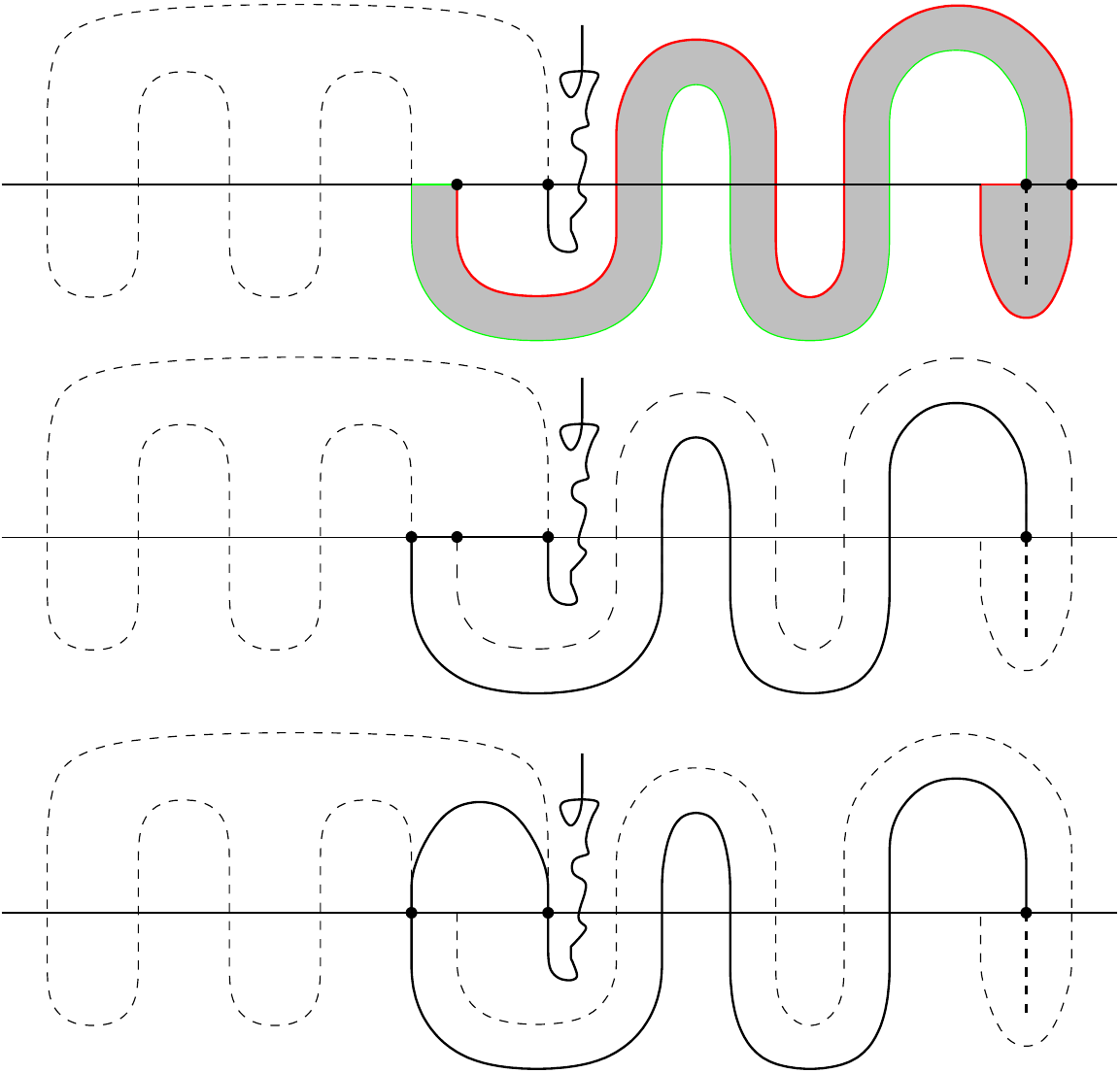}%
\end{picture}%
\setlength{\unitlength}{2960sp}%
\begingroup\makeatletter\ifx\SetFigFont\undefined%
\gdef\SetFigFont#1#2#3#4#5{%
  \reset@font\fontsize{#1}{#2pt}%
  \fontfamily{#3}\fontseries{#4}\fontshape{#5}%
  \selectfont}%
\fi\endgroup%
\begin{picture}(7374,7052)(1189,-12212)
\put(1276,-7261){\makebox(0,0)[lb]{\smash{{\SetFigFont{9}{10.8}{\rmdefault}{\mddefault}{\updefault}{\color[rgb]{0,0,0}(a)}%
}}}}
\put(1276,-12136){\makebox(0,0)[lb]{\smash{{\SetFigFont{9}{10.8}{\rmdefault}{\mddefault}{\updefault}{\color[rgb]{0,0,0}(c)}%
}}}}
\put(1276,-9661){\makebox(0,0)[lb]{\smash{{\SetFigFont{9}{10.8}{\rmdefault}{\mddefault}{\updefault}{\color[rgb]{0,0,0}(b)}%
}}}}
\put(6901,-5386){\makebox(0,0)[rb]{\smash{{\SetFigFont{9}{10.8}{\rmdefault}{\mddefault}{\updefault}{\color[rgb]{0,0,0}$Q_A$}%
}}}}
\put(3901,-6961){\makebox(0,0)[rb]{\smash{{\SetFigFont{9}{10.8}{\rmdefault}{\mddefault}{\updefault}{\color[rgb]{0,0,0}$Q_T$}%
}}}}
\end{picture}%
\caption{(a) The annuli $Q_A$ and $Q_T$ together bound a solid torus (b) The result of unfurling along $Q = Q_A \cup Q_T$\qua  (c) The result of a further final isotopy}
\label{fig:unfurl3}
\end{figure}

After unfurling along $R$ and then along $Q$, the annulus $A_{(m,\,n)}$ has been repositioned as $T_0$, and the arcs $K_{(m,\,n)}$ and $K_{(t-n+1,\, t-m+1)}$ now lie as transverse arcs of $T_0$.  See \fullref{fig:unfurl3}(b).   The rest of $K$ is as it was (up to perhaps height-preserving isotopies near $a_1$, $a_m$, and $a_n$).  By a slight further isotopy of the arcs $K_{(m,\,n)}$ and $K_{(t-n+1,\, t-m+1)}$ into $X_0$, we obtain another Morse presentation of $K$ of width not greater than previously.  See \fullref{fig:unfurl3}(c).  Since $K$ was in thin position before these isotopies, it must be in another thin position now.  Hence $n = m+1$ and $x_0$ is now a leaf of $\T_*$.  Furthermore this implies that $x_0$ must have been a leaf of $\T$ connected to $x_1 \in \T_*$ before the isotopies. 

If $\T$ had a leaf not in $\T_*$ of opposite sign from $x_0$, then after these isotopies there would be a high disk and low disk pair contradicting \fullref{highdisklowdisk}.  If $\T$ had a nonleaf vertex that was not in $\T_*$, then $\T$ has a penultimate leaf not in $\T_*$.  Since $x_0$ was a leaf of $\T$ connected to $\T_*$ by an edge, it could not be a leaf adjoined to this penultimate leaf.  Thus there would have been a long disk contradicting \fullref{longdisk}.  The conclusion of this lemma follows.
\end{proof}

\begin{lemma}\label{lem:movingx0}
If $x_0 \not \in \T_*$, then
there is an isotopy of $K$ into another thin position which takes each annulus $A_i$ to itself for $i < k$ and $A_k$ to another annulus such that the resulting trees have $x_0 \in \T_*$.
\end{lemma}

\begin{proof}
The proof of this lemma is basically contained within the proof of \fullref{unfurling2}.  

Assume $x_0 \not \in \T_*$.  By \fullref{unfurling2}, $\bdry X_0$ consists of the annulus $A_k$ and the annulus, say, $T_k \in \hatT \cut A$ bounded by the curves $a_k$ and $a_{k-1}$.  Let $X_2 \in \mathcal{X}$ be the other solid torus that contains $A_k$ in its boundary.  Since $k \geq 3$, $X_2$ cannot correspond to a leaf of $\T$.  \fullref{unfurling2} implies that the interior of $X_2$ must intersect $K$.

The isotopy employed in the proof of \fullref{unfurling2} has the effect of rearranging $S$ in $N(\bigcup_{i=1}^k A_i \cup T_k)$ so that $A_k$ is moved from one side of $X_0$ to a slight push off $T_k$ on the other side and each $A_i$ for $1\leq i <k$ is taken to itself.  For $\T$, this is tantamount to moving the label $x_0$ to the vertex in $\T_*$ that previously corresponded to $X_2$.
\end{proof}

\subsection{The least extreme critical points}

Recall that in the beginning of this \fullref{sec:annuliandtrees} we let $\{A_2, \dots, A_k\}$ be the annuli and $A_1$ be the M\"obius band or complex associated to the face of an extended Scharlemann cycle of order $2$ or $3$.  Let $A_{k+1}$ be the complex formed from identifying the common corners of the bigon and trigon at the ``ends'' of the face bounded by a forked extended Scharlemann cycle that contains the extended Scharlemann cycle of order $2$ provided such a forked extended Scharlemann cycle exists.  Otherwise, set $A_{k+1} = \emptyset$.
\newpage

\begin{Lemma}\label{highestmin}
After perhaps a width-preserving isotopy,
the highest minimum (resp.\ lowest maximum) below (resp.\ above) $\hatT$ lies on an arc of $K$ that together with an arc on $\hatT$ bounds a low disk (resp.\ high disk) with interior disjoint from $ \bigcup_{i=1}^{k+1} A_i$ and from $K$.  Furthermore, if the boundary of this low disk (resp.\ high disk) intersects $\bigcup_{i=1}^{k+1} \bdry A_i$ then either it intersects only $\bdry A_{k+1} - a_k$ if $A_{k+1} \neq \emptyset$ or it intersects only $a_k$ if $A_{k+1} = \emptyset$.
\end{Lemma}

\begin{proof}
We work with the highest minimum below $\hatT$.  The other case follows in the same manner.

Let $m$ be the first critical value of $h(K)$ below $0$, and let $p_m$ be the corresponding critical point.  Since $\hatT = h^{-1}(0)$ is a thick level, $p_m$ is a minimum.  For a suitably small $\epsilon >0$, $K \cap h^{-1}[m-\epsilon,0]$ is a collection of arcs of $K$ transverse to the induced product structure on $h^{-1}[m-\epsilon,0]$ together with the arc $\kappa$ containing $p_m$ that has both endpoints on $\hatT$.  Thus there exists a low disk for $\kappa$.  

Let $D$ be a low disk for $\kappa$ so that its interior is disjoint from $K$ and $\bdry D$ consists of $\kappa$ and an arc on $\hatT$.  If $\kappa$ is contained in $A_i$ for some $i$, we may assume that $N(\bdry \kappa) \cap (D-\kappa) \cap A_i = \emptyset$ at the expense of perhaps increasing $|(\bdry D - \kappa) \cap A_i|$.  
Let us further assume $D$ is transverse to $\smash{\bigcup_{i=1}^{k+1}} A_i$ and that $D$ has been chosen among all such disks (with $N(\bdry \kappa) \cap (D-\kappa) \cap A_i = \emptyset$) so that $|D \cap \smash{\bigcup_{i=1}^{k+1}} A_i|$ is minimized.

Assume $(D-\kappa) \cap \Int A_i \neq \emptyset$ for some $i \in \{1, \dots, k+1\}$.  Since $D$ is disjoint from $K$ except along $\kappa$, $(D-\kappa) \cap A_i = D \cap (A_i \cut K)$.  If the intersection contained any simple closed curves (that do not intersect $\hatT$ in the case of $A_1$), then standard innermost disk arguments would imply a contradiction either to the minimality of $|\Int D \cap \bigcup_{i=1}^{k+1} A_i|$ or to the incompressibility of $S$.  Therefore $D \cap A_i$ is a collection of arcs for each $i$.  Since each component of $A_i \cut K$ is a bigon or trigon of $G_S$, each arc of $D \cap A_i$ either bounds a subdisk of $A_i \cut K$ with an arc of $\hatT$ or bounds a rectangle of $A_i$ with two arcs of $\hatT$ and an arc of $K$.  Standard outermost arc arguments show that we may assume that only arcs of the second type occur.

Let $\alpha$ be an arc of $D \cap (A_i \cut K)$, outermost on $D$.  Thus $\alpha$ cuts off a subdisk $D_\alpha \subseteq D$ disjoint from $\kappa$.  Since $\alpha$ must be an arc of the second type, let $R \subseteq A_i$ be a rectangle with boundary composed of $\alpha$, two arcs of $\bdry A_i \subseteq \hatT$, and an arc, say $\kappa'$, of $K \cap A_i$.  Since the interior of $D_\alpha$ is disjoint from $A_i$, $D_\alpha$ intersects $R$ only along $\alpha$.  Thus by a slight isotopy of the disk $R \cup D_\alpha$, we may form a low disk $D'$ with boundary composed of $\kappa'$ and an arc on $\hatT$ such that $|\Int D' \cap \bigcup_{i=1}^{k+1} A_i| = 0$.  

If $\kappa' = \kappa$ then the existence of $D'$ contradicts the minimality of $|\Int D \cap \smash{\bigcup_{i=1}^{k+1}} A_i|$.  Thus either $|\Int D \cap \smash{\bigcup_{i=1}^{k+1}} A_i| = 0$ satisfying the lemma or $\kappa' \neq \kappa$.

If $\kappa' \neq \kappa$ then $D'$ guides a width-preserving isotopy of $K$ so that the minimum of $\kappa'$ is higher than the minimum of $\kappa$. 

Given the disk $D$ with $|\Int D \cap \bigcup_{i=1}^{k+1} A_i| = 0$, the final conclusion of the lemma is immediate.
\end{proof}

If $A_{k+1} \neq \emptyset$, then each arc of $A_{k+1} \cap K$ bounds a disk with an arc of $\hatT$ contained in $N(\bdry A_{k+1})$ whose interior is disjoint from $K$ and $a_k$.  These disks may be constructed similarly to the construction of the disk $D_g$ in \fullref{thinningtwoforks}. 

\begin{Lemma}\label{uniqueleaf}
If a leaf of $\T$ is not in $\T_*$, then it is the only such leaf.  Furthermore $A_k$ is contained in the boundary of the solid torus corresponding to this leaf.
\end{Lemma}

\begin{proof}
Assume $k=2$.  If there are two leaves of $\T$ not in $\T_*$, then $x_* = x_1 = \T_*$. Hence $t=4$ contradicting that $t \geq 6$.  Thus we assume $k \geq 3$.

Let $x_m$ be a leaf of $\T$ not in $\T_*$.  Without loss of generality we may assume its sign is positive.  By \fullref{unfurling}, $x_1 \in \T_*$.  \fullref{lem:movingx0} implies that we may assume $x_0 \in \T_*$.  Thus $x_m$ is neither $x_0$ nor $x_1$.

  Let $A_m$ be the annulus of $A \cut \hatT$ which together with an annulus $T_m \in \hatT \cut A$ cuts off the solid torus $X_m \in \mathcal{X}$ corresponding to the leaf $x_m$.   The arcs of $K \cap A_m = K_{(m-1,\,m)} \cup K_{(t-m+1,\,t-m+2)}$ each with an arc $\tau_{(m-1,\,m)}$ or $\tau_{(t-m+1,\,t-m+2)}$ of $T_m$ bound a meridional disk of $X_m$.  These two meridional disks are high disks and may be assumed to be disjoint.

By \fullref{highestmin} the arc $\kappa$ of $K \cut \hatT$ containing the highest minimum below $\hatT$ together with an arc $\tau_\kappa$ of $\hatT$ bounds a low disk $D$ with interior disjoint from $\bigcup_{i=1}^k A_i$.  If the interior of $\tau_\kappa$ is disjoint from the interior of either $\tau_{(m-1,\,m)}$ or $\tau_{(t-m+1,\,t-m+2)}$, then by \fullref{highdisklowdisk} $K$ is not in thin position.  Hence $\tau_\kappa$ must intersect the interior of both.  Thus $\tau_\kappa \cap \Int T_m \neq \emptyset$.  For this to occur, either $\bdry \kappa \cap \Int T_m \neq \emptyset$, $\bdry \kappa \subseteq \bdry T_m$, or $\Int \tau_\kappa \cap \bdry T_m \neq \emptyset$.

Since the interior of $X_m$ is disjoint from $K$, $\bdry \kappa \cap \Int T_m = \emptyset$.  If $\bdry \kappa \subseteq \bdry T_m$ then $\kappa$ must join the two arcs of $K \cap A_m$.  Hence $a_k \subseteq \bdry T_m$.  If $\Int \tau_\kappa \cap \bdry T_m \neq \emptyset$ then, as in \fullref{highestmin}, $\Int \tau_\kappa$ may only intersect $a_k$.  Hence $a_k \subseteq \bdry T_m$.  Thus $A_k \subseteq \bdry X_m$. 

Assume $x_n$ is another leaf of $\T$ not in $\T_*$.  \fullref{unfurling2} implies that $x_m$ and $x_n$ have the same sign.  We may also assume $x_n$ is neither $x_0$ nor $x_1$.  The same argument above thus applies for the leaf $x_n$.  Hence if $X_n \in \mathcal{X}$ is the solid torus corresponding to $x_n$,  then $A_k \subseteq \bdry X_n$.  Therefore $A_k$ separates $X^+$ into the solid tori $X_m$ and $X_n$ of $\mathcal{X}$.  In other words, $X^+ = X_m \cup A_k \cup X_n$.  

Since $k \geq 3$, $A_{k-2}$ exists.  Furthermore $A_{k-2}$ must be contained in $X^+$ since $A_k$ is.  Yet since $A_k$ is the only subannulus of $A = \bigcup_{i=2}^k A_i$ contained in $X^+$, it must be the case that $k=3$ and $A_1 \subseteq X^+$.  But because neither $x_m = x_1$ nor $x_n = x_1$, $A_1 \not \subseteq X^+$.  This is a contradiction.
\end{proof}

\begin{lemma}\label{lem:lonevertex}
There may be at most one vertex of $\T$ not in $\T_*$.
\end{lemma}

\begin{proof}
This follows from \fullref{unfurling2} and \fullref{uniqueleaf}.
\end{proof}


\subsection{Bounds on extended Scharlemann cycles}
We may now get good estimates on how many labels are accounted for by the face of an extended $S2$ or $S3$ cycle.  Recall that we are assuming $t \geq 6$.

\begin{Lemma}\label{bounded}
For any $x \in \mathbf{t}$ the face bounded by an extended $S2$ or $S3$ cycle on $\smash{G_S^x}$ may account for at most $t/2 +1$ of the labels if $4 {\nmid} t$ and at most $t/2$ of the labels if $4|t$.
\end{Lemma}

\begin{proof}
Let $\sigma'$ be an extended $S2$ or $S3$ cycle of $\smash{G_S^x}$ for some $x \in \mathbf{t}$ bounding the disk face $f'$.  
Let $\sigma$ be the Scharlemann cycle contained in $f'$.  Assume we have relabeled so that $\sigma$ has label pair $\{t,1\}$ and $\sigma'$ has label pair $\{t-k+1, k\}$.  Thus $f'$ accounts for $2k$ labels.  If $k=1$ then $\sigma = \sigma'$ and the lemma holds.  Therefore we may assume $k \geq 2$.

From $f'$ form the complex $A_1$, the annuli $A_i$ for $i=2, \dots, k$, the annulus $A=\smash{\bigcup_{i=2}^{k}} A_i$, and the corresponding trees.  \fullref{lem:lonevertex} implies that there can be at most one vertex of $\T$ not contained in $\T_*$. 

Since $\T$ has $k$ edges, $\T_*$ has at least $k-1$ edges.  By \fullref{evenintersections}, the interior of the arc $K_{(k,\,t-k+1)}$ intersects $\hatT$ at least $2(k-1)$ times.  Hence at least $2(k-1)$ labels are not accounted for by $f'$.  Therefore $t \geq 2k + 2(k-1) = 4k - 2$, ie\ $t/2 +1 \geq 2k$.  Since $k$ must be an integer, if $t$ is divisible by $4$ then $t/2 \geq 2k$.
\end{proof}

\begin{lemma}\label{boundedforked}
For any $x \in \mathbf{t}$ the trigon bounded by a forked extended $S2$ cycle on $\smash{G_S^x}$ may account for at most $t/2$ of the labels if $4 {\nmid} t$ and at most $t/2-1$ of the labels if $4|t$.
\end{lemma}

\begin{proof}
Let $\sigma''$ be a forked extended $S2$ cycle bounding the disk face $f''$.  Let $\sigma'$ be the outermost extended $S2$ cycle contained in $f''$, and let $f'$ be the disk face it bounds.  As in the above proof of \fullref{bounded}, let $\sigma$ be the $S2$ cycle contained in $f'$.  Assume we have relabeled so that $\sigma$ has label pair $\{t,1\}$ and $\sigma'$ has label pair $\{t-k+1, k\}$.  Since $f'$ accounts for $2k$ labels, $f''$ accounts for $2k+1$ labels.  If $k=1$ then $\sigma = \sigma'$, $f''$ accounts for $3$ labels, and the lemma holds.  Therefore we may assume $k\geq 2$.

Again form the complex $A_1$, the annuli $A_i$ for $i=2, \dots,k$, the annulus $A$, and the corresponding trees; and again note that $\T$ can have at most one vertex not in $\T_*$.  Also form the complex $A_{k+1}$ from the bigon and trigon of $f'' \cut f'$.  Note that $A_{k+1} \cap A = a_k$.

Let $T_k$ and $T_k'$ be the two annuli of $\hatT \cut A$ with $a_k$ as a boundary component.  Since $A_{k+1} \cap K = K_{(k,\,k+1)} \cup K_{(t-k,\,t-k+1)}$, either $T_k$ or $T_k'$ contains both the end points $K_{k+1}$ and $K_{t-k}$ (because $a_k$ contains the end points $K_k$ and $K_{t-k+1}$).  Say $T_k$ contains these end points.  Since by \fullref{lem:lonevertex} at most one vertex of $\T$ is not in $\T_*$, either (a) $T_k'$ corresponds to the edge $t_k'$ of $\T$ incident to the leaf of $\T$ not in $\T_*$ or (b) the interior of $K_{(k+1,\,t-k)}$ intersects $T_k'$ nontrivially.

{\bf Case (a)}\qua  Assume $T_k'$ corresponds to an edge $t_k'$ of $\T$ incident to the leaf of $\T$ not in $\T_*$.  It must be the annulus $A_k$ that together with $T_k'$ cuts off the solid torus corresponding to this leaf.  By \fullref{lem:movingx0} we may assume this leaf is not $x_0$.  Furthermore, by paying attention to the isotopy of \fullref{lem:movingx0} (described in the proof of \fullref{unfurling2}) one may note that $A_{k+1}$ is preserved.
Since we may assume this leaf is not $x_0$ so that $A_k$ and $T_k'$ are parallel, the arcs of $K \cap A_k$ together with arcs contained in $T_k'$ then bound high (or low) disks.  Furthermore, as in the construction of the high disks of \fullref{thinningtwoforks} shown in \fullref{fig:twoforkshighlowdisks}, the arcs of $K \cap A_{k+1}$ together with arcs of $T_k$ bound low (or high) disks.  By \fullref{highdisklowdisk}, this contradicts the thinness of $K$.

{\bf Case (b)}\qua  Assume the interior of $K_{(k+1,\,t-k)}$ intersects $T_k'$ nontrivially.

  If $\T \neq \T_*$, then there is a leaf of $\T$ not in $\T_*$.  By \fullref{uniqueleaf} the solid torus corresponding to this leaf contains $A_k$ and hence $a_k$ in its boundary.  Thus it must contain $T_k$ or $T_k'$.  This contradicts our assumption that neither $T_k$ nor $T_k'$ has interior disjoint from $K$.  Thus $\T = \T_*$.

Since $\T$ has $k$ edges, $\T_*$ has $k$ edges.  By \fullref{evenintersections}, each edge of $\T_*$ corresponds to at least two intersections of the interior of the arc $K_{(k,\,t-k+1)}$ with each of the $k$ annuli $\hatT \cut A$.  Since $\T$ is a tree and the points $K_{k+1}$ and $K_{t-k}$ are contained in $T_k$, the interior of the arc $K_{(k+1,\,t-k)}$ must intersect $T_k$ twice more so that $K$ may intersect the interior of $T_k'$.  Hence $K$ intersects the interior of $T_k$ at least $4$ times.

 It follows that the interior of the arc $K_{(k,\,t-k+1)}$ intersects $\hatT$ at least $2k+2$ times.  Therefore $t \geq 2k + 2k+2$, ie\ $t/2 -1 \geq 2k$.  Since $k$ must be an integer, if $t$ is divisible by $4$ then $t/2-2 \geq 2k$.  Since $f''$ accounts for one more label than $f'$, it may account for at most $t/2$ labels if $4{\nmid}t$ and at most $t/2-1$ labels if $4|t$.
\end{proof}

\section{Two Scharlemann cycles}\label{sec:twoschcycles}

\subsection{Accounting for labels}

Let $\mathcal{F}_{\{x,\, x{+}1\}}$ be the collection of bigons and trigons of $G_S^y$ for every $y \in \mathbf{t}$ that contain a Scharlemann cycle (of order $2$ or $3$) of $G_S$ with label pair $\{x, x{+}1\}$.  Let $\Lambda_{\{x,\, x{+}1\}}$ be the subset of labels in $\mathbf{t}$ accounted by the faces of $\mathcal{F}_{\{x,\, x{+}1\}}$.  

\begin{lemma} \label{labelaccount}
There exists an $F \in \mathcal{F}_{\{x,\, x{+}1\}}$ such that $F$ accounts for every label in $\Lambda_{\{x,\, x{+}1\}}$.
\vspace{-5pt}
\end{lemma}

\begin{proof}
For $F_1, F_2 \in \mathcal{F}_{\{x,\, x{+}1\}}$, say $F_1 \leq F_2$ if the set of labels for which $F_1$ accounts is a subset of the set of labels for which $F_2$ accounts.  It is clear that $\leq$ is a partial ordering.  We claim that $\leq$ indeed totally orders the finite set $\mathcal{F}_{\{x,\, x{+}1\}}$.  A maximal element of $\mathcal{F}_{\{x,\, x{+}1\}}$ then accounts for every label in $\Lambda_{\{x,\, x{+}1\}}$.

Assume neither $F_1$ nor $F_2 \in \mathcal{F}_{\{x,\, x{+}1\}}$ are bounded by forked extended $S2$ cycles.   Now $F_i$ accounts for the labels $\{x-k_i, x-k_i+1, \dots, x, x+1, \dots, x+k_i, x+k_i+1\}$ for some integer $k_i \geq 0$ for each $i=1,2$, since $F_i$ is an extended $S2$ or $S3$ cycle.  Thus either $F_1 \leq F_2$ or $F_2 \leq F_1$ since either $k_1 \leq k_2$ or $k_2 \leq k_1$ respectively.

By \fullref{prop:nosaladforks}, if $F_1, F_2 \in \mathcal{F}_{\{x,\, x{+}1\}}$ are two trigons bounded by forked extended $S2$ cycles containing outermost extended $S2$ cycles with the same label pair, then they must account for the same set of labels.  Thus $F_1 = F_2$ with respect to the ordering $\leq$.

Since a trigon bounded by a forked extended $S2$ cycle accounts for an odd number of labels whereas an extended $S2$ or $S3$ cycle accounts for an even number of labels, it follows that for any $F_1, F_2 \in \mathcal{F}_{\{x,\, x{+}1\}}$, either $F_1 \leq F_2$ or $F_2 \leq F_1$.  Thus $\leq$ is a total ordering on $\mathcal{F}_{\{x,\, x{+}1\}}$.
\end{proof}

\begin{prop}\label{disjointlabelpairs}
There exist two Scharlemann cycles on $G_S$ each of order $2$ or $3$ with disjoint label pairs.  Furthermore, there cannot be a third Scharlemann cycle of order $2$ or $3$ with label pair distinct from the other two.
\vspace{-5pt}
\end{prop}

\begin{proof}
Assume any Scharlemann cycle of $G_S$ of order $2$ or $3$ has label pair $\{x, x{+}1\}$.  By \fullref{labelaccount} there exists an $F \in \mathcal{F}_{\{x,\, x{+}1\}}$ that accounts for the labels $\Lambda_{\{x,\, x{+}1\}}$.   Then by \fullref{bounded} and \fullref{boundedforked}, $|\Lambda_{\{x,\, x{+}1\}}| \leq t/2 +1$.  Since $t \geq 6$, we have $\Lambda_{\{x,\, x{+}1\}} \neq \mathbf{t}$.  Therefore there must be a second $S2$ or $S3$ cycle with label pair distinct from $\{x, x{+}1\}$.

Assume there are two Scharlemann cycles $\sigma$ and $\sigma'$ of $G_S$ of order $2$ or $3$ with label pairs $\{x-1, x\}$ and $\{x, x{+}1\}$ respectively.  Again, by \fullref{labelaccount} there exists $F \in \mathcal{F}_{\{x-1,\, x\}}$ that accounts for the labels $\Lambda_{\{x-1,\, x\}}$ and $F' \in \mathcal{F}_{\{x,\, x{+}1\}}$ that accounts for the labels $\Lambda_{\{x,\, x{+}1\}}$.  Also, by \fullref{bounded} and \fullref{boundedforked}, $|\Lambda_{\lambda}| \leq t/2 +1$ for each $\lambda = \{x{-}1, x\}, \{x, x{+}1\}$.  Because for each $\{x{-}1, x\}$ and $\{x, x{+}1\}$ the sequence of labels $\lambda$ contains the middle of the set $\Lambda_{\lambda}$ (when ordered sequentially), $|\Lambda_{\{x-1,\, x\}} \cup \Lambda_{\{x,\, x{+}1\}}| \leq t/2 + 1 + 1$.  Since $t\geq 6$, $\Lambda_{\{x-1,\, x\}} \cup \Lambda_{\{x,\, x{+}1\}} \neq \mathbf{t}$.  Therefore there must be a third $S2$ or $S3$ cycle $\sigma''$ in $G_S$ with label pair distinct from the label pair of each $\sigma$ and $\sigma'$.

Assume there are three Scharlemann cycles of order $2$ or $3$ with mutually distinct label pairs.  Then two of them have disjoint label pairs and must bound faces on the same side of $\hatT$.  This contradicts \fullref{oppositesides}.
\end{proof}

Pulling \fullref{labelaccount} and the proof of  \fullref{disjointlabelpairs} together, there exists two bigons or trigons, say $F^-$ and $F^+$ of the subgraphs $\smash{G_S^x}$ and $G_S^y$ respectively for some $x, y \in \mathbf{t}$, that account for all the labels $\mathbf{t}$.  Interior to $F^-$ and $F^+$ are order $2$ or $3$ Scharlemann cycles $\sigma^-$ and $\sigma^+$ respectively with disjoint label pairs.  These Scharlemann cycles bound faces $f^-$ and $f^+$ respectively of $G_S$.  Let us assume we have labeled the intersections of $K \cap \hatT$ so that the Scharlemann cycle $\sigma^-$ within $F^-$ has label pair $\{t, 1\}$.

\begin{Lemma} \label{tis61014}
$4{\nmid} t$
\end{Lemma}
\begin{proof} 
Assume $4|t$.  Then by \fullref{bounded} and \fullref{boundedforked} the maximum number of labels for which each $F^-$ and $F^+$ may account is $t/2$.  In order for $F^-$ and $F^+$ to account for all $t$ labels together, their label sets must each realize this maximum and must be disjoint from one another.  Since $\sigma^-$ has label pair $\{t, 1\}$, $F^-$ accounts for the labels $\{t-t/4+1, \dots, t, 1, \dots, t/4\}$.  Therefore $F^+$ accounts for the labels $\{t/4+1, \dots, t/2, t/2+1, \dots, t-t/4\}$.  This however implies that $\sigma^+$ has label pair $\{t/2, t/2+1\}$.  Thus $f_{\sigma^-}$ and $f_{\sigma^+}$ lie on the same side of $\hatT$ contradicting \fullref{oppositesides}.
\end{proof}

\begin{Lemma} \label{oneconfig}
$F^+$ and $F^-$ are each bounded by extended Scharlemann cycles of order $2$ or $3$.   
\end{Lemma}
\begin{proof}
Given that $t \geq 6$ by assumption and $4{\nmid} t$ by \fullref{tis61014}, then it is direct to check that the bounds of \fullref{bounded} and \fullref{boundedforked} give us only two cases that we must consider.  For both cases it turns out that $\sigma^+$ has label pair $\{t/2, t/2+1\}$.  Without loss of generality, these two cases are:
\begin{enumerate}
\item \label{twoschcycles}
\vspace{10pt}
\begin{itemize}
\item $F^-$ is bounded by an extended Scharlemann cycle of order $2$ or $3$ that accounts for the labels $\{t-(t+2)/4+1\, \dots, t, 1, \dots, (t+2)/4\}$, and
\item $F^+$ is bounded by an extended Scharlemann cycle of order $2$ or $3$ that accounts for the labels $\{(t+2)/4+1, \dots, t/2, t/2+1, \dots, t-(t+2)/4\}$.
\end{itemize} 
\item \label{twoforks}
\vspace{10pt}
\begin{itemize}
\item $F^-$ is bounded by a forked extended $S2$ cycle that accounts for the labels $\{t-(t+2)/4+2, \dots, t, 1, \dots, (t+2)/4-1, (t+2)/4\}$, and
\item $F^+$ is bounded by a forked extended $S2$ cycle that accounts for the labels $\{(t+2)/4+1, \dots, t/2, t/2+1, \dots, t-(t+2)/4, t-(t+2)/4+1\}$.
\end{itemize}
\end{enumerate}
  
Case~\eqref{twoschcycles} satisfies the lemma.  Note also that since $F^+$ only accounts for $t/2 - 1$ labels it could be interior to the face of another extended Scharlemann cycle or a forked extended $S2$ cycle.

In Case~\eqref{twoforks}, notice that $F^-$ is a trigon of $G_{\smash{S}}^{(t+2)/4}$ bounded by a forked extended $S2$ cycle with the arc $K_{\smash{((t+2)/4-1, (t+2)/4)}}$ on its boundary.  Similarly, $F^+$ is a trigon of $G_S^{\smash{t-(t+2)/4+1}}$ which is bounded by a forked extended $S2$ cycle with the arc $K_{\smash{(t-(t+2)/4, t-(t+2)/4+1)}}$ on its boundary.  Furthermore, since $(t+2)/4-1$ and $t-(t+2)/4$ have opposite parity, these two arcs of $K$ lie on opposite sides of $\hatT$.  A high disk and a low disk may then be constructed from $F^+$ and $F^-$ as in \fullref{thinningtwoforks}.  By construction these disks may be isotoped to be disjoint.  This contradicts \fullref{highdisklowdisk}.
\end{proof}

We will henceforth assume $F^+$ and $F^-$ are as in Case~\eqref{twoschcycles} in the proof of \fullref{oneconfig}.

\subsection{Two trees}  

As in \fullref{sec:annuliandtrees}, from each $F^-$ and $F^+$ we may construct corresponding complexes $$A^- = \tsty\bigcup_{i=1}^{(t+2)/4-1} A_i^- \quad \text{and} \quad A^+ = \tsty\bigcup_{i=1}^{(t+2)/4-1} A_i^+$$ respectively where $A_i^\pm$ are annuli for $i = 2, \dots, (t+2)/4-1$ and $A_1^\pm$ are M\"obius bands or complexes.
Notice that $A^-$ and $A^+$ are disjoint.  Define the curves $\smash{a_i^\pm} = \bdry \smash{A_i^+} -\bdry \smash{A_{i-1}^\pm}$ for $2 \leq i \leq (t+2)/4-1$ and $\smash{a_1^\pm = \bdry A_1^\pm}$.  We may further define the annulus $\smash{A_{(t+2)/4}^-}$ and the curve $\smash{a_{(t+2)/4}^-} = \bdry \smash{A_{(t+2)/4}^- - a_{(t+2)/4-1}^-}$.  Note that $\smash{A_{(t+2)/4}^-}$ is disjoint from $A^+$ too.

We may also define trees $\T^+$, $\smash{\T_*^+}$, $\T^-$, and $\smash{\T_*^-}$ associated to $A^+$ and $A^-$, as well as the corresponding labelings of vertices $\smash{x_1^+}$, $\smash{x_*^+}$, $\smash{x_0^+}$, $\smash{x_1^-}$, $\smash{x_*^-}$, and $\smash{x_0^-}$.  Since the solid torus corresponding to $\smash{x_0^+}$ must contain $\smash{A_1^-}$ and the solid torus corresponding to $\smash{x_0^-}$ must contain $\smash{A_1^+}$, $\smash{x_0^- \in \T_*^-}$ and $\smash{x_0^+ \in \T_*^+}$.

\begin{Lemma}\label{highestmin2}
After perhaps a width-preserving isotopy, the highest minimum (resp.\ lowest maximum) below (resp.\ above) $\hatT$ lies on an arc of $K$ that together with an arc on $\hatT$ bounds a low (resp.\ high) disk with interior disjoint from $A^- \cup A_{(t+2)/4}^- \cup A^+$.  
\end{Lemma}

\begin{proof}
Because $A^- \cup A_{(t+2)/4}^-$ and $A^+$ are disjoint, this lemma quickly follows from the proof of \fullref{highestmin}.
\end{proof}

\begin{Lemma}\label{oneannulusintersectingK}
The interiors of the arcs of $K \cut (A^- \cup A^+)$ intersect $\hatT$ in a single annulus of $\hatT \cut (A^- \cup A^+)$.
\end{Lemma}

\begin{proof}
The interior of the arcs of $K \cut (A^- \cup A^+)$ intersect $\hatT$ in only two points, $K_{(t+2)/4}$ and $K_{t-(t+2)/4 +1}$.  Since $A_{(t+2)/4}^-$ is disjoint from $A^+$, $a_{(t+2)/4}^-$ is contained in a component of $\hatT \cut (A^- \cup A^+)$.  The points $K_{(t+2)/4}$ and $K_{t-(t+2)/4+1}$ are contained in $a_{(t+2)/4}^-$.
\end{proof}

Define $T_{(t+2)/4}$ to be this annulus of $\hatT \cut (A^- \cup A^+)$ containing $a_{(t+2)/4}^-$.

\begin{Lemma}\label{differentboundarycurves}
Assume $t \geq 10$ and $4 {\nmid} t$.  
For any annulus $R \in \hatT \cut (A^- \cup A^+)$, we have $\bdry R = a_i^- \cup a_j^+$ for some $i, j \in \{1, \dots, (t+2)/4-1\}$.
\end{Lemma}
\begin{proof}
Assume otherwise.  Since $|A^- \cap \hatT| = |A^+ \cap \hatT| = (t+2)/4-2 \geq 2$, there must be two annuli, say $R^+$ and $R^-$, of $\hatT \cut (A^- \cup A^+)$ such that $\bdry R^+ = a_i^+ \cup a_j^+$ and $\bdry R^- = a_k^- \cup a_l^-$ for integers $1 \leq i < j \leq (t+2)/4-1$ and $1 \leq k < l \leq (t+2)/4-1$.  By \fullref{oneannulusintersectingK}, at least one of these two annuli has its interior disjoint from $K$ and is thus not $T_{(t+2)/4}$. 

Consider the annulus $\bigcup_{i=2}^k A_i^-$ where $k=(t+2)/4$ and its associated trees $\T$ and $\T_*$.  As noted in \fullref{evenintersections} since $k > t/4$, there exists a vertex of $\T$ not in $\T_*$.  By \fullref{uniqueleaf} there is a single leaf of $\T$ not contained in $\T_*$.  Furthermore, $A_k^- = A_{(t+2)/4}^-$ is contained in the boundary of the solid torus $X_{(t+2)/4}$ corresponding to this leaf.  Thus $\bdry X_{(t+2)/4} \cut A_{(t+2)/4}^-$ is one component of $T_{(t+2)/4} \cut a_{(t+2)/4}$.  Therefore one component of $\bdry T_{(t+2)/4}$ is $a_{(t+2)/4-1}^-$.

Since the other component of $T_{(t+2)/4} \cut a_{(t+2)/4}$ has interior disjoint from $K$, the other component of $\bdry T_{(t+2)/4}$ cannot be $a_n^+$ for any $1 \leq n \leq (t+2)/4-2$ without contradicting \fullref{uniqueleaf}.  Thus neither $R^+$ nor $R^-$ is $T_{(t+2)/4}$.

Now consider the annuli $A^+$ and $A^-$ and their corresponding trees.

If $i$ and $j$ have the same parity, then the torus $R^+ \cup \smash{\bigcup_{s=i+1}^j A_s^+}$ separates $\smash{A_1^+}$ from $K \cut A^+$.  Thus the vertex $x_1^+ \in \T^+$ is not contained in $\T_*^+$.  \fullref{unfurling} then applies contradicting the thinness of $K$.  Similarly $k$ and $l$ cannot have the same parity.

Since $i$ and $j$ have opposite parity, then $R^+$ corresponds to an edge of $\T^+$ that separates vertices of $\T^+$ from $\T_*^+$.  By \fullref{uniqueleaf} $R^+$ separates a single leaf of $\T^+$ from $\T_*^+$.  Thus $i+1 = j = (t+2)/4-1$.  This leaf corresponds to a solid torus $V^+$ with interior disjoint from $K$ that contains $\smash{A_{(t+2)/4-1}^+}$ in its boundary.  Since $V^+$ contains neither $\smash{A_1^+}$ nor $\smash{A_1^-}$, the core of $R^+$ is a longitudinal curve of $V^+$.  It follows that there are meridional disks of $V^+$ which form low (or high) disks for the arcs of $A_{(t+2)/4-1}^+ \cap K$.  Let $D^+$ be one of these disks.

Similarly since $k$ and $l$ have different parity, then $k+1 = l = (t+2)/4$ and there exists a solid torus $V^-$ such that $\bdry V^- = R^- \cup A_{(t+2)/4}^-$.  Furthermore there are meridional disks of $V^-$ which form high (or low) disks for the arcs of $A_{(t+2)/4-1}^- \cap K$.  Let $D^-$ be one of these disks.

Note that the annuli $A_{\smash[b]{(t+2)/4-1}}^+$ and $A_{\smash[b]{(t+2)/4-1}}^-$ lie on opposite sides of $\hatT$.  Thus $V^+$ and $V^-$ lie on opposite sides of $\hatT$.  Therefore $D^+$ and $D^-$ form a pair of disjoint high and low disks for $K$.  By \fullref{highdisklowdisk} this contradicts the thinness of $K$.
\end{proof}

\begin{lemma}\label{alternatingcurves}
 We have $\bdry T_{\smash[b]{(t+2)/4}} = a_{\smash[b]{(t+2)/4-1}}^+ \cup a_{\smash[b]{(t+2)/4-1}}^-$.  Furthermore,  except for $T_{(t+2)/4}$ and the annulus (other than $T_{(t+2)/4}$ if $t=6$)  bounded by $a_1^+ \cup a_1^-$, all of the other annuli of $\hatT \cut (A^+ \cup A^-)$ are bounded by either $a_{i-1}^+ \cup a_{i}^-$ or $a_{i-1}^- \cup a_{i}^+$ for $i \in \{2, \dots, (t+2)/4-1\}$.
\end{lemma}

\begin{proof}
 \fullref{tis61014} implies that $4{\nmid}t$.  Recall that we are assuming $t \geq 6$.
  
The lemma follows immediately for $t=6$ since in this case $(t+2)/4-1 = 1$.  We have the complexes $A^- = \smash{A_1^-}$ and $A^+ = \smash{A_1^+}$ and the annulus $\smash{A_2^-}$ which give the three curves $\smash{a_1^-}$, $\smash{a_1^+}$, and $\smash{a_2^-}$ on $\hatT$.  There are only two annuli of $\hatT \cut (A^- \cup A^+)$---one of which is $T_{(t+2)/4}=T_2$---and both have boundary $a_1^+ \cup a_1^-$.

\fullref{differentboundarycurves} which applies for $t \geq 10$ implies that each annulus of $\hatT \cut (A^+ \cup A^-)$ has boundary $a_i^+ \cup a_j^-$ for some $i$ and some $j \in \{1, \dots, (t+2)/4-1\}$.   

Consider the edges of the trees $\T^+$ and $\T^-$.  Since each edge of $\T^+$ (resp.\ $\T^-$) corresponds to an annulus of $\hatT \cut A^+$ (resp.\ $\hatT \cut A^-$), such an annulus must intersect $A^-$ (resp.\ $A^+$) exactly once.  Therefore the annulus corresponding to an edge incident to a leaf of $\T^+$ (resp.\ $\T^-$) must intersect $A^-$ (resp.\ $A^+$) either in $a_1^-$ or $\smash{a_{(t+2)/4-1}^-}$ (resp.\ $\smash{a_1^+}$ or $\smash{a_{(t+2)/4-1}^+}$).  Thus each tree $\T^+$ and $\T^-$ may only have two leaves and is thus homeomorphic to a line segment.  

The lemma now follows for $t=10$; see \fullref{fig:tis10A+-}(a) for the case $t = 10$.  Hence we may assume $t\geq14$.

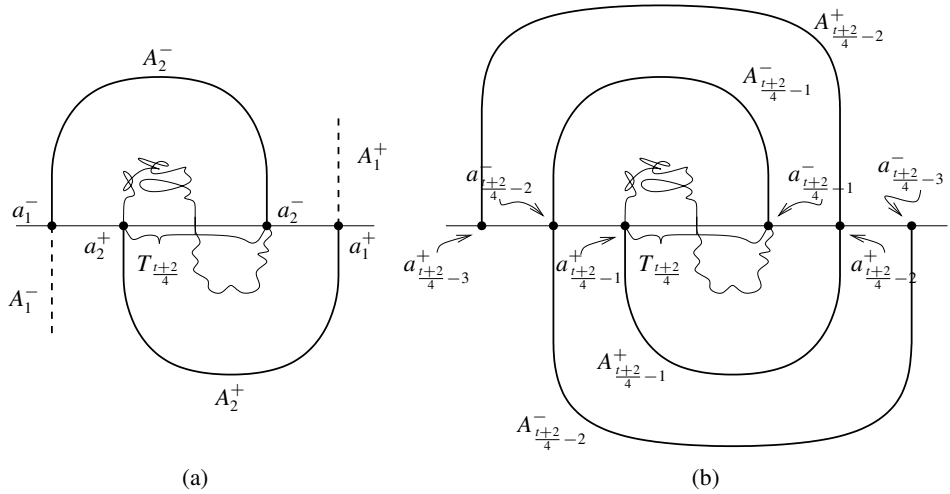
\begin{figure}[ht!]
\centering
\begin{picture}(0,0)%
\includegraphics{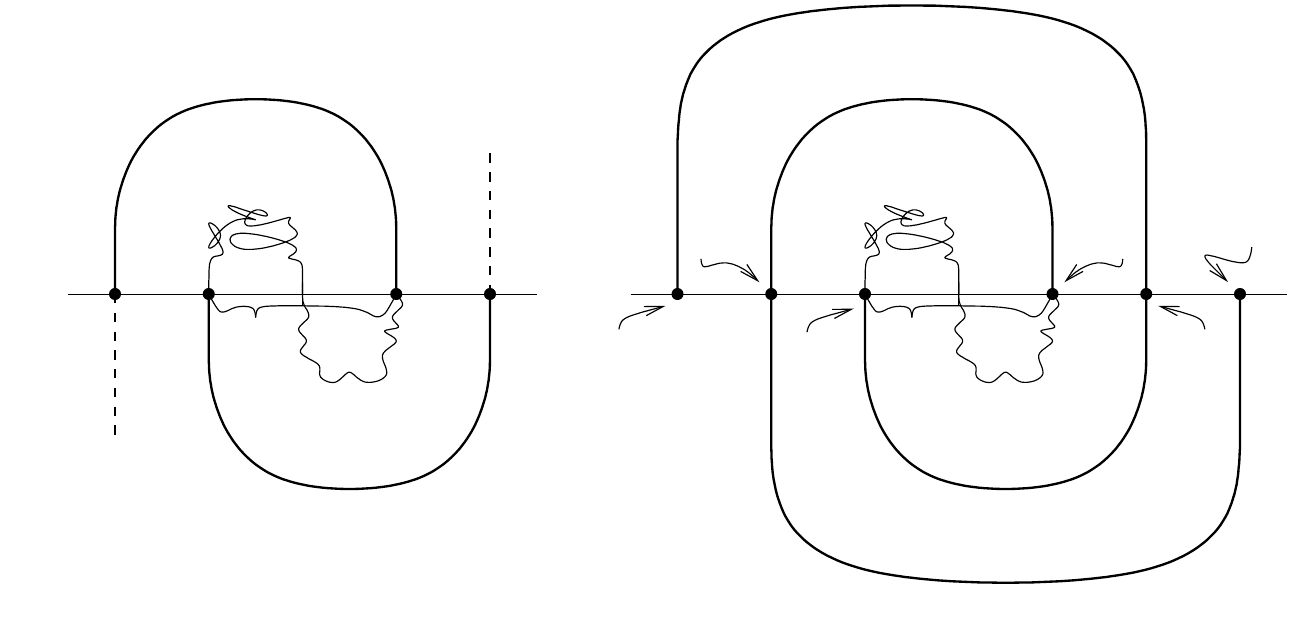}%
\end{picture}%
\setlength{\unitlength}{2960sp}%
\begingroup\makeatletter
\gdef\SetFigFont#1#2#3#4#5{%
  \reset@font\fontsize{#1}{#2pt}%
  \fontfamily{#3}\fontseries{#4}\fontshape{#5}%
  \selectfont}%
\endgroup%
\begin{picture}(8249,4098)(2789,-8590)
\put(10726,-5911){\makebox(0,0)[b]{\smash{{\SetFigFont{9}{10.8}{\rmdefault}{\mddefault}{\updefault}{\color[rgb]{0,0,0}$a_{\frac{t+2}{4}-3}^-$}%
}}}}
\put(4726,-8536){\makebox(0,0)[b]{\smash{{\SetFigFont{9}{10.8}{\rmdefault}{\mddefault}{\updefault}{\color[rgb]{0,0,0}(a)}%
}}}}
\put(5026,-7861){\makebox(0,0)[b]{\smash{{\SetFigFont{9}{10.8}{\rmdefault}{\mddefault}{\updefault}{\color[rgb]{0,0,0} $A_2^+$}%
}}}}
\put(4426,-5011){\makebox(0,0)[b]{\smash{{\SetFigFont{9}{10.8}{\rmdefault}{\mddefault}{\updefault}{\color[rgb]{0,0,0} $A_2^-$}%
}}}}
\put(6076,-5836){\makebox(0,0)[lb]{\smash{{\SetFigFont{9}{10.8}{\rmdefault}{\mddefault}{\updefault}{\color[rgb]{0,0,0}$\smash{A_1^+}$}%
}}}}
\put(3451,-6286){\makebox(0,0)[rb]{\smash{{\SetFigFont{9}{10.8}{\rmdefault}{\mddefault}{\updefault}{\color[rgb]{0,0,0}$a_1^-$}%
}}}}
\put(4051,-6586){\makebox(0,0)[rb]{\smash{{\SetFigFont{9}{10.8}{\rmdefault}{\mddefault}{\updefault}{\color[rgb]{0,0,0}$a_2^+$}%
}}}}
\put(6001,-6586){\makebox(0,0)[lb]{\smash{{\SetFigFont{9}{10.8}{\rmdefault}{\mddefault}{\updefault}{\color[rgb]{0,0,0}$a_1^+$}%
}}}}
\put(5401,-6286){\makebox(0,0)[lb]{\smash{{\SetFigFont{9}{10.8}{\rmdefault}{\mddefault}{\updefault}{\color[rgb]{0,0,0}$a_2^-$}%
}}}}
\put(4426,-6736){\makebox(0,0)[b]{\smash{{\SetFigFont{9}{10.8}{\rmdefault}{\mddefault}{\updefault}{\color[rgb]{0,0,0}$T_{\frac{t+2}{4}}$}%
}}}}
\put(3451,-7036){\makebox(0,0)[rb]{\smash{{\SetFigFont{9}{10.8}{\rmdefault}{\mddefault}{\updefault}{\color[rgb]{0,0,0}$\smash{A_1^-}$}%
}}}}
\put(8626,-6736){\makebox(0,0)[b]{\smash{{\SetFigFont{9}{10.8}{\rmdefault}{\mddefault}{\updefault}{\color[rgb]{0,0,0}$T_{\frac{t+2}{4}}$}%
}}}}
\put(9901,-4711){\makebox(0,0)[lb]{\smash{{\SetFigFont{9}{10.8}{\rmdefault}{\mddefault}{\updefault}{\color[rgb]{0,0,0}$A_{\frac{t+2}{4}-2}^+$}%
}}}}
\put(9001,-8536){\makebox(0,0)[b]{\smash{{\SetFigFont{9}{10.8}{\rmdefault}{\mddefault}{\updefault}{\color[rgb]{0,0,0}(b)}%
}}}}
\put(9301,-5161){\makebox(0,0)[lb]{\smash{{\SetFigFont{9}{10.8}{\rmdefault}{\mddefault}{\updefault}{\color[rgb]{0,0,0} $A_{\frac{t+2}{4}-1}^-$}%
}}}}
\put(10201,-6736){\makebox(0,0)[lb]{\smash{{\SetFigFont{9}{10.8}{\rmdefault}{\mddefault}{\updefault}{\color[rgb]{0,0,0}$a_{\frac{t+2}{4}-2}^+$}%
}}}}
\put(7051,-6736){\makebox(0,0)[rb]{\smash{{\SetFigFont{9}{10.8}{\rmdefault}{\mddefault}{\updefault}{\color[rgb]{0,0,0}$a_{\frac{t+2}{4}-3}^+$}%
}}}}
\put(8701,-7561){\makebox(0,0)[rb]{\smash{{\SetFigFont{9}{10.8}{\rmdefault}{\mddefault}{\updefault}{\color[rgb]{0,0,0} $A_{\frac{t+2}{4}-1}^+$}%
}}}}
\put(8026,-8086){\makebox(0,0)[rb]{\smash{{\SetFigFont{9}{10.8}{\rmdefault}{\mddefault}{\updefault}{\color[rgb]{0,0,0}$A_{\frac{t+2}{4}-2}^-$}%
}}}}
\put(9676,-5986){\makebox(0,0)[lb]{\smash{{\SetFigFont{9}{10.8}{\rmdefault}{\mddefault}{\updefault}{\color[rgb]{0,0,0}$a_{\frac{t+2}{4}-1}^-$}%
}}}}
\put(8326,-6736){\makebox(0,0)[rb]{\smash{{\SetFigFont{9}{10.8}{\rmdefault}{\mddefault}{\updefault}{\color[rgb]{0,0,0}$a_{\frac{t+2}{4}-1}^+$}%
}}}}
\put(7576,-5986){\makebox(0,0)[rb]{\smash{{\SetFigFont{9}{10.8}{\rmdefault}{\mddefault}{\updefault}{\color[rgb]{0,0,0}$a_{\frac{t+2}{4}-2}^-$}%
}}}}
\end{picture}%
\caption{(a) $A^+$ and $A^-$ for $t=10$\qua  (b) The ends of $A^+$ and $A^-$ for $t > 10$}
\label{fig:tis10A+-}
\end{figure}

Since each tree is homeomorphic to a line segment, the two annuli of $\hatT \cut A^-$ (resp.\ $\hatT \cut A^+$) with boundary containing the curve $\smash{a_{(t+2)/4-1}^-}$ (resp.\ $a_{(t+2)/4-1}$) also contain either $\smash{a_{(t+2)/4-2}^-}$ or $\smash{a_{(t+2)/4-3}^-}$ (resp.\ $\smash{a_{(t+2)/4-2}^+}$ or $\smash{a_{(t+2)/4-3}^+}$).  Note that one of these two annuli must contain $T_{(t+2)/4}$.  

By \fullref{differentboundarycurves} each of these two annuli must intersect $A^+$ (resp.\ $A^-$) in exactly one curve.  Indeed, since the annulus with boundary $a_{(t+2)/4-2}^-  \cup a_{(t+2)/4-1}^-$ (resp.\ $\smash{a_{(t+2)/4-2}^+  \cup a_{(t+2)/4-1}^+}$) corresponds to an edge of $\T^+$ which is incident to a leaf, it must intersect either $\smash{a_1^+}$ or $\smash{a_{(t+2)/4-1}^+}$ (resp.\ $\smash{a_1^-}$ or $\smash{a_{(t+2)/4-1}^-}$).  

If it intersects $a_1^+$ (resp.\ $a_1^-$), then the other annulus must contain $T_{(t+2)/4}$.  In this case, $$\bdry T_{(t+2)/4} = a_{(t+2)/4-1}^- \cup a_2^+ \quad \text{(resp. }a_{(t+2)/4-1}^+ \cup a_2^-).$$  Moreover, $T_{(t+2)/4}$ must thus be contained in the annulus of $\hatT \cut A^+$ (resp.\ $\hatT \cut A^-$) with boundary $a_1^- \cup a_2^-$ (resp.\ $a_1^+ \cup a_2^+$).  This however forces $a_2^- = a_{(t+2)/4-1}^-$ (resp.\ $\smash{a_2^+ = a_{(t+2)/4-1}^+}$.)  Hence $t = 10$ contrary to our assumption.  
This implies $$\bdry T_{(t+2)/4} = a_{(t+2)/4-1}^+ \cup a_{(t+2)/4-1}^-.$$  Because each $\T^+$ and $\T^-$ is homeomorphic to a line segment, $A^+$ and $A^-$ spiral around one another as in \fullref{fig:tis10A+-}(b).  The remainder of the lemma then follows. 
\end{proof}

\begin{thm}\label{doubleunfurl}
$t \leq 6$
\end{thm}

\begin{proof}
This proof is phrased and its figures drawn so that $A_{(t+2)/4-1}^-$ is above $\hatT$.  \fullref{tis61014} implies that $4{\nmid}t$.  Hence we assume $t \geq 10$.  If $t \equiv 2 \mod 8$, this coincides with our convention that $K_{(t, 1)}$ is below $\hatT$.  If $t \equiv 6 \mod 8$, then one ought to flip the figures and make the appropriate corresponding changes to the language.

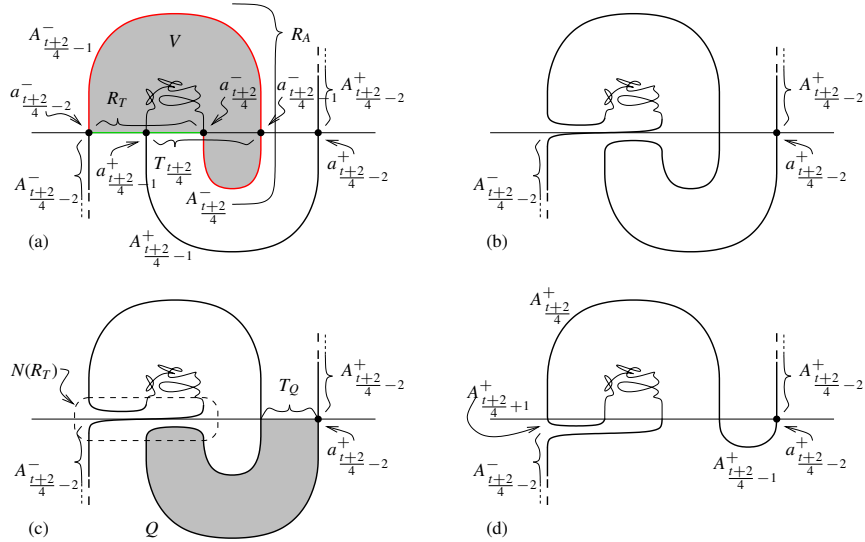
\begin{figure}[ht!]
\centering
\begin{picture}(0,0)%
\includegraphics{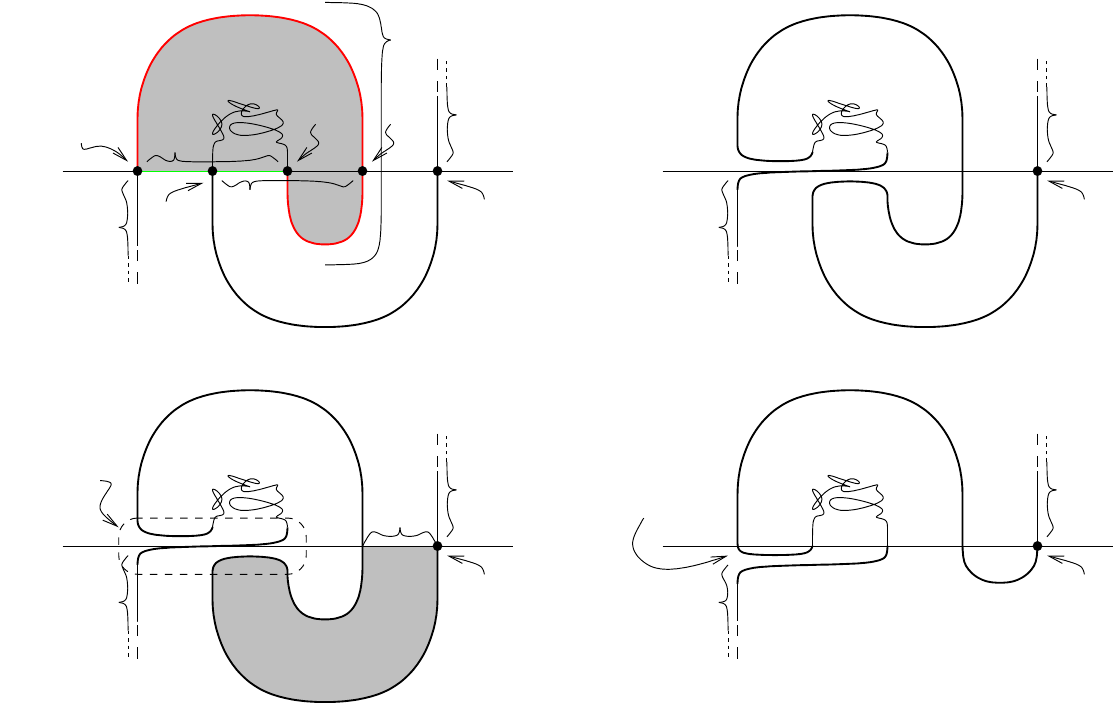}%
\end{picture}%
\setlength{\unitlength}{2368sp}%
\begingroup\makeatletter
\gdef\SetFigFont#1#2#3#4#5{%
  \reset@font\fontsize{#1}{#2pt}%
  \fontfamily{#3}\fontseries{#4}\fontshape{#5}%
  \selectfont}%
\endgroup%
\begin{picture}(8912,5631)(2501,-11830)
\put(10876,-7936){\makebox(0,0)[lb]{\smash{{\SetFigFont{7}{8.4}{\rmdefault}{\mddefault}{\updefault}{\color[rgb]{0,0,0}$a_{\frac{t+2}{4}-2}^+$}%
}}}}
\put(11026,-7111){\makebox(0,0)[lb]{\smash{{\SetFigFont{7}{8.4}{\rmdefault}{\mddefault}{\updefault}{\color[rgb]{0,0,0}$A_{\frac{t+2}{4}-2}^+$}%
}}}}
\put(8326,-8161){\makebox(0,0)[rb]{\smash{{\SetFigFont{7}{8.4}{\rmdefault}{\mddefault}{\updefault}{\color[rgb]{0,0,0}$A_{\frac{t+2}{4}-2}^-$}%
}}}}
\put(7876,-8761){\makebox(0,0)[b]{\smash{{\SetFigFont{7}{8.4}{\rmdefault}{\mddefault}{\updefault}{\color[rgb]{0,0,0}(b)}%
}}}}
\put(10876,-10936){\makebox(0,0)[lb]{\smash{{\SetFigFont{7}{8.4}{\rmdefault}{\mddefault}{\updefault}{\color[rgb]{0,0,0}$a_{\frac{t+2}{4}-2}^+$}%
}}}}
\put(11026,-10111){\makebox(0,0)[lb]{\smash{{\SetFigFont{7}{8.4}{\rmdefault}{\mddefault}{\updefault}{\color[rgb]{0,0,0}$A_{\frac{t+2}{4}-2}^+$}%
}}}}
\put(8326,-11161){\makebox(0,0)[rb]{\smash{{\SetFigFont{7}{8.4}{\rmdefault}{\mddefault}{\updefault}{\color[rgb]{0,0,0}$A_{\frac{t+2}{4}-2}^-$}%
}}}}
\put(8701,-9361){\makebox(0,0)[rb]{\smash{{\SetFigFont{7}{8.4}{\rmdefault}{\mddefault}{\updefault}{\color[rgb]{0,0,0}$A_{\frac{t+2}{4}}^+$}%
}}}}
\put(8251,-10336){\makebox(0,0)[rb]{\smash{{\SetFigFont{7}{8.4}{\rmdefault}{\mddefault}{\updefault}{\color[rgb]{0,0,0}$A_{\frac{t+2}{4}+1}^+$}%
}}}}
\put(10126,-11086){\makebox(0,0)[lb]{\smash{{\SetFigFont{7}{8.4}{\rmdefault}{\mddefault}{\updefault}{\color[rgb]{0,0,0}$A_{\frac{t+2}{4}-1}^+$}%
}}}}
\put(7876,-11761){\makebox(0,0)[b]{\smash{{\SetFigFont{7}{8.4}{\rmdefault}{\mddefault}{\updefault}{\color[rgb]{0,0,0}(d)}%
}}}}
\put(6076,-7936){\makebox(0,0)[lb]{\smash{{\SetFigFont{7}{8.4}{\rmdefault}{\mddefault}{\updefault}{\color[rgb]{0,0,0}$a_{\frac{t+2}{4}-2}^+$}%
}}}}
\put(3451,-7186){\makebox(0,0)[rb]{\smash{{\SetFigFont{7}{8.4}{\rmdefault}{\mddefault}{\updefault}{\color[rgb]{0,0,0}$a_{\frac{t+2}{4}-2}^-$}%
}}}}
\put(3901,-7261){\makebox(0,0)[b]{\smash{{\SetFigFont{7}{8.4}{\rmdefault}{\mddefault}{\updefault}{\color[rgb]{0,0,0}$R_T$}%
}}}}
\put(4501,-6661){\makebox(0,0)[b]{\smash{{\SetFigFont{7}{8.4}{\rmdefault}{\mddefault}{\updefault}{\color[rgb]{0,0,0}$V$}%
}}}}
\put(6076,-10936){\makebox(0,0)[lb]{\smash{{\SetFigFont{7}{8.4}{\rmdefault}{\mddefault}{\updefault}{\color[rgb]{0,0,0}$a_{\frac{t+2}{4}-2}^+$}%
}}}}
\put(3676,-6586){\makebox(0,0)[rb]{\smash{{\SetFigFont{7}{8.4}{\rmdefault}{\mddefault}{\updefault}{\color[rgb]{0,0,0} $A_{\frac{t+2}{4}-1}^-$}%
}}}}
\put(5701,-10261){\makebox(0,0)[b]{\smash{{\SetFigFont{7}{8.4}{\rmdefault}{\mddefault}{\updefault}{\color[rgb]{0,0,0}$T_Q$}%
}}}}
\put(4351,-11761){\makebox(0,0)[rb]{\smash{{\SetFigFont{7}{8.4}{\rmdefault}{\mddefault}{\updefault}{\color[rgb]{0,0,0}$Q$}%
}}}}
\put(4726,-8761){\makebox(0,0)[rb]{\smash{{\SetFigFont{7}{8.4}{\rmdefault}{\mddefault}{\updefault}{\color[rgb]{0,0,0} $A_{\frac{t+2}{4}-1}^+$}%
}}}}
\put(5701,-6586){\makebox(0,0)[lb]{\smash{{\SetFigFont{7}{8.4}{\rmdefault}{\mddefault}{\updefault}{\color[rgb]{0,0,0}$R_A$}%
}}}}
\put(3301,-10111){\makebox(0,0)[rb]{\smash{{\SetFigFont{7}{8.4}{\rmdefault}{\mddefault}{\updefault}{\color[rgb]{0,0,0}$N(R_T)$}%
}}}}
\put(4501,-7936){\makebox(0,0)[b]{\smash{{\SetFigFont{7}{8.4}{\rmdefault}{\mddefault}{\updefault}{\color[rgb]{0,0,0}$T_{\frac{t+2}{4}}$}%
}}}}
\put(4351,-8011){\makebox(0,0)[rb]{\smash{{\SetFigFont{7}{8.4}{\rmdefault}{\mddefault}{\updefault}{\color[rgb]{0,0,0}$a_{\frac{t+2}{4}-1}^+$}%
}}}}
\put(4951,-7111){\makebox(0,0)[lb]{\smash{{\SetFigFont{7}{8.4}{\rmdefault}{\mddefault}{\updefault}{\color[rgb]{0,0,0}$a_{\frac{t+2}{4}}^-$}%
}}}}
\put(5551,-7111){\makebox(0,0)[lb]{\smash{{\SetFigFont{7}{8.4}{\rmdefault}{\mddefault}{\updefault}{\color[rgb]{0,0,0}$a_{\frac{t+2}{4}-1}^-$}%
}}}}
\put(6226,-7111){\makebox(0,0)[lb]{\smash{{\SetFigFont{7}{8.4}{\rmdefault}{\mddefault}{\updefault}{\color[rgb]{0,0,0}$A_{\frac{t+2}{4}-2}^+$}%
}}}}
\put(6226,-10111){\makebox(0,0)[lb]{\smash{{\SetFigFont{7}{8.4}{\rmdefault}{\mddefault}{\updefault}{\color[rgb]{0,0,0}$A_{\frac{t+2}{4}-2}^+$}%
}}}}
\put(3526,-8161){\makebox(0,0)[rb]{\smash{{\SetFigFont{7}{8.4}{\rmdefault}{\mddefault}{\updefault}{\color[rgb]{0,0,0}$A_{\frac{t+2}{4}-2}^-$}%
}}}}
\put(3526,-11161){\makebox(0,0)[rb]{\smash{{\SetFigFont{7}{8.4}{\rmdefault}{\mddefault}{\updefault}{\color[rgb]{0,0,0}$A_{\frac{t+2}{4}-2}^-$}%
}}}}
\put(5101,-8311){\makebox(0,0)[rb]{\smash{{\SetFigFont{7}{8.4}{\rmdefault}{\mddefault}{\updefault}{\color[rgb]{0,0,0}$A_{\frac{t+2}{4}}^-$}%
}}}}
\put(3076,-11761){\makebox(0,0)[b]{\smash{{\SetFigFont{7}{8.4}{\rmdefault}{\mddefault}{\updefault}{\color[rgb]{0,0,0}(c)}%
}}}}
\put(3076,-8761){\makebox(0,0)[b]{\smash{{\SetFigFont{7}{8.4}{\rmdefault}{\mddefault}{\updefault}{\color[rgb]{0,0,0}(a)}%
}}}}
\end{picture}%
\caption{(a) The ``ends'' of $A^+$ and $A^-$, $A_{(t+2)/4}^-$, $A_R$, $T_R$, and $V$\qua  (b) The result of an unfurling isotopy along $R = A_R \cup T_R$\qua  (c) The subannulus $Q \subseteq \Theta_{2\pi}(A^+)$ and the subannulus $T_Q \subseteq \hatT$\qua  (d) An isotopy of $Q$ towards $T_Q$ and a further slight isotopy in $N(T_R)$}
\label{fig:doubleunfurl}
\end{figure}

By \fullref{alternatingcurves}, the ``ends'' of $A^+$ and $A^-$ are as in \fullref{fig:tis10A+-}.  Let $T_R$ be the annulus on $\hatT$ bounded by $a_{(t+2)/4}^-$ and $a_{(t+2)/4-2}^-$ that contains $a_{(t+2)/4-1}^+$.  Let $A_R$ be the union $\smash{A_{(t+2)/4-1}^- \cup A_{(t+2)/4}^-}$.  By \fullref{dividingtorus}, the torus $R = T_R \cup A_R$ bounds a solid torus $V$.  Note that $V$ contains the two arcs $K - (A^+ \cup A^- \cup \smash{A_{(t+2)/4}^-})$ in its interior, the arcs $K \cap (\smash{A_{(t+2)/4-1}^- \cup A_{(t+2)/4}^-})$ on its boundary, as well as the points $K \cap a_{(t+2)/4-1}^+$ on its boundary.  See \fullref{fig:doubleunfurl}(a).  

\vspace{5pt}

Let $\delta$ and $\delta'$ be parallel simple closed curves on $R$ that transversely cross $T_R$ and $A_R$ just once so that $(\delta \cup \delta') \cap A_R = K \cap A_R$ and $(\delta \cup \delta') \cap a_{(t+2)/4-1}^+ \subseteq K \cap a_{(t+2)/4-1}^+$.  Orient $\delta$ so that it crosses from $T_R$ into $A_{(t+2)/4}^-$.
Although $K$ intersects the interior of $T_R$ in two points, we may perform a Dehn twist $\Theta_u$ along the torus $R$ in the direction of $\delta$ analogous to the unfurling isotopy of \fullref{main-unfurl}.  We must mind the effect of the isotopy on the arcs $K \cap N(a_{(t+2)/4-1}^+)$.

\vspace{5pt}

The annulus $\Theta_{2 \pi}(A_R)$ lies as $T_R$.  A slight further isotopy makes this annulus and the arcs of $\Theta_{2 \pi}(K)$ on it transverse to the height function.   These arcs of $K$ now have no critical points; at least four fewer than before the isotopy.  The critical points of $K$ that were on $A_{(t+2)/4-1}^-$ and $A_{(t+2)/4}^-$ have been removed.  See \fullref{fig:doubleunfurl}(b).

The annulus $A_{(t+2)/4-1}^+$ (slightly extended through $N(R)$) near $a_{(t+2)/4-1}^+$ is spun once around the meridian of $R$ by $\Theta_u$.  The resulting annulus may be regarded as the double-curve sum of $A_{(t+2)/4-1}^+$ with $R$ (suitably oriented).  The arcs of $\Theta_{2\pi}(K)$ on this annulus may similarly be regarded as resulting from this double-curve sum of $K \cap N(A_{(t+2)/4-1}^+)$ with the two curves $\delta$ and $\delta'$ on $R$.  A further slight isotopy with support in $N(T_R)$ makes this resulting annulus and arcs of $\Theta_{2\pi}(K)$ on it transverse to the height function.  We may further assume that after these isotopies the arcs of $\Theta_{2\pi}(K \cap N(a_{(t+2)/4-1}^+))$ in $N(T_R)$ each have just one critical point.  See again \fullref{fig:doubleunfurl}(b).  Indeed, outside of $N(T_R)$, we have $K = \Theta_{2\pi}(K)$ and $A^- \cup A_{(t+2)/4}^- \cup A^+ = \Theta_{2\pi}(A^- \cup A_{(t+2)/4}^- \cup A^+)$.  

After these isotopies there are four ``new'' critical points of $\Theta_{2\pi}(K)$:  two minima just above $\hatT$ and two maxima just below $\hatT$.  Due to these critical points, the width of $\Theta_{2\pi}(K)$ is currently greater than the width of $K$.

Consider the subannulus $Q$ of $\Theta_{2 \pi}(A^+)$ lying between the curves formerly labeled $a_{\smash{(t+2)/4-1}}^-$ and $a_{(t+2)/4-2}^+$.  Let $T_Q$ be the subannulus of $\hatT$ bounded by $\bdry Q$ with interior disjoint from $\Theta_{2\pi}(K)$.  Note that $Q$ is parallel to $T_Q$ through the solid torus that they together bound.  This solid torus and the annuli $Q$ and $T_Q$ are shown in \fullref{fig:doubleunfurl}(c).  We may isotop $Q$ along with the arcs of $\Theta_{2\pi}(K) \cap Q$ towards $T_Q$ so that the arcs  $\Theta_{2\pi}(K) \cap Q$ after this isotopy each have only one minimum just below $\hatT$.  We may also slightly isotop the other annuli and arcs of the knot in $N(T_R)$ downwards so that the minima that were just above $\hatT$ are now just below $\hatT$.  The result of these isotopies may be seen in \fullref{fig:doubleunfurl}(d).  Notice the collection of annuli and arcs of $K$ now above $\hatT$ is indistinguishable from the collection of annuli and arcs of $K$ above $\hatT$ prior to the isotopies.  Furthermore, notice that (mainly due to the isotopy of $Q$) the width of the resulting position of the knot is at most the width of $K$.  Since $K$ was originally in thin position, these widths must be equal.

Let us continue to refer to the annuli $A_i^\pm$ that are unaffected by these isotopies by their former labels.  The annulus $A_{(t+2)/4-2}^-$ has been elongated (in a height-preserving manner) with $\Theta_{2\pi}(A_R)$; we shall also refer to this elongated annulus by $\smash{A_{(t+2)/4-2}^-}$.  After the isotopies, $\Theta_{2\pi}(\smash{A_{(t+2)/4-1}^+})$ is cut into three annuli by $\hatT$.  Let us use $\smash{A_{(t+2)/4-1}^+}$ to refer to the isotoped annulus $Q$ and the subsequent two annuli as $\smash{A_{(t+2)/4}^+}$ and $\smash{A_{(t+2)/4+1}^+}$.  Note that the post-isotopies annulus $\smash{A_{(t+2)/4}^+}$ coincides with the pre-isotopies annulus $\smash{A_{(t+2)/4-1}^-}$.  Again, see \fullref{fig:doubleunfurl}(d).  Note also that $\hatT$ remains a thick level.
\newpage
Prior to any isotopies in this proof, \fullref{highestmin2} implies that the lowest maximum of $K$ above $\hatT$ lies on an arc $\kappa$ of $K \cut \hatT$ that together with an arc $\tau$ of $\hatT$ bounds a high disk $\Delta$ with interior disjoint from $A^+ \cup A^-$.  Since the interior of $\Delta$ is disjoint from $A_{(t+2)/4-1}^-$, either $\tau \subseteq T_R \cup T_{(t+2)/4}$ or $\tau \subseteq \hatT - \Int(T_R \cup T_{(t+2)/4})$.  Because the annuli and arcs of $K$ above $\hatT$ before the isotopies are indistinguishable from the annuli and arcs of $K$ above $\hatT$ after the isotopies, the high disk $\Delta$ exists after these isotopies too.

After the isotopies, an arc of $K$ on the annulus $\smash{A_{(t+2)/4+1}^+}$ bounds a low disk $\Delta_{+1}$ with an arc of $T_R \cut T_{(t+2)/4}$, and an arc of $K$ on the annulus $A_{(t+2)/4-1}^+$ bounds a low disk $\Delta_{-1}$ with an arc of $T_Q$.  Note that the interior of $\tau$ is disjoint from at least one of $\Delta_{+1} \cap \hatT$ or $\Delta_{-1} \cap \hatT$.  Since the position of $K$ after the isotopies has width at most that of $K$ before the isotopies, $\Delta$ together with either $\Delta_{+1}$ or $\Delta_{-1}$ forms a pair of high and low disks that by \fullref{highdisklowdisk} contradict the thinness of the original thin position of $K$.
\end{proof}

\section[The case t=6]{The case $t=6$}\label{sec:tis6}

In this section we show that if $K$ is in thin position, then $t \neq 6$.  Assuming $t=6$, our approach is to isotop $\Int S$ fixing $K$ so that we may gain an understanding of the resulting pieces of $S \cut T$ on each side of $\hatT$.  Then since there are $3s$ arcs of $S \cap T$ and $s \geq 4g-1$, we obtain conflicting estimates on the Euler characteristic of $S$.

\subsection[Isotopies of S]{Isotopies of $S$}

By \fullref{oneconfig} we may assume that among the graphs $\smash{G_S^x}$ for $x \in \mathbf{t}$ there is a face $F$ bounded by an extended $S2$ cycle or an extended $S3$ cycle accounting for the labels $\{1, 2, 5, 6\}$ and a face $g$ bounded by an $S2$ cycle or an $S3$ cycle accounting for the labels $\{3, 4\}$.  Let $f$ be the face of $G_S$ in $F$ bounded by the Scharlemann cycle.

\begin{prop} \label{prop:isotopS} After perhaps reassigning labels of $\mathbf{t}$ by $x \mapsto 7-x$, one may isotop the interior of $S$ so that for some integer $n \geq 0$ the following holds: 
\begin{enumerate}
\item \vspace{10pt}\begin{itemize}
\item The resulting arcs of $S \cap T$ are essential in $S$ and in $T$, and
\item either $|S \cap T|$ remains minimized over isotopies of $\Int S$ or there exists annuli in $S \cut T$.
\end{itemize}
\item \vspace{10pt}\begin{itemize}
\item For each label pair $\lambda = \{1, 6\}, \{2, 5\}, \{3, 4\}$ there are $s-n$ edges of $G_S$ with label pair $\lambda$, and they lie in an essential annulus in $\hatT$,
\item for each label pair $\lambda = \{1, 4\}, \{2, 3\}, \{5, 6\}$ there are $n$ edges of $G_S$ with label pair $\lambda$ and
\item for each label pair $\lambda = \{1, 2\}, \{3, 6\}, \{4, 5\}$ there are no edges of $G_S$ with label pair $\lambda$.
\end{itemize}

\item \vspace{10pt}\begin{itemize}
\item If $f$ is a bigon, then there are $(s-n)/2$ faces of $S \cap X^-$ parallel to $f$, and
\item if $f$ is a trigon, then there are either $(s-n)/3$ or $(s-2n)/3$ faces of $S \cap X^-$ parallel to $f$.
\end{itemize}
\item \vspace{10pt}\begin{itemize}
\item If $g$ is a bigon then $S \cap X^+$ is a collection of $3(s-n)/2$ bigons and $n$ trigons, and
\item if $g$ is a trigon then $S \cap X^+$ is either a collection of $s-2n$ bigons and $(s+4n)/3$ trigons or a collection of $s-n$ bigons and $(s+2n)/3$ trigons.
\end{itemize}
\end{enumerate}
\end{prop} 

Let $V$ be the solid torus component of $X^+ \cut (\bar{N}(K) \cup F)$.  Let $W$ be the solid torus $X^+ \cut (V \cup \bar{N}(K) \cup g)$.  Let $T_V$ be the annulus of $\hatT \cap V$ so that $T_V$ is bounded by edges of $G_T$ with label pairs $\{1,6\}$ and $\{2, 5\}$.  Observe that $\bdry V \cut T_V$ is parallel through $V$ onto $T_V$.  Let $T_W$ and $T_W'$ be the two annuli of $\hatT \cap W$ where $T_W$ is bounded by edges of $G_T$ with label pairs $\{3, 4\}$ and $\{1, 6\}$ and $T_W'$ is bounded by edges of $G_T$ with label pairs $\{3, 4\}$ and $\{2, 5\}$.  See \fullref{fig:toriVandW}.
\def\SetFigFont#1#2#3#4#5{\small}%
\begin{figure}[ht!]
\centering
\begin{picture}(0,0)%
\includegraphics{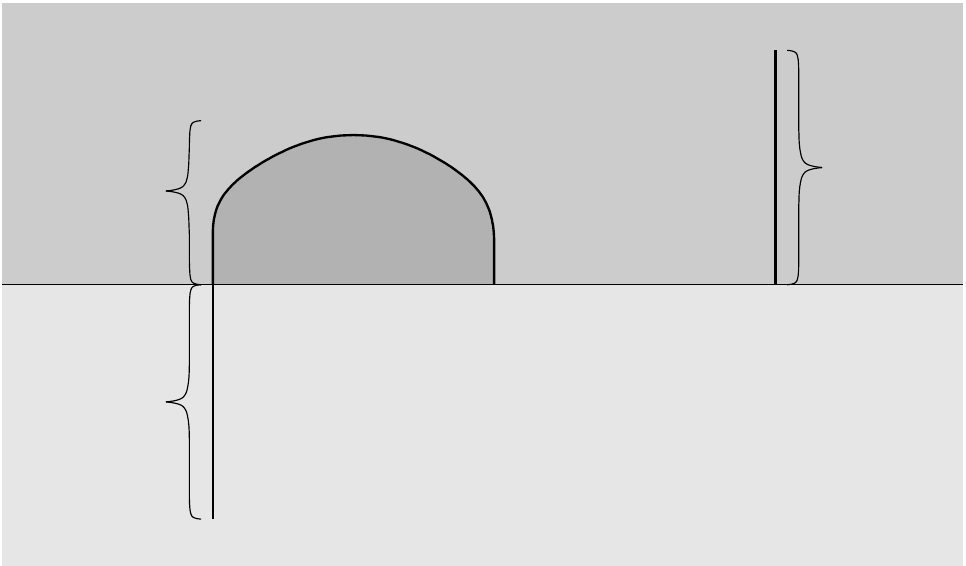}%
\end{picture}%
\setlength{\unitlength}{2960sp}%
\begingroup\makeatletter\ifx\SetFigFont\undefined%
\gdef\SetFigFont#1#2#3#4#5{%
  \reset@font\fontsize{#1}{#2pt}%
  \fontfamily{#3}\fontseries{#4}\fontshape{#5}%
  \selectfont}%
\fi\endgroup%
\begin{picture}(6174,3602)(1039,-6962)
\put(3301,-4861){\makebox(0,0)[b]{\smash{{\SetFigFont{9}{10.8}{\rmdefault}{\mddefault}{\updefault}{\color[rgb]{0,0,0}$V$}%
}}}}
\put(4951,-4111){\makebox(0,0)[b]{\smash{{\SetFigFont{9}{10.8}{\rmdefault}{\mddefault}{\updefault}{\color[rgb]{0,0,0}$W$}%
}}}}
\put(2026,-4636){\makebox(0,0)[rb]{\smash{{\SetFigFont{9}{10.8}{\rmdefault}{\mddefault}{\updefault}{\color[rgb]{0,0,0}$F\cut f$}%
}}}}
\put(2026,-5986){\makebox(0,0)[rb]{\smash{{\SetFigFont{9}{10.8}{\rmdefault}{\mddefault}{\updefault}{\color[rgb]{0,0,0}$f$}%
}}}}
\put(6376,-4486){\makebox(0,0)[lb]{\smash{{\SetFigFont{9}{10.8}{\rmdefault}{\mddefault}{\updefault}{\color[rgb]{0,0,0}$g$}%
}}}}
\put(4951,-6136){\makebox(0,0)[b]{\smash{{\SetFigFont{9}{10.8}{\rmdefault}{\mddefault}{\updefault}{\color[rgb]{0,0,0}$X^- \cut f$}%
}}}}
\end{picture}%
\caption{The solid tori $V$, $W$, and $X^-\cut f$ schematically}
\label{fig:toriVandW}
\end{figure}
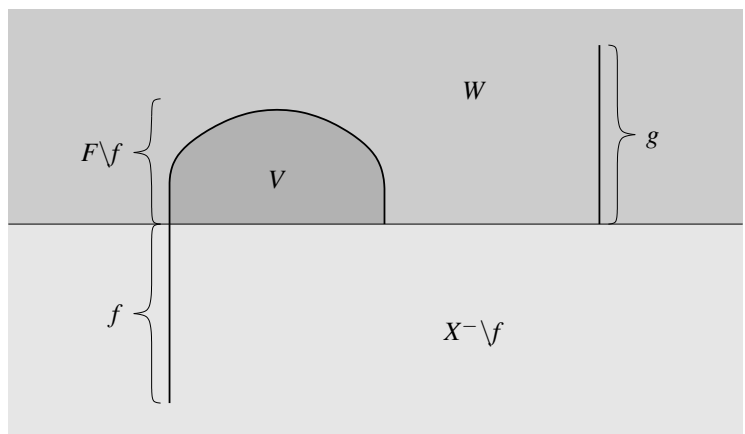

To prove \fullref{prop:isotopS} we begin in $V$ and push what we can of $S$ along $F$ through $X^-$ into $W$ and back into $X^-$.  When in $W$ we also use $g$ to guide the isotopy.  This puts $S \cap X^+$ into a position that we can understand and count.  The unknown parts of $S \cap X^-$ are dealt with in subsequent subsections.


\subsubsection{Pushing $S$ along bigons and trigons}

We say two edges $e$ and $e'$ of $G_T$ are {\em parallel\/} if they bound an embedded bigon on $T$.   
We say two faces $R$ and $R'$ of $S \cut T$ are {\em parallel\/} if $R$ and $R'$ cobound a product region in $(X-N(K)) \cut T$.
\vspace{-4pt}

Similarly we say two edges of $G_T$ or two faces of $S \cut T$ are {\em adjacent\/} if they are parallel and no other edge or face respectively lies between them. 
\vspace{-4pt}

\begin{lemma}\label{lem:paralleldisks}
Let $R$ and $R'$ be two disks of $S \cut T$ such that $\bdry R$ is parallel to $\bdry R'$ on $\bdry (X^\pm - N(K))$.  Then $R$ is parallel to $R'$.  Furthermore, every face of $S \cut T$ in the product region between $R$ and $R'$ is a disk parallel to $R$ and $R'$. 
\vspace{-4pt}
\end{lemma}

\begin{proof}
Since $\bdry R$ is parallel to $\bdry R'$, there is an annulus $A$ on $\bdry (X^\pm - N(K))$ connecting them.  Thus $R \cup A \cup R' \cong S^2$.  Since $X^\pm - N(K)$ is irreducible, $R \cup A \cup R'$ is the boundary of a solid ball.  This ball is the requisite product region $R \times [0,1]$.
\vspace{-4pt}

If $P \in S \cut T$ is contained in the ball bounded by $R \cup A \cup R'$, then each component of $\bdry P$ is an essential curve in $A$.  Since $P$ must be incompressible, $P$ is a disk.  It follows that $P$ is parallel to both $R$ and $R'$.
\end{proof}

\begin{lemma} \label{lem:parallelbigons}
If $B$ and $B'$ are two bigons of $G_S$ on the same side of $\hatT$ with edges $e_1 \subseteq B$ and $\smash{e_1'} \subseteq B'$ that are parallel as edges of $G_T$, then $B$ and $B'$ are parallel in $X$.
\vspace{-4pt}
\end{lemma}

\begin{proof}
We have the edges $e_1$ and $\smash{e_1'}$ of $B$ and $B'$ respectively which are parallel on $T$.  Let $e_2$ and $\smash{e_2'}$ be the other edges of $B$ and $B'$ respectively.  Note that $e_2$ and $\smash{e_2'}$ have the same label pair.  
\vspace{-4pt}

Assume $B$ and $B'$ are not parallel.  By \fullref{lem:paralleldisks} $B$ and $B'$ the edges $e_2$ and $\smash{e_2'}$ lie in an essential annulus $A$.  Due to the existence of $f$ and $g$, \fullref{GT:L2.1} implies that the core of $A$ does not bound a disk in $X$.  
\vspace{-4pt}

Join $B$ and $B'$ together along their corners and their parallel edges $e_1$ and $\smash{e_1'}$.  This forms a disk in $X$ whose boundary is an essential curve in $A$ giving a contradiction.
\end{proof}
\vspace{-4pt}

\begin{lemma}\label{lem:isotopacrossdisk}
Let $R$ be a region of $S \cut T$ in $X^+$. 
Let $D$ be a boundary compressing disk for $R$ in $X^+ - N(K)$.
Let $D \cap R = \alpha$ and $D \cap T = \tau$.  Assume that $\alpha$ is not parallel on $R$ to a corner of $R$.
Then there exists an isotopy of $\Int S$ with support in $N(D)$ such that after the isotopy
\begin{itemize}
\item in $X^+$, $R$ is replaced by the result of boundary compression of $R$ along $D$,
\item in $X^-$, the regions of $S \cap X^-$ incident to $\tau$ are joined by surgery (restricted to $X^-$) along $\tau$, 
\item each arc of $S \cap T$ is essential in both $S$ and $T$ and
\item if $\tau$ connects distinct components of $R \cap T$ then such components are arcs and $|S \cap T|$ remains minimized.
\end{itemize}
The same holds with $X^+$ and $X^-$ interchanged.
\end{lemma}

\begin{remark}
This is the type $A$ isotopy described by Jaco \cite{jaco3mfldtop}.
\end{remark}

\begin{proof}
Consider the isotopy of $S$ with support in $N(D)$ that pushes the arc $\alpha \subseteq R$ through $D$.  The lemma follows directly.  

 In $X^+$, this isotopy is effectively a boundary compression of $R$ along $D$.  Such a compression could only create a monogon if $\alpha$ were parallel to a corner of $R$.  We however explicitly avoid this case.  In $X^-$, this isotopy is effectively surgery of $S$ along the arc $\tau$ restricted to $X^-$.  Such a surgery cannot create a monogon.  See \fullref{fig:isotopyalongdisk}.

\begin{figure}[ht!]
\centering
\begin{picture}(0,0)%
\includegraphics{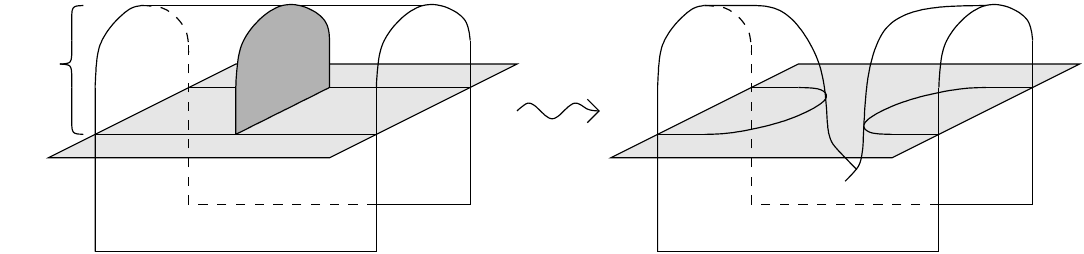}%
\end{picture}%
\setlength{\unitlength}{2960sp}%
\begingroup\makeatletter\ifx\SetFigFont\undefined%
\gdef\SetFigFont#1#2#3#4#5{%
  \reset@font\fontsize{#1}{#2pt}%
  \fontfamily{#3}\fontseries{#4}\fontshape{#5}%
  \selectfont}%
\fi\endgroup%
\begin{picture}(6920,1605)(2093,-5623)
\put(5251,-4336){\makebox(0,0)[lb]{\smash{{\SetFigFont{9}{10.8}{\rmdefault}{\mddefault}{\updefault}{\color[rgb]{0,0,0}$\hatT$}%
}}}}
\put(8851,-4336){\makebox(0,0)[lb]{\smash{{\SetFigFont{9}{10.8}{\rmdefault}{\mddefault}{\updefault}{\color[rgb]{0,0,0}$\hatT$}%
}}}}
\put(3901,-4486){\makebox(0,0)[b]{\smash{{\SetFigFont{9}{10.8}{\rmdefault}{\mddefault}{\updefault}{\color[rgb]{0,0,0}$D$}%
}}}}
\put(2401,-4486){\makebox(0,0)[rb]{\smash{{\SetFigFont{9}{10.8}{\rmdefault}{\mddefault}{\updefault}{\color[rgb]{0,0,0}$R$}%
}}}}
\end{picture}%
\caption{The isotopy of $R$ through $D$}
\label{fig:isotopyalongdisk}
\end{figure}
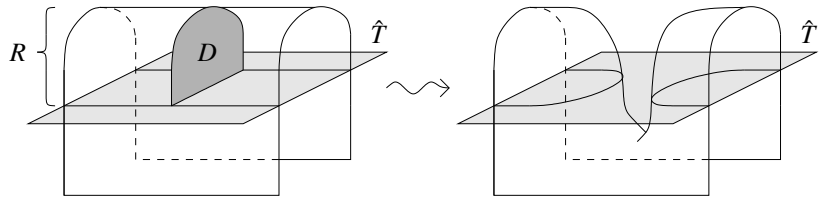

Indeed, restricting our view to just $S \cap T$ on $T$, this isotopy has the effect of surgering $R \cap T$ along $\tau$.  If $\tau$ has both end points on the same component of $\gamma \in R \cap T$, then there is an arc $\gamma_0 \subseteq \gamma$ connecting the end points of $\tau$.  The isotopy will then render $\gamma$ into the two components $\gamma_0 \cup \tau'$ and $(\gamma - \gamma_0) \cup \tau''$ where $\tau'$ and $\tau''$ are suitable pushoffs of $\tau$.

If $\tau$ has its endpoints on distinct arc components of $R \cap T$, then $|S \cap T|$ remains unchanged after the isotopy.  However, if a simple closed curve of $R \cap T$ contains just one end point of $\tau$, then such an isotopy would reduce $|S \cap T|$ contradicting the assumed minimality of $|S \cap T|$.  
\end{proof}

\begin{lemma}\label{lem:pushalongbigon}
Let $B$ and $R$ be regions of $S \cut T$ on the same side of $\hatT$ such that $B$ is a bigon and there is an edge of $B$ adjacent to an edge of $R$.
Then either
\begin{itemize}
\item  $R$ is a bigon parallel to $B$,
\item  $R$ is a trigon or
\item  there is an isotopy of $\Int S$ supported in $N(B \cup R)$ such that after the isotopy $R$ is replaced by a bigon $B'$ parallel to $B$ and a region $R'$ with two fewer corners than $R$ and two edges of regions (or of a single region) of $S \cut T$ on the opposite side of $\hatT$ are joined by surgery along an arc on $T$.
\end{itemize}
Furthermore in the last case the isotopy preserves the property that each arc of $S \cap T$ is essential in both $S$ and $T$.  If the two edges are distinct, then $|S \cap T|$ remains minimized over all isotopies of $\Int S$.  If the two edges are not distinct, then $R'$ is not a disk.
\end{lemma}
 
\begin{proof}
We will construct a disk $D$ that gives a boundary compression for $R$ in cut exterior $(X-N(K)) \cut T$ and then apply \fullref{lem:isotopacrossdisk}.  The desired disk $D$ is a pushoff of $B$ with two corners and an edge on $R$ and remaining edge on $T$.  

Let $e_B$ and $e_R$ be the adjacent edges of $B$ and $R$.  Together they bound a bigon $\delta$ of $T \cut S$.  Let $e_B'$ be the other edge of $B$.  Assume $e_B$ and $e_R$ have label pairs $\{x, y{+}1\}$ so that, by the parity rule, the corners of $B$ are on $\bdry H_{(x,\, x{+}1)}$ and $\bdry H_{(y,\,y+1)}$.  It may be the case that $x=y$.  Note that the corners of $R$ adjacent to $e_R$ bound rectangles $\rho_{(x,\, x{+}1)} \subseteq \bdry H_{(x,\, x{+}1)}$ and $\rho_{(y,\,y+1)} \subseteq \bdry H_{(y,\,y+1)}$ with the corners of $\delta$, the corners of $B$, and arcs on the vertices $U_{x+1}$ and $U_{y}$.  Since $\rho_{(x,\, x{+}1)}$ and $\rho_{(y,\,y+1)}$ have their interiors disjoint from $S$, we form a disk $D$ that is a slight pushoff of the disk $B \cup \delta \cup \rho_{(x,\, x{+}1)} \cup \rho_{(y,\,y+1)}$.  Notice that $\bdry D$ is composed of an arc $\alpha$ of $R$ and an arc $\tau$ of $T \cut S$.  The arc $\alpha$ is a slight pushoff of the edge $e_R$ and the two corners of $R$ to which it is incident.   

Assume $R$ is neither a bigon nor a trigon.  Then $\alpha$ is not parallel on $R$ either to a corner or into an edge of $R$.  Thus we may apply the isotopy of \fullref{lem:isotopacrossdisk}.  If the endpoints of $\tau$ lie on a single component of $R \cap T$ then the component of $\bdry R$ containing $e_R$ has just two edges.  Since $R$ is not a bigon, it cannot be a disk.  Thus after the isotopy $R'$ is not a disk.

If $R$ is a bigon, then \fullref{lem:parallelbigons} implies that $R$ is parallel to $B$.
\end{proof}

The above lemma applies directly if $B$ is a face of $G_S$ bounded by an $S2$ cycle.  If $B$ is a face of $G_S$ bounded by an $S3$ cycle, we may obtain a similar statement.

\begin{lemma}\label{lem:pushalongS3cycle}
Let $B$ be a face of $G_S$ bounded by an $S3$ cycle.  Let $R$ be a region of $S \cut T$ on the same side of $\hatT$ as $B$ such that there is an edge of $B$ adjacent to an edge of $R$.  Further assume this edge of $R$ is not between the two parallel edges of $B$.  Then either
\begin{itemize}
\item $R$ is a trigon parallel to $B$,
\item $R$ is a tetragon or
\item there is an isotopy of $\Int S$ supported in $N(B \cup R)$ such that after the isotopy $R$ is replaced by a trigon $B'$ parallel to $B$ and a region $R'$ with three fewer corners than $R$ and two edges of regions (or of a single region) of $S \cut T$ on the opposite side of $\hatT$ are joined by surgery along an arc on $\hatT$.
\end{itemize}
Furthermore in the last case the isotopy preserves the property that each arc of $S \cap T$ is essential in both $S$.  If the two edges are distinct, then $|S \cap T|$ remains minimized over all isotopies of $S$.  If the two edges are not distinct, then $R'$ is not a disk.
\end{lemma}

\begin{proof}
Again, we will construct a disk $D$ that gives a boundary compression for $R$ in $(X-N(K)) \cut T$ and then apply \fullref{lem:isotopacrossdisk}.  The desired disk $D$ is a pushoff of $B$ with three corners and two edges on $R$ and remaining edge on $T$.

Let $e_B$ and $e_R$ be the adjacent edges of $B$ and $R$.  Together they bound a bigon $\delta$ of $T \cut S$.  Assume $e_B$ and $e_R$ have label pairs $\{x, x{+}1\}$ so that the corners of $B$ are on $\bdry H_{(x,\, x{+}1)}$.  Note that the corners of $R$ adjacent to $e_R$ bound rectangles $\rho, \rho' \subseteq \bdry H_{(x,\, x{+}1)}$ with the corners of $\delta$, two corners of $B$, and arcs on the vertices $U_{x}$ and $U_{x+1}$.  

Two of the three edges of $B$ are parallel on $T$.  These two edges bound a bigon $\delta_B$ on $T$.  By assumption, $e_R \not \subseteq \delta_B$.  Nevertheless, following around the component of $\bdry R$ containing $e_R$, one of the edges before or after $e_R$ is contained in $\delta_b$ and is adjacent to an edge of $B$ other than $e_B$.  Let $e_B'$ and $e_R'$ be these edges of $B$ and $R$ respectively.  Let $\delta'$ be the bigon of $T \cut S$ bounded by $e_B'$ and $e_R'$.

 We may assume that the rectangle $\rho'$ is bounded by the corner of $B$ that connects the edges $e_B$ to $e_B'$, the corner of $R$ that connects the edges $e_R$ to $e_R'$, and a corner of each $\delta$ and $\delta'$.  The next corner of $R$ then bounds a rectangle $\rho'' \subseteq \bdry H_{(x,\, x{+}1)}$ with the next corner of $B$, a corner of $\delta'$ and an arc on one of the vertices $U_{x}$ or $U_{x+1}$.

Since $\rho$, $\rho'$, and $\rho''$ have their interiors disjoint from $S$, we form a disk $D$ that is a slight pushoff of the disk $B \cup \delta \cut \delta' \cup \rho \cup \rho' \cup \rho''$.  Notice that $\bdry D$ is composed of an arc $\alpha$ of $R$ and an arc $\tau$ of $T \cut S$.  The arc $\alpha$ is a slight pushoff of the edges $e_R$ and $e_R'$, the corner between them and the two corners surrounding them. 

The remainder of this proof follows completely analogously to the above proof of \fullref{lem:pushalongbigon}.
\smallskip

Assume $R$ is neither a trigon nor a tetragon.  Then $\alpha$ is not parallel on $R$ either to a corner or into an edge of $R$.  Thus we may apply the isotopy of \fullref{lem:isotopacrossdisk}.  If the endpoints of $\tau$ lie on a single component of $R \cap T$ then the component of $\bdry R$ containing $e_R$ has just three edges.  Since $R$ is not a trigon, it cannot be a disk.  Thus after the isotopy $R'$ is not a disk.

If $R$ is a trigon, then a proof analogous to that of \fullref{lem:parallelbigons} implies that $R$ is parallel to $B$.
\end{proof}


\subsubsection{Arranging $S$ in the solid torus $V$}

Before performing any isotopies of $S$ we need to see how it may presently be positioned in $V$.  

\begin{lemma} \label{lem:bigonsdisksandannuliinsolidtorus}
Each region $R \in S \cap V$ is either a bigon parallel to a bigon of $F \cut f$, a meridional disk of $V$ with an odd number of at least $3$ corners, or an annulus with two corners.
\end{lemma}

\begin{proof}
Choose an orientation on $S$.  Such a choice induces an orientation on each $R \in S \cap V$ which in turn induces an orientation on each boundary component of $R$.  Let us then consider the collection of oriented simple closed curves $\mathcal{C} = \{ \bdry R | R \in S \cap V\}$ as they lie on the torus $\bdry V$.  The collection $\mathcal{C}$ may contain both trivial and essential curves on $\bdry V$. 

Let $\gamma$ be an (arbitrarily oriented) essential circle on $\bdry V$ disjoint from $\hatT$ intersecting each corner of each $R \in S \cap V$ (ie\ each arc of $C \cap \bdry \bar{N}(K)$ for each $C \in \mathcal{C}$) exactly once.  Since a corner of an $R \in S \cap V$ is an arc on either $\bdry H_{(1, 2)}$ or $\bdry H_{(5, 6)}$, $\gamma$ may be divided into two arcs, say $\gamma_+$ and $\gamma_-$, according to the direction in which a curve $C \in \mathcal{C}$ crosses the arc.

If $C \in \mathcal{C}$ is a trivial curve on $\bdry V$, then $C$ must be the boundary of a disk of $S \cap V$.  If $C \cap \gamma = \emptyset$, then $C$ is a trivial curve on $T$.  Since $C$ bounds a disk in $S$ and a disk in $T$, these disks together form a sphere in $V$ which bounds a solid ball in $V$.  This implies that there exists an isotopy of $S$ that reduces $|S \cap T|$ contradicting the assumption that $|S \cap T|$ is minimized.  Thus $C \cap \gamma \neq \emptyset$.  

Since $C$ bounds a disk on $\bdry V$, $\gamma$ must alternately cross in and out of the disk that $C$ bounds.  Therefore the direction in which $C$ crosses $\gamma$ alternates around $\gamma$.  Because $\gamma$ is divided into the two arcs $\gamma_+$ and $\gamma_-$ that dictate the direction in which $C$ may cross, $C$ must intersect $\gamma$ only twice.  Hence $C$ is the boundary of a bigon.  Such a bigon has one corner on $K_{(1, 2)}$ and one corner on $K_{(5, 6)}$.  Since $C$ is a trivial curve on $\bdry V$, the bigon it bounds must be parallel to one of the bigons of $F \cut f$.

If $C \in \mathcal{C}$ is a meridional curve, then it intersects $\gamma$ algebraically once.  Hence $|C \cap \gamma|$ is odd.  Furthermore, $C$ must bound a disk region in $S$ since otherwise there would be a compression of $S$.  Since there are no monogons of $S \cut T$, we know that $|C \cap \gamma| > 1$.   Thus any meridional disk of $S \cap V$ has an odd number of at least three corners.

Assume $C \in \mathcal{C}$ is an essential nonmeridional curve on $\bdry V$.  Since $C \subseteq \bdry R$ for some $R \in S \cap V$, then $\bdry R \subseteq \mathcal{C}$ is a collection of essential nonmeridional curves on $\bdry V$ that is nullhomologous in $V$.  Since $R$ is incompressible, $R$ must be an annulus.  Let $C'$ be the other component of $\bdry R$.

Let $A$ be an annulus on $\bdry V$ between $C$ and $C'$ oriented so that the orientations of $C$ and $C'$ respect the induced boundary orientation.  In order for the curves $C$ and $C'$ to intersect each $\gamma_+$ and $\gamma_-$ in the prescribed directions either $\gamma$ crosses $A$ transversely once, $\gamma$ intersects $A$ once but is disjoint from one component of $\bdry A$, or $\bdry A$ is disjoint from $\gamma$.  Thus the annulus $R$ has respectively one corner and one edge on each boundary component, one boundary component in $\hatT$ and one with two corners and two edges, or both boundary components in $\hatT$.  However, if both boundary components of $R$ are in $\hatT$, then $R$ is isotopic onto $A$ within $V$ (since the annulus $\bdry V \cut T_V$ is parallel to $A$ through $V$).  Thus there is an isotopy of $S$ that reduces $|S \cap T|$ contradicting the assumption that $|S \cap T|$ is minimized.  Therefore an annulus of $S \cap V$ must have exactly two corners.
\end{proof}

\begin{lemma}\label{lem:meridionaldisks}
Assume $P \in S \cap V$ is a trigon and $R \in S \cap V$ is a meridional disk of $V$ distinct from $P$.  Then $R$ has an edge parallel on $\hatT$ to an edge of a bigon of $S \cap V$ such that no edge of $P$ lies between these edges.  Furthermore, if $R$ is also a trigon, it is parallel to $P$.
\end{lemma}

\begin{proof}
We continue with the setup of the proof of \fullref{lem:bigonsdisksandannuliinsolidtorus}.

As $\bdry P \subseteq \bdry V$ must cross $\gamma$ geometrically three times and algebraically once, cut torus $\bdry V \cut (\bdry P \cup \gamma)$ consists of three components.  Two of the components are disks, say $D_1$ and $D_2$, each with boundary composed of an arc of $\bdry P$ and an arc of $\gamma$.  The third component is  a disk with boundary composed alternately of four arcs of $\bdry P$ and four arcs of $\gamma$.  Note that each disk $D_1$ and $D_2$ corresponds to a parallelism on $T$ between an edge of $P$ and an edge of a bigon of $S \cap V$.  

The two disks $D_1$ and $D_2$ meet at a point $p$ of $\bdry P \cap \gamma$ where $\bdry P$ crosses $\gamma$ in a direction opposite that of the other two points.  Without loss of generality, we may assume $p \in \gamma_+$ and $(\bdry P \cap \gamma)-p \in \gamma_-$.  Thus $\gamma_+ \subseteq (D_1 \cup D_2) \cap \gamma$.  
Since $\bdry R$ must intersect $\gamma_+$, $\bdry R \cap D_i \neq \emptyset$ for either $i=1$ or $2$.  An arc of $\bdry R \cap D_i$ then cuts off a disk disjoint from $\bdry P$ that corresponds to a parallelism between an edge of $R$ and an edge of a bigon of $S \cap V$.  This parallelism does not contain an edge of $\bdry P$.

Furthermore, if $R$ is a trigon then $|\bdry R \cap \gamma|=3$.  The end points of an arc of $\bdry R \cap D_i$ for $i=1$ or $2$ account for two of these intersections.  In the complement of $\bdry P$, there is only one way to complete this arc into an essential simple closed curve that crosses $\gamma$ just once more.  The rectangle on $D_i$ between $\bdry P \cap D_i$ and $\bdry R \cap D_i$ extends to an annulus between $\bdry P$ and $\bdry R$.  This annulus is cut into three rectangles by $\gamma$.  These rectangles imply each edge of $R$ is parallel to an edge of $P$ and hence $\bdry R$ is parallel to $\bdry P$ on $\bdry (X^+ - N(K))$.   Thus by \fullref{lem:paralleldisks} $R$ is parallel to $P$.
\end{proof}

\begin{prop}\label{prop:groomingV}
Let $B_1$ and $B_2$ be the two bigons of $F \cut f$ on $\bdry V$.  Then one may isotop $\Int S$ so that
\begin{itemize}
\item $S \cap V$ is a collection of bigons parallel to either $B_1$ or $B_2$ and a collection of mutually parallel trigons, 
\item each arc of $S \cap T$ is essential in both $S$ and $T$ and
\item $|S \cap T|$ remains minimized.
\end{itemize}
\end{prop}

\begin{proof}
By \fullref{lem:bigonsdisksandannuliinsolidtorus}, a region $R \in S \cap V$ is either a bigon parallel to either $B_1$ or $B_2$, a meridional disk of $V$ with an odd number of edges, or an annulus with two corners.

We first note that $S \cap V$ cannot simultaneously contain meridional disks of $V$ and annuli.  If an annulus $Q$ of $S \cap V$ were to exist with a meridional disk of $S \cap V$, then the curves $\bdry Q$ and its core curve must all be meridional curves.  Thus the core curve of $Q$ must bound a disk in $V$.  This disk implies $S$ is compressible contradicting that it is incompressible.  Hence we consider annuli and meridional disks of $S \cap V$ separately.

{\bf Case 1}\qua  $S \cap V$ contains annuli.

Let $\mathcal{Q}$ be the collection of annuli of $S \cap V$.  By the proof of \fullref{lem:bigonsdisksandannuliinsolidtorus}, all the annuli in $\mathcal{Q}$ either (a) have one boundary component contained in $T_V$ or (b) have each boundary component intersecting the annulus $T_V$ in a single spanning arc.  We will refer to such annuli as type (a) or type (b) accordingly.

Since properly embedded annuli in a solid torus are parallel into the boundary torus, each $Q \in \mathcal{Q}$ is isotopic to an embedded annulus $Q'$ in $\bdry V$ bounded by $\bdry Q$.  Let $V_Q \subseteq V$ be the solid torus through which $Q$ is parallel to $Q'$.  If $Q_1 \in \mathcal{Q}$ is contained in $V_Q$, then $Q_1$ is parallel to an annulus $Q_1' \subseteq Q'$.  If $Q_2 \in \mathcal{Q}$ is not contained in $V_Q$, then it is parallel to an annulus $Q_2' \subseteq \bdry V$ such that either $Q_2' \cap Q' = \emptyset$ or $Q' \subseteq Q_2'$.  Let $\mathcal{Q}'$ be a collection of annuli in $\bdry V$ to which the annuli of $\mathcal{Q}$ are parallel.  We may assume $\mathcal{Q}'$ has been chosen so that it is partially ordered by nesting.

Choose a meridional curve $m \subseteq \bdry V$ for $V$ such that for each $Q' \in \mathcal{Q}'$ $m \cap Q' \subseteq \hatT$ is a collection of spanning arcs of $Q'$.  Thus there exists a meridional disk $D$ with $\bdry D = m$ such that $D$ intersects each annulus $Q \in \mathcal{Q}$ in transverse arcs that are isotopic through a subdisk of $D$ onto $\hatT$.  Though there may be several, let us  choose such a subdisk $D_Q$ of $D$ for each $Q \in \mathcal{Q}$ (and corresponding $Q' \in \mathcal{Q}'$).  Since the $Q' \in \mathcal{Q}'$ are nested on $\bdry V$, the disks $D_Q$ are nested on $D$.  Consider a disk $D_Q$ outermost on $D$ and its corresponding annulus $Q$.  Note that $\Int D_Q \cap S = \emptyset$ and $\bdry D_Q \cap S \subseteq Q$ is not parallel to any corner of $Q$.  Thus there is an isotopy of $S$ through $D_Q$ with support in $N(D_Q)$ as in \fullref{lem:isotopacrossdisk}.

Assume $Q \in \mathcal{Q}$ is an annulus of type (a).  Since one component of $\bdry Q$ before the isotopy is wholly contained in $\hatT$, after the isotopy the components of $\bdry Q \cap \hatT$ that contained $\bdry \delta$ are joined.  Thus $|S \cap T|$ is reduced contradicting the original minimality assumption on $|S \cap T|$.  Hence $\mathcal{Q}$ contains no annuli of type (a).  

Assume $Q \in \mathcal{Q}$ is an annulus of type (b).  Since $\bdry Q \cap \hatT$ is composed of two arcs before the isotopy, the isotopy exchanges them for two new arcs.  Moreover, viewing the isotopy of $Q$ as a boundary compression along $D_Q$, it follows that $Q$ becomes a bigon parallel to a bigon of $F \cut f$.  A similar isotopy may be then performed for each remaining annulus of $\mathcal{Q}$.

{\bf Case 2}\qua $S \cap V$ contains meridional disks.

A meridional disk $R \in S \cap V$ must have an edge parallel to an edge of $B_1$ or $B_2$.  If not, each edge of $\bdry R$ would necessarily have the same label pair and must lie in the annulus $T_V$.  Thus if $R$ is an $n$--gon, $\bdry R$ would intersect a longitude of $V$ minimally $n$ times on $\bdry V$.  Since $R$ cannot be a monogon, this contradicts that $R$ is a meridional disk.

Recall that by \fullref{lem:bigonsdisksandannuliinsolidtorus} each bigon of $S \cap V$ is parallel to either $B_1$ or $B_2$.  Let $R$ be a nonbigon face of $S \cap V$ such that $R$ has an edge $e_R$ adjacent to an edge of a bigon $B$ of $S \cap V$.  If $R$ is not a trigon, then by \fullref{lem:pushalongbigon} there exists an isotopy of $\Int S$ in $N(B \cup R)$ yielding in $V$ a bigon parallel to either $B_1$ or $B_2$ and a meridional disk of $V$ with two fewer corners than $R$.  Perform such isotopies until each nonbigon face of $S \cap V$ with an edge adjacent to an edge of a bigon of $S \cap V$ is a trigon.  

Now after all the isotopies thus far, assume there exists a nonbigon, nontrigon face $R \in S \cap V$.  Then since by assumption $S \cap V$ contains no annuli, \fullref{lem:bigonsdisksandannuliinsolidtorus} implies that $R$ is a meridional disk of $V$ with at least $5$ corners.  By \fullref{lem:meridionaldisks} the trigons of $S \cap V$ are mutually parallel.  The same lemma then also implies that $R$ has an edge $e_R$ parallel to an edge $e_B$ of $B_1$ or $B_2$ such that no edge of a trigon of $S \cap V$ lies between $e_R$ and $e_B$.  But then in between $e_R$ and $e_B$ there must be an edge of a nonbigon, nontrigon face of $S \cap V$ that is adjacent to an edge of a bigon of $S \cap V$.  This contradict that we have isotoped $S$ so that each nonbigon face of $S \cap V$ with an edge adjacent to an edge of a bigon of $S \cap V$ is a trigon. 
\end{proof}


\subsubsection{Arranging $S$ near $f$ in $X^-$}

Perform the isotopy of \fullref{prop:groomingV}.  By relabeling if necessary, we may assume that the edges of the trigons of $S \cap V$ that are not parallel to edges of bigons of $S \cap V$ connect the vertices $U_5$ and $U_6$.  Thus because the edges of $g$ (which is an $S2$ cycle or an $S3$ cycle with label pair $\{3, 4\}$) obstruct them, there are no edges of $G_T$ connecting $U_1$ to $U_2$ in $\hatT$ .  Also, since the trigons of $S \cap V$ are all mutually parallel, there are at most three parallelism edge classes of $G_T$ in $T_V$ incident to $U_6$.

\begin{prop}\label{prop:groomingf}
Fixing $S \cap V$, one may isotop $\Int S$ so that
\begin{itemize}
\item if $R \in S \cap X^-$ has an edge with label pair $\{1, 6\}$ that does not lie between two parallel edges of $f$, then $R$ is parallel to $f$,
\item if $R \in S \cap X^-$ is not parallel to $f$ but has an edge with label pair $\{1, 6\}$, then the corners of $R$ incident to this edge are also incident to edges with label pairs $\{1, 4\}$ and $\{5, 6\}$,
\item there are no edges of $G_T$ with label pairs $\{1, 2\}$ or $\{3, 6\}$,
\item each arc of $S \cap T$ is essential in both $S$ and $T$ and
\item either $|S \cap T|$ remains minimized or $S \cap T$ contains essential simple closed curves in the annulus on $T$ between edges with label pairs $\{1, 6\}$ and $\{3, 4\}$.
\end{itemize}
\end{prop}

\begin{proof}
Recall that $f$ is either a bigon or a trigon bounded by an $S2$ cycle or an $S3$ cycle respectively.  Note that $f$ has two parallel edges only if $f$ is a trigon.

Let $\mathcal{E}$ be the collection of edges of $G_T$ with label pair $\{1, 6\}$ that lie in $T_V$, do not lie between two parallel edges of $f$, and are not an edge of a face of $S \cap X^-$ parallel to $f$.  If $\mathcal{E} \neq \emptyset$ then let $R$ be a face of $S \cap X^-$ that has an edge $e_R \in \mathcal{E}$ adjacent to an edge of a face of $S \cap X^-$ that is parallel to $f$.

Since $f$ is a $p$--gon bounded by an $Sp$ cycle for $p = 2$ or $3$, we claim that $R$ cannot be a $(p{+}1)$--gon.  To the contrary, assume $R$ is a $(p{+}1)$--gon.  The corners of $f$ divide $\bdry H_{(6, 1)}$ into $p$ rectangles.  These rectangles are joined cyclically by edges of $G_T$ with label pair $\{1, 6\}$.  Since $R$ is a $(p{+}1)$--gon, every edge of $R$ cannot have label pair $\{1, 6\}$.  Therefore $R$ must have an edge with label pair $\{1, 4\}$, $\{5, 6\}$ or $\{3, 6\}$.  Thus it must have a second edge with one of these three label pairs.  Moving in both directions along $\bdry R$ from $e_R$, counting a second edge if $f$ is a trigon, the final two edges of $R$ must emanate from $U_1$ and $U_6$ and lie in $\hatT - T_V$.  Thus the final two edges of $R$ must have label pairs $\{1, 4\}$ and $\{3, 6\}$.  Since $H_{(3, 4)} \subseteq X^+$, $R$ cannot have a corner connecting these edges.  Thus $R$ is not a $(p{+}1)$--gon.

Since $R$ is not a $(p{+}1)$--gon, by either \fullref{lem:pushalongbigon} or \fullref{lem:pushalongS3cycle}, one may isotop $R$ to produce a new face $B'$ of $S \cap X^-$ which has $e_R$ as an edge that is parallel in $X-N(K)$ to $f$ and a new face $R'$ of $S \cap X^-$ with two fewer corners than $R$ while joining faces of $S \cap X^+$ along an arc $\tau \subseteq \hatT \cut S$.  By construction, $\tau$ connects an edge of $R$ incident to $U_1$ not in $T_V$ to an edge of $R$ incident to $U_6$ not in $T_V$.  Therefore any arc of $S \cap T$ created by this isotopy is not contained in $T_V$.  Recall that because edges of $F$ and $f$ delineate the boundary of $T_V$ there may be edges of $G_T$ parallel on $T$ to edges of $f$ which are not contained in $T_V$.  Furthermore, since $e_R$ is now an edge of the face $B'$ which is parallel to $f$, it is no longer included in $\mathcal{E}$.  Thus the isotopy strictly decreases $|\mathcal{E}|$.  Since such an isotopy creates no monogons, each arc of $S \cap T$ is essential in each $S$ and $T$. 

The isotopy may only increase $|S \cap T|$ if the arc $\tau$ connects an edge of $R$ to itself.  In this case neither $R$ nor $R'$ is a disk.  The region $R'$ has a boundary component that lies as an essential simple closed curve in the annulus of $\hatT$ between the edges of $G_T$ with label pairs $\{1, 6\}$ and those with label pairs $\{3, 4\}$.

Repeatedly perform such isotopies until $\mathcal{E} = \emptyset$.  If for some $R \in S \cap X^-$ not parallel to $f$ has an edge $e_R$ with label pair $\{1, 6\}$, then $e_R$ must lie between two parallel edges of $f$.  Otherwise, since $e_R$ is not contained in $T_V$, moving in both directions along $\bdry R$ away from $e_R$, picking up a second edge that lies between two parallel edges of $f$ if $f$ is a trigon, we arrive at two edges of $R$ that lie in $T_V$.  One edge has an end point on $U_6$ and the other has an end point on $U_1$.  Since the edges of $G_T$ in $T_V$ must have label pair $\{1, 6\}$ or $\{5, 6\}$,  the latter edge must be in $\mathcal{E}$.  This contradicts that $\mathcal{E} = \emptyset$.

Assume $e_R$ is an edge of a face $R \in S \cap X^-$ not parallel to $f$ such that $e_R$ lies between two parallel edges of $f$.  Then the two corners of $R$ incident to $e_R$ are incident to edges each with an endpoint incident to either $U_1$ or $U_6$.  Since these two edges cannot lie between two parallel edges of $f$, they cannot have label pair $\{1, 6\}$.  Thus they have the label pairs $\{1, 4\}$ and $\{1, 6\}$.

It similarly follows that there are no edges with label pair $\{3, 6\}$ due to the fact that there are no edges with label pair $\{1, 2\}$.
\end{proof}


\subsubsection{Arranging $S$ in $X^+ \cut V$}

Recall that $g$ is a $p$--gon of $S \cap X^+$ bounded by an $Sp$ cycle for $p = 2$ or $3$.  Also recall that $W$ is the solid torus $X^+ \cut (V \cup \bar{N}(K) \cup g)$.

\begin{prop}\label{prop:groomingW} 
One may isotop $S$ so that each face of $S \cap W$ is either
\begin{itemize}
\item a bigon parallel to a bigon of $F \cut f$,
\item a $p$--gon parallel to the $p$--gon $g$,
\item a meridional disk of $W$ with $2$ or $3$ corners or 
\item an annulus with a boundary component in each component $T_W$ and $T_W'$ of $T \cap W$.
\end{itemize}
Furthermore, this isotopy creates no monogons and preserves the results of the isotopies of \fullref{prop:groomingV} and \fullref{prop:groomingf}.
\vspace{-10pt}
\end{prop}

\begin{proof}
   Let $\mathcal{E}$ be the collection of edges of $G_T$ with label pair either $\{3, 4\}$ or $\{1, 6\}$ that lie in $T_W$, do not lie between two parallel edges of $g$, and are not an edge of a face of $S \cap W$ parallel to either $g$ or a bigon of $F \cut f$.  Let $\mathcal{E}'$ be the collection of edges of $G_T$ with label pair either $\{3, 4\}$ or $\{2, 5\}$ that lie in $T_W'$, do not lie between two parallel edges of $g$, and are not an edge of a face of $S \cap W$ parallel to either $g$ or a bigon of $F \cut f$.

If $\mathcal{E} \neq \emptyset$ then let $R$ be a face of $S \cap W$ that has an edge $e_R \in \mathcal{E}$ adjacent to an edge of a face $S \cap W$ that is parallel to either $g$ or a bigon of $F \cut f$.

 If $e_R$ has label pair $\{1, 6\}$ then $R$ cannot be a trigon.  If $g$ is a $p$--gon and $e_R$ has label pair $\{3, 4\}$ then $R$ cannot be a $(p{+}1)$--gon.  This may be seen by attempting to connect a set of $3$ or $4$ corners in $\bdry W$ with the edges of $R$.

Therefore, by either \fullref{lem:pushalongbigon} or \fullref{lem:pushalongS3cycle}, one may isotop $R$ to produce a new face $B'$ of $S \cap W$ which has $e_R$ as an edge that is parallel to either $g$ or a bigon of $F \cut f$ and a new face $R'$ of $S \cap W$ with two fewer corners than $R$ while joining faces of $S \cap X^-$ along an arc $\tau \subseteq T_W'$.  Since $e_R$ is now an edge of the face $B'$ which is parallel to either $g$ or a bigon of $F$, it is no longer included in $\mathcal{E}$.  Thus the isotopy strictly decreases $|\mathcal{E}|$.  Such an isotopy creates no monogons, and hence each arc of $S \cap T$ is essential in each $S$ and $T$.  Repeatedly perform such isotopies until $\mathcal{E} = \emptyset$.  Because the only arcs of $S \cap T$ affected are ones in $T_W'$, the results of these isotopies do not have any effect upon those of \fullref{prop:groomingV} and \fullref{prop:groomingf}.

By \fullref{prop:groomingf} there are no edges with label pair $\{3, 6\}$.  Thus $\mathcal{E} = \emptyset$ implies $\mathcal{E}' = \emptyset$.  If otherwise, then there exists an edge $e_R \in \mathcal{E}'$ of a region $R \in S \cap W$ that is not parallel to either $g$ or a bigon of $F \cut f$.  If $e_R$ has label pair $\{2, 5\}$, then it is incident to a corner of $R$ on $\bdry H_{(5, 6)}$.  Therefore the next edge of $R$ must lie in $T_W$ and have an end point labeled $6$.  Since there are no edges of $G_T$ with label pair $\{3, 6\}$, this edge must have label pair $\{1, 6\}$ and is thus an edge of $\mathcal{E}$ contradicting that $\mathcal{E} = \emptyset$.  If $e_R$ has label pair $\{3, 4\}$ then after the corner of $R$ incident to the end point of $e_R$ labeled $4$ the next edge must have label pair $\{3, 4\}$.  If this edge is between two parallel edges of $g$, then the next edge too must have label pair $\{3, 4\}$.  Regardless, one of these edges lies in $T_W$ and is thus an edge of $\mathcal{E}$ contradicting that $\mathcal{E} = \emptyset$.

Similarly, it follows that there can be no edges with label pair $\{4, 5\}$.  An edge with label pair $\{4, 5\}$ must be an edge of a face $S \cap W$ with a corner on $\bdry H_{(5, 6)}$ followed by an edge incident to $U_6$.  This edge must then have label pair $\{1, 6\}$ and thus be an edge of $\mathcal{E}$ contradicting that $\mathcal{E} = \emptyset$.

Assume a face $R$ of $S \cap W$ has an edge $e_R$ with label pair $\{1, 4\}$.  Since there are no edges with label pair $\{3, 6\}$, no edge of $S \cap W$ can transversely cross $T_W$ in the direction opposite from $e_R$.  Thus $R$ must be a meridional disk with one edge in $T_W$ and one edge in $T_W'$.  If $g$ is a trigon, $R$ may have a third edge between two parallel edges of $g$.  Thus $R$ is a meridional disk of $W$ with either two or three edges and two or three corners respectively.  Note that because there may be two distinct edges classes of $G_T$ with label pair $\{1, 4\}$ there may be two parallelism classes of disks of $S \cap W$ that are meridional disks of $W$.  Nevertheless, all such disks have the same number of corners.  Furthermore, the number of edges with label pair $\{1, 4\}$ equals the number of edges with label pair $\{2, 5\}$.

If $S \cap T_W$ or $S \cap T_W'$ contains a simple closed curve component $c$, then there can be no edges with label pair $\{1, 4\}$ or $\{2, 5\}$.  Thus an edge of $G_T$ in $T_W$ or $T_W'$ must be an edge of a face of $S \cap W$ parallel to either $g$ or a bigon of $F \cut f$.  Let $R \in S \cap W$ be the face with $c \in \bdry R$.  Since $c$ is a longitudinal curve of $W$, $R$ must be an annulus with its other boundary component contained in whichever of $T_W$ or $T_W'$ that does not contain $c$.  Otherwise, if these two boundary components are contained in the same annulus, then $R$ is parallel into $T_W$ or $T_W'$.  An isotopy of $S$ will push $R$ (and any other annuli of $S \cap W$ between $R$ and $\hatT$) into $X^-$.  This will violate the minimality of $|S \cap T|$.  Hence we may assume $S$ has been isotoped so that any annulus of $S \cap W$ has one boundary component in $T_W$ and its other one in $T_W'$.
\end{proof} 

\subsubsection[Proof of \ref{prop:isotopS}]{Proof of
\fullref{prop:isotopS}}\label{sec:proofofpropisotopS}
\vspace{-10pt}

\begin{proof}  Perform the isotopies of $\Int S$ first in \fullref{prop:groomingV}, then in \fullref{prop:groomingf}, and finally in \fullref{prop:groomingW}.  This arranges the edges of $G_T$ into the edge classes described in the proposition.  The isotopies also maintain that the arcs of $S \cap T$ are essential in both $S$ and $T$.
\vspace{-5pt}

Recall that on $G_T$ there are $s$ edges around each vertex. 
Since there are $n$ edges with each label pair $\{1, 4\}$, $\{2, 3\}$, and $\{5, 6\}$ and no edges with label pair $\{1, 2\}$, $\{3, 6\}$, or $\{4, 5\}$, there must be a total of $s-n$ edges with each label pair $\{1, 6\}$, $\{2, 5\}$, and $\{3, 4\}$. Each collection of $s-n$ edges with label pair $\{1, 6\}$, $\{2, 5\}$, and $\{3, 4\}$ lie in an essential annulus.
\vspace{-5pt}

{\bf The faces of $S \cap X^+$}\qua
The proof of Case 2 of \fullref{prop:groomingV} implies each edge of $G_T$ with label pair $\{5, 6\}$ is an edge of a trigon in $V$ whose other two edges have label pairs $\{1, 6\}$ and $\{2, 5\}$.  Since there are $n$ edges with label pair $\{5, 6\}$, there must be $n$ trigons contained in $V$.  

Together the proofs of \fullref{prop:groomingV} and \fullref{prop:groomingW} imply that an edge with label pair $\{1, 6\}$ or $\{2, 5\}$ that is not an edge of a trigon in $V$ is an edge of a bigon that is parallel to a bigon of $F \cut f$.  Since there are $s-n$ edges with label pair $\{1, 6\}$, $n$ of which are edges of trigons in $V$, there are $s-2n$ edges with label pair $\{1, 6\}$ that are edges of bigons parallel to bigons of $F \cut f$.  Since a bigon of $F \cut f$ has only one edge with label pair $\{1, 6\}$, each edge with label pair $\{1, 6\}$ that is not an edge of a trigon is an edge of a distinct bigon.  Thus there are $s-2n$ bigons parallel to a bigon of $F \cut f$.  Furthermore, of the $s-n$ edges with label pair $\{2, 5\}$, $n$ belong to a trigon in $V$ and $s-2n$ belong to a bigon parallel to a bigon of $F \cut f$.
\vspace{-2pt}

By \fullref{prop:groomingW} an edge with label pair $\{2, 3\}$ is an edge of a meridional disk of $W$ with $2$ or $3$ corners and hence $2$ or $3$ edges respectively.  One corner runs along $\bdry H_{(1, 2)}$ from the edge with label pair $\{2, 3\}$ to an edge with label pair $\{1, 4\}$.  If $g$ is a bigon, then the meridional disk of $W$ is a bigon itself with its final corner on $\bdry H_{(3, 4)}$.  If $g$ is a trigon, then two edges of $g$ are parallel.  Also, the corners of $g$ divide $\bdry H_{(3, 4)}$ into three rectangles $\rho$, $\rho'$, and $\rho''$.  The bigon $\delta$ on $\hatT$ between the two parallel edges of $g$ connects two of these rectangles, say $\rho$ and $\rho'$.  Either the edges with label pair $\{2, 3\}$ and $\{1, 4\}$ are incident to $\rho \cup \delta \cup \rho'$ or they are incident to $\rho''$.  If they are incident to $\rho \cup \delta \cup \rho'$, then a meridional disk of $W$ is a trigon with a corner on each of $\rho$ and $\rho'$ and an edge that lies in $\delta$.  If the edges are incident to $\rho''$, then a meridional disk of $W$ is a bigon with its final corner on $\rho''$.  
\vspace{-2pt}

Since there are just $n$ edges with label pair $\{2, 3\}$, there are $n$ such meridional disks of $W$.  They are either all bigons or all trigons.  In either case, each has one edge with label pair $\{2, 3\}$ and one edge with label pair $\{1, 4\}$.  If the disks are trigons, then their third edge has label pair $\{3, 4\}$.  If $g$ is a bigon, then each edge with label pair $\{3, 4\}$ is an edge of a bigon parallel to $g$.  Since there are $s-n$ edges with label pair $\{3, 4\}$, there are $(s-n)/2$ bigons parallel to $g$.  If $g$ is a trigon, then each edge with label pair $\{3, 4\}$ that is not an edge of a meridional disk of $W$ is an edge of a trigon parallel to $g$.  If the meridional disks of $W$ are bigons, then they have no edge with label pair $\{3, 4\}$.  Hence there are $(s-n)/3$ trigons parallel to $g$.  If the meridional disks of $W$ are trigons, then they each have one edge with label pair $\{3, 4\}$.  Hence there are $(s-2n)/3$ trigons parallel to $g$.
\vspace{-2pt}

If $S \cap X^+$ contains no meridional disks of $W$, then $n=0$.  Similarly, if there are annuli of $S \cap W$ then $n=0$.  
\vspace{-2pt}

Let us now summarize the collection $S \cap X^+$.
If $g$ is a bigon, then $S \cap X^+$ is composed of $n$ trigons contained in $V$, $s-2n$ bigons parallel to bigons of $F \cut f$, $n$ bigons that are meridional disks of $W$, $(s-n)/2$ bigons parallel to $g$, and annuli if $n = 0$.  Thus if $g$ is a bigon, then $S \cap X^+$ contains $(s-2n) + n + (s-n)/2 = 3(s-n)/2$ bigons, $n$ trigons, and possibly annuli if $n=0$.

If $g$ is a trigon, then $S \cap X^+$ is composed of $n$ trigons contained in $V$, $s-2n$ bigons parallel to bigons of $F \cut f$, either $n$ bigons that are meridional disks of $W$ and $(s-n)/3$ trigons parallel to $g$ or $n$ trigons that are meridional disks of $W$ and $(s-2n)/3$ trigons parallel to $g$, and annuli if $n=0$.  Thus if $g$ is a trigon, then $S \cap X^+$ contains either $(s-2n) + n= s-n$ bigons and $n + (s-n)/3 = (s+2n)/3$ trigons or $s-2n$ bigons and $n + n + (s-2n)/3 = (s-4n)/3$ trigons and possibly annuli if $n=0$.
\vspace{-2pt}

{\bf The faces of $S \cap X^-$}\qua
By \fullref{prop:groomingf}, an edge with label pair $\{1, 6\}$ that does not lie between two parallel edges of $f$ is an edge of a face parallel to $f$.  
\vspace{-2pt}

If $f$ is a bigon, then no edges of $f$ are parallel.  Hence all $s-n$ edges with label pair $\{1, 6\}$ are edges of bigons parallel to $f$.  Thus there are $(s-n)/2$ bigons parallel to $f$.  
\vspace{-2pt}

If $f$ is a trigon, then each of the $n$ edges with label pair $\{5, 6\}$ are incident to a corner on $\bdry H_{(6, 1)}$ with are in turn either incident to an edge with label pair $\{1, 4\}$ or incident to an edge with label pair $\{1, 6\}$ that lies between the parallel edges of $f$.  In the former case the $s-n$ edges with label pair $\{1, 6\}$ are all edges of trigons parallel to $f$;  hence there are $(s-n)/3$ trigons of $S \cap X^-$ parallel to $f$.  In the latter case the remaining $s-2n$ edges with label pair $\{1, 6\}$ are edges of trigons parallel to $f$; hence there are $(s-2n)/3$ trigons of $S \cap X^-$ parallel to $f$.
\end{proof}
\vspace{-2pt}


\subsection{Euler characteristic estimates}

The arcs and simple closed curves of $S \cap T$ break $S$ into faces which lie either in $X^+$ or $X^-$.  The arcs of $S \cap T$ form the edges of $G_S$.  Because we are assuming $t = 6$, $G_S$ has a total of $3s$ edges.

\begin{lemma}\label{lem:upperbound} Both $\displaystyle\sum_{R \in S \cap X^+} \chi (R) \leq 3s/2\ \ $ and $\displaystyle\sum_{R \in S \cap X^-} \chi (R) \leq 3s/2$.
\end{lemma}

\begin{proof} 
Since $\displaystyle \sum_{R \in S \cap X^+} \chi (R) \leq \#(\mbox{disks in } S \cap X^+)$
\begin{align*} 
 2 \cdot \#(\mbox{disks in } S \cap X^+) &\leq \#(\mbox{edges of disks in } S \cap X^+) \tag*{\hbox{and}}\\ 
&\leq \#(\mbox{edges of faces in } S \cap X^+) = 3s,\\
\sum_{R \in S \cap X^+} \chi (R) &\leq 3s/2. \tag*{\hbox{we have}}\end{align*}
Replacing $X^+$ with $X^-$ achieves the other case.
\end{proof}

Note that equality is realized in \fullref{lem:upperbound} if $S \cap X^+$ (or $S \cap X^-$) is a collection of bigons.  Utilizing \fullref{prop:isotopS} we can obtain an exact count for $\displaystyle\sum_{R \in S \cap X^+} \chi (R)$.

\begin{lemma}\label{lem:betterupperbound}  Assume $S$ has been isotoped in accordance with \fullref{prop:isotopS}.  

If $g$ is a bigon then $\displaystyle \sum_{R \in S \cap X^+} \chi (R) = 3s/2 - n/2$.

If $g$ is a trigon then
$\displaystyle\sum_{R \in S \cap X^+} \chi (R) = 
\begin{cases}
4s/3 - 2n/3 \\
 \mbox{or} \\ 
4s/3 - n/3. 
\end{cases} $
\end{lemma}

\begin{proof}
If $g$ is a bigon then $S \cap X^+$ is a collection of $3(s-n)/2$ bigons and $n$ trigons.   
\[ \sum_{R \in S \cap X^+} \chi (R) = 3(s-n)/2 + n = 3s/2 - n/2.\leqno{\hbox{Thus}}\]
If $g$ is a trigon then $S \cap X^+$ is either a collection of $s-2n$ bigons and $(s+4n)/3$ trigons or a collection of $s-n$ bigons and $(s+2n)/3$ trigons.
\[ \sum_{R \in S \cap X^+} \chi (R) = 
\begin{cases}
s-2n + (s+4n)/3 = 4s/3 - 2n/3 \\
 \mbox{or} \\
s-n + (s+2n)/3 = 4s/3 - n/3.
\end{cases}\]\vskip-32pt\proved
\end{proof}

\begin{lemma} \label{lem:lowerfacesbound} 
If $g$ is a bigon then
\[s + 1/2 + n/2 \leq \sum_{R \in S \cap X^-} \chi (R).\]
If $g$ is trigon then either
\begin{align*}7s/6 + 1/2 + 2n/3 &\leq \sum_{R \in S \cap X^-} \chi (R)\\
7s/6 + 1/2 + n/3 &\leq \sum_{R \in S \cap X^-} \chi (R).\tag*{\hbox{or}}\end{align*}
In any case, $s < \displaystyle\sum_{R \in S \cap X^-} \chi (R)$.
 \end{lemma}

\begin{proof} 
Since $S$ is a once punctured orientable surface of genus $g$, we have
\[ \chi(S) = 1 - 2g = -3s + \sum_{R \in S \cap X^-} \chi (R) +  \sum_{R \in S \cap X^+} \chi (R). \]
Because $s \geq 4g-1$, $(s+1)/2 \geq 2g$.  Thus
\begin{align*}1 - (s+1)/2 \leq -3s + &\sum_{R \in S \cap X^-} \chi (R) +
\sum_{R \in S \cap X^+} \chi (R)\qquad\\[-1ex]
5s/2 +1/2 - \sum_{R \in S \cap X^+} \chi (R) \leq  &\sum_{R \in S \cap X^-} \chi (R).\tag*{\hbox{and so}}\end{align*}
Due to \fullref{lem:betterupperbound}, if $g$ is a bigon
\[5s/2 + 1/2 - (3s/2 - n/2) = s + 1/2 + n/2 \leq  \sum_{R \in S \cap X^-}
\chi (R).\]\vspace{-5pt}
If $g$ is a trigon then one of following occurs:
\[\eqalignbot{5s/2 + 1/2 - (4s/3 - 2n/3) &= 7s/6 +1/2 + 2n/3 \leq  \sum_{R \in S \cap X^-} \chi (R),\cr
5s/2 + 1/2 - (4s/3 - n/3) &= 7s/6 + 1/2 + n/3 \leq  \sum_{R \in S \cap X^-} \chi (R).}\proved\]
\end{proof}
\vspace{-10pt}

\subsubsection{The disks of $S \cap X^-$}
\vspace{-5pt}
The remainder of this section is devoted to showing $$\displaystyle \sum_{R \in S \cap X^-} \chi (R) \leq s.$$  To do so we must better understand the disks of $S \cap X^-$.  
\vspace{-5pt}


\begin{lemma}\label{lem:nobigonortrigonon1456}
There are no bigons or trigons of $S \cap X^-$ that have an edge with label pair $\{1, 4\}$ or $\{5, 6\}$.  If there exist edges of $G_T$ with label pair $\{1, 6\}$ that are not edges of faces of $S \cap X^-$ parallel to $f$, then there are also no tetragons of $S \cap X^-$ with an edge having label pair $\{1, 4\}$ or $\{5, 6\}$. 
\end{lemma}
\vspace{-10pt}

\begin{proof}
In light of \fullref{prop:groomingf}, if a face $\kreis{E} \in S \cap X^-$ has an edge with one of the label pairs $\{1, 4\}$ and $\{5, 6\}$, then it also has an edge with the other label pair on the same component of $\bdry \kreis{E}$ connected either (a) by a single $(6, 1)$ corner or (b) by a sequence of a $(6, 1)$ corner, an edge with label pair $\{1, 6\}$, and another $(6, 1)$ corner.  If $\kreis{E} \in S \cap X^-$ is not parallel to $f$ but has an edge with label pair $\{1,6\}$, then this edge must appear both on $\bdry \kreis{E}$ as in situation (b) and on $\hatT$ between two parallel edges of $f$ as arranged by \fullref{prop:groomingf}.  Note that in either situation, in the other arc of $\bdry \kreis{E}$ between the two edges with label pairs $\{1, 4\}$ and $\{5, 6\}$ there cannot be just one other edge.  This is because such an edge would necessarily have label pair $\{4,5\}$ which is prohibited by \fullref{prop:isotopS}. Thus in situation (a) $\kreis{E}$ cannot be a trigon, and in situation (b) $\kreis{E}$ cannot be a tetragon (nor a bigon).  Therefore we consider situation (a) where $\kreis{E}$ is a bigon and situation (b) where $\kreis{E}$ is a trigon.
\vspace{-5pt}

If $n=0$, then there are no edges with label pair $\{1, 4\}$ or $\{5, 6\}$, and the lemma is trivial.  Therefore assume $n > 0$.  By \fullref{prop:groomingV} there exists a trigon $\kreis{D}$ of $S \cap V$ that shares an edge with $\kreis{E}$.  The three edges of $\kreis{D}$ have label pairs $\{1, 6\}$, $\{2, 5\}$, and $\{5, 6\}$.  
\vspace{-5pt}

Let $\kreis{A}$ be a bigon of $S \cap V$ with an edge parallel and closest to the edge of $\kreis{D}$ with label pair $\{2, 5\}$.  The other edge of $\kreis{A}$ has label pair $\{1, 6\}$.  Let $\Delta$ be the bigon on $T$ between these edges of $\kreis{A}$ and $\kreis{D}$ with label pair $\{2, 5\}$.  Let $\rho_{(1, 2)}$ and $\rho_{(5, 6)}$ be the rectangles on $\bdry H_{(1, 2)}$ and $\bdry H_{(5, 6)}$ respectively that contain the corners of $\kreis{A}$, $\kreis{D}$, and $\Delta$.  We may assemble $\kreis{A}$, $\kreis{D}$, $\Delta$, $\rho_{(1, 2)}$, and $\rho_{(5, 6)}$ into a high disk for $K_{(5, 6)}$.  See \fullref{fig:startinghighdisk}.
\vspace{-5pt}

\begin{figure}[ht!]
\centering
\begin{picture}(0,0)%
\includegraphics{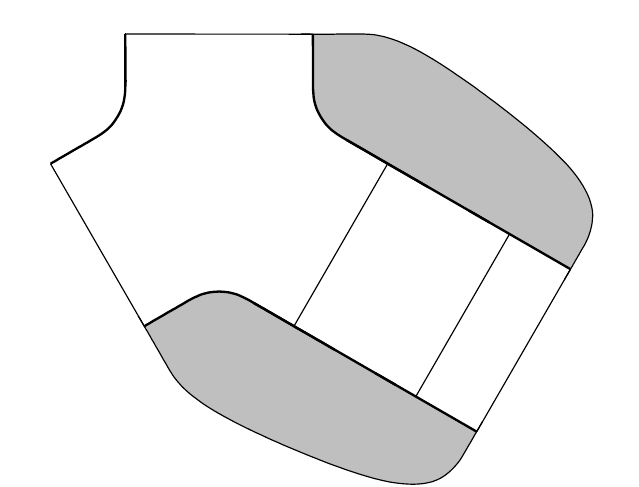}%
\end{picture}%
\setlength{\unitlength}{2960sp}%
\begingroup\makeatletter\ifx\SetFigFont\undefined%
\gdef\SetFigFont#1#2#3#4#5{%
  \reset@font\fontsize{#1}{#2pt}%
  \fontfamily{#3}\fontseries{#4}\fontshape{#5}%
  \selectfont}%
\fi\endgroup%
\begin{picture}(4077,3100)(3999,-5354)
\put(4276,-3211){\makebox(0,0)[rb]{\smash{{\SetFigFont{9}{10.8}{\rmdefault}{\mddefault}{\updefault}{\color[rgb]{0,0,0}$5$}%
}}}}
\put(6526,-3211){\makebox(0,0)[lb]{\smash{{\SetFigFont{9}{10.8}{\rmdefault}{\mddefault}{\updefault}{\color[rgb]{0,0,0}$2$}%
}}}}
\put(5851,-4561){\makebox(0,0)[rb]{\smash{{\SetFigFont{9}{10.8}{\rmdefault}{\mddefault}{\updefault}{\color[rgb]{0,0,0}$5$}%
}}}}
\put(5401,-3436){\makebox(0,0)[b]{\smash{{\SetFigFont{9}{10.8}{\rmdefault}{\mddefault}{\updefault}{\color[rgb]{0,0,0}$\bigkreis{D}$}%
}}}}
\put(6601,-5011){\makebox(0,0)[rb]{\smash{{\SetFigFont{9}{10.8}{\rmdefault}{\mddefault}{\updefault}{\color[rgb]{0,0,0}$5$}%
}}}}
\put(7276,-3586){\makebox(0,0)[lb]{\smash{{\SetFigFont{9}{10.8}{\rmdefault}{\mddefault}{\updefault}{\color[rgb]{0,0,0}$2$}%
}}}}
\put(4726,-2536){\makebox(0,0)[rb]{\smash{{\SetFigFont{9}{10.8}{\rmdefault}{\mddefault}{\updefault}{\color[rgb]{0,0,0}$6$}%
}}}}
\put(7201,-4411){\makebox(0,0)[b]{\smash{{\SetFigFont{9}{10.8}{\rmdefault}{\mddefault}{\updefault}{\color[rgb]{0,0,0}$\bigkreis{A}$}%
}}}}
\put(6076,-2386){\makebox(0,0)[lb]{\smash{{\SetFigFont{9}{10.8}{\rmdefault}{\mddefault}{\updefault}{\color[rgb]{0,0,0}$1$}%
}}}}
\put(7800,-3961){\makebox(0,0)[lb]{\smash{{\SetFigFont{9}{10.8}{\rmdefault}{\mddefault}{\updefault}{\color[rgb]{0,0,0}$1$}%
}}}}
\put(7126,-5161){\makebox(0,0)[lb]{\smash{{\SetFigFont{9}{10.8}{\rmdefault}{\mddefault}{\updefault}{\color[rgb]{0,0,0}$6$}%
}}}}
\put(4801,-4486){\makebox(0,0)[rb]{\smash{{\SetFigFont{9}{10.8}{\rmdefault}{\mddefault}{\updefault}{\color[rgb]{0,0,0}$6$}%
}}}}
\put(6526,-4111){\makebox(0,0)[b]{\smash{{\SetFigFont{9}{10.8}{\rmdefault}{\mddefault}{\updefault}{\color[rgb]{0,0,0}$\Delta$}%
}}}}
\put(7351,-4711){\makebox(0,0)[lb]{\smash{{\SetFigFont{9}{10.8}{\rmdefault}{\mddefault}{\updefault}{\color[rgb]{0,0,0}$e_A$}%
}}}}
\end{picture}%
\caption{The construction of a high disk for $K_{(5, 6)}$}
\label{fig:startinghighdisk}
\end{figure}
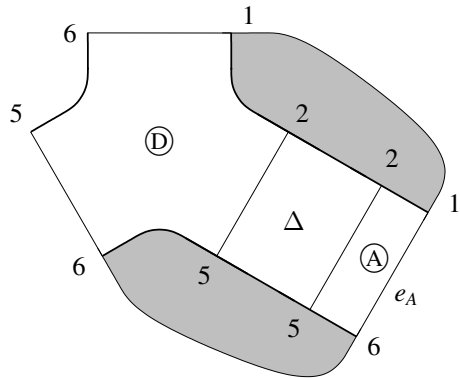
\vspace{-5pt}

{\bf Situation (a)}\qua $\kreis{E}$ is a bigon.  In this situation, $f$ may be either a bigon or a trigon.
\vspace{-5pt}

 Let $\kreis{B}$ and $\kreis{C}$ be two bigons or trigons in $X^- - N(K)$ parallel to $f$ and disjoint from both $f$ and $\kreis{E}$ so that $\kreis{A} \cap \kreis{B}$ and $\kreis{C} \cap \kreis{D}$ are each edges with label pair $\{1, 6\}$.  Note that $\kreis{B}$ and $\kreis{C}$ may actually be faces of $S \cap X^-$.  

Let $\rho_{(6, 1)}$ be the rectangle on $\bdry H_{(6, 1)}$ with boundary containing corners of $\kreis{B}$, $\kreis{E}$, and $\rho_{(5, 6)}$.  

If $f$ is a bigon, then let $\smash{\rho_{(6, 1)}'}$ be the rectangle on $\bdry H_{(6, 1)}$ with boundary containing corners of $\kreis{B}$, $\kreis{C}$, and $\rho_{(1, 2)}$.  

If $f$ is a trigon, then there is an edge of $\kreis{B}$ and an edge of $\kreis{C}$ that lies between the two parallel edges of $f$.  Let $\delta$ be the bigon on $\hatT$ bounded by these two edges of $\kreis{B}$ and $\kreis{C}$.  Then let $\smash{\rho_{(6, 1)}'}$ be the rectangle on $\bdry H_{(6, 1)}$ bounded by corners of $\kreis{B}$, $\kreis{C}$, $\rho_{(1, 2)}$, and $\delta \cap U_6$.  Let $\rho_{(6, 1)}''$ be the rectangle on $\bdry H_{(6, 1)}$ with boundary containing corners of $\kreis{B}$, $\kreis{C}$, and $\delta \cap U_1$.

Assemble $\kreis{A}, \kreis{B}, \kreis{C}, \kreis{D}, \kreis{E}, \Delta$, $\rho_{(1, 2)}$, $\rho_{(5, 6)}$, $\rho_{(6, 1)}$, and $\smash{\rho_{(6, 1)}'}$ (and $\delta$ and $\smash{\rho_{(6, 1)}''}$ if $f$ is a trigon) to form the embedded disk $D$ as shown in \fullref{fig:twolongdisks}(a) and (b).  The boundary of $D$ may be slightly extended into $H_{(4, 1)}$ so that it is composed of the arc $K_{(4, 1)}$ and an arc on $\hatT$.  Thus $D$ is a long disk.  By \fullref{longdisk}, such a disk cannot exist.

\begin{figure}[t!]
\centering
\begin{picture}(0,0)%
\includegraphics{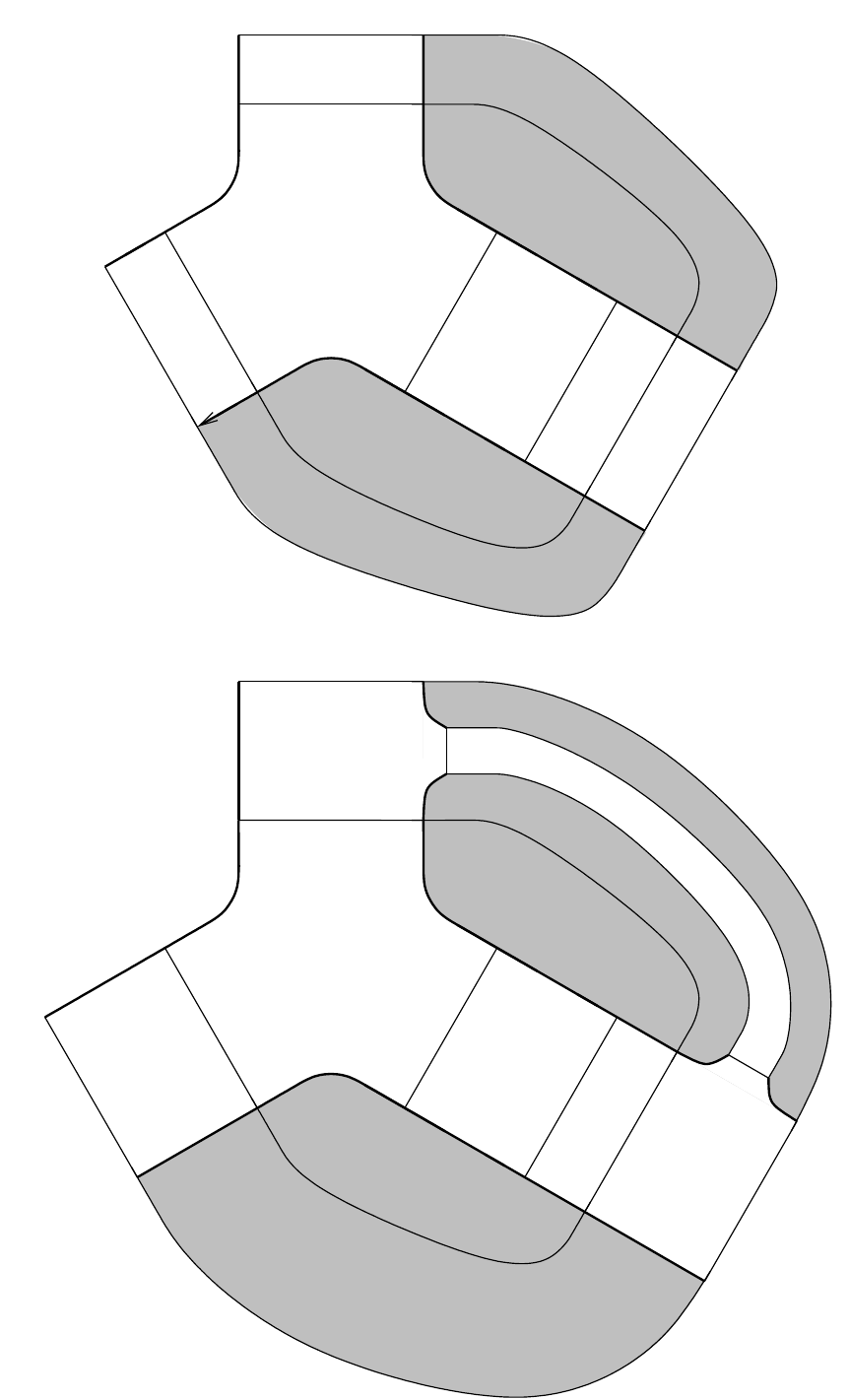}%
\end{picture}%
\setlength{\unitlength}{2960sp}%
\begingroup\makeatletter\ifx\SetFigFont\undefined%
\gdef\SetFigFont#1#2#3#4#5{%
  \reset@font\fontsize{#1}{#2pt}%
  \fontfamily{#3}\fontseries{#4}\fontshape{#5}%
  \selectfont}%
\fi\endgroup%
\begin{picture}(5503,9067)(3549,-10871)
\put(4126,-3436){\makebox(0,0)[rb]{\smash{{\SetFigFont{9}{10.8}{\rmdefault}{\mddefault}{\updefault}{\color[rgb]{0,0,0}$4$}%
}}}}
\put(4801,-5386){\makebox(0,0)[b]{\smash{{\SetFigFont{9}{10.8}{\rmdefault}{\mddefault}{\updefault}{\color[rgb]{0,0,0}(a)}%
}}}}
\put(4576,-3211){\makebox(0,0)[rb]{\smash{{\SetFigFont{9}{10.8}{\rmdefault}{\mddefault}{\updefault}{\color[rgb]{0,0,0}$5$}%
}}}}
\put(6826,-3211){\makebox(0,0)[lb]{\smash{{\SetFigFont{9}{10.8}{\rmdefault}{\mddefault}{\updefault}{\color[rgb]{0,0,0}$2$}%
}}}}
\put(6151,-4561){\makebox(0,0)[rb]{\smash{{\SetFigFont{9}{10.8}{\rmdefault}{\mddefault}{\updefault}{\color[rgb]{0,0,0}$5$}%
}}}}
\put(5701,-2311){\makebox(0,0)[b]{\smash{{\SetFigFont{9}{10.8}{\rmdefault}{\mddefault}{\updefault}{\color[rgb]{0,0,0}$\bigkreis{C}$}%
}}}}
\put(5701,-3436){\makebox(0,0)[b]{\smash{{\SetFigFont{9}{10.8}{\rmdefault}{\mddefault}{\updefault}{\color[rgb]{0,0,0}$\bigkreis{D}$}%
}}}}
\put(6901,-5011){\makebox(0,0)[rb]{\smash{{\SetFigFont{9}{10.8}{\rmdefault}{\mddefault}{\updefault}{\color[rgb]{0,0,0}$5$}%
}}}}
\put(7576,-3586){\makebox(0,0)[lb]{\smash{{\SetFigFont{9}{10.8}{\rmdefault}{\mddefault}{\updefault}{\color[rgb]{0,0,0}$2$}%
}}}}
\put(5026,-2536){\makebox(0,0)[rb]{\smash{{\SetFigFont{9}{10.8}{\rmdefault}{\mddefault}{\updefault}{\color[rgb]{0,0,0}$6$}%
}}}}
\put(7501,-5011){\makebox(0,0)[lb]{\smash{{\SetFigFont{9}{10.8}{\rmdefault}{\mddefault}{\updefault}{\color[rgb]{0,0,0}$e_2$}%
}}}}
\put(7051,-4786){\makebox(0,0)[lb]{\smash{{\SetFigFont{9}{10.8}{\rmdefault}{\mddefault}{\updefault}{\color[rgb]{0,0,0}$e_1$}%
}}}}
\put(7501,-4411){\makebox(0,0)[b]{\smash{{\SetFigFont{9}{10.8}{\rmdefault}{\mddefault}{\updefault}{\color[rgb]{0,0,0}$\bigkreis{A}$}%
}}}}
\put(7801,-4711){\makebox(0,0)[b]{\smash{{\SetFigFont{9}{10.8}{\rmdefault}{\mddefault}{\updefault}{\color[rgb]{0,0,0}$\bigkreis{B}$}%
}}}}
\put(6376,-1936){\makebox(0,0)[lb]{\smash{{\SetFigFont{9}{10.8}{\rmdefault}{\mddefault}{\updefault}{\color[rgb]{0,0,0}$6$}%
}}}}
\put(6376,-2386){\makebox(0,0)[lb]{\smash{{\SetFigFont{9}{10.8}{\rmdefault}{\mddefault}{\updefault}{\color[rgb]{0,0,0}$1$}%
}}}}
\put(8476,-4186){\makebox(0,0)[lb]{\smash{{\SetFigFont{9}{10.8}{\rmdefault}{\mddefault}{\updefault}{\color[rgb]{0,0,0}$6$}%
}}}}
\put(8101,-3961){\makebox(0,0)[lb]{\smash{{\SetFigFont{9}{10.8}{\rmdefault}{\mddefault}{\updefault}{\color[rgb]{0,0,0}$1$}%
}}}}
\put(7201,-5161){\makebox(0,0)[rb]{\smash{{\SetFigFont{9}{10.8}{\rmdefault}{\mddefault}{\updefault}{\color[rgb]{0,0,0}$6$}%
}}}}
\put(4801,-4786){\makebox(0,0)[rb]{\smash{{\SetFigFont{9}{10.8}{\rmdefault}{\mddefault}{\updefault}{\color[rgb]{0,0,0}$1$}%
}}}}
\put(7726,-5461){\makebox(0,0)[lb]{\smash{{\SetFigFont{9}{10.8}{\rmdefault}{\mddefault}{\updefault}{\color[rgb]{0,0,0}$1$}%
}}}}
\put(5401,-4486){\makebox(0,0)[lb]{\smash{{\SetFigFont{9}{10.8}{\rmdefault}{\mddefault}{\updefault}{\color[rgb]{0,0,0}$6$}%
}}}}
\put(4726,-4036){\makebox(0,0)[b]{\smash{{\SetFigFont{9}{10.8}{\rmdefault}{\mddefault}{\updefault}{\color[rgb]{0,0,0}$\bigkreis{E}$}%
}}}}
\put(5026,-2086){\makebox(0,0)[rb]{\smash{{\SetFigFont{9}{10.8}{\rmdefault}{\mddefault}{\updefault}{\color[rgb]{0,0,0}$1$}%
}}}}
\put(4801,-10486){\makebox(0,0)[b]{\smash{{\SetFigFont{9}{10.8}{\rmdefault}{\mddefault}{\updefault}{\color[rgb]{0,0,0}(b)}%
}}}}
\put(6826,-7861){\makebox(0,0)[lb]{\smash{{\SetFigFont{9}{10.8}{\rmdefault}{\mddefault}{\updefault}{\color[rgb]{0,0,0}$2$}%
}}}}
\put(6151,-9211){\makebox(0,0)[rb]{\smash{{\SetFigFont{9}{10.8}{\rmdefault}{\mddefault}{\updefault}{\color[rgb]{0,0,0}$5$}%
}}}}
\put(5701,-8086){\makebox(0,0)[b]{\smash{{\SetFigFont{9}{10.8}{\rmdefault}{\mddefault}{\updefault}{\color[rgb]{0,0,0}$\bigkreis{D}$}%
}}}}
\put(5026,-7186){\makebox(0,0)[rb]{\smash{{\SetFigFont{9}{10.8}{\rmdefault}{\mddefault}{\updefault}{\color[rgb]{0,0,0}$6$}%
}}}}
\put(7501,-9061){\makebox(0,0)[b]{\smash{{\SetFigFont{9}{10.8}{\rmdefault}{\mddefault}{\updefault}{\color[rgb]{0,0,0}$\bigkreis{A}$}%
}}}}
\put(4651,-7861){\makebox(0,0)[rb]{\smash{{\SetFigFont{9}{10.8}{\rmdefault}{\mddefault}{\updefault}{\color[rgb]{0,0,0}$5$}%
}}}}
\put(5026,-6286){\makebox(0,0)[rb]{\smash{{\SetFigFont{9}{10.8}{\rmdefault}{\mddefault}{\updefault}{\color[rgb]{0,0,0}$1$}%
}}}}
\put(3826,-8311){\makebox(0,0)[rb]{\smash{{\SetFigFont{9}{10.8}{\rmdefault}{\mddefault}{\updefault}{\color[rgb]{0,0,0}$4$}%
}}}}
\put(5701,-6736){\makebox(0,0)[b]{\smash{{\SetFigFont{9}{10.8}{\rmdefault}{\mddefault}{\updefault}{\color[rgb]{0,0,0}$\bigkreis{C}$}%
}}}}
\put(4501,-8686){\makebox(0,0)[b]{\smash{{\SetFigFont{9}{10.8}{\rmdefault}{\mddefault}{\updefault}{\color[rgb]{0,0,0}$\bigkreis{E}$}%
}}}}
\put(8101,-9361){\makebox(0,0)[b]{\smash{{\SetFigFont{9}{10.8}{\rmdefault}{\mddefault}{\updefault}{\color[rgb]{0,0,0}$\bigkreis{B}$}%
}}}}
\put(5251,-9286){\makebox(0,0)[rb]{\smash{{\SetFigFont{9}{10.8}{\rmdefault}{\mddefault}{\updefault}{\color[rgb]{0,0,0}$6$}%
}}}}
\put(6976,-9661){\makebox(0,0)[rb]{\smash{{\SetFigFont{9}{10.8}{\rmdefault}{\mddefault}{\updefault}{\color[rgb]{0,0,0}$5$}%
}}}}
\put(7276,-9961){\makebox(0,0)[lb]{\smash{{\SetFigFont{9}{10.8}{\rmdefault}{\mddefault}{\updefault}{\color[rgb]{0,0,0}$6$}%
}}}}
\put(7501,-8311){\makebox(0,0)[lb]{\smash{{\SetFigFont{9}{10.8}{\rmdefault}{\mddefault}{\updefault}{\color[rgb]{0,0,0}$2$}%
}}}}
\put(7801,-8461){\makebox(0,0)[lb]{\smash{{\SetFigFont{9}{10.8}{\rmdefault}{\mddefault}{\updefault}{\color[rgb]{0,0,0}$1$}%
}}}}
\put(6376,-7261){\makebox(0,0)[lb]{\smash{{\SetFigFont{9}{10.8}{\rmdefault}{\mddefault}{\updefault}{\color[rgb]{0,0,0}$1$}%
}}}}
\put(4426,-9661){\makebox(0,0)[rb]{\smash{{\SetFigFont{9}{10.8}{\rmdefault}{\mddefault}{\updefault}{\color[rgb]{0,0,0}$1$}%
}}}}
\put(6451,-6961){\makebox(0,0)[lb]{\smash{{\SetFigFont{9}{10.8}{\rmdefault}{\mddefault}{\updefault}{\color[rgb]{0,0,0}$6$}%
}}}}
\put(8326,-8536){\makebox(0,0)[rb]{\smash{{\SetFigFont{9}{10.8}{\rmdefault}{\mddefault}{\updefault}{\color[rgb]{0,0,0}$6$}%
}}}}
\put(6301,-6136){\makebox(0,0)[lb]{\smash{{\SetFigFont{9}{10.8}{\rmdefault}{\mddefault}{\updefault}{\color[rgb]{0,0,0}$6$}%
}}}}
\put(8776,-9211){\makebox(0,0)[lb]{\smash{{\SetFigFont{9}{10.8}{\rmdefault}{\mddefault}{\updefault}{\color[rgb]{0,0,0}$6$}%
}}}}
\put(8626,-8836){\makebox(0,0)[lb]{\smash{{\SetFigFont{9}{10.8}{\rmdefault}{\mddefault}{\updefault}{\color[rgb]{0,0,0}$1$}%
}}}}
\put(8176,-10261){\makebox(0,0)[lb]{\smash{{\SetFigFont{9}{10.8}{\rmdefault}{\mddefault}{\updefault}{\color[rgb]{0,0,0}$1$}%
}}}}
\put(6451,-6436){\makebox(0,0)[lb]{\smash{{\SetFigFont{9}{10.8}{\rmdefault}{\mddefault}{\updefault}{\color[rgb]{0,0,0}$1$}%
}}}}
\put(7201,-7186){\makebox(0,0)[b]{\smash{{\SetFigFont{9}{10.8}{\rmdefault}{\mddefault}{\updefault}{\color[rgb]{0,0,0}$\rho_{(6, 1)}'$}%
}}}}
\put(7501,-6661){\makebox(0,0)[b]{\smash{{\SetFigFont{9}{10.8}{\rmdefault}{\mddefault}{\updefault}{\color[rgb]{0,0,0}$\rho_{(6, 1)}''$}%
}}}}
\put(7126,-2386){\makebox(0,0)[b]{\smash{{\SetFigFont{9}{10.8}{\rmdefault}{\mddefault}{\updefault}{\color[rgb]{0,0,0}$\rho_{(6, 1)}'$}%
}}}}
\put(5626,-5161){\makebox(0,0)[b]{\smash{{\SetFigFont{9}{10.8}{\rmdefault}{\mddefault}{\updefault}{\color[rgb]{0,0,0}$\rho_{(6, 1)}$}%
}}}}
\put(6826,-8761){\makebox(0,0)[b]{\smash{{\SetFigFont{9}{10.8}{\rmdefault}{\mddefault}{\updefault}{\color[rgb]{0,0,0}$\Delta$}%
}}}}
\put(6226,-4861){\makebox(0,0)[b]{\smash{{\SetFigFont{9}{10.8}{\rmdefault}{\mddefault}{\updefault}{\color[rgb]{0,0,0}$\rho_{(5, 6)}$}%
}}}}
\put(6826,-4111){\makebox(0,0)[b]{\smash{{\SetFigFont{9}{10.8}{\rmdefault}{\mddefault}{\updefault}{\color[rgb]{0,0,0}$\Delta$}%
}}}}
\put(6826,-2836){\makebox(0,0)[b]{\smash{{\SetFigFont{9}{10.8}{\rmdefault}{\mddefault}{\updefault}{\color[rgb]{0,0,0}$\rho_{(1, 2)}$}%
}}}}
\put(5701,-10186){\makebox(0,0)[b]{\smash{{\SetFigFont{9}{10.8}{\rmdefault}{\mddefault}{\updefault}{\color[rgb]{0,0,0}$\rho_{(6, 1)}$}%
}}}}
\put(7651,-7111){\makebox(0,0)[b]{\smash{{\SetFigFont{9}{10.8}{\rmdefault}{\mddefault}{\updefault}{\color[rgb]{0,0,0}$\delta$}%
}}}}
\put(6226,-9511){\makebox(0,0)[b]{\smash{{\SetFigFont{9}{10.8}{\rmdefault}{\mddefault}{\updefault}{\color[rgb]{0,0,0}$\rho_{(5, 6)}$}%
}}}}
\put(6826,-7561){\makebox(0,0)[b]{\smash{{\SetFigFont{9}{10.8}{\rmdefault}{\mddefault}{\updefault}{\color[rgb]{0,0,0}$\rho_{(1, 2)}$}%
}}}}
\end{picture}%
\caption{The construction of a long disk if $f$ is (a) a bigon and (b) a trigon}
\label{fig:twolongdisks}
\end{figure}
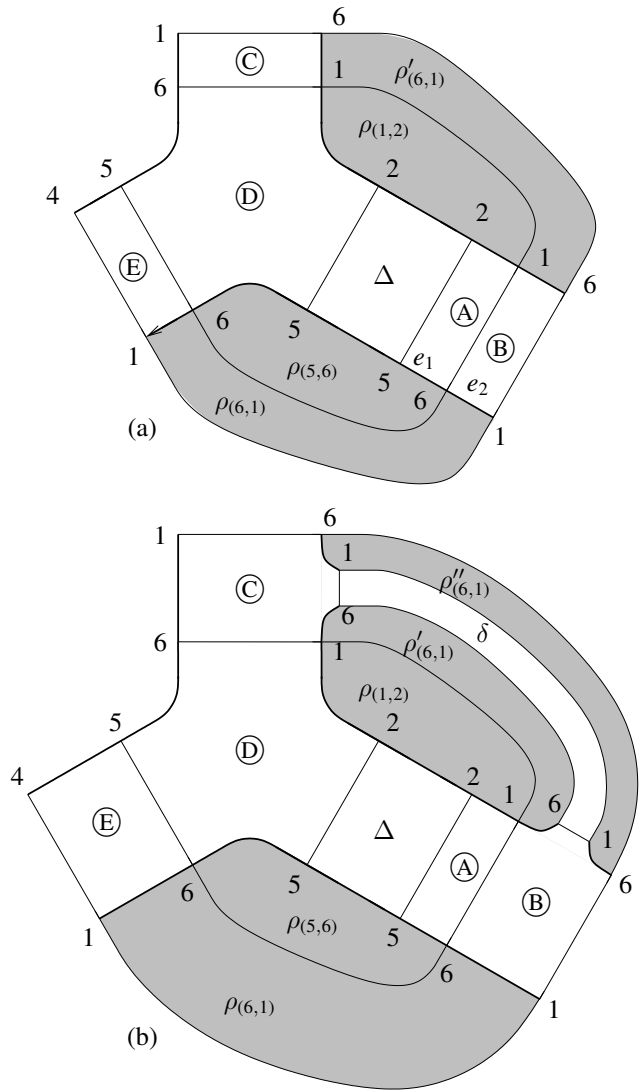

{\bf Situation (b)}\qua $\kreis{E}$ is a trigon.  Here $f$ is necessarily a trigon.

Again, let $\kreis{B}$ and $\kreis{C}$ be two bigons or trigons in $X^- - N(K)$ parallel to $f$ and disjoint from both $f$ and $\kreis{E}$ so that $\kreis{A} \cap \kreis{B}$ and $\kreis{C} \cap \kreis{D}$ are each edges with label pair $\{1, 6\}$.  Note that $\kreis{B}$ and $\kreis{C}$ may actually be faces of $S \cap X^-$.  

 Let $\rho_{(6, 1)}$ be the rectangle on $\bdry H_{(6, 1)}$ bounded by corners of $\kreis{B}$, $\kreis{C}$, and $\rho_{(1, 2)}$.  

Since $f$ is a trigon, there is an edge of $\kreis{B}$ and an edge of $\kreis{E}$ that lies between the two parallel edges of $f$.  Let $\delta$ be the bigon on $T$ bounded by these two edges of $\kreis{B}$ and $\kreis{E}$.  Note that there is an edge of $\kreis{C}$ between the two parallel edges of $f$, but it is not contained in $\delta$.  

Let $\smash{\rho_{(6, 1)}'}$ be the rectangle on $\bdry H_{(6, 1)}$ bounded by corners of $\kreis{B}$, $\kreis{E}$, $\rho_{(5, 6)}$, and $\delta \cap U_1$.  Let $\rho_{(6, 1)}''$ be the rectangle on $\bdry H_{(6, 1)}$ with boundary containing corners of $\kreis{B}$, $\kreis{E}$, and $\delta \cap U_6$.

Assemble $\kreis{A}, \kreis{B}, \kreis{C}, \kreis{D}, \kreis{E}, \Delta$, $\rho_{(1, 2)}$, $\rho_{(5, 6)}$, $\rho_{(6, 1)}$, $\smash{\rho_{(6, 1)}'}$, and $\smash{\rho_{(6, 1)}''}$ to form the embedded disk $D$, a ``lopsided bigon,'' as shown in \fullref{fig:lopsidedbigon}.  The boundary of $D$ may be slightly extended into $H_{(4, 1)}$ so that it is composed of the arc $K_{(4, 1)}$, the arc $K_{(6, 1)}$, and two arcs, say $\tau_1$ and $\tau_2$, on $\hatT$.  

\begin{figure}[ht!]
\centering
\begin{picture}(0,0)%
\includegraphics{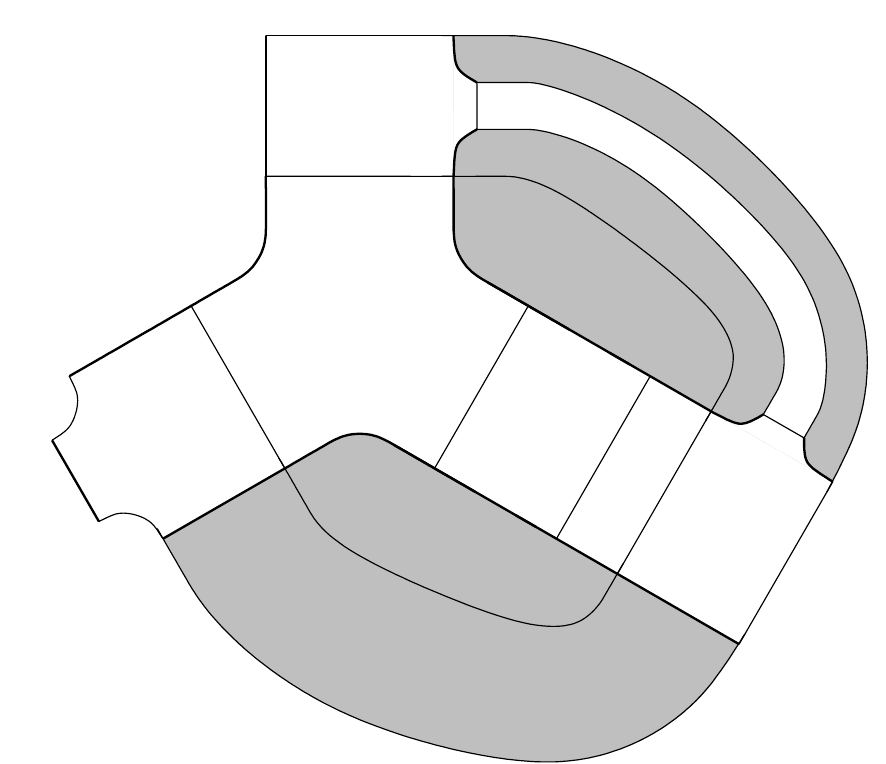}%
\end{picture}%
\setlength{\unitlength}{2960sp}%
\begingroup\makeatletter\ifx\SetFigFont\undefined%
\gdef\SetFigFont#1#2#3#4#5{%
  \reset@font\fontsize{#1}{#2pt}%
  \fontfamily{#3}\fontseries{#4}\fontshape{#5}%
  \selectfont}%
\fi\endgroup%
\begin{picture}(5653,4867)(2499,-6371)
\put(3076,-4936){\makebox(0,0)[rb]{\smash{{\SetFigFont{9}{10.8}{\rmdefault}{\mddefault}{\updefault}{\color[rgb]{0,0,0}$6$}%
}}}}
\put(4802,-2236){\makebox(0,0)[b]{\smash{{\SetFigFont{9}{10.8}{\rmdefault}{\mddefault}{\updefault}{\color[rgb]{0,0,0}$\bigkreis{C}$}%
}}}}
\put(3601,-4187){\makebox(0,0)[b]{\smash{{\SetFigFont{9}{10.8}{\rmdefault}{\mddefault}{\updefault}{\color[rgb]{0,0,0}$\bigkreis{E}$}%
}}}}
\put(7201,-4861){\makebox(0,0)[b]{\smash{{\SetFigFont{9}{10.8}{\rmdefault}{\mddefault}{\updefault}{\color[rgb]{0,0,0}$\bigkreis{B}$}%
}}}}
\put(4351,-4786){\makebox(0,0)[rb]{\smash{{\SetFigFont{9}{10.8}{\rmdefault}{\mddefault}{\updefault}{\color[rgb]{0,0,0}$6$}%
}}}}
\put(6076,-5161){\makebox(0,0)[rb]{\smash{{\SetFigFont{9}{10.8}{\rmdefault}{\mddefault}{\updefault}{\color[rgb]{0,0,0}$5$}%
}}}}
\put(6376,-5461){\makebox(0,0)[lb]{\smash{{\SetFigFont{9}{10.8}{\rmdefault}{\mddefault}{\updefault}{\color[rgb]{0,0,0}$6$}%
}}}}
\put(6601,-3811){\makebox(0,0)[lb]{\smash{{\SetFigFont{9}{10.8}{\rmdefault}{\mddefault}{\updefault}{\color[rgb]{0,0,0}$2$}%
}}}}
\put(6901,-3961){\makebox(0,0)[lb]{\smash{{\SetFigFont{9}{10.8}{\rmdefault}{\mddefault}{\updefault}{\color[rgb]{0,0,0}$1$}%
}}}}
\put(5476,-2761){\makebox(0,0)[lb]{\smash{{\SetFigFont{9}{10.8}{\rmdefault}{\mddefault}{\updefault}{\color[rgb]{0,0,0}$1$}%
}}}}
\put(3526,-5161){\makebox(0,0)[rb]{\smash{{\SetFigFont{9}{10.8}{\rmdefault}{\mddefault}{\updefault}{\color[rgb]{0,0,0}$1$}%
}}}}
\put(5552,-2461){\makebox(0,0)[lb]{\smash{{\SetFigFont{9}{10.8}{\rmdefault}{\mddefault}{\updefault}{\color[rgb]{0,0,0}$6$}%
}}}}
\put(7426,-4036){\makebox(0,0)[rb]{\smash{{\SetFigFont{9}{10.8}{\rmdefault}{\mddefault}{\updefault}{\color[rgb]{0,0,0}$6$}%
}}}}
\put(5401,-1636){\makebox(0,0)[lb]{\smash{{\SetFigFont{9}{10.8}{\rmdefault}{\mddefault}{\updefault}{\color[rgb]{0,0,0}$6$}%
}}}}
\put(7876,-4711){\makebox(0,0)[lb]{\smash{{\SetFigFont{9}{10.8}{\rmdefault}{\mddefault}{\updefault}{\color[rgb]{0,0,0}$6$}%
}}}}
\put(7726,-4336){\makebox(0,0)[lb]{\smash{{\SetFigFont{9}{10.8}{\rmdefault}{\mddefault}{\updefault}{\color[rgb]{0,0,0}$1$}%
}}}}
\put(7276,-5761){\makebox(0,0)[lb]{\smash{{\SetFigFont{9}{10.8}{\rmdefault}{\mddefault}{\updefault}{\color[rgb]{0,0,0}$1$}%
}}}}
\put(5551,-1936){\makebox(0,0)[lb]{\smash{{\SetFigFont{9}{10.8}{\rmdefault}{\mddefault}{\updefault}{\color[rgb]{0,0,0}$1$}%
}}}}
\put(2776,-4411){\makebox(0,0)[rb]{\smash{{\SetFigFont{9}{10.8}{\rmdefault}{\mddefault}{\updefault}{\color[rgb]{0,0,0}$1$}%
}}}}
\put(5926,-3361){\makebox(0,0)[lb]{\smash{{\SetFigFont{9}{10.8}{\rmdefault}{\mddefault}{\updefault}{\color[rgb]{0,0,0}$2$}%
}}}}
\put(5251,-4711){\makebox(0,0)[rb]{\smash{{\SetFigFont{9}{10.8}{\rmdefault}{\mddefault}{\updefault}{\color[rgb]{0,0,0}$5$}%
}}}}
\put(4801,-3586){\makebox(0,0)[b]{\smash{{\SetFigFont{9}{10.8}{\rmdefault}{\mddefault}{\updefault}{\color[rgb]{0,0,0}$\bigkreis{D}$}%
}}}}
\put(4126,-2686){\makebox(0,0)[rb]{\smash{{\SetFigFont{9}{10.8}{\rmdefault}{\mddefault}{\updefault}{\color[rgb]{0,0,0}$6$}%
}}}}
\put(6601,-4561){\makebox(0,0)[b]{\smash{{\SetFigFont{9}{10.8}{\rmdefault}{\mddefault}{\updefault}{\color[rgb]{0,0,0}$\bigkreis{A}$}%
}}}}
\put(3751,-3361){\makebox(0,0)[rb]{\smash{{\SetFigFont{9}{10.8}{\rmdefault}{\mddefault}{\updefault}{\color[rgb]{0,0,0}$5$}%
}}}}
\put(4126,-1786){\makebox(0,0)[rb]{\smash{{\SetFigFont{9}{10.8}{\rmdefault}{\mddefault}{\updefault}{\color[rgb]{0,0,0}$1$}%
}}}}
\put(2926,-3811){\makebox(0,0)[rb]{\smash{{\SetFigFont{9}{10.8}{\rmdefault}{\mddefault}{\updefault}{\color[rgb]{0,0,0}$4$}%
}}}}
\put(5926,-4261){\makebox(0,0)[b]{\smash{{\SetFigFont{9}{10.8}{\rmdefault}{\mddefault}{\updefault}{\color[rgb]{0,0,0}$\Delta$}%
}}}}
\put(6751,-2611){\makebox(0,0)[b]{\smash{{\SetFigFont{9}{10.8}{\rmdefault}{\mddefault}{\updefault}{\color[rgb]{0,0,0}$\delta$}%
}}}}
\end{picture}%
\caption{The construction of the ``lopsided bigon''}
\label{fig:lopsidedbigon}
\end{figure}
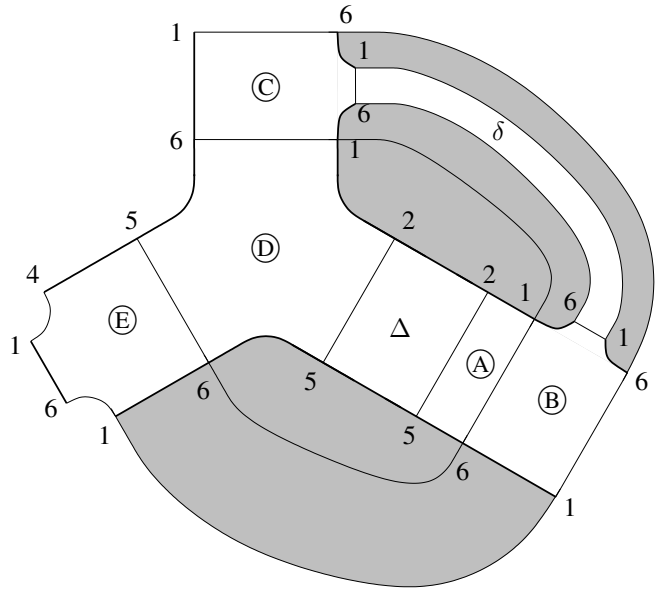

Isotop the arc $K_{(4, 1)}$ across $D$ so that the decomposition $K = K_{(1, 4)} \cup K_{(4, 1)}$ becomes $K = K_{(1, 4)} \cup \tau_1 \cup K_{(6, 1)} \cup \tau_2$.  Since the ends of $K_{(6, 1)}$ are in $X^-$ while the ends of $K_{(1, 4)}$ are in $X^+$, a further slight isotopy of $K$ in $N(\tau_1 \cup \tau_2)$ puts $K$ in Morse position without introducing new critical points.  Since the extrema of $K$ are in a $1-1$ correspondence with a proper subset of the former extrema of $K$ (with the correspondence taking maxima to maxima and minima to minima), the width of $K$ has been reduced.  In other words, the width of the two arcs $K_{(1, 4)} \cup K_{(6, 1)}$ is less than the width of the knot $K_{(1, 4)} \cup K_{(4, 1)}$.  This contradicts that $K$ is in thin position.
\end{proof}

\begin{lemma}\label{lem:bigonlabels}
A bigon of $S \cap X^-$ that is not parallel to $f$ must have edges with label pairs $\{2, 5\}$ and $\{3, 4\}$.
\end{lemma}

\begin{proof}
Due to \fullref{lem:nobigonortrigonon1456}, a bigon of $S \cap X^-$ cannot have an edge with label pair $\{1, 4\}$ or $\{5, 6\}$.  Thus an edge of a bigon of $S \cap X^-$ not parallel to $f$ may have label pair $\{2, 5\}$, $\{3, 4\}$, or $\{2, 3\}$. 

If an bigon of $S \cap X^-$ has label pair $\{2, 3\}$, then both edges must have this label pair.  Such a bigon is thus bounded by an $S2$ cycle.  Since $f$ is in $X^-$ too, this contradicts \fullref{oppositesides}.  

Because $H_{(3, 4)} \subseteq X^+$, the conclusion of the lemma follows.
\end{proof}

\begin{lemma}\label{lem:existenceofbigon}
There must exist a bigon of $S \cap X^-$ whose edges have label pairs $\{2, 5\}$ and $\{3, 4\}$.
\end{lemma}

\begin{proof}
Assume every bigon of $S \cap X^-$ is parallel to $f$.

{\bf Case 1}\qua $f$ is a bigon.\qua

By \fullref{prop:isotopS}, $(s-n)/2$ faces of $S \cap X^-$ are bigons parallel to $f$.  Since by \fullref{lem:nobigonortrigonon1456} the edges with label pair $\{1, 4\}$ and $\{5, 6\}$ cannot be edges of bigons or trigons, at worst they are edges of tetragons.  Since each edge with label pair $\{1, 4\}$ is connected to an edge with label pair $\{5, 6\}$ by a corner on $\bdry H_{(6, 1)}$ (due to \fullref{lem:nobigonortrigonon1456} because the edges of the bigon $f$ are not parallel on $\hatT$) and each such edge class has $n$ edges, together with $2n$ edges that have label pairs $\{2, 5\}$, $\{3, 4\}$, or $\{2, 3\}$, these may bound at most $n$ tetragons.  

The remaining $3s - 2 \cdot (s-n)/2 - 4 \cdot n = 2s - 3n$ edges all have label pair $\{2, 5\}$, $\{3, 4\}$, or $\{2, 3\}$.  Any trigon with all its edge among these, must have one edge with each of these three label pairs.  Since there are only $n$ edges with label pair $\{2, 3\}$, there may be at most $n$ such trigons.  The remaining $2s-3n - 3 \cdot n = 2s-6n$ edges may at worst bound tetragons.  
  Thus 
\[ \sum_{R \in S \cap X^-} \chi (R)  \leq (s-n)/2 + n + n + (2s-6n)/4 = s.\]

{\bf Case 2}\qua $f$ is a trigon.\qua

 By \fullref{prop:isotopS}, either $(s-n)/3$ or $(s-2n)/3$ of the faces of $S \cap X^-$ are trigons parallel to $f$.

{\bf Case 2a}\qua
If there are $(s-n)/3$ trigons parallel to $f$, then each edge with label pair $\{1, 6\}$ is an edge of a trigon parallel to $f$.  Then by \fullref{lem:nobigonortrigonon1456} the edges with label pair $\{1, 4\}$ and $\{5, 6\}$ cannot be edges of bigons or trigons.  At worst they are edges of tetragons.  Since each edge with label pair $\{1, 4\}$ is connected to an edge with label pair $\{5, 6\}$ by a corner on $\bdry H_{(6, 1)}$  (because in order to have $(s-n)/3$ trigons parallel to $f$, each edge of $G_T$ between the parallel edges of $f$ must belong to one of these trigons---see the last paragraph of the proof of \fullref{prop:isotopS} in \fullref{sec:proofofpropisotopS}) and each such edge class has $n$ edges, then together with $2n$ of the edges that have label pairs $\{2, 5\}$, $\{3, 4\}$, or $\{2, 3\}$, these may bound at most $n$ tetragons.  

The remaining $3s - 3 \cdot (s-n)/3 - 4 \cdot n = 2s - 3n$ edges all have label pair $\{2, 5\}$, $\{3, 4\}$, or $\{2, 3\}$.  Any trigon with all its edges among these, must have one edge with each of these three label pairs.  Since there are only $n$ edges with label pair $\{2, 3\}$, there may be at most $n$ such trigons.  The remaining $2s-3n - 3 \cdot n = 2s-6n$ edges may at worst bound $(2s-6n)/4$ tetragons.   
\[ \sum_{R \in S \cap X^-} \chi (R)  \leq (s-n)/3 + n + n + (2s-6n)/4 = 5s/6 + n/6.\leqno{\hbox{Thus}}\]
{\bf Case 2b}\qua
If there are $(s-2n)/3$ trigons parallel to $f$ (and $n \neq 0$), then there exists an edge with label pair $\{1, 6\}$ that is not an edge of a trigon parallel to $f$ (indeed in $G_T$ between the parallel edges of $f$ there are $n$ edges that do not belong to faces of $S \cap X^-$ that are trigons parallel to $f$---again, see the last paragraph of the proof of \fullref{prop:isotopS} in \fullref{sec:proofofpropisotopS}).  Then by \fullref{lem:nobigonortrigonon1456} the edges with label pair $\{1, 4\}$ and $\{5, 6\}$ cannot be edges of bigons, trigons, or tetragons.  At worst they are edges of pentagons.  Since each edge with label pair $\{1, 4\}$ is connected to an edge with label pair $\{5, 6\}$ by a two corners on $\bdry H_{(6, 1)}$ and an edge with label pair $\{1, 6\}$ and each such edge class has $n$ edges, together with $2n$ edges that have label pairs $\{2, 5\}$, $\{3, 4\}$, or $\{2, 3\}$, these may bound at most $n$ pentagons.

The remaining $3s - 3 \cdot (s-2n)/3 - 5 \cdot n = 2s - 3n$ edges all have label pair $\{2, 5\}$, $\{3, 4\}$, or $\{2, 3\}$.  Any trigon with all its edges among these, must have one edge with each of these three label pairs.  Since there are only $n$ edges with label pair $\{2, 3\}$, there may be at most $n$ such trigons.  The remaining $2s-3n - 3 \cdot n = 2s-6n$ edges may at worst bound $(2s-6n)/4$ tetragons.
\[ \sum_{R \in S \cap X^-} \chi (R)  \leq (s-2n)/3 + n + n + (2s-6n)/4 = 5s/6 - n/6.\leqno{\hbox{Thus}}\]
Since necessarily $n<s$, we have $5s/6 - n/6 < 5s/6 + n/6 < s$.  Thus 
\[ \sum_{R \in S \cap X^-} \chi (R)  \leq s.\]
However, by \fullref{lem:lowerfacesbound}
\[s < \sum_{R \in S \cap X^-} \chi (R)\]
 regardless of whether $g$ is a bigon or trigon.  This, in all cases, is a contradiction.  Thus there must be a bigon of $S \cap X^-$ that is not parallel to $f$.  By \fullref{lem:bigonlabels}, the edges of such a bigon must have label pairs $\{2, 5\}$ and $\{3, 4\}$.
\end{proof}

\begin{lemma}\label{lem:arrangingedges2345}

Let $B$ be a bigon of $S \cap X^-$ whose edges $b_0$ and $b_1$ have label pairs $\{2, 5\}$ and $\{3, 4\}$ respectively.  Let $R$ be an $m$--gon of $S \cap X^-$ whose edges have label pairs among $\{2, 5\}$, $\{3, 4\}$, and $\{2, 3\}$.  For each $i = 0,1$ let $p_i$ be the number of edges of $R$ in the same edge class of $G_T$ as $b_i$, and let $\bar{p}_i$ be the number of edges of $R$ with the same label pair as $b_i$ but not in the same edge class as $b_i$.

Assume that the edges of $B$ and $R$ together lie in an essential annulus $A$ on $\hatT$.  If $M = \bar{N}(A \cup H_{(2, 3)} \cup H_{(4, 5)} \cup B \cup R)$ is a solid torus such that the core of $A$ is a longitude, then $R$ has exactly one edge with label pair $\{2, 3\}$.  Thus for each $i = 0, 1$ we have $p_i + \bar{p}_i  = (m-1)/2$.  Furthermore either
\begin{itemize}
\item $\bar{p}_i = 0$ for $i = 0,1$, or
\item $p_0 = \bar{p}_1$ and $p_1 = \bar{p}_0$.
\end{itemize}
\end{lemma}

One may care to compare this lemma in the case that $R$ is a trigon with \fullref{lem:funnytrigon}.

\begin{figure}[ht!]
\centering
\begin{picture}(0,0)%
\includegraphics{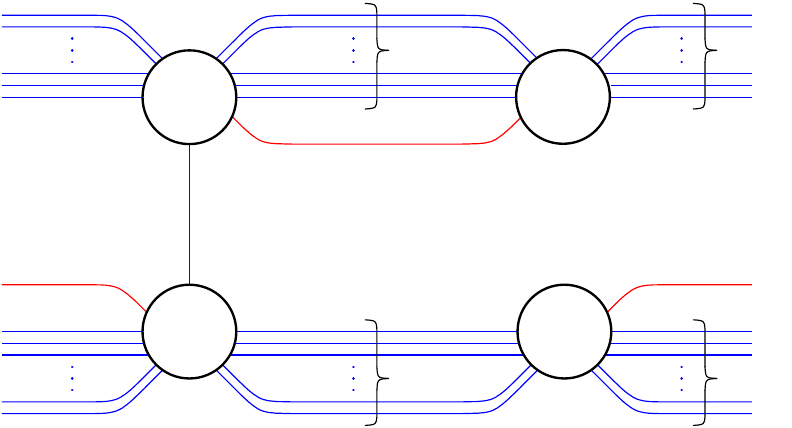}%
\end{picture}%
\setlength{\unitlength}{2960sp}%
\begingroup\makeatletter\ifx\SetFigFont\undefined%
\gdef\SetFigFont#1#2#3#4#5{%
  \reset@font\fontsize{#1}{#2pt}%
  \fontfamily{#3}\fontseries{#4}\fontshape{#5}%
  \selectfont}%
\fi\endgroup%
\begin{picture}(5048,2724)(4639,-7873)
\put(5101,-6886){\makebox(0,0)[b]{\smash{{\SetFigFont{9}{10.8}{\rmdefault}{\mddefault}{\updefault}{\color[rgb]{0,0,0}$b_1$}%
}}}}
\put(7051,-6286){\makebox(0,0)[b]{\smash{{\SetFigFont{9}{10.8}{\rmdefault}{\mddefault}{\updefault}{\color[rgb]{0,0,0}$b_0$}%
}}}}
\put(7126,-5536){\makebox(0,0)[lb]{\smash{{\SetFigFont{9}{10.8}{\rmdefault}{\mddefault}{\updefault}{\color[rgb]{0,0,0}$p_0$}%
}}}}
\put(9226,-5536){\makebox(0,0)[lb]{\smash{{\SetFigFont{9}{10.8}{\rmdefault}{\mddefault}{\updefault}{\color[rgb]{0,0,0}$\bar{p}_0$}%
}}}}
\put(9226,-7636){\makebox(0,0)[lb]{\smash{{\SetFigFont{9}{10.8}{\rmdefault}{\mddefault}{\updefault}{\color[rgb]{0,0,0}$p_1$}%
}}}}
\put(7126,-7636){\makebox(0,0)[lb]{\smash{{\SetFigFont{9}{10.8}{\rmdefault}{\mddefault}{\updefault}{\color[rgb]{0,0,0}$\bar{p}_1$}%
}}}}
\put(9001,-6886){\makebox(0,0)[b]{\smash{{\SetFigFont{9}{10.8}{\rmdefault}{\mddefault}{\updefault}{\color[rgb]{0,0,0}$b_1$}%
}}}}
\put(8251,-5836){\makebox(0,0)[b]{\smash{{\SetFigFont{9}{10.8}{\rmdefault}{\mddefault}{\updefault}{\color[rgb]{0,0,0}$5$}%
}}}}
\put(5851,-5836){\makebox(0,0)[b]{\smash{{\SetFigFont{9}{10.8}{\rmdefault}{\mddefault}{\updefault}{\color[rgb]{0,0,0}$2$}%
}}}}
\put(5851,-7336){\makebox(0,0)[b]{\smash{{\SetFigFont{9}{10.8}{\rmdefault}{\mddefault}{\updefault}{\color[rgb]{0,0,0}$3$}%
}}}}
\put(8251,-7336){\makebox(0,0)[b]{\smash{{\SetFigFont{9}{10.8}{\rmdefault}{\mddefault}{\updefault}{\color[rgb]{0,0,0}$4$}%
}}}}
\end{picture}%
\caption{The edges of $\Gamma$ on $A$ when $\bar{p}_0 \neq 0$ and $\bar{p}_1 \neq 0$}
\label{fig:edgesofGammaonA}
\end{figure}
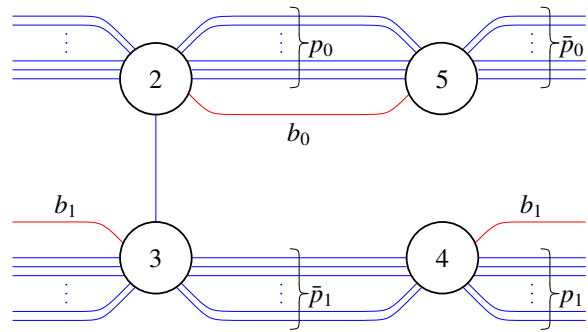

\begin{proof}
Consider the genus $3$ handlebody $\bar{N}(A \cup H_{(2, 3)} \cup H_{(4, 5)})$ which we view as the torus $A \times [0,1]$ with the $1$--handles $H_{(2, 3)}$ and $H_{(4, 5)}$ attached to $A \times \{1\}$.   We then obtain $M$ by attaching to $A \times [0,1] \cup \bdry H_{(2, 3)} \cup \bdry H_{(4, 5)}$ the $2$--handles $\smash{\bar{N}(B)}$ and $\smash{\bar{N}(R)}$.  Assume $M$ is a solid torus such that the core of $A$ is a longitude.

Let $E$ be a disk.  By attaching the $2$--handle $\bar{N}(E)$ to $A \times [0,1]$ along the core of $A \times \{0\}$ we form a $3$--ball $W = \smash{\bar{N}(A \cup E)}$.  Let $E^+$ and $E^-$ be the two sides of $\bar{N}(E)$ on $\bdry W$.  Notice that $\bdry W - \Int(E^+ \cup E^-)$ may be identified with the annulus $A$.  

To the ball $W$ we attach the ``$1$--handle'' $H=\smash{H_{(2, 3)} \cup \bar{N}(B) \cup H_{(4, 5)}}$ forming the solid torus $M' = \smash{\bar{N}(A \cup H_{(2, 3)} \cup H_{(4, 5)} \cup B \cup E)}$.  Note $M' \cup \smash{\bar{N}(R) = M \cup \bar{N}(E)}$ and $M \cup \smash{\bar{N}(E)} \cong B^3$.  Therefore on $\bdry M'$ the curve $\bdry R$ must be isotopic to a curve that algebraically intersects the cocore of $H$ just once.  Because $R$ has no edges with label pair $\{4, 5\}$, it must have exactly one edge with label pair $\{2, 3\}$.  It then follows that $p_i + \bar{p}_i = (m-1)/2$ for each $i = 0, 1$. 

Let $\Gamma$ be the subgraph of $G_T$ consisting of the edges of $B$ and $R$ and the vertices $U_2$, $U_3$, $U_4$, and $U_5$.  Since $\Gamma \subseteq A$, $\Gamma \subseteq \bdry W$.  The positions of $E^+$ and $E^-$ relative to $\Gamma$ define how $\Gamma$ is contained in $A$.

Unless $E^+$ and $E^-$ are contained in distinct bigons of $\bdry W \cut \Gamma$ whose edges lie in two distinct edge classes of $\Gamma$ on $\bdry W$, either $\bar{p}_0=0$ or $\bar{p}_1=0$.  Thus we assume both $\bar{p}_0 \neq 0$ and $\bar{p}_1 \neq 0$.  This forces $\Gamma$ to appear on $A$ as in \fullref{fig:edgesofGammaonA}.  It is then clear (eg\ by sliding over $B$ each of the $p_i$ edges of $R$ that are parallel to $b_i$) that in order for $M$ to be a solid torus, $p_0 = \bar{p}_1$ and $p_1 = \bar{p}_0$. 
\end{proof} 

\begin{lemma} \label{lem:onlybigons}
The only disks of $S \cap X^-$ whose edges have label pairs among $\{2, 5\}$, $\{3, 4\}$, and $\{2, 3\}$ are bigons.  Furthermore all such bigons are parallel.
\end{lemma}

\begin{proof}
\fullref{lem:existenceofbigon} implies there exists a bigon $B$ of $S \cap X^-$ whose edges have label pairs $\{2, 5\}$ and $\{3, 4\}$.  By \fullref{lem:parallelbigons} any bigon with edges having label pairs among $\{2, 5\}$, $\{3, 4\}$, and $\{2, 3\}$ is parallel to $B$.

Assume $n=0$.  Then there are no edges with label pair $\{2,3\}$.  Thus, due to the existence of $B$, any boundary component of a face of $S \cap X^-$ with edges having all its label pairs among $\{2, 5\}$ and $\{3, 4\}$ may only have two edges.  

To see this, let $A_{\{2,5\}}$ and $A_{\{3,4\}}$ be narrow essential annuli in $\hatT$ which contain all edges of $G_T$ with label pair $\{2,5\}$ and $\{3, 4\}$ respectively and the pairs of fat vertices $\{U_2, U_5\}$ and $\{U_3, U_4\}$ respectively.   Abstractly cap off the boundary components of $A_{\{2,5\}}$ and $A_{\{3,4\}}$ and delete the interiors of the fat vertices to form annuli $\hatA_{\{2,5\}}$ and $\hatA_{\{3,4\}}$ in which the edges of $G_T$ they contain run from one boundary component to the other.  Then joining $\hatA_{\{2,5\}}$ and $\hatA_{\{3,4\}}$ with the annuli $\bdry H_{(2,3)} \cut (U_2 \cup U_3)$ and $\bdry H_{(4,5)} \cut (U_4 \cup U_5)$ (with the corners of faces of $S \cap X^-$ on them) forms a torus.  The edges and corners of $S \cap X^-$ on this torus form a $1$--manifold in which the closed components all meet each $\hatA_{\{2,5\}}$, $\hatA_{\{3,4\}}$, $\bdry H_{(2,3)}$, and $\bdry H_{(4,5)}$ the same number of times.  Since $\bdry B$ is one such component, any other such component must meet  $\hatA_{\{2,5\}}$ and $\hatA_{\{3,4\}}$ each just once as well, ie\ must have just two edges.  Therefore any boundary component of a face of $S \cap X^-$ with edges having all its label pairs among $\{2, 5\}$ and $\{3, 4\}$ is such a closed component of this torus and hence must have just two edges.    
\vspace{-2pt}

Therefore a disk of $S \cap X^-$ whose edges have label pairs among $\{2, 5\}$ and $\{3, 4\}$ must be a bigon.  Since \fullref{lem:parallelbigons} implies this bigon is parallel to $B$, the lemma at hand is satisfied.
\vspace{-2pt}

Assume $n>0$.  Construct a high disk $D^+$ for the arc $K_{(5,6)}$ as in \fullref{lem:nobigonortrigonon1456}.  See \fullref{fig:startinghighdisk}.  Assume $R$ is an $m$--gon of $S \cap X^-$ for $m>2$ with edges having label pairs among $\{2, 3\}$, $\{2, 5\}$, and $\{3, 4\}$.  
\vspace{-2pt}

{\bf Case 1}\qua  The edges of $B$ and $R$ lie in an essential annulus $A$ in $\hatT$.  
\vspace{-2pt}

By taking $A$ smaller if necessary we may assume the interior of the arc $D^+ \cap \hatT$ is disjoint from $A$.   Let $M = \bar{N}(A \cup H_{(2, 3)} \cup H_{(4, 5)} \cup B \cup R)$.

Due to the existence of $B$, we may assume $R$ has at least one edge with label pair $\{2, 3\}$ or else the argument assuming $n=0$ applies.  Therefore $\bdry R$ is a nonseparating curve on $\bdry \bar{N}(A \cup H_{(2, 3)} \cup H_{(4, 5)} \cup B)$, and $\bdry M$ is a torus.  Since $M$ is contained in the solid torus $X^- - N(f \cup H_{(1,6)})$ in which the core of $A$ is longitudinal, $M$ must also be a solid torus such that the core of $A$ is a longitude in $M$.

Let $b_i$, $p_i$, and $\bar{p}_i$ for $i = 0,1$ be as in \fullref{lem:arrangingedges2345}.  Then by \fullref{lem:arrangingedges2345} $R$ has exactly one edge with label pair $\{2, 3\}$ and either $\bar{p}_0=0$, $\bar{p}_1=0$, or both $p_0 = \bar{p}_1$ and $p_1 = \bar{p}_0$.  In particular, either $p_0 \neq 0$ or $p_1 \neq 0$ and always $\bar{p}_0 \leq p_1$.  Note that the edges of $R$ parallel to $b_0$ lie on the opposite side of $B$ from the edges parallel to $b_1$.

Beginning with the edge of $R$ adjacent to $b_0$ label the $p_0$ edges of $R$ parallel to $b_0$ as $e_j$ for $j = 1, \dots, p_0$.  Then beginning with the edge of $R$ adjacent to $b_1$ label the first $\bar{p}_0$ edges that are parallel to $b_1$ as $e_j$ for $j = p_0 +1, \dots, (m-1)/2$.  If either $p_0 = 0$ or $p_1 = 0$ then all edges $e_j$ for $j = 1, \dots, (m-1)/2$ are of the second type or first type respectively.   Also for each $j$ let $B_j$ be a copy of $B$.  

For $j = 1, \dots, p_0$, let $\Delta_j$ be the bigon on $\hatT$ bounded by $e_j$ and $b_0$.  Let $\rho_j'$ and $\rho_j''$ be the rectangles on $\bdry H_{(2, 3)}$ and $\bdry H_{(4, 5)}$ respectively bounded by the corners of $R$ incident to $e_j$, the corners of $\Delta_j$, the corners of $B_j$, and the appropriate final arcs of $\bdry U_3$ and $\bdry U_4$.  Name these final arcs $\gamma_j'$ and $\gamma_j''$ of $\rho_j'$ and $\rho_j''$ respectively.

For $j = p_0+1, \dots, (m-1)/2$, let $\Delta_j$ be the bigon on $\hatT$ bounded by $e_j$ and $b_1$.  Let $\rho_j'$ and $\rho_j''$ be the rectangles on $\bdry H_{(2, 3)}$ and $\bdry H_{(4, 5)}$ respectively bounded by the corners of $R$ incident to $e_j$, the corners of $\Delta_j$, the corners of $B_j$, and the appropriate final arcs of $\bdry U_2$ and $\bdry U_5$.  Name these final arcs $\gamma_j'$ and $\gamma_j''$ of $\rho_j'$ and $\rho_j''$ respectively.

Form the disks $D_j = \Delta_j \cup B_j \cup \rho_j' \cup \rho_j''$ for each $j = 1, \dots, (m-1)/2$.  Sequentially attach each disk $D_j$ to $R$ along $e_j$ and the corners incident to it and slightly isotop the current attached disk off of the remaining disks.  The resulting complex is a low disk $D^-$ for $K_{(2, 3)}$.  (Alternatively, one may conceive of $D^-$ as the result of ``sliding'' each edge $e_j$ over the bigon $B$.)  By construction, the arcs $\gamma_j''$ are disjoint from the edges of $G_T$ with label pair $\{5, 6\}$.  It follows that $D^+ \cap \hatT$ and $D^- \cap \hatT$ have disjoint interiors.  Thus by \fullref{highdisklowdisk} the existence of the pair of disks $D^+$ and $D^-$ contradicts the thinness of $K$.

{\bf Case 2}\qua  The edges of $B$ and $R$ lie in a disk $D$ in $\hatT$.

Consider the solid torus $M = \bar{N}(D \cup H_{(2, 3)} \cup H_{(4, 5)} \cup B)$.  If $R$ has more than one edge with label pair $\{2, 3\}$ then $M \cup \bar{N}(R)$ is a punctured lens space contained in a solid torus which cannot occur.  If $R$ has no edges with label pair $\{2, 3\}$ then the argument assuming $n=0$ applies.  Thus we may assume $R$ has just one edge with label pair $\{2, 3\}$.

The procedure of Case 1 of constructing a low disk $D^-$ that is disjoint from $D^+$ now applies with $p_0 = p_1 = (m-1)/2$ and $\bar{p}_0 = \bar{p}_1 = 0$.
\end{proof}

\begin{prop}\label{prop:tnot6}
$t \neq 6$
\end{prop}

\begin{proof} 
Assuming $t=6$, we shall find contradictions to \fullref{lem:lowerfacesbound}.

If a disk of $S \cap X^-$ is neither a bigon nor parallel to $f$, then by \fullref{lem:onlybigons} it must have an edge with label pair $\{1, 4\}$ and an edge with label pair $\{5, 6\}$.  Since there are only $n$ edges with label pair $\{1, 4\}$ (and only $n$ edges with label pair $\{5, 6\}$), there may be at most $n$ disks that are neither bigons nor parallel to $f$.

If a disk of $S \cap X^-$ is a bigon yet not parallel to $f$, then again by \fullref{lem:onlybigons} its edges must have label pairs $\{2, 5\}$ and $\{3, 4\}$, and all such bigons are parallel.  Thus there may be at most as many bigons of $S \cap X^-$ not parallel to $f$ as the number of edges in an edge (parallelism) class with label pair $\{3, 4\}$.  As a result of the isotopies of \fullref{prop:isotopS}, if $g$ is a bigon then each edge class with label pair $\{3,4\}$ contains $(s-n)/2$ edges.  If $g$ is a trigon then either there  are $(s-n)/3$ trigons parallel to $g$ whose edges account for all edges with label pair $\{3,4\}$ or there are $(s-2n)/3$ trigons parallel to $g$ whose edges along with an edge from each of the $n$ trigons that are meridional disks of $W$ account for all edges with label pair $\{3, 4\}$.  In the former case, there are at most $2(s-n)/3$ edges in an edge class with label pair $\{3,4\}$.  In the latter case, there are at most $2(s-2n)/3+n$ edges in an edge class with label pair $\{3,4\}$.  Therefore if $g$ is a bigon there are at most $(s-n)/2$ bigons of $S \cap X^-$ not parallel to $f$, and if $g$ is a trigon there are at most $2(s-2n)/3+n$ bigons of $S \cap X^-$ not parallel to $f$.

Any other disk of $S \cap X^-$ must be parallel to $f$.

{\bf Case 1}\qua $g$ is a bigon.\qua
By \fullref{lem:lowerfacesbound}, 
\[s+ \frac{1}{2} + \frac{n}{2} \leq \sum_{R \in S \cap X^-} \chi(R). \]
{\bf Case 1a}\qua $f$ is a bigon.\qua
Since $f$ is a bigon then by \fullref{prop:isotopS} there are $(s-n)/2$ faces of $S \cap X^-$ parallel to $f$.  Since $g$ is a bigon then there are at most $(s-n)/2$ bigons of $S \cap X^-$ not parallel to $f$.  Since there are at most $n$ disks of $S \cap X^-$ that are not bigons (parallel to $f$ or otherwise), 
\begin{align*}
 \sum_{R \in S \cap X^-} \chi(R) &= \sum_{\mbox{{\scriptsize disks} }R \in S \cap X^-} \chi(R) + \sum_{\mbox{{\scriptsize nondisks} }R \in S \cap X^-} \chi(R) \\
    &\leq \Big(\frac{s-n}{2} + \frac{s-n}{2} + n\Big) + 0\\
    & = s.
\end{align*}
But $s+ \frac{1}{2} + \frac{n}{2} \leq s$ is a contradiction.

{\bf Case 1b}\qua $f$ is a trigon.\qua
Since $f$ is a trigon then by \fullref{prop:isotopS} there are at most $(s-n)/3$ faces of $S \cap X^-$ parallel to $f$.  Since $g$ is a bigon then there are at most $(s-n)/2$ bigons of $S \cap X^-$ not parallel to $f$.  Since there are at most $n$ disks of $S \cap X^-$ that are neither bigons nor trigons parallel to $f$,
\begin{align*}
 \sum_{R \in S \cap X^-} \chi(R) &= \sum_{\mbox{{\scriptsize disks} }R \in S \cap X^-} \chi(R) + \sum_{\mbox{{\scriptsize nondisks} }R \in S \cap X^-} \chi(R) \\
    &\leq \Big(\frac{s-n}{3} + \frac{s-n}{2} + n\Big) + 0\\
    & = \frac{5}{6}s+\frac{n}{6}n.
\end{align*}
But $s+ \frac{1}{2} + \frac{n}{2} \leq \frac{5}{6}s+\frac{n}{6}$ is a contradiction.

{\bf Case 2}\qua $g$ is a trigon.\qua
By \fullref{lem:lowerfacesbound}, 
\[\frac{7}{6}s+ \frac{1}{2} + \frac{n}{3} \leq \sum_{R \in S \cap X^-} \chi(R). \]
{\bf Case 2a}\qua $f$ is a bigon.\qua
Since $f$ is a bigon then by \fullref{prop:isotopS} there are $(s-n)/2$ faces of $S \cap X^-$ parallel to $f$.  Since $g$ is a trigon then there are at most $2(s-2n)/3+n$ bigons of $S \cap X^-$ not parallel to $f$.  Since there are at most $n$ disks of $S \cap X^-$ that are not bigons (parallel to $f$ or otherwise), 
\begin{align*}
 \sum_{R \in S \cap X^-} \chi(R) &= \sum_{\mbox{{\scriptsize disks} }R \in S \cap X^-} \chi(R) + \sum_{\mbox{{\scriptsize nondisks} }R \in S \cap X^-} \chi(R) \\
    &\leq \Big(\frac{s-n}{2} + 2\frac{s-2n}{3}+n + n\Big) + 0\\
    & = \frac{7}{6}s+\frac{n}{6}.
\end{align*}
But $\frac{7}{6}s+ \frac{1}{2} + \frac{n}{3} \leq \frac{7}{6}s+\frac{n}{6}$ is a contradiction.

{\bf Case 2b}\qua $f$ is a trigon.\qua
Since $f$ is a trigon then by \fullref{prop:isotopS} there are at most $(s-n)/3$ faces of $S \cap X^-$ parallel to $f$.  Since $g$ is a trigon then there are at most $2(s-2n)/3+n$ bigons of $S \cap X^-$ not parallel to $f$.  Since there are at most $n$ disks of $S \cap X^-$ that are not bigons (parallel to $f$ or otherwise), 
\begin{align*}
 \sum_{R \in S \cap X^-} \chi(R) &= \sum_{\mbox{{\scriptsize disks} }R \in S \cap X^-} \chi(R) + \sum_{\mbox{{\scriptsize nondisks} }R \in S \cap X^-} \chi(R) \\
    &\leq \Big(\frac{s-n}{3} + 2\frac{s-2n}{3}+n + n\Big) + 0\\
    & = s+\frac{n}{3}.
\end{align*}
But $\frac{7}{6}s+ \frac{1}{2} + \frac{n}{3} \leq s+\frac{n}{3}$ is a contradiction.
\end{proof}

\section[The case t=4]{The case $t=4$}\label{sec:bridgeposition}

We continue to assume that $K$ is in thin position and that $r \geq 4$.  In light of \fullref{prop:tnot6}, we have that $|K \cap \hatT| = t \leq 4$.  If $t=4$ we will show that $K$ must also be in bridge position.  Then since each graph $\smash{G_S^x}$ for $x = 1, 2, 3, 4$ must have a bigon or a trigon by \fullref{musthavebigons}, we will exhibit a width reducing isotopy of $K$ that contradicts its supposed thinness.

The following lemma relating bridge position and thin position is an immediate consequence of definitions.
\begin{lemma}\label{thintobridge}
If a knot $K$ in thin position has no thin levels, then it is in bridge position.
\end{lemma}

To show that thin position of $K$ is bridge position, we must consider the existence of a second thick level.  If $\hatT'$ is another thick level, let $T' = \hatT' - N(K)$.  We may assume that $S$ has been isotoped (with support outside a neighborhood of $\hatT$) so that $S$ and $\hatT'$ intersect transversely and every arc of $S \cap T'$ is essential in both $S$ and $T'$.  We then form the fat vertexed graphs $G_T'$ and $G_S'$.  Since $\hatT$ and $\hatT'$ are disjoint, the edges of the graphs $G_S$ and $G_S'$ are disjoint.  Hence the graphs $G_S$ and $G_S'$  may be considered simultaneously on $\what{S}$.
\vspace{-5pt}

\begin{lemma}\label{nextthicklevel}
Assume $f$ is the face of an $Sp$ cycle for $p=2$ or $3$ with label pair $\{x, x{+}1\}$ so that $\bdry f \cap K = K_{(x,\,x+1)}$.  Further assume that $f$ is below (resp.\ above) $\hatT$.  If $\hatT'$ is another thick level below (resp.\ above) $\hatT$ then $\hatT' \cap K_{(x,\,x+1)} \neq \emptyset$.  In particular $f$ contains the face of an $Sp$ cycle of $G_S'$.
\end{lemma}
\vspace{-5pt}

\begin{proof}
We show the case that $f$ is below $\hatT$.  
\vspace{-2pt}

By \fullref{GT:L2.1} the edges of $f$ lie in an essential annulus $A$ on $\hatT$ and the core of $A$ runs $p$ times longitudinally in $X^-$.  
\vspace{-2pt}

If $\hatT'$ is disjoint from $f$, then $f$ is disjoint from $\hatT_{-\infty}$.   Thus the curve $\hatT_{-\infty}$ must be isotopic to a curve on $\hatT$ that is disjoint from $A$.  This contradicts that $\hatT_{-\infty}$ runs just once in the longitudinal direction of $X^-$.  Therefore $f \cap \what{T'} \neq \emptyset$.
\vspace{-2pt}

If $\hatT'$ is disjoint from $K_{(x,\, x{+}1)}$, then $\hatT'$ intersects $f$ in simple closed curves.  A simple closed curve of $\hatT' \cap f$ innermost on $f$ bounds a disk.  This contradicts \fullref{circlesofintersection}.  Thus $\hatT' \cap K_{(x,\, x{+}1)} \neq \emptyset$. 
\vspace{-2pt}

Since the each arc of $S \cap T'$ is essential in $S$, each edge of $G_S'$ on $f$ must be parallel on $f$ to one of the edges of $\bdry f$.  Note that $\hatT'$ intersects $K_{(x,\, x{+}1)}$ an even number of times.  Each of these intersections appear on each corner of $f$.  Since $p = 2$ or $3$, it is clear that $f$ contains an $Sp$ cycle of $G_S'$.
\end{proof}
\vspace{-5pt}

\begin{thm}\label{thm:t=2}
$t=2$
\end{thm}
\vspace{-5pt}

\begin{figure}[ht!]
\centering
\begin{picture}(0,0)%
\includegraphics{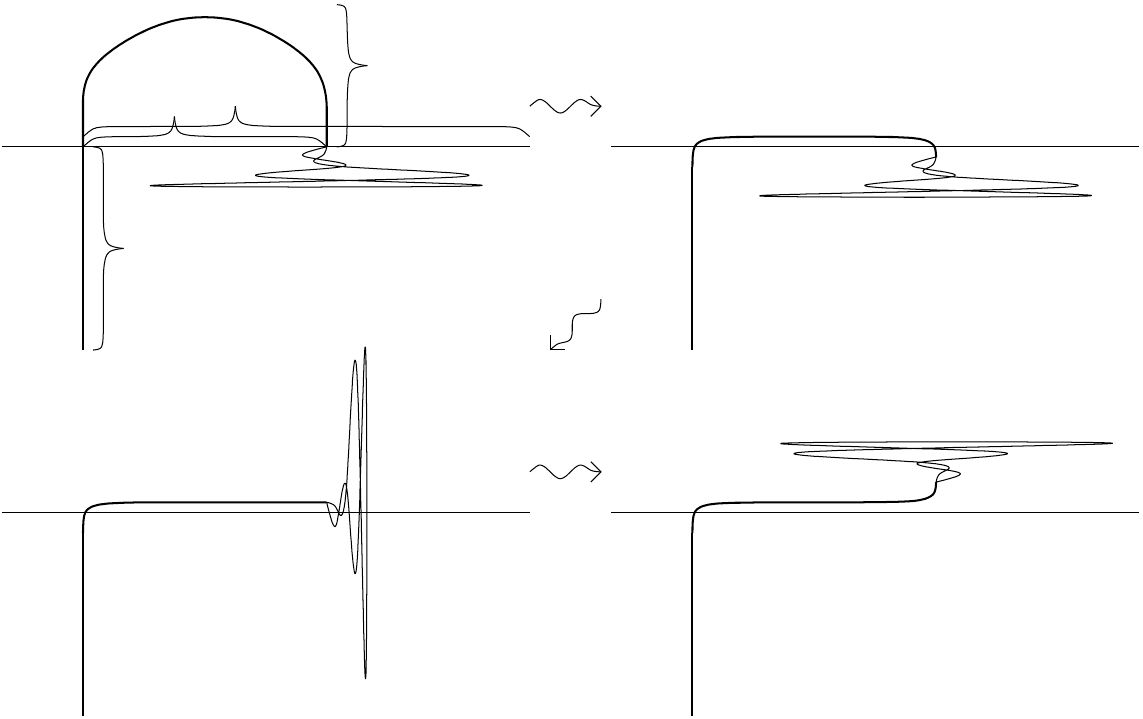}%
\end{picture}%
\setlength{\unitlength}{2565sp}%
\begingroup\makeatletter\ifx\SetFigFont\undefined%
\gdef\SetFigFont#1#2#3#4#5{%
  \reset@font\fontsize{#1}{#2pt}%
  \fontfamily{#3}\fontseries{#4}\fontshape{#5}%
  \selectfont}%
\fi\endgroup%
\begin{picture}(8424,5284)(1489,-9383)
\put(2701,-4936){\makebox(0,0)[lb]{\smash{{\SetFigFont{8}{9.6}{\rmdefault}{\mddefault}{\updefault}{\color[rgb]{0,0,0}$R$}%
}}}}
\put(3151,-4861){\makebox(0,0)[lb]{\smash{{\SetFigFont{8}{9.6}{\rmdefault}{\mddefault}{\updefault}{\color[rgb]{0,0,0}$R'$}%
}}}}
\put(4426,-5761){\makebox(0,0)[b]{\smash{{\SetFigFont{8}{9.6}{\rmdefault}{\mddefault}{\updefault}{\color[rgb]{0,0,0}$K_{(2, 3)}$}%
}}}}
\put(4276,-4636){\makebox(0,0)[lb]{\smash{{\SetFigFont{8}{9.6}{\rmdefault}{\mddefault}{\updefault}{\color[rgb]{0,0,0}$K_{(1, 2)}, K_{(3, 4)} \subset F\cut f$}%
}}}}
\put(1651,-5086){\makebox(0,0)[b]{\smash{{\SetFigFont{8}{9.6}{\rmdefault}{\mddefault}{\updefault}{\color[rgb]{0,0,0}$\hatT$}%
}}}}
\put(2476,-5986){\makebox(0,0)[lb]{\smash{{\SetFigFont{8}{9.6}{\rmdefault}{\mddefault}{\updefault}{\color[rgb]{0,0,0}$K_{(4, 1)} \subset f$}%
}}}}
\end{picture}%
\caption{The isotopy of $K$ from a $2$--bridge position to a $1$--bridge position}
\label{fig:2bridgeisotopy}
\end{figure}
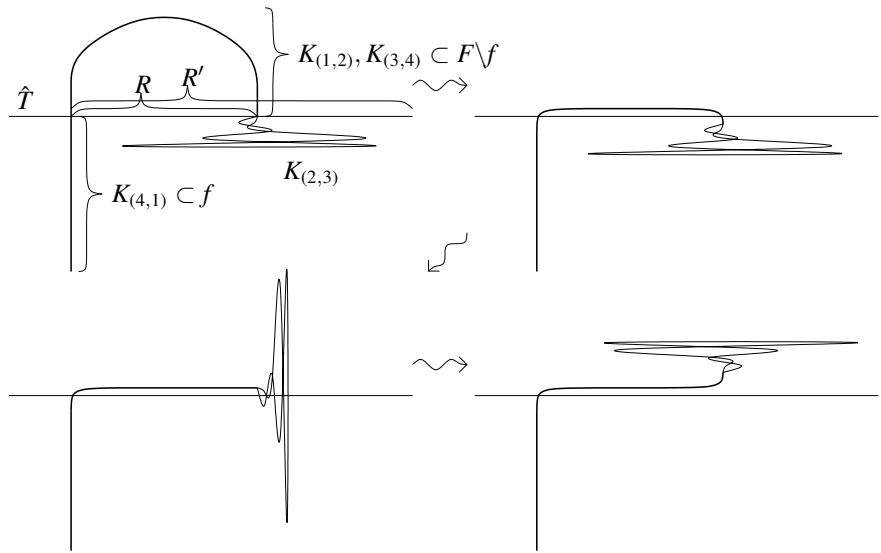

\begin{proof}
By \fullref{prop:tnot6}, $t \leq 4$.  Assume $t=4$.
\vspace{-2pt}

By \fullref{musthavebigons} for each of the labels $y \in \mathbf{t} = \{1, 2, 3, 4\}$ the graph $G_S^y$ must have a bigon or a trigon.  
\vspace{-2pt}

Due to \fullref{oppositesides}, since $t=4$ any two distinct label pairs of $S2$ and $S3$ cycles in $G_S$ must intersect.  Therefore the set of label pairs of $S2$ and $S3$ cycles in $G_S$ has cardinality at most two.  Hence the $S2$ and $S3$ cycles in $G_S$ account for at most three of the four labels.  Thus for some $y \in \mathbf{t}$ the graph $G_S^y$ contains an extended $S2$ cycle, an extended $S3$ cycle, or a forked extended $S2$ cycle.

Let $F$ be the face of this extended $S2$ cycle, extended $S3$ cycle, or forked extended $S2$ cycle.  Let $\sigma$ be the Scharlemann cycle contained in the interior of $F$ and let $f$ be the face of $G_S$ that $\sigma$ bounds.  By relabeling we may assume that $f$ lies below $\hatT$ and that $\sigma$ has label pair $\{1, 4\}$.  Since $r \geq 4$, by \fullref{GT:L2.1} the edges of $\sigma$ lie in an essential annulus in $\hatT$.  The arc $K_{(2, 3)}$ is the only arc of $K$ not appearing among the corners of $F$. 

If $F$ is bounded by an extended $S2$ cycle or extended $S3$ cycle, then the arcs $K_{(1, 2)}$ and $K_{(3, 4)}$ lie on an annulus formed from bigons of $F \cut f$.  This annulus is necessarily parallel onto $\hatT$ through one of the two solid tori it cuts off from $X^+$.  Since $K$ is in thin position, the arcs $K_{(1, 2)}$ and $K_{(3, 4)}$ each have only one critical point (a maximum) in their interiors.

If $F$ is bounded by a forked extended $S2$ cycle, then the arcs $K_{(1, 2)}$ and $K_{(3, 4)}$ lie on the complex formed from the bigon and trigon of $F \cut f$.  As in \fullref{thinningtwoforks}, we may construct a high disk $D_g$ for one of the two arcs of $K$.  Since the two arcs of $K \cap X^+$ are joined by a bigon of $G_S$, it follows that they each have only one critical point in their interiors.

Regardless of whether $F$ is bounded by an extended Scharlemann cycle or a forked extended Scharlemann cycle, the two arcs $K_{(1, 2)}$ and $K_{(3, 4)}$ each have high disks that intersect $\hatT$ in an annulus $R$ that contains edges of $f$ in just one of its boundaries.

\fullref{highestmin} applies to show that the arc of $K$ containing the highest minimum below $\hatT$ bounds a low disk $D$ with an arc of $\hatT$ such that the interior of $D$ is disjoint from $F$ and $K$.

By \fullref{nextthicklevel}, if $K_{(4, 1)}$ contains the highest minimum then there cannot be another thick level below $\hatT$.  Thus $K$ is in bridge position.  Therefore the arc $K_{(2, 3)}$ has a low disk $D$ with interior disjoint from $f$.   Alternately, if we assume the arc $K_{(2, 3)}$ contains the highest minimum, then by \fullref{highestmin} $K_{(2, 3)}$ again has a low disk $D$ with interior disjoint from $f$.  

Since the high disks of the arcs $K_{(1, 2)}$ and $K_{(3, 4)}$ intersect $\hatT$ in $R$ and the interior of the low disk $D$ for the arc $K_{(2, 3)}$ is disjoint from $f$, these three arcs are contained within a solid torus $R' \times [-\epsilon, \epsilon]$ for some annulus $R' \subseteq \hatT$ that extends $R$ to contain $D \cap \hatT$ and some $\epsilon > 0$.  As $R$ contains edges of $f$ in just one of its boundaries, so does $R'$.  Hence $R' \times [-\epsilon, \epsilon]$ may be assumed to have its interior disjoint from $f$.  

Via the two high disks, isotop the arcs $K_{(1, 2)}$ and $K_{(3, 4)}$ onto $R$.  Then, along the component of $\bdry R$ in $\Int R'$, pivot $D$ in $R' \times [-\epsilon, \epsilon]$ $180^\circ$ through $R' \cut R$ so that it is contained in $X^+$.  Finally, a further small isotopy to make $K_{(1, 2)}$ and $K_{(3, 4)}$ transverse to the height function will reduce the width of $K$.  See \fullref{fig:2bridgeisotopy}.  This contradicts the thinness of $K$.  Hence $t \neq 4$.
\end{proof}

\begin{lemma}\label{lem:t=2}
If $t=2$ then $K$ is at most $1$--bridge.
\end{lemma}

\begin{proof}
Since the first critical point of $K$ above $\hatT$ is a maximum, there can be no other critical points of $K$ above $\hatT$.  Similarly, the arc of $K$ below $\hatT$ may have no critical points other than a minimum.  Hence there are no thin levels.  By \fullref{thintobridge} $K$ is in bridge position.  Since $t=2$, $b(K) \leq 1$. 
\end{proof}

\section[The case r=3]{The case $r=3$}\label{sec:genus1case}

In this section we assume $r=3$.  Since $r \geq s \geq 4g-1$, either $g=0$ or  $r = s$ and $g=1$.  If $g=0$ then $N(S) \cup N(K)$ is a punctured lens space of order $3$ and hence $K$ is $0$--bridge.  Thus we assume $g=1$.  Goda and Teragaito have shown the following.
\begin{thm}[Goda--Teragaito~\cite{gt:dsokwylsagok}]
No Dehn surgery on a genus one, hyperbolic knot in $S^3$ gives a lens space.
\end{thm}
We adapt their proof of this theorem for our case that $r=3$ and $g=1$.

\begin{thm} \label{thm:r=3}
If $r=3$ then $b(K) \leq 1$.
\end{thm} 

To prove this, we follow Section~5 of \cite{gt:dsokwylsagok}.  First we need to adapt some lemmas from Section~3 of \cite{gt:dsokwylsagok}.
\vspace{-8pt}

\begin{lemma}[cf {{\cite[Lemma~3.1]{gt:dsokwylsagok}}}]
\label{lem:GT3.1}
Let $\{e_1, e_2, \dots, e_t\}$ be mutually parallel edges in $G_S$ numbered successively.  If $r$ is odd, then $\{e_{t/2}, e_{t/2+1}\}$ is an $S2$ cycle.
\end{lemma}
\vspace{-8pt}

\begin{proof}
As in the first paragraph of the proof of Lemma~3.1 of \cite{gt:dsokwylsagok} we assume $e_i$ has the label $i$ at one endpoint for $1 \leq i \leq t$ and that $e_2j$ has the label $1$ at its other endpoint for some $j < t/2$.  Thus we obtain two $S2$ cycles $\sigma_1$ and $\sigma_2$ with disjoint label pairs.  Similarly we let $f_i$ be the face of $G_S$ bounded by $\sigma_i$ for $i=1,2$.  
\vspace{-2pt}

Since $r \neq 2$, by \fullref{GT:L2.3} the edges of each $\sigma_1$ and $\sigma_2$ do not lie in a disk.  Thus they lie in an essential annulus.  The proof of \fullref{oppositesides} then carries through to show that $f_1$ and $f_2$ lie on opposite sides of $\hatT$.  Furthermore, by disk exchanges outside of $N(\sigma_1 \cup \sigma_2)$, we may construct a Heegaard torus $\hatT'$ from $\hatT$ so that the edges of $\sigma_i$ lie in an essential annulus in $\hatT'$ and $\Int f_i \cap \hatT' = \emptyset$ for $i=1,2$.
\vspace{-2pt}

We then construct two disjoint M\"obius bands $B_1$ and $B_2$ on either side of $\hatT'$ as in the proof of Lemma~2.5 of \cite{gt:dsokwylsagok} and the beginning of our \fullref{sec:annuliandtrees}.  Since $\bdry B_1$ and $\bdry B_2$ are parallel on $\hatT'$, they divide $\hatT'$ into two annuli.  Let $A$ be one of these annuli.  Then $B_1 \cup A \cup B_2$ is an embedded Klein bottle in our lens space $X$.  Thus the order of $X$ is even.  This contradicts that $r$ is odd.
\vspace{-2pt}

Therefore the edge $e_t$ has the label $1$ at its other endpoint and $\{e_{t/2}, e_{t/2+1}\}$ is an $S2$ cycle.
\end{proof} 
\vspace{-8pt}

\begin{lemma}{\rm (cf\ \cite[Lemma~3.2]{gt:dsokwylsagok})}\qua\label{lem:GT3.2}
If $r$ is odd, then $G_S$ does not contain more than $t$ mutually parallel edges.
\end{lemma}
\vspace{-8pt}

\begin{proof}
In the proof of Lemma~3.2 of \cite{gt:dsokwylsagok}, replace each occurrence of Lemma~3.1 with our \fullref{lem:GT3.1}.
\end{proof}
\vspace{-8pt}

\begin{lemma}[cf {{\cite[Lemma~5.2]{gt:dsokwylsagok}}}]
\label{lem:GT5.2} 
If $r$ is odd, then $G_S$ cannot have more than $t/2$ mutually parallel edges.
\end{lemma}

\begin{proof}
In the proof of Lemma~5.2 of \cite{gt:dsokwylsagok}, replace Lemma~3.1 with our \fullref{lem:GT3.1} and Lemma~3.2 with our \fullref{lem:GT3.2}.
\end{proof}

\begin{proof}[Proof of \fullref{thm:r=3}]
If $t=2$ then \fullref{lem:t=2} implies $b(K) \leq 1$.  Thus we assume $t \geq 4$.  Also assume the interior of $S$ has been isotoped to minimize $|S \cap T|$.

Following the proof of Lemma~5.3 of \cite{gt:dsokwylsagok}:  The vertex of $G_S$ has valency $3t$, and there are a total of $3t/2$ edges of $G_S$.  By \fullref{lem:GT5.2}, $G_S$ consists of three families of mutually parallel edges each containing exactly $t/2$ edges.  Thus there is no $S2$ cycle in $G_S$, but there are two $S3$ cycles $\sigma_1$ and $\sigma_2$ which we may assume to have label pairs $\{t, 1\}$ and $\{t/2, t/2+1\}$ respectively by an appropriate choice of relabeling.  Let $g_i$ be the face of $G_S$ bounded by $\sigma_i$ for $i=1,2$.  

Note that each graph $\smash{G_S^x}$ for $x \in \mathbf{t}$ consists of three edges with label pair $\{x, t-x+1\}$.  Each such triple of edges $\smash{G_S^x}$ is a ``twice'' (extended) $S3$ cycle:  to one side it bounds a trigon containing $g_1$ and to the other side it bounds a trigon containing $g_1$.  

The proof of Claim~5.5 of \cite{gt:dsokwylsagok} carries through without alteration to show each of these extended $S3$ cycles lie in essential annuli in $\hatT$.  Thus $S \cap T$ cannot contain a simple closed curve that bounds a disk in $\hatT$ without bounding a disk in $T$.  Because every face of $G_S$ is a disk, every simple closed curve of $S \cap T$ must bound a disk on $S$.  Thus each simple closed curve of $S \cap T$ that is trivial on $\hatT$ is also trivial on $S$.  

Let $D_T$ be an disk on $\hatT$ bounded by a simple closed curve of $S \cap T$ that is both trivial and innermost on $\hatT$.  Let $D_S$ be the disk in $S$ bounded by $\bdry D_T$.  Since $D_T$ must be disjoint from $K$, $D_S \cup D_T$ forms a $2$--sphere that is disjoint from $K$.  Since lens spaces are irreducible, $D_S \cup D_T$ bounds a ball $B$.  Because $K$ is not nullhomologous in our lens space $X$, $K \not \subseteq B$.  Thus there is an isotopy of $\Int S$ with support in a neighborhood of $B$ pushing $D_S$ past $D_T$ thereby reducing $|S \cap T|$.  This contradicts our minimality assumption.  Thus any simple closed curve of $S \cap T$ is essential on $\hatT$.  Furthermore, any simple closed curve of $S \cap T$ innermost on $S$ must bound a meridional disk of the same solid torus on one side of $\hatT$.  

For each $i=1,2$, let $A_i$ be a narrow annulus in $\hatT$ in which $\sigma_i$ lies that does not contain any simple closed curves of $S \cap T$.  If $S \cap T$ does indeed contain simple closed curves, then let $\delta$ be an innermost simple closed curve on $S$ bounding the disk $D$.  The cores of the $A_i$ are parallel on $\hatT$ to $\delta$, and $\delta \subseteq \hatT \cut (A_1 \cup A_2)$.  For each $i=1,2$, let $A_i'$ be the annulus on $\hatT$ between $A_i$ and $\delta$, and let $D_i' = A_i \cup A_i' \cup D$ slightly pushed off $D$ so that $D_1'$ and $D_2'$ are disjoint from each other and from $D$.  If $\Int(g_1 \cup g_2) \cap (A_1' \cup A_2') \neq \emptyset$ then perform disk exchanges on $g_1 \cup g_2$ with $D_1' \cup D_2'$ to produce trigons $g_i'$ from $g_i$ for each $i=1,2$ so that $\Int(g_1' \cup g_2') \cap (A_1' \cup A_2') = \emptyset$.  Note that $\bdry g_i' = \bdry g_i$.  Then $\bar{N}(D_1' \cup H_{(t,1)} \cup g_1')$ and $\bar{N}(D_2' \cup H_{(t/2, t/2+1)} \cup g_2')$ gives two disjoint punctured lens spaces in $X$, which is absurd.  Thus $S \cap T$ contains no simple closed curves.

The proof of Claim~5.4 of \cite{gt:dsokwylsagok} (from which the preceding paragraph takes inspiration) may now be used to show that $g_1$ and $g_2$ lie on opposite sides.  \fullref{GT:L2.1} then implies that the cores of the $A_i$ are not meridional curves for the solid tori on either side of $\hatT$.  Thus $t/2$ is odd, and so $t \geq 6$.

\fullref{sec:annuliandtrees} now applies to the extended $S3$ cycle $\sigma_2$ considered as bounding the trigon $F$ containing $g_1$.  \fullref{bounded} implies that $F$ may account for at most $t/2+1$ labels.  This contradicts that $F$ accounts for all $t$ labels.
\end{proof}

\bibliographystyle{gtart}
\bibliography{link}

\end{document}